*Set Linear Algebra and*
*Set Fuzzy Linear Algebra*

# SET LINEAR ALGEBRA AND SET FUZZY LINEAR ALGEBRA

**W. B. Vasantha Kandasamy**
e-mail: vasanthakandasamy@gmail.com
web: http://mat.iitm.ac.in/~wbv
www.vasantha.net

**Florentin Smarandache**
e-mail: smarand@unm.edu

**K Ilanthenral**
e-mail: ilanthenral@gmail.com

# CONTENTS

# PREFACE

In this book, the authors define the new notion of set vector spaces which is the most generalized form of vector spaces. Set vector spaces make use of the least number of algebraic operations, therefore, even a non-mathematician is comfortable working with it. It is with the passage of time, that we can think of set linear algebras as a paradigm shift from linear algebras. Here, the authors have also given the fuzzy parallels of these new classes of set linear algebras.

This book abounds with examples to enable the reader to understand these new concepts easily. Laborious theorems and proofs are avoided to make this book approachable for non-mathematicians.

The concepts introduced in this book can be easily put to use by coding theorists, cryptologists, computer scientists, and socio-scientists.

Another special feature of this book is the final chapter containing 304 problems. The authors have suggested so many problems to make the students and researchers obtain a better grasp of the subject.

This book is divided into seven chapters. The first chapter briefly recalls some of the basic concepts in order to make this book self-contained. Chapter two introduces the notion of set vector spaces which is the most generalized concept of vector spaces. Set vector spaces lends itself to define new classes of vector spaces like semigroup vector spaces and group vector spaces. These are also generalization of vector spaces. The fuzzy analogue of these concepts are given in Chapter three.

In Chapter four, set vector spaces are generalized to biset bivector spaces and not set vector spaces. This is done taking into account the advanced information technology age in which we live. As mathematicians, we have to realize that our computer-dominated world needs special types of sets and algebraic structures.

Set n-vector spaces and their generalizations are carried out in Chapter five. Fuzzy n-set vector spaces are introduced in the sixth chapter. The seventh chapter suggests more than three hundred problems. When a researcher sets forth to solve them, she/he will certainly gain a deeper understanding of these new notions.

Our thanks are due to Dr. K. Kandasamy for proof-reading this book. We also acknowledge our gratitude to Kama and Meena for their help with the corrections and layout.

W.B.VASANTHA KANDASAMY
FLORENTIN SMARANDACHE
K.ILANTHENRAL

## Chapter One

# BASIC CONCEPTS

This chapter has three sections. In section one a brief introduction to linear algebra is given. Section two gives the basic notions in bilinear algebra and the final section gives the definition of fuzzy vector spaces. For more about these concepts, please refer [48, 60].

## 1.1 Definition of Linear Algebra and its Properties

In this section we just recall the definition of linear algebra and enumerate some of its basic properties. We expect the reader to be well versed with the concepts of groups, rings, fields and matrices. For these concepts will not be recalled in this section.

Throughout this section, V will denote the vector space over F where F is any field of characteristic zero.

**DEFINITION 1.1.1:** *A vector space or a linear space consists of the following:*

*i.* *a field F of scalars.*
*ii.* *a set V of objects called vectors.*

iii. *a rule (or operation) called vector addition; which associates with each pair of vectors* $\alpha, \beta \in V$*;* $\alpha + \beta$ *in V, called the sum of* $\alpha$ *and* $\beta$ *in such a way that*

a. *addition is commutative* $\alpha + \beta = \beta + \alpha$.

b. *addition is associative* $\alpha + (\beta + \gamma) = (\alpha + \beta) + \gamma$.

c. *there is a unique vector 0 in V, called the zero vector, such that*

$$\alpha + 0 = \alpha$$

*for all* $\alpha$ *in V.*

d. *for each vector* $\alpha$ *in V there is a unique vector* $-\alpha$ *in V such that*

$$\alpha + (-\alpha) = 0.$$

e. *a rule (or operation), called scalar multiplication, which associates with each scalar c in F and a vector* $\alpha$ *in V, a vector* $c \bullet \alpha$ *in V, called the product of c and* $\alpha$*, in such a way that*

1. $1 \bullet \alpha = \alpha$ *for every* $\alpha$ *in V.*
2. $(c_1 \bullet c_2) \bullet \alpha = c_1 \bullet (c_2 \bullet \alpha)$.
3. $c \bullet (\alpha + \beta) = c \bullet \alpha + c \bullet \beta$.
4. $(c_1 + c_2) \bullet \alpha = c_1 \bullet \alpha + c_2 \bullet \alpha$.

*for* $\alpha, \beta \in V$ *and* $c, c_1 \in F$.

It is important to note as the definition states that a vector space is a composite object consisting of a field, a set of 'vectors' and two operations with certain special properties. The same set of vectors may be part of a number of distinct vectors.

We simply by default of notation just say V a vector space over the field F and call elements of V as vectors only as matter of convenience for the vectors in V may not bear much resemblance to any pre-assigned concept of vector, which the reader has.

***Example 1.1.1:*** Let R be the field of reals. R[x] the ring of polynomials. R[x] is a vector space over R. R[x] is also a vector space over the field of rationals Q.

***Example 1.1.2:*** Let Q[x] be the ring of polynomials over the rational field Q. Q[x] is a vector space over Q, but Q[x] is clearly not a vector space over the field of reals R or the complex field **C**.

***Example 1.1.3:*** Consider the set V = R × R × R. V is a vector space over R. V is also a vector space over Q but V is not a vector space over **C**.

***Example 1.1.4:*** Let $M_{m \times n} = \{(a_{ij}) \mid a_{ij} \in Q\}$ be the collection of all m × n matrices with entries from Q. $M_{m \times n}$ is a vector space over Q but $M_{m \times n}$ is not a vector space over R or **C**.

***Example 1.1.5:*** Let

$$P_{3 \times 3} = \left\{ \begin{pmatrix} a_{11} & a_{12} & a_{13} \\ a_{21} & a_{22} & a_{23} \\ a_{31} & a_{32} & a_{33} \end{pmatrix} \middle| \ a_{ij} \in Q, 1 \leq i \leq 3, \ 1 \leq j \leq 3 \right\}.$$

$P_{3 \times 3}$ is a vector space over Q.

***Example 1.1.6:*** Let Q be the field of rationals and G any group. The group ring, QG is a vector space over Q.

***<u>Remark</u>*****:** All group rings KG of any group G over any field K are vector spaces over the field K.

We just recall the notions of linear combination of vectors in a vector space V over a field F. A vector β in V is said to be a linear combination of vectors $\nu_1, \ldots, \nu_n$ in V provided there exists scalars $c_1, \ldots, c_n$ in F such that

$$\beta = c_1\nu_1 + \ldots + c_n\nu_n = \sum_{i=1}^{n} c_i \; \nu_i \; .$$

Now we proceed on to recall the definition of subspace of a vector space and illustrate it with examples.

**DEFINITION 1.1.2:** *Let V be a vector space over the field F. A subspace of V is a subset W of V which is itself a vector space over F with the operations of vector addition and scalar multiplication on V.*

We have the following nice characterization theorem for subspaces; the proof of which is left as an exercise for the reader to prove.

**THEOREM 1.1.1:** *A non-empty subset W of a vector V over the field F; V is a subspace of V if and only if for each pair $\alpha$, $\beta$ in W and each scalar c in F the vector $c\alpha + \beta$ is again in W.*

***Example 1.1.7:*** Let $M_{n \times n} = \{(a_{ij}) \mid a_{ij} \in Q\}$ be the vector space over Q. Let $D_{n \times n} = \{(a_{ii}) \mid a_{ii} \in Q\}$ be the set of all diagonal matrices with entries from Q. $D_{n \times n}$ is a subspace of $M_{n \times n}$.

***Example 1.1.8:*** Let $V = Q \times Q \times Q$ be a vector space over Q. P $= Q \times \{0\} \times Q$ is a subspace of V.

***Example 1.1.9:*** Let V = R[x] be a polynomial ring, R[x] is a vector space over Q. Take $W = Q[x] \subset R[x]$; W is a subspace of R[x].

It is well known results in algebraic structures. The analogous result for vector spaces is:

**THEOREM 1.1.2:** *Let V be a vector space over a field F. The intersection of any collection of subspaces of V is a subspace of V.*

*Proof***:** This is left as an exercise for the reader.

**DEFINITION 1.1.3:** *Let P be a set of vectors of a vector space V over the field F. The subspace spanned by W is defined to be the intersection of W of all subspaces of V which contains P, when P is a finite set of vectors, P =* $\{\alpha_1, \ldots, \alpha_m\}$ *we shall simply call W the subspace spanned by the vectors* $\alpha_1, \alpha_2, \ldots, \alpha_m$.

**THEOREM 1.1.3:** *The subspace spanned by a non-empty subset P of a vector space V is the set of all linear combinations of vectors in P.*

*Proof*: Direct by the very definition.

**DEFINITION 1.1.4:** *Let* $P_1, \ldots, P_k$ *be subsets of a vector space V, the set of all sums* $\alpha_1 + \ldots + \alpha_k$ *of vectors* $\alpha_i \in P_i$ *is called the sum of subsets of* $P_1, P_2, \ldots, P_k$ *and is denoted by* $P_1 + \ldots + P_k$ *or by* $\sum_{i=1}^{k} P_i$.

*If* $U_1, U_2, \ldots, U_k$ *are subspaces of V, then the sum*

$$U = U_1 + U_2 + \ldots + U_k$$

*is easily seen to be a subspace of V which contains each of the subspaces* $U_i$.

Now we proceed on to recall the definition of basis and dimension.

*Let V be a vector space over F. A subset P of V is said to be linearly dependent (or simply dependent) if there exists distinct vectors,* $\alpha_1, \ldots, \alpha_t$ *in P and scalars* $c_1, \ldots, c_k$ *in F not all of which are 0 such that* $c_1\alpha_1 + c_2\alpha_2 + \ldots + c_k\alpha_k = 0$.

*A set which is not linearly dependent is called independent. If the set P contains only finitely many vectors* $\alpha_1, \ldots, \alpha_k$ *we sometimes say that* $\alpha_1, \ldots, \alpha_k$ *are dependent (or independent) instead of saying P is dependent or independent.*

i. *A subset of a linearly independent set is linearly independent.*
ii. *Any set which contains a linearly dependent set is linearly dependent.*
iii. *Any set which contains the 0 vector is linear by dependent for 1.0 = 0.*
iv. *A set P of vectors is linearly independent if and only if each finite subset of P is linearly independent i.e. if and only if for any distinct vectors $\alpha_1, \ldots, \alpha_k$ of P, $c_1\alpha_1 + \ldots + c_k \alpha_k = 0$ implies each $c_i = 0$.*

*For a vector space V over the field F, the basis for V is a linearly independent set of vectors in V, which spans the space V. The space V is finite dimensional if it has a finite basis.*

We will only state several of the theorems without proofs as results and the reader is expected to supply the proof.

**Result 1.1.1:** Let V be a vector space over F which is spanned by a finite set of vectors $\beta_1, \ldots, \beta_t$. Then any independent set of vectors in V is finite and contains no more than t vectors.

**Result 1.1.2:** If V is a finite dimensional vector space then any two bases of V have the same number of elements.

**Result 1.1.3:** Let V be a finite dimensional vector space and let n = dim V. Then

i. any subset of V which contains more than n vectors is linearly dependent.
ii. no subset of V which contains less than n vectors can span V.

**Result 1.1.4:** If W is a subspace of a finite dimensional vector space V, every linearly independent subset of W is finite, and is part of a (finite) basis for W.

**Result 1.1.5:** If W is a proper subspace of a finite dimensional vector space V, then W is finite dimensional and dim W < dim V.

**Result 1.1.6:** In a finite dimensional vector space V every non-empty linearly independent set of vectors is part of a basis.

**Result 1.1.7:** Let A be a n × n matrix over a field F and suppose that row vectors of A form a linearly independent set of vectors; then A is invertible.

**Result 1.1.8:** If $W_1$ and $W_2$ are finite dimensional subspaces of a vector space V then $W_1 + W_2$ is finite dimensional and dim $W_1$ + dim $W_2$ = dim ($W_1 \cap W_2$) + dim ($W_1 + W_2$). We say $\alpha_1, \ldots, \alpha_t$ are linearly dependent if there exists scalars $c_1, c_2, \ldots, c_t$ not all zero such that $c_1 \alpha_1 + \ldots + c_t \alpha_t = 0$.

***Example 1.1.10:*** Let $V = M_{2 \times 2} = \{(a_{ij}) \mid a_{ij} \in Q\}$ be a vector space over Q. A basis of V is

$$\left\{\begin{pmatrix} 0 & 1 \\ 0 & 0 \end{pmatrix}, \begin{pmatrix} 0 & 0 \\ 1 & 0 \end{pmatrix}, \begin{pmatrix} 1 & 0 \\ 0 & 0 \end{pmatrix}, \begin{pmatrix} 0 & 0 \\ 0 & 1 \end{pmatrix}\right\}.$$

***Example 1.1.11:*** Let V = R × R × R be a vector space over R. Then {(1, 0, 0), (0, 1, 0), (0, 0, 1)} is a basis of V.

If V = R × R × R is a vector space over Q, V is not finite dimensional.

***Example 1.1.12:*** Let V = R[x] be a vector space over R. V = R[x] is an infinite dimensional vector spaces. A basis of V is $\{1, x, x^2, \ldots, x^n, \ldots\}$.

***Example 1.1.13:*** Let $P_{3 \times 2} = \{(a_{ij}) \mid a_{ij} \in R\}$ be a vector space over R. A basis for $P_{3 \times 2}$ is

$$\left\{ \begin{pmatrix} 1 & 0 \\ 0 & 0 \\ 0 & 0 \end{pmatrix}, \begin{pmatrix} 0 & 1 \\ 0 & 0 \\ 0 & 0 \end{pmatrix}, \begin{pmatrix} 0 & 0 \\ 1 & 0 \\ 0 & 0 \end{pmatrix}, \begin{pmatrix} 0 & 0 \\ 0 & 1 \\ 0 & 0 \end{pmatrix}, \begin{pmatrix} 0 & 0 \\ 0 & 0 \\ 1 & 0 \end{pmatrix}, \begin{pmatrix} 0 & 0 \\ 0 & 0 \\ 0 & 1 \end{pmatrix} \right\}.$$

Now we just proceed on to recall the definition of linear algebra.

**DEFINITION 1.1.5:** *Let F be a field. A linear algebra over the field F is a vector space A over F with an additional operation called multiplication of vectors which associates with each pair of vectors $\alpha$, $\beta$ in A a vector $\alpha\beta$ in A called the product of $\alpha$ and $\beta$ in such a way that*

i. *multiplication is associative $\alpha(\beta\gamma) = (\alpha\beta)\gamma$.*
ii. *multiplication is distributive with respect to addition*

$$\alpha(\beta + \gamma) = \alpha\beta + \alpha\gamma$$
$$(\alpha + \beta)\gamma = \alpha\gamma + \beta\gamma.$$

iii. *for each scalar c in F, $c(\alpha\beta) = (c\alpha)\beta = \alpha(c\beta)$.*

*If there is an element 1 in A such that $1\alpha = \alpha 1 = \alpha$ for each $\alpha$ in A we call $\alpha$ a linear algebra with identity over F and call 1 the identity of A. The algebra A is called commutative if $\alpha\beta = \beta\alpha$ for all $\alpha$ and $\beta$ in A.*

***Example 1.1.14:*** F[x] be a polynomial ring with coefficients from F. F[x] is a commutative linear algebra over F.

***Example 1.1.15:*** Let $M_{5 \times 5} = \{(a_{ij}) \mid a_{ij} \in Q\}$; $M_{5 \times 5}$ is a linear algebra over Q which is not a commutative linear algebra.

All vector spaces are not linear algebras for we have got the following example.

***Example 1.1.16:*** Let $P_{5 \times 7} = \{(a_{ij}) \mid a_{ij} \in R\}$; $P_{5 \times 7}$ is a vector space over R but $P_{5 \times 7}$ is not a linear algebra.

It is worthwhile to mention that by the very definition of linear algebra all linear algebras are vector spaces and not conversely.

## 1.2 Basic Properties of Linear Bialgebra

In this section we for the first time introduce the notion of linear bialgebra, prove several interesting results and illustrate them also with example.

**DEFINITION 1.2.1:** *Let $V = V_1 \cup V_2$ be a bigroup. If $V_1$ and $V_2$ are linear algebras over the same field F then we say V is a linear bialgebra over the field F.*

*If both $V_1$ and $V_2$ are of infinite dimension vector spaces over F then we say V is an infinite dimensional linear bialgebra over F. Even if one of $V_1$ or $V_2$ is infinite dimension then we say V is an infinite dimensional linear bialgebra. If both $V_1$ and $V_2$ are finite dimensional linear algebra over F then we say $V = V_1 \cup V_2$ is a finite dimensional linear bialgebra.*

***Examples 1.2.1:*** Let $V = V_1 \cup V_2$ where $V_1$ = {set of all $n \times n$ matrices with entries from Q} and $V_2$ be the polynomial ring Q [x]. $V = V_1 \cup V_2$ is a linear bialgebra over Q and the linear bialgebra is an infinite dimensional linear bialgebra.

***Example 1.2.2:*** Let $V = V_1 \cup V_2$ where $V_1 = Q \times Q \times Q$ abelian group under '+', $V_2$ = {set of all $3 \times 3$ matrices with entries from Q} then $V = V_1 \cup V_2$ is a bigroup. Clearly V is a linear bialgebra over Q. Further dimension of V is 12 V is a 12 dimensional linear bialgebra over Q.

The standard basis is {(0 1 0), (1 0 0), (0 0 1)} $\cup$

$$\left\{\begin{array}{l}\begin{pmatrix}1&0&0\\0&0&0\\0&0&0\end{pmatrix},\begin{pmatrix}0&1&0\\0&0&0\\0&0&0\end{pmatrix},\begin{pmatrix}0&0&1\\0&0&0\\0&0&0\end{pmatrix},\begin{pmatrix}0&0&0\\1&0&0\\0&0&0\end{pmatrix},\begin{pmatrix}0&0&0\\0&1&0\\0&0&0\end{pmatrix},\\ \begin{pmatrix}0&0&0\\0&0&1\\0&0&0\end{pmatrix},\begin{pmatrix}0&0&0\\0&0&0\\1&0&0\end{pmatrix},\begin{pmatrix}0&0&0\\0&0&0\\0&1&0\end{pmatrix},\begin{pmatrix}0&0&0\\0&0&0\\0&0&1\end{pmatrix}\end{array}\right\}$$

***Example 1.2.3:*** Let $V = V_1 \cup V_2$ where $V_1$ is a collection of all $8 \times 8$ matrices over Q and $V_2$ = {the collection of all $3 \times 2$ matrices over Q}. Clearly V is a bivector space of dimension 70 and is not a linear bialgebra.

From this example it is evident that there exists bivector spaces which are not linear bialgebras.

Now if in a bivector space $V = V_1 \cup V_2$ one of $V_1$ or $V_2$ is a linear algebra then we call V as a semi linear bialgebra. So the vector space given in example 1.2.3 is a semi linear bialgebra.

We have the following interesting theorem.

**THEOREM 1.2.1:** *Every linear bialgebra is a semi linear bialgebra. But a semi linear bialgebra in general need not be a linear bialgebra.*

*Proof:* The fact that every linear bialgebra is a semi linear bialgebra is clear from the definition of linear bialgebra and semi linear bialgebra.

To prove that a semi linear bialgebra need not in general be a linear bialgebra. We consider an example. Let $V = V_1 \cup V_2$ where $V_1 = Q \times Q$ and $V_2$ the collection of all $3 \times 2$ matrices with entries from Q clearly $V = Q \times Q$ is a linear algebra of dimension 2 over Q and $V_2$ is not a linear algebra but only a vector space of dimension 6, as in $V_2$ we cannot define matrix multiplication. Thus $V = V_1 \cup V_2$ is not a linear bialgebra but only a semi linear bialgebra.

Now we have another interesting result.

**THEOREM 1.2.2:** *Every semi linear bialgebra over Q is a bivector space over Q but a bivector space in general is not a semi linear bialgebra.*

*Proof:* Every semi linear bialgebra over Q is clearly by the very definition a bivector space over Q. But to show a bivector space over Q in general is not a semi linear bialgebra we give an example. Let $V = V_1 \cup V_2$ where $V_1$ = {the set of all $2 \times 5$ matrices with entries from Q} and $V_2$ = {all polynomials of degree less than or equal to 5 with entries from Q}. Clearly both $V_1$ and $V_2$ are only vector spaces over Q and none of them are linear algebra. Hence $V = V_1 \cup V_2$ is only a bivector space and not a semi linear bialgebra over Q.

Hence the claim.

Now we define some more types of linear bialgebra. Let V $= V_1 \cup V_2$ be a bigroup. Suppose $V_1$ is a vector space over Q $\left(\sqrt{2}\right)$ and $V_2$ is a vector space over Q $\left(\sqrt{3}\right)$ (Q $\left(\sqrt{2}\right)$ and Q $\left(\sqrt{3}\right)$ are fields).

Then V is said to be a strong bivector space over the bifield Q $\left(\sqrt{3}\right) \cup$ Q $\left(\sqrt{2}\right)$. Similarly if $V = V_1 \cup V_2$ be a bigroup and if V is a linear algebra over F and $V_2$ is a linear algebra over K, $K \neq F$ $K \cap F \neq F$ or K. i.e. if $K \cup F$ is a bifield then we say V is a strong linear bialgebra over the bifield.

Thus now we systematically give the definitions of strong bivector space and strong linear bialgebra.

**DEFINITION 1.2.2:** *Let $V = V_1 \cup V_2$ be a bigroup. $F = F_1 \cup F_2$ be a bifield. If $V_1$ is a vector space over $F_1$ and $V_2$ is a vector space over $F_2$ then $V = V_1 \cup V_2$ is called the strong bivector space over the bifield $F = F_1 \cup F_2$. If $V = V_1 \cup V_2$ is a bigroup and if $F = F_1 \cup F_2$ is a bifield. If $V_1$ is a linear algebra over $F_1$* *and $V_2$ is a linear algebra over $F_2$. Then we say $V = V_1 \cup V_2$ is a strong linear bialgebra over the field $F = F_1 \cup F_2$.*

***Example 1.2.4:*** Consider the bigroup $V = V_1 \cup V_2$ where $V_1 = Q(\sqrt{2}) \times Q(\sqrt{2})$ and $V_2$ = {set of all 3 × 3 matrices with entries from $Q(\sqrt{3})$}.

Let $F = Q(\sqrt{2}) \cup Q(\sqrt{3})$, Clearly F is a bifield. $V_1$ is a linear algebra over $Q(\sqrt{2})$ and $V_2$ is a linear algebra over $Q(\sqrt{3})$. So $V = V_1 \cup V_2$ is a strong linear bialgebra over the bifield $F = Q(\sqrt{2}) \cup Q(\sqrt{3})$.

It is interesting to note that we do have the notion of weak linear bialgebra and weak bivector space.

***Example 1.2.5:*** Let $V = V_1 \cup V_2$ where $V_1 = Q \times Q \times Q \times Q$ be the group under '+' and $V_2$ = The set of all polynomials over the field $Q(\sqrt{2})$. Now $V_1$ is a linear algebra over Q and $V_2$ is a linear algebra over $Q(\sqrt{2})$. Clearly $Q \cup Q(\sqrt{2})$ is a not a bifield as $Q \subseteq Q(\sqrt{2})$. Thus $V = V_1 \cup V_2$ is a linear bialgebra over Q we call $V = V_1 \cup V_2$ to be a weak linear bialgebra over $Q \cup Q(\sqrt{2})$. For $Q \times Q \times Q \times Q = V_1$ is not a linear algebra over $Q(\sqrt{2})$. It is a linear algebra only over Q.

Based on this now we give the definition of weak linear bialgebra over $F = F_1 \cup F_2$. where $F_1$ and $F_2$ are fields and $F_1 \cup F_2$ is not a bifield.

**DEFINITION 1.2.3:** *Let $V = V_1 \cup V_2$ be a bigroup. Let $F = F_1 \cup F_2$. Clearly F is not a bifield (For $F_1 \underset{\neq}{\subset} F_2$ or $F_2 \underset{\neq}{\subset} F_1$) but $F_1$* *and $F_2$ are fields. If $V_1$ is a linear algebra over $F_1$ and $V_2$ is a linear algebra over $F_2$, then we call $V_1 \cup V_2$ a weak linear bialgebra over $F_1 \cup F_2$. One of $V_1$ or $V_2$ is not a linear algebra over $F_2$ or $F_1$ respectively.*

On similar lines we can define weak bivector space.

***Example 1.2.6:*** Let $V = V_1 \cup V_2$ be a bigroup. Let $F = Q \cup Q\left(\sqrt{2},\sqrt{3}\right)$ be a union of two fields, for $Q \subset Q\left(\sqrt{2},\sqrt{3}\right)$ so $Q \cup Q\left(\sqrt{2},\sqrt{3}\right) = Q\left(\sqrt{2},\sqrt{3}\right)$ so is a field. This cannot be always claimed, for instance if $F = Q\left(\sqrt{2}\right) \cup Q\left(\sqrt{3}\right)$ is not a field only a bifield. Let $V_1 = Q \times Q \times Q$ and $V_2 = Q\left(\sqrt{2},\sqrt{3}\right)[x]$. $V_1$ is a vector space, in fact a linear algebra over Q but $V_1$ is not a vector space over the field $Q\left(\sqrt{2},\sqrt{3}\right)$. $V_2$ is a vector space or linear algebra over Q or $Q\left(\sqrt{2},\sqrt{3}\right)$. Since $Q \subset Q\left(\sqrt{2},\sqrt{3}\right)$, we see V is a weak bivector space over $Q \cup Q\left(\sqrt{2},\sqrt{3}\right)$ infact V is a weak linear bialgebra over $Q \cup Q\left(\sqrt{2},\sqrt{3}\right)$.

***Example 1.2.7:*** Let $V = V_1 \cup V_2$ be a bigroup. Let $F = Q\left(\sqrt{2},\sqrt{3}\right) \cup Q\left(\sqrt{5},\sqrt{7},\sqrt{11}\right)$ be a bifield. Suppose $V_1 = Q\left(\sqrt{2},\sqrt{3}\right)[x]$ be a linear algebra over $Q\left(\sqrt{2},\sqrt{3}\right)$ and $V_2 = Q\left(\sqrt{5},\sqrt{7},\sqrt{11}\right) \times Q\left(\sqrt{5},\sqrt{7},\sqrt{11}\right)$ be a linear algebra over $Q\left(\sqrt{5},\sqrt{7},\sqrt{11}\right)$. Clearly $V = V_1 \cup V_2$ is a strong linear bialgebra over the bifield $Q\left(\sqrt{2},\sqrt{3}\right) \cup Q\left(\sqrt{5},\sqrt{7},\sqrt{11}\right)$.

Now we have the following results.

**THEOREM 1.2.3:** *Let $V = V_1 \cup V_2$ be a bigroup and V be a strong linear bialgebra over the bifield $F = F_1 \cup F_2$. V is not a linear bialgebra over F.*

*Proof:* Now we analyze the definition of strong linear bialgebra and the linear bialgebra. Clearly the strong linear bialgebra has no relation with the linear bialgebra or a linear bialgebra has no relation with strong linear bialgebra for linear bialgebra is defined over a field where as the strong linear bialgebra is defined over a bifield, hence no relation can ever be derived. In the similar means one cannot derive any form of relation between the weak linear bialgebra and linear bialgebra.

All the three notions, weak linear bialgebra, linear bialgebra and strong linear bialgebra for a weak linear bialgebra is defined over union of fields $F = F_1 \cup F_2$ where $F_1 \subset F_2$ or $F_2 \subset F_1$; $F_1$ and $F_2$ are fields; linear bialgebras are defined over the same field where as the strong linear bialgebras are defined over bifields. Thus these three concepts are not fully related.

It is important to mention here that analogous to weak linear bialgebra we can define weak bivector space and analogous to strong linear bialgebra we have the notion of strong bivector spaces.

***Example 1.2.8:*** Let $V = V_1 \cup V_2$ where $V_1$ = {set of all linear transformation of a n dimensional vector space over Q to a m dimensional vector space W over Q} and $V_2$ = {All polynomials of degree ≤ 6 with coefficients from R}.

Clearly $V = V_1 \cup V_2$ is a bigroup. V is a weak bivector space over $Q \cup R$.

***Example 1.2.9:*** Let $V = V_1 \cup V_2$ be a bigroup. $V_1$ = {set of all polynomials of degree less than or equal to 7 over $Q\left(\sqrt{2}\right)$} and $V_2$ = {set of all 5 × 2 matrices with entries from $Q\left(\sqrt{3}, \sqrt{7}\right)$}. $V = V_1 \cup V_2$ is a strong bivector space over the bifield $F = Q\left(\sqrt{2}\right)$ $\cup$ $Q\left(\sqrt{3},\sqrt{7}\right)$. Clearly V = $V_1 \cup V_2$ is not a strong linear bialgebra over F.

Now we proceed on to define linear subbialgebra and subbivector space.

**DEFINITION 1.2.4:** *Let V = $V_1 \cup V_2$ be a bigroup. Suppose V is a linear bialgebra over F. A non empty proper subset W of V is said to be a linear subbialgebra of V over F if*

1. *W = $W_1 \cup W_2$ is a subbigroup of V = $V_1 \cup V_2$.*
2. *$W_1$ is a linear subalgebra over F.*
3. *$W_2$ is a linear subalgebra over F.*

***Example 1.2.10:*** Let V = $V_1 \cup V_2$ where $V_1 = Q\times Q\times Q\times Q$ and $V_2$ ={set of all 4 × 4 matrices with entries from Q}. V = $V_1 \cup V_2$ is a bigroup under '+' V is a linear bialgebra over Q.

Now consider W = $W_1 \cup W_2$ where $W_1 = Q \times \{0\} \times Q \times \{0\}$ and $W_2$ = {collection of all upper triangular matrices with entries from Q}. W = $W_1 \cup W_2$ is a subbigroup of V = $V_1 \cup V_2$. Clearly W is a linear subbialgebra of V over Q.

**DEFINITION 1.2.5:** *Let V = $V_1 \cup V_2$ be a bigroup. F = $F_1 \cup F_2$ be a bifield. Let V be a strong linear bialgebra over F. A non empty subset W = $W_1 \cup W_2$ is said to be a strong linear subbialgebra of V over F if*

1. *W = $W_1 \cup W_2$ is a subbigroup of V = $V_1 \cup V_2$.*
2. *$W_1$ is a linear algebra over $F_1$ and*
3. *$W_2$ is a linear algebra over $F_2$.*

***Example 1.2.11:*** Let V = $V_1 \cup V_2$ where $V_1$ = {Set of all 2 × 2 matrices with entries from $Q\left(\sqrt{2}\right)$} be a group under matrix addition and $V_2$ = {collection of all polynomials in the variable x over $Q\left(\sqrt{5},\sqrt{3}\right)$}. $V_2$ under polynomial addition is a group.

Thus $V = V_1 \cup V_2$ is a strong linear bialgebra over the bifield F $= Q\left(\sqrt{2}\right) \cup Q\left(\sqrt{3}, \sqrt{5}\right)$.

Consider $W = W_1 \cup W_2$ where

$$W_1 = \left\{ \begin{pmatrix} 0 & 0 \\ 0 & a \end{pmatrix} \middle| a \in Q\left(\sqrt{2}\right) \right\}$$

and $W_2$ = {all polynomials of even degree i.e. p (x) = $\sum_{i=1}^{n} p_i x^{2i} n \in Q(\sqrt{3}, \sqrt{5})$. Clearly, $W = W_1 \cup W_2$ is a subbigroup of $V = V_1 \cup V_2$. Clearly W is a strong linear subbialgebra of V over the bifield $F = Q\left(\sqrt{2}\right) \cup Q\left(\sqrt{5}, \sqrt{3}\right)$.

**DEFINITION 1.2.6:** *Let $V = V_1 \cup V_2$ be a bigroup. Let V be a weak linear bialgebra over the union of fields $F_1 \cup F_2$ (i.e. $F_1$ is a subfield of $F_2$ or $F_2$ is a subfield of $F_1$) where $V_1$ is a linear algebra over $F_1$ and $V_2$ is a linear algebra over $F_2$. A non empty subset $W = W_1 \cup W_2$ is a weak linear sub bialgebra of V over $F_1 \cup F_2$ if*

1. *W is a subbigroup of V.*
2. *$W_1$ is a linear subalgebra of $V_1$.*
3. *$W_2$ is a linear subalgebra of $V_2$.*

***Example 1.2.12:*** Let $V = V_1 \cup V_2$ be a bigroup. Let $F = Q \cup Q\left(\sqrt{3}, \sqrt{7}\right)$ be the union of fields. $V_1$ be a linear algebra over Q and $V_2$ be a linear algebra over $Q\left(\sqrt{3}, \sqrt{7}\right)$; where $V_1$ = {all polynomials in the variable x with coefficient from Q} and $V_2$ = {3 × 3 matrices with entries from $Q\left(\sqrt{2}, \sqrt{3}\right)$}. $V = V_1 \cup V_2$ is a weak linear bialgebra over $Q \cup Q\left(\sqrt{3}, \sqrt{7}\right)$. Consider $W = W_1 \cup W_2$, $W_1$ = {all polynomial of only even degree in the variable x with coefficients from Q} and

$$W_2 = \left\{ \begin{pmatrix} a & b & 0 \\ e & d & 0 \\ 0 & 0 & 0 \end{pmatrix} \middle| a, b, c, d \in Q\left(\sqrt{3}, \sqrt{7}\right) \right\}.$$

Clearly $W = W_1 \cup W_2$ is a subbigroup of $V = V_1 \cup V_2$. Further $W = W_1 \cup W_2$ is a weak linear subbialgebra of V over $F = Q \cup Q\left(\sqrt{2}, \sqrt{3}\right)$.

It is very interesting to note that at times a weak linear bialgebra can have a subbigroup which happens to be just a linear subbialgebra over $F_1$ (if $F_1 \subset F_2$ i.e., $F_1$ is a subfield of $F_2$) and not a weak linear subbialgebra.

We define such structures as special weak linear bialgebra.

**DEFINITION 1.2.7:** *Let $V = V_1 \cup V_2$ be a bigroup. $F = F_1 \cup F_2$ (where either $F_1$ is a subfield of $F_2$ or $F_2$ is a subfield of $F_1$) be union of fields and V is a weak linear bialgebra over F. Suppose $W = W_1 \cup W_2$ is a subbigroup of $V = V_1 \cup V_2$ and if W is a linear bialgebra only over $F_1$ (or $F_2$) (which ever is the subfield of the other). Then we call $V = V_1 \cup V_2$ a special weak linear bialgebra.*

***Example 1.2.13:*** Let $V = V_1 \cup V_2$ be a bigroup where $V_1 = Q \times Q \times Q$ and $V_2 = \{2 \times 2$ matrices with entries from $Q\left(\sqrt{2}\right)\}$. $V = V_1 \cup V_2$ is a weak linear bialgebra over $Q \cup Q\left(\sqrt{2}\right)$.

Consider $W = W_1 \cup W_2$ where $W_1 = Q \times \{0\} \times Q$ and $W_2 = \{$set of all $2 \times 2$ matrices over $Q\}$. Then $W = W_1 \cup W_2$ is a linear bialgebra over Q. Clearly $W_2$ is note linear algebra over $Q\left(\sqrt{2}\right)$. We call $V = V_1 \cup V_2$ a special weak linear algebra over $Q \cup Q\left(\sqrt{2}\right)$.

**DEFINITION 1.2.8:** *Let $V = V_1 \cup V_2$ be a bigroup, $F = F_1 \cup F_2$ be a union fields (i.e., $F_1 = F$ or $F_2 = F$), V be the weak linear bialgebra over F. We say V is a strong special weak linear bialgebra if V has a proper subset W such that $W = W_1 \cup W_2$ is a subbigroup of V and there exists subfields of $F_1$ and $F_2$ say $K_1$ and $K_2$ respectively such that $K = K_1 \cup K_2$ is a bifield and W is a strong linear bialgebra over K.*

*Note:* Every weak linear bialgebra need not in general be a strong special weak linear bialgebra. We illustrate both by examples.

***Example 1.2.14:*** Let $V = V_1 \cup V_2$ where $V_1 = Q \times Q$ and $V_2$ = {Collection of all $5 \times 5$ matrices with entries from $Q\left(\sqrt{2}\right)$}.V is a weak linear bialgebra over $F = Q \cup Q\left(\sqrt{2}\right)$. Clearly V is never a strong special weak linear bialgebra.

Thus all weak linear bialgebras need not in general be a strong special weak linear bialgebras.

Thus all weak linear bialgebras need not in general be a strong special weak linear bialgebras. We will prove a theorem to this effect.

***Example 1.2.15:*** Let $V = V_1 \cup V_2$ be a bigroup. $V_1 = Q\left(\sqrt{2},\sqrt{3}\right)[x]$ be the polynomial ring over $Q\left(\sqrt{2},\sqrt{3}\right) = F_1$. Clearly $V_1$ is a abelian group under addition. Let $Q\left(\sqrt{2},\sqrt{3},\sqrt{7},\sqrt{11}\right) = F_2$ be the field. Let $V_2$ = {the set of all $3 \times 3$ matrices with entries from $F_2$}. $V_2$ is an abelian group under matrix addition. $V = V_1 \cup V_2$ is a weak linear bialgebra over $F_1 \cup F_2$.

Now consider the subbigroup, $W = W_1 \cup W_2$ where $W_1$ = {Set of all polynomials in x with coefficients from the field $Q\left(\sqrt{2}\right)$} and $W_2$ = {set of all $3 \times 3$ matrices with entries from $Q\left(\sqrt{3},\sqrt{7}\right)$}. $F = Q\left(\sqrt{2}\right) \cup Q\left(\sqrt{3},\sqrt{7}\right)$ is a bifield and the subbigroup $W = W_1 \cup W_2$ is a strong linear bialgebra over the bifield $F = Q(\sqrt{2}) \cup Q(\sqrt{3}, \sqrt{7})$. Thus V has a subset W such that W is a strong linear bialgebra so V is a strong special weak linear bialgebra.

Now we give condition for a weak linear bialgebra to be strong special weak linear bialgebra.

**THEOREM 1.2.4:** *Let $V = V_1 \cup V_2$ be a bigroup which has a proper subbigroup. Suppose V is a weak linear bialgebra over $F = F_1 \cup F_2$. V is not a strong special weak linear bialgebra if and only if one of $F_1$ or $F_2$ is a prime field.*

*Proof:* Let $V = V_1 \cup V_2$ be a bigroup and V be a weak linear bialgebra over the field $F = F_1 \cup F_2$. Suppose we assume $F_1$ is a prime field then clearly $F_2$ has $F_1$ to be its subfield. So $F = F_1 \cup F_2$ has no subset say $K = K_1 \cup K_2$ such that K is a bifield. Since F has no subset $K = K_1 \cup K_2$ such that K is a bifield, we see for no subbigroup W can be a strong linear bialgebra.

Conversely if $F = F_1 \cup F_2$ and if one of $F_1$ or $F_2$ is a prime field then there does not exist $K \subset F$, $K = K_1 \cup K_2$ such that K is a bifield, then V cannot have a subbigroup which is a strong linear bialgebra.

Now we give the conditions for the weak linear bialgebra to be a strong special weak linear bialgebra.

1. $V = V_1 \cup V_2$ should have proper subbigroups.
2. In $F = F_1 \cup F_2$ ($F = F_1$ or $F = F_2$) the subfield must not be a prime field and the extension field must contain some other non prime subfield other than the already mentioned subfield i.e. V should have subset which is a bifield.
3. V should have a subbigroup which is a strong linear bialgebra over the bifield.

**THEOREM 1.2.5:** *Let $V = V_1 \cup V_2$ be a bigroup. $F = F_1 \cup F_2$ be the union of field. V be a weak linear bialgebra over F. If V has no proper subbigroup then*

1. *V is not a strong special weak linear bialgebra over F.*
2. *V is not a special weak linear bialgebra.*

The proof of the above theorem is left as an exercise for the reader.

## 1.3 Fuzzy Vector Spaces

In this section we recall the definition and properties of fuzzy vector spaces. The study of fuzzy vector spaces started as early as 1977. For more refer [55].

Throughout this section V denotes a vector space over a field F. A fuzzy subset of a non-empty set S is a function from S into [0, 1]. Let A denote a fuzzy subspace of V over a fuzzy subfield K of F and let X denote a fuzzy subset of V such that $X \subseteq A$. Let $\langle X \rangle$ denote the intersection of all fuzzy subspaces of V over K that contain X and are contained in A.

**DEFINITION 1.3.1:**

1. *A fuzzy subset K of F is a fuzzy subfield of F, if K (1) = 1 and for all c, d $\in$ F, $K(c - d) \geq \min\{K(c), K(d)\}$ and $K(cd^{-1}) \geq \min\{K(c), K(d)\}$ where $d \neq 0$.*

2. *A fuzzy subset A of V is a fuzzy subspace over a fuzzy subfield K of F, if $A(0) > 0$ and for all x, y $\in$ V and for all c $\in$ F, $A(x - y) \geq \min\{A(x), A(y)\}$ and $A(cx) \geq \min\{K(c), A(x)\}$. If K is a fuzzy subfield of F and if $x \in F$, $x \neq 0$, then $K(0) = K(1) \geq K(x) = K(-x) = K(+ x^{-1})$.*

In the following we let L denote the set of all fuzzy subfields of F and let $A_L$ denote the set of all fuzzy subspaces of V over K $\in$ L.

If A and B are fuzzy subsets of V then $A \subset B$ means $A(x) \leq B(x)$ for all $x \in A$. For $0 \leq t \leq 1$, let $A_t = \{x \in V \mid A(x) \geq t\}$.

Now we proceed onto recall the concept of fuzzy spanning.

**DEFINITION 1.3.2:** *Let $A_1, A_2, \ldots, A_n$ be fuzzy subsets of V and let K be a fuzzy subset of F.*

1. *Define the fuzzy subset $A_1 + \ldots + A_n$ on V by the following: for all $x \in V$, $(A_1 + \ldots + A_n)(x) = \sup\{\min\{A_1(x_1), \ldots, A_n(x_n)\} / x = x_1 + \ldots + x_n, x_i \in V\}$.*
2. *Define the fuzzy subset K o A of V by for all $x \in V$, $(K \circ A)(x) = \sup\{\min\{K(c), A(y)\} \mid c \in F, y \in V, x = cy\}$.*

**DEFINITION 1.3.3:** *Let S be a set $x \in S$ and $0 \leq \lambda \leq 1$. Define the fuzzy subset $x_\lambda$ of S by $x_\lambda(y) = \lambda$ if $y = x$ and $x_\lambda(y) = 0$ if $y \neq x$. $x_\lambda$ is called a fuzzy singleton.*

**DEFINITION 1.3.4:** *A fuzzy vector space $(V, \eta)$ or $\eta_V$ is an ordinary vector space V with a map $\eta : V \to [0, 1]$ satisfying the following conditions.*

1. *$\eta(a + b) \geq \min\{\eta(a), \eta(b)\}$.*
2. *$\eta(-a) = \eta(a)$.*
3. *$\eta(0) = 1$.*
4. *$\eta(ra) \geq \eta(a)$ for all $a, b \in V$ and $r \in F$ where F is a field.*

**DEFINITION 1.3.5:** *For an arbitrary fuzzy vector space $\eta_V$ and its vector subspace $\eta_W$, the fuzzy vector space $(V/W, \hat{\eta})$ or $\eta_{VW}$ determined by*

$$\hat{\eta}(\upsilon + W) = \begin{cases} 1 & \text{if } \upsilon \in W \\ \sup\limits_{\omega \in W} \eta(\upsilon + \omega) & \text{otherwise} \end{cases}$$

*is called the fuzzy quotient vector space, $\eta_V$ by $\eta_W$.*

**DEFINITION 1.3.6:** *For an arbitrary fuzzy vector space $\eta_V$ and its fuzzy vector subspace $\eta_W$, the fuzzy quotient space of $\eta_V$ by $\eta_W$ is determined by*

$$\bar{\eta}(\nu + W) = \begin{cases} 1 & \upsilon \in W \\ \inf\limits_{\omega \in W} \eta(\upsilon + \omega) & \upsilon \notin W \end{cases}$$

*It is denoted by* $\bar{\eta}_{v/w}$.

## Chapter Two

# SET VECTOR SPACES

In this chapter for the first time we introduce the new concept of set vector spaces and set linear algebra. This new notion would help even non-mathematicians to use the tools of set vector spaces and set linear algebra over sets in fields like computer science, cryptology and sociology. Since the abstract algebraic concepts of abelian groups, fields and vector spaces are not easily understood by non-mathematicians more so by technical experts, engineers and social scientists we in this book try to introduce the notion of set vector spaces and set linear algebra over sets.

This chapter has five sections. In section one we introduce the new notion of set vector space over a set and illustrate it with examples. In section two set linear transformation of set vector spaces is introduced. Section three introduces the notion of set linear algebra.

In the fourth section the new notion of semigroup vector spaces are defined and a few of its important properties analyzed. The final section for the first time introduces the notion of group vector spaces.

## 2.1 Set Vector Spaces and their Properties

In this section we introduce for the first time the notion of set vector spaces over sets and illustrate them with examples and give some interesting properties about them.

**DEFINITION 2.1.1:** *Let S be a set. V another set. We say V is a set vector space over the set S if for all v $\in$ V and for all s $\in$ S; vs and sv $\in$ V.*

***Example 2.1.1:*** Let V = {1, 2, ..., ∞} be the set of positive integers. S = {2, 4, 6, ..., ∞} the set of positive even integers. V is a set vector space over S. This is clear for sv = vs $\in$ V for all s $\in$ S and v $\in$ V.

It is interesting to note that any two sets in general may not be a set vector space over the other. Further even if V is a set vector space over S then S in general need not be a set vector space over V.

For from the above example 2.1.1 we see V is a set vector space over S but S is also a set vector space over V for we see for every s $\in$ S and v $\in$ V, v.s = s.v $\in$ S. Hence the above example is both important and interesting as one set V is a set vector space another set S and vice versa also hold good inspite of the fact S $\neq$ V.

Now we illustrate the situation when the set V is a set vector space over the set S. We see V is a set vector space over the set S and S is not a set vector space over V.

***Example 2.1.2:*** Let V = {$Q^+$ the set of all positive rationals} and S = {2, 4, 6, 8, …, ∞}, the set of all even integers. It is easily verified that V is a set vector space over S but S is not a set vector space over V, for $\frac{7}{3} \in$ V and 2 $\in$ S but $\frac{7}{3}.2 \notin$ S. Hence the claim.

Now we give some more examples so that the reader becomes familiar with these concepts.

***Example 2.1.3:*** Let $\left\{ M_{2\times 2} = \begin{pmatrix} a & b \\ c & d \end{pmatrix} \right|$ a, b, c, d $\in$ {set of all positive integers together with zero} be the set of $2 \times 2$ matrices with entries from N. Take S = {0, 2, 4, 6, 8, …, ∞}. Clearly $M_{2x2}$ is a set vector space over the set S.

***Example 2.1.4:*** Let V = $\{Z^+ \times Z^+ \times Z^+\}$ such that $Z^+$ = {set of positive integers}. S = {2, 4, 6, …, ∞}. Clearly V is a set vector space over S.

***Example 2.1.5:*** Let V = $\{Z^+ \times Z^+ \times Z^+\}$ such that $Z^+$ is the set of positive integers. S = {3, 6, 9, 12, …, ∞}. V is a set vector space over S.

***Example 2.1.6:*** Let $Z^+$ be the set of positive integers. $pZ^+ = S$, p any prime. $Z^+$ is a set vector space over S.

**Note:** Even if p is replaced by any positive integer n still $Z^+$ is a set vector space over $nZ^+$. Further $nZ^+$ is also a set vector space over $Z^+$. This is a collection of special set vector spaces over the set.

***Example 2.1.7:*** Let Q[x] be the set of all polynomials with coefficients from Q, the field of rationals. Let S = {0, 2, 4, …, ∞}. Q[x] is a set vector space over S. Further S is not a set vector space over Q[x].

Thus we see all set vector spaces V over the set S need be such that S is a set vector space over V.

***Example 2.1.8:*** Let R be the set of reals. R is a set vector space over the set S where S = {0, 1, 2, …, ∞}. Clearly S is not a set vector space over R.

***Example 2.1.9:*** $\mathbb{C}$ be the collection of all complex numbers. Let $Z^+ \cup \{0\} = \{0, 1, 2, \ldots, \infty\} = S$. $\mathbb{C}$ is a set vector space over S but S is not a set vector space over $\mathbb{C}$**.**

At this point we propose a problem.

Characterize all set vector spaces V over the set S which are such that S is a set vector space over the set V.

Clearly $Z^+ \cup \{0\} = V$ is a set vector space over $S = p\, Z^+ \cup \{0\}$, p any positive integer; need not necessarily be prime. Further S is also a set vector space over V.

***Example 2.1.10:*** Let $V = Z_{12} = \{0, 1, 2, \ldots, 11\}$ be the set of integers modulo 12. Take $S = \{0, 2, 4, 6, 8, 10\}$. V is a set vector space over S. For $s \in S$ and $v \in V$, $sv \equiv vs \pmod{12}$.

***Example 2.1.11:*** Let $V = Z_p = \{0, 1, 2, \ldots, p - 1\}$ be the set of integers modulo p (p a prime). $S = \{1, 2, \ldots, p - 1\}$ be the set. V is a set vector space over S.

In fact V in example 2.1.11 is a set vector space over any proper subset of S also.

This is not always true, yet if V is a set vector space over S. V need not in general be a set vector space over every proper subset of S. Infact in case of $V = Z_p$ every proper subset $S_i$ of S $= \{1, 2, \ldots, p - 1\}$ is such that V is a set vector space over $S_i$.

**Note:** It is important to note that we do not have the notion of linear combination of set vectors but only the concept of linear magnification or linear shrinking or linear annulling. i.e., if $v \in V$ and s is an element of the set S where V is the set vector space defined over it and if sv makes v larger we say it is linear magnification.

***Example 2.1.12:*** If $V = Z^+ \cup \{0\} = \{0, 1, 2, \ldots, \infty\}$ and $S = \{0, 2, 4, \ldots, \infty\}$ we say for s = 10 and v = 21, sv = 10.21 is a linear magnification.

Now for v = 5 and s = 0, sv = 0.5 = 0 is a linear annulling. Here we do not have the notion of linear shrinking.

So in case of any modeling where the researcher needs only linear magnification or annulling of the data he can use the set vector space V over the set S.

*Notation:* We call the elements of the set vector space V over the set S as set vectors and the scalars in S as set scalars.

***Example 2.1.13:*** Suppose $Q^+ \cup \{0\}$ = {all positive rationals} = $Q^*$ = V and

$$S = \left\{0, 1, \frac{1}{2}, \frac{1}{2^2}, \frac{1}{2^3}, \ldots\right\}.$$

Then V is a set vector space over S. We see for $6 \in Q^* = V$ and $\frac{1}{2} \in S$ $\frac{1}{2}.6 = 3$. This is an instance of linear shrinking. Suppose s = 0 and $v = \frac{21}{3}$ then s.v = 0. $\frac{21}{3} = 0$ is an instance of linear annulling.

Now we see this vector space is such that there is no method to get a linearly magnified element or the possibility of linear magnification.

Now we have got a peculiar instance for $1 \in S$ and 1.v = v for every $v \in V$ linear neutrality or linearly neutral element. The set may or may not contain the linearly neutral element. We illustrate yet by an example in which the set vector space which has linear magnifying, linear shrinking, linear annulling and linearly neutral elements.

***Example 2.1.14:*** Let $Q^+ \cup \{0\}$ = V = {set of all positive rationals with zero}.

$$S = \left\{0, 1, \frac{1}{3}, \frac{1}{3^2}, \ldots, 2, 4, 6, \ldots, \infty\right\}$$

be the set. V is a set vector space over the set S. We see V is such that there are elements in S which linearly magnify some elements in V; for instance if v = 32 and s = 10 then s.v = 10.32 = 320 is an instance of linear magnification. Consider v = 30 and s = $\frac{1}{3}$ then s.v. = $\frac{1}{3} \times 30 = 10 \in$ V. This is an instance of linear shrinking. Thus we have certain elements in V which are linearly shrunk by elements of the set S.

Now we have $0 \in S$ such that 0 v = 0 for all $v \in V$ which is an instance of linear annulling. Finally we see $1 \in S$ is such that 1.v = v for all $v \in V$; which is an instance of linearly neutral. Thus we see in this set vector space V over the set S given in example 2.1.14, all the four properties hold good. Thus if a researcher needs all the properties to hold in the model he can take them without any hesitation. Now we define yet another notion called linearly normalizing element of the set vector space V. Suppose $v \in V$ is a set vector and $s \in S$ is a set scalar and $1 \in V$ which is such that 1.v = v.1 = v and s.1 = 1.s = s for all $v \in V$ and for all $s \in S$; we call the element scalar s in S to be a linearly normalizing element of S if we have $v \in V$ such that s.v = v.s = 1.

For example in the example 2.1.14, we have for $3 \in V$ and $\frac{1}{3} \in S$; $\frac{1}{3} . 3 = 1 \in S$. Thus $\frac{1}{3}$ is a linearly normalizing element of S.

It is important to note that as in case of linearly annulling or linearly neutral element the scalar need not linearly normalize every element of the set vector space V. In most cases an element can linearly normalize only one element.

Having seen all these notions we now proceed on to define the new notion of set vector subspace of a set vector space V.

**DEFINITION 2.1.2:** *Let V be a set vector space over the set S. Let $W \subset V$ be a proper subset of V. If W is also a set vector space over S then we call W a set vector subspace of V over the set S.*

We illustrate this by a few examples so that the reader becomes familiar with this concept.

***Example 2.1.15:*** Let V = {0, 1, 2, …, ∞} be the set of positive integers together with zero. Let S = {0, 2, 4, 6, …, ∞}, the set of positive even integers with zero. V is a set vector space over S. Take W = {0, 3, 6, 9, 12, …, ∞} set of all multiples of 3 with zero. $W \subseteq V$ ; W is also a set vector space over S. Thus W is a set vector subspace of V over S.

***Example 2.1.16:*** Let Q[x] be the set of all polynomials with coefficients from Q; the set of rationals. Q[x] is a set vector space over the set S = {0, 2, 4, …, ∞}. Take W = {0, 1, …, ∞} the set of positive integers with zero. W is a set vector space over the set S. Now $W \subseteq Q \subseteq Q[x]$; so W is a set vector subspace of V over the set S.

***Example 2.1.17:*** Let $V = Z^+ \times Z^+ \times Z^+ \times Z^+$ a set of vector space over the set $Z^+$ = {0, 1, 2, …, ∞}. Let $W = Z^+ \times Z^+ \times \{0\} \times \{0\}$, a proper subset of V. W is also a set vector space over $Z^+$, i.e., W is a set vector subspace of V.

***Example 2.1.18:*** Let $V = 2Z^+ \times 3Z^+ \times Z^+$ be a set, V is a set vector space over the set S = {0, 2, 4, …, ∞}. Now take $W = 2Z^+ \times \{0\} \times 2Z^+ \subseteq V$; W is a set vector subspace of V over the set S.

***Example 2.1.19:*** Let

$$V = M_{3\times 2} = \left\{ \begin{pmatrix} a & d \\ b & e \\ c & f \end{pmatrix} \middle| \; a, b, c, d, e, f \in Z^+ \cup \{0\} \right\}$$

be the set of all 3 × 2 matrices with entries from the set of positive integers together with zero. V is a set vector space over the set $S = Z^+ \cup \{0\}$. Take

$$W = \left\{ \begin{pmatrix} a & 0 \\ b & 0 \\ c & 0 \end{pmatrix} \middle| \; a, b, c \in Z^+ \cup \{0\} \right\} \subseteq V;$$

W is a set vector subspace of V over the set S.

Now having defined set vector subspaces, we proceed on to define the notion of zero space of a set vector space V over the set S. We as in the case of usual vector spaces cannot define set zero vector space at all times. The set zero vector space of a set vector space V exists if and only if the set vector space V over S has {0} in V i.e., {0} is the linearly annulling element of S or $0 \in V$ and $0 \notin S$ in either of the two cases we have the set zero subspace of the set vector space V.

***Example 2.1.20:*** Let $Z^+ = V = \{1, 2, \ldots, \infty\}$ be a set vector space over the set $S = \{2, 4, 6, \ldots, \infty\}$; V is a set vector space but $0 \notin V$, so V does not have a set vector zero subspace.

It is interesting to mention here that we can always adjoin the zero element to the set vector space V over the set S and this does not destroy the existing structure. Thus the element {0} can always be added to make the set vector space V to contain a set zero subspace of V.

We leave it for the reader to prove the following theorem.

**THEOREM 2.1.1:** *Let V be a set vector space over the set S. Let $W_1, \ldots, W_n$ be n proper set vector subspaces of V over S. Then $\bigcap_{i=1}^{n} W_i$ is a set vector subspace of V over S. Further $\bigcap_{i=1}^{n} W_i = \phi$* *can also occur, if even for a pair of set vector subspaces $W_i$ and $W_j$ of V we have $W_i \cap W_j = \phi$ ($i \neq j$).*

*Clearly even if $0 \in V$ then also we cannot say $\bigcap_{i=1}^{n} W_i = \{0\}$ as 0 need not be present in every set vector subspace of V.*

We illustrate the situation by the following example.

***Example 2.1.21:*** Let $Z^+ = \{1, 2, \ldots, \infty\}$ be a set vector space over the set $S = \{2, 4, 6, \ldots, \infty\}$. Take

$$W_1 = \{2, \ldots, \infty\},$$
$$W_2 = \{3, 6, \ldots, \infty\}, \ldots,$$

and

$$W_p = \{p, 2p, 3p, \ldots, \infty\}.$$

We see $W_i \cap W_j \neq \phi$ for every i, j ($i \neq j$).

Will $\bigcap W_i = \phi$ if $i = 1$ to $\infty$ ?

Will $\bigcap_{i=1}^{n} W_i \neq \phi$ if $i = 1, 2, \ldots, n$; $n < \infty$?

***Example 2.1.22:*** Let $V = \{1, 2, \ldots, \infty\}$ be a set vector space over $S = \{2, , \ldots, \infty\}$. $W = \{2, 2^2, \ldots, 2^n \ldots\}$ is set vector subspace of V over the set S.

$W_1 = \{3, 3^2, 3^3, \ldots, \infty\}$ is a proper subset of V but $W_1$ is not a set vector subspace of V as $W_1$ is not a set vector space over S as $2.3 = 6 \notin W_1$ for $2 \in S$ and $3 \in W_1$. Thus we by this example show that in general every proper subset of the set V need not be a set vector space over the set S.

Now having seen the set vector subspaces of a set vector space we now proceed to define yet another new notion about set vector spaces.

**DEFINITION 2.1.3:** *Let V be a set vector space over the set S. Let T be a proper subset of S and W a proper subset of V. If W is a set vector space over T then we call W to be a subset vector subspace of V over T.*

We first illustrate this situation by examples before we proceed on to give more properties about them.

***Example 2.1.23:*** Let V = $Z^+ \cup \{0\} = \{0, 1, 2, \ldots, \infty\}$ be a set vector space over the set S = $\{2, 4, 6, \ldots, \infty\}$. Consider W = $\{3, 6, 9, \ldots\} \in$ V.

Take T = $\{6, 12, 18, 24, \ldots\} \subseteq$ S. Clearly W is a set vector space over T; so W is a subset vector subspace over T.

***Example 2.1.24:*** Let V = $Z^+ \times Z^+ \times Z^+$ be a set vector space over $Z^+$ = S. Take W = $Z^+ \times \{0\} \times \{0\}$ a proper subset of V and T = $\{2, , 6, \ldots, \infty\} \subseteq Z^+$ = S. Clearly W is a set vector space over T i.e., W is a subset vector subspace over T.

Now we show all proper subsets of a set vector space need not be a subset vector subspace over every proper subset of S. We illustrate this situation by the following examples.

***Example 2.1.25:*** Let V = $Z^+ \times Z^+ \times Z^+$ be a set vector space over $Z^+$ = S = $\{0, 1, 2, ..., \infty\}$. Let W = $\{3, 3^2, \ldots, \infty\} \times \{5, 5^2, \ldots\} \times \{0\} \subseteq$ V. Take T = $\{2, 4, 6, \ldots\} \subseteq$ S.

Clearly W is not a set vector space over T. That is W is not a subset vector subspace of V over T. Thus we see every subset of a set vector space need not in general be a subset vector subspace of the set vector space V over any proper subset T of the set S.

We illustrate this concept with some more examples.

***Example 2.1.26:*** Let V = $2Z^+ \times 3Z^+ \times 5Z^+ = \{(2n, 3m, 5t) \mid n, m, t \in Z^+\}$; V is a set vector space over the set S = $Z^+ = \{1, 2, \ldots, \infty\}$. Take W = $\{2, 2^2, \ldots, \infty\} \times \{3, 3^2, \ldots \} \times \{5, 5^2, \ldots\} \subseteq$

V. W is not a subset vector subspace over the subset T = {2, 4, …, ∞} ⊆ $Z^+$. Further W is not even a set vector subspace over $Z^+$ = S. Take W = {2, $2^2$, …} × {0} × {0} a proper subset of V. Choose T = {2, $2^3$, $2^5$, $2^7$}. W is a subset vector subspace of V over the subset T ⊆ S.

***Example 2.1.27:*** Let

$$V = \left\{ \begin{pmatrix} a & b \\ c & d \end{pmatrix} \middle| \; a, b, c, d \in Z^+ \cup \{0\} = \{0, 1, 2, \ldots\} \right\}$$

be the set of all 2 × 2 matrices with entries from $Z^+ \cup \{0\}$. V is a set vector space over the set S = $Z^+$ = {1, 2, …, ∞}. Take

$$W = \left\{ \begin{pmatrix} x & y \\ z & w \end{pmatrix} \middle| \; x, y, z, w \in 2Z^+ = \{2, 4, 6, \ldots\} \right\} \subseteq V.$$

W is a subset vector subspace over the subset T = {2, 4, ..., ∞}.

***Example 2.1.28:*** Let

$$V = \left\{ \begin{pmatrix} a & b & c & d \\ e & f & g & h \end{pmatrix} \middle| \; a, b, c, d, e, f, g \text{ and } h \in Z \right\}$$

be the set of all 2 × 4 matrices with entries from the set of integers. V is a set vector space over the set S = $Z^+$ = {1, 2, …, ∞}. Take

$$W = \left\{ \begin{pmatrix} a & b & c & d \\ e & f & g & h \end{pmatrix} \middle| \; a, b, c, d, e, f, g \text{ and } h \in Z^+ \right\} \subseteq V.$$

Clearly W is a subset vector subspace over T = {2, 4, …, ∞} ⊆ S = $Z^+$.

Now as in case of vector spaces we can have the following theorem.

**THEOREM 2.1.2:** *Let V be a set vector space over the set S. Suppose $W_1, \ldots, W_n$ be a set of n subset vector subspaces of V over the same subset T of S. Then either $\bigcap_{i=1}^{n} W_i = \phi$ or $\bigcap_{i=1}^{n} W_i$ is a subset vector subspace over T or if each $W_i$ contains 0 then $\bigcap_{i=1}^{n} W_i = \{0\}$, the subset vector zero subspace of V over T.*

The proof is left as an exercise for the reader as the proof involves only simple set theoretic techniques.

*Note:* We cannot say anything when the subset vector subspaces of V are defined over different subsets of S. We illustrate this situation by some more examples.

***Example 2.1.29:*** Let $V = Z^+ \times Z^+ \times Z^+$ be a set vector space over a set $S = Z^+ = \{1, 2, \ldots, \infty\}$. $W = \{2, 4, 6, 8, \ldots\} \times \{\phi\} \times \{\phi\}$ is a subset vector subspace of V over the subset T = {2, 4, 8, 16, …}. $W_1 = \{3, 6, 9, \ldots\} \times \phi \times \phi \subset V$, is a subset subvector space over the subset $T_1 = \{3, 3^2, \ldots\}$. Clearly we cannot define $W \cap W_1$ for we do not have even a common subset over which it can be defined as $T \cap T_1 = \phi$.

From this example a very natural question is if $T \cap T_1 \neq \phi$ and if $W_1 \cap W$ is not empty can we define some new structure. For this we make the following definition.

**DEFINITION 2.1.4:** *Let V be a set vector space over the set S. Suppose $W_i$ is a subset vector subspace defined over the subset $T_i$ of S for i = 1, 2, …, n; $n < \infty$ and if $\cap W_i \neq \phi$ and $T = \cap T_i \neq \phi$; then we call W to be a sectional subset vector sectional subspace of V over T.*

*Note:* We call it a sectional for every subset vector subspace contributes to it.

We give illustration of the same.

***Example 2.1.30:*** Let $V = Z^+ \times Z^+ \times Z^+ \times Z^+$ be a set vector space over the set $Z^+ = \{1, 2, \ldots, \infty\}$. Let $W_1 = \{2, 4, 6, \ldots\} \times \{2, 4, 6, \ldots\} \times Z^+ \times Z^+$ be a subset vector subspace over $T_1 = \{2, 4, \ldots, \infty\}$. $W_2 = \{Z^+\} \times \{2, 4, 6, ...\} \times Z^+ \times Z^+$ be a subset vector subspace over $T_2 = \{2, 2^2, \ldots, \infty\}$. Let $W_3 = \{2, 2^2, \ldots, \infty\} \times \{Z^+\} \times \{2, 4, 6, \ldots\} \times \{2, 2^2, \ldots, \infty\}$ be a subset vector subspace over $T_3 = \{2, 2^3, 2^5, \ldots, \infty\}$. Consider $W_1 \cap W_2 \cap W_3 = \{2, 2^2, ..., \infty\} \times \{2, 4, 6, ...\} \times \{2, 4, 6, ..., \infty\} \times \{2, 2^2, ..., \infty\} = W$. Now $T = T_1 \cap T_2 \cap T_3 = \{2, 2^3, 2^5, \ldots\}$; W is a sectional subset vector sectional subspace of V over T.

We see a sectional subset vector sectional subspace is a subset vector subspace but a sectional subset vector sectional subspace in general is not a subset vector subspace.

We prove the following interesting theorem.

**THEOREM 2.1.3:** *Every sectional subset vector sectional subspace W of set vector space V over the set S is a subset vector subspace of a subset of S but not conversely.*

*Proof:* Let V be the set vector space over the set S. $W = \bigcap_{i=1}^{n} W_i$ be a sectional subset vector sectional subspace over $T = \bigcap_{i=1}^{n} T_i$ where each $W_i$ is a subset vector subspace over $T_i$ for $i = 1, 2, \ldots, n$ with $W_i \neq W_j$ and $T_i \neq T_j$; for $i \neq j$, $1 \leq i, j \leq n$. We see W is a sectional subset vector subspace over $T_i$ ; $i = 1, 2, \ldots, n$.

We illustrate the converse by an example.

***Example 2.1.31:*** Let $V = \{2, 2^2, \ldots\}$ be a set vector space over the set $S = \{2, 2^2\}$. $W_1 = \{2^2, 2^4, \ldots\}$ is a subset vector subspace of V over the set $T_1 = \{2\}$ and $W_2 = \{2^3, 2^6, \ldots, \infty\}$ is a subset vector subspace of over the set $T_2 = \{2^2\}$. Now $W_1 \cap W_2 \neq \phi$ but $T_1 \cap T_2 = \phi$. We do not and cannot make $W_1$ or $W_2$ as sectional subset vector sectional subspace.

That is why is general for every set vector space V over the set S we cannot say for every subset vector subspace $W(\subset V)$ over the subset $T(\subset S)$ we can find atleast a subset vector subspace $W_1(\subset V)$ over the subset $T_1(\subset S)$ such that $W_1 \cap W_2 \neq \phi$ and $T \cap T_1 \neq \phi$. Hence the theorem 2.1.3.

Now before we define the notion of basis, dimension of a set S we mention certain important facts about the set vector spaces.

1. All vector spaces are set vector spaces but not conversely.

The converse is proved by giving counter examples.

Take $V = Z^+ = \{0, 1, 2, \ldots, \infty\}$, V is a set vector space over the set $S = \{2, 4, 6, \ldots, \infty\}$. We see V is not an abelian group under addition and S is not a field so V can never be a vector space over S. But if we have V to be a vector space over the field F we see V is a set vector space over the set F as for every $c \in F$ and $v \in V$ we have $cv \in V$. Thus every vector space is a set vector space and not conversely.

2. All semivector spaces over the semifield F is a set vector space but not conversely.

We see if V is a semivector space over the semifield F then V is a set and F is a set and for every $v \in V$ and $a \in F$; $av \in V$ hence V is trivially a set vector space over the set F.

However take $\{-1, 0, 1, 2, \ldots, \infty\} = V$ and $S = \{0, 1\}$ V is a set vector space over S but V is not even closed with respect to '+' so V is not a group. Further the set $S = \{0, 1\}$ is not a semifield so V is not a semivector space over S. Hence the claim.

Thus we see the class of all set vector spaces contains both the collection of all vector spaces and the collection of all semivector spaces. Thus set vector spaces happen to be the most generalized concept.

Now we proceed on to define the notion of generating set of a set vector space V over the set S.

**DEFINITION 2.1.5:** *Let V be a set vector space over the set S. Let $B \subseteq V$ be a proper subset of V, we say B generates V if every element v of V can be got as sb for some $s \in S$ and $b \in B$. B is called the generating set of V over S.*

***Example 2.1.32:*** Let $Z^+ \cup \{0\} = V$ be a set vector space over the set S = {1, 2, 4, …, ∞}. The generating set of V is B = {0, 1, 3, 5, 7, 9, 11, 13, …, 2n + 1, …} B is unique. Clearly the cardinality of B is infinite.

***Examples 2.1.33:*** Let $V = Z^+ \cup \{0\}$ be a set vector space over the set $S = \{1, 3, 3^2, \ldots, \infty\}$. B = {0, 1, 2, 4, 5, 6, 7, 8, 10, 11, 12, 13, 14, 15, 16, 17, 18, 19, 20, 21, 22, 23, 24, 25, 26, 28, …, 80, 82, …, 242, 244, …, $3^6 - 1$, …, $3^6 + 1$, …, ∞} ⊆ V is the generating set of V.

***Example 2.1.34:*** Let $V = \{0, 3, 3^2, \ldots, 3^n, \ldots\}$ be a set vector space over S = {0, 1, 3}. B = {0, 3} is a generating set for 3.0 = 0, 3.1 = 3, $3.3 = 3^2$, $3.3^2 = 3^3$ so on.

***Example 2.1.35:*** Let

$$V = \left\{ \begin{pmatrix} a & b \\ c & d \end{pmatrix} \middle| \; a, b, c, d \in Z^+ = \{1, 2, \ldots, \infty\} \right\}$$

be the set vector space over the set $S = Z^+$. The generating set of V is infinite.

Thus we see unlike in vector space or semivector spaces finding the generating set of a set vector space V is very difficult.

When the generating set is finite for a set vector space V we say the set vector space is finite set hence finite cardinality or finite dimension, otherwise infinite or infinite dimension.

***Example 2.1.36:*** Let $Z^+ = V = \{1, 2, \ldots, \infty\}$ be a set vector space over the set $S = \{1\}$. The dimension of V is infinite.

***Example 2.1.37:*** Let $V = \{1, 2, \ldots, \infty\}$ be a set vector space over the set $S = \{1, 2, \ldots, \infty\} = V$. Then dimension of V is 1 and B is uniquely generated by {1}. No other element generates V.

Thus we see in case of set vector space we may have only one generating set. It is still an open problem to study does every set vector space have one and only one generating subset?

***Example 2.1.38:*** Let $V = \{2, 4, 6, \ldots\}$ be a set vector space over the set $S = \{1, 2, 3, \ldots\}$. $B = \{2\}$ is the generating set of V and dimension of V is one.

Thus we have an important property enjoyed by set vector spaces. We have a vector space V over a field F of dimension one only if V = F, but we see in case of set vector space V, dimension V is one even if $V \neq S$. The example 2.1.38 gives a set vector space of dimension one where $V \neq S$.

***Example 2.1.39:*** Let $V = \{2, 4, 6, \ldots, \infty\}$ be a set vector space over the set $S = \{2, 2^2, \ldots, \infty\}$. $B = \{2, 6, 10, 12, 14, 18, 20, 22, 24, \ldots, 30, 34, \ldots\}$ is a generating set of V.

Now in case of examples 2.1.38 and 2.1.39 we see $V = \{2, 4, \ldots, \infty\}$ but only the set over which they are defined are different so as in case of vector spaces whose dimension is dependent on the field over which it is defined are different so as in case of vector spaces whose dimension is dependent on the field over which it is defined so also the cardinality of the generating set of a vector space V depends on the set over which V is defined. This is clear from examples 2.1.38 and 2.1.39.

Now having defined cardinality of set vector spaces we define the notion of linearly dependent and linearly independent set of a set vector space V over the set S.

**DEFINITION 2.1.6:** *Let V be a set vector space over the set S. B a proper subset of V is said to be a linearly independent set if x, y $\in$ B then x $\neq$ sy or y $\neq$ s'x for any s and s' in S. If the set B is not linearly independent then we say B is a linearly dependent set.*

We now illustrate the situation by the following examples.

***Example 2.1.40:*** Let $Z^+$ = V = {1, 2, ..., ∞} be a set vector space over the set S = {2, 4, 6, …, ∞}. Take B = {2, 6, 12} ⊂ V; B is a linearly dependent subset for 12 = 6.2, for 6 ∈ S. B = {1, 3} is a linearly independent subset of V.

***Example 2.1.41:*** Let V = {0, 1, 2, …, ∞} be a set vector space over the set S = {3, $3^2$, …}. B = {1, 2, 4, 8, 16, …} is a linearly independent subset of V. As in case of vector spaces we can in case of set vector spaces also say a set B which is the largest linearly dependent subset of V? A linearly independent subset B of V which can generate V, then we say B is a set basis of V or the generating subset of V and cardinality of B gives the dimension of V.

## 2.2 Set Linear Transformation of Set Vector Spaces

In this section we proceed on to define the notion of set linear transformation of set vector spaces. As in case of vector spaces we can define set linear transformation of set vector spaces only if the set vector spaces are defined over the same set S.

**DEFINITION 2.2.1:** *Let V and W be two set vector spaces defined over the set S. A map T from V to W is said to be a set linear transformation if*

$$T(v) = w$$

*and*

$$T(sv) = sw = sT(v)$$

*for all v ∈ V, s ∈ S and w ∈ W.*

***Example 2.2.1:*** Let V = $Z^+$ = {0, 1, 2, …, ∞} and W = {0, 2, 4, …, ∞} be set vector spaces over the set S = {0, 2, $2^2$, $2^3$, …}.

$$
\begin{aligned}
&T : V \to W \\
&T(0) = 0, \\
&T(1) = 2, \\
&T(2) = 4, \\
&T(3) = 6
\end{aligned}
$$

and so on. T ($2^n$p) = $2^n$ (2p) for all p ∈ V. Thus T is a set linear transformation from V to W. As in case of vector spaces we will not be always in a position to define the notion of null space of T in case of set vector spaces. For only if 0 ∈ V as well as W; we will be in a position to define null set of a set linear transformation T.

Now we proceed on to define the notion of set linear operator of a set vector space V over the set S.

**DEFINITION 2.2.2:** *Let V be a set vector space over the set S. A set linear transformation T from V to V is called the set linear operator of V.*

***Example 2.2.2:*** Let V = $Z^+$ = {1, 2, …, ∞} be a set vector space over the set S = {1, 3, $3^2$, …, ∞}. Define T from V to V by T(x) = 2x for every x ∈ V. T is a set linear operator on V.

Now it is easy to state that if V is a set vector space over the set S and if $O_S(V)$ denotes the set of all set linear operators on V then $O_S(V)$ is also a set vector space over the same set S. Similarly if V and W are set vector spaces over the set S and $T_S$ (V, W) denotes the set of all linear transformation from V to W then $T_S$ (V,W) is also a set vector space over the set S. For if we want to prove some set V is a set vector space over a set S it is enough if we show for every s ∈ S and v ∈ V; sv ∈ V. Now in case of $O_S$ (V) = {set of all set linear operators from V to V}. $O_S$ (V) is a set and clearly for every s ∈ S and for every T ∈ $O_S(V)$, sT ∈ $O_s$ (V) i.e., sT is again a set linear operator of V. Hence $O_S$ (V) is a set vector space over the set S.

Likewise if we consider the set $T_S$ (V, W) = {set of all set linear transformations from V to W}, then the set $T_S$ (V, W) is again a set vector space over S. For if s $\in$ S and T $\in$ $T_S$ (V, W) we see sT is also a set linear operator from V to W. Hence the claim. Now we can talk about the notion of invertible set linear transformation of the set vector spaces V and W defined over the set S.

Let T be a set linear transformation from V into W. We say T is set invertible if there exists a set linear transformation U from W into V such that UT and TU are set identity set maps on V and W respectively. If T is set invertible, the map U is called the set inverse of T and is unique and is denoted by $T^{-1}$.

Further more T is set invertible if and only if (1) T is a one to one set map that is $T\alpha = T\beta$ implies $\alpha = \beta$ (2) T is onto that is range of T is all of W.

We have the following interesting theorem.

**THEOREM 2.2.1:** *Let V and W be two set vector spaces over the set S and T be a set linear transformation from V into W. If T is invertible the inverse map $T^{-1}$ is a set linear transformation from W onto V.*

*Proof:* Given V and W are set vector spaces over the set S. T is a set linear transformation from V into W. When T is a one to one onto map, there is a uniquely determined set inverse map $T^{-1}$ which set maps W onto V such that $T^{-1}T$ is the identity map on V and $TT^{-1}$ is the identity function on W.

Now what we want to prove is that if a set linear transformation T is set invertible then $T^{-1}$ is also set linear.

Let x be a set vector in W and c a set scalar from S. To show

$$T^{-1} (cx) = cT^{-1} (x).$$

Let $y = T^{-1}(x)$ that is y be the unique set vector in V such that Ty = x. Since T is set linear T(cy) = cTy. Then cy is the unique set vector in V which is set by T into cx and so

$$T^{-1}(cx) = cy = c\,T^{-1}(x)$$

and $T^{-1}$ is set linear !

Suppose that we have a set invertible set linear transformation T from V onto W and a set invertible set linear transformation U from W onto Z. Then UT is set invertible and $(UT)^{-1} = T^{-1} U^{-1}$.

To obtain this it is enough if we verify $T^{-1} U^{-1}$ is both left and right set inverse of UT.

Thus we can say in case of set linear transformation; T is one to one if and only if $T\alpha = T\beta$ if and only if $\alpha = \beta$.

Since in case of set linear transformation we will not always be in a position to have zero to be an element of V we cannot define nullity T or rank T. We can only say if $0 \in V$, V a set vector space over S and $0 \in W$, W also a set vector space over S then we can define the notion of rank T and nullity T, where T is a set linear transformation from V into W.

We may or may not be in a position to have results of linear transformation from vector spaces.

Also the method of representing every vector space V over a field F of dimension n as $V \cong \underbrace{F \times \ldots \times F}_{n-\text{times}}$ may not be feasible in case of set vector spaces. Further the concept of representation of set linear transformation as a matrix is also not possible for all set vector spaces. So at this state we make a note of the inability of this structure to be always represented in this nice form.

Next we proceed on to define the new notion of set linear functionals of a set vector space V over a set S.

**DEFINITION 2.2.3:** *Let V be a set vector space over the set S. A set linear transformation from V onto the set S is called a set linear functional on V, if*

$$f : V \to S$$

$$f(c\alpha) = cf(\alpha)$$

*for all* $c \in S$ *and* $\alpha \in V$.

***Example 2.2.3:*** Let V = {0, 1, 2, ..., ∞} and S = {0, 2, 4, …, ∞}. f : V → S defined by f(0) = 0, f(1) = 2, f(2) = 4, ..., is a set linear functional on V.

***Example 2.2.4:*** Let $V = Z^+ \times Z^+ \times Z^+$ be a set vector space over $Z^+$. A set linear functional $f : V \to Z^+$ defined by f(x, y, z) = x + y + z.

Now we proceed onto define the new notion of set dual space of a set vector space V.

**DEFINITION 2.2.4:** *Let V be a set vector space over the set S. Let L(V,S) denote the set of all linear functionals from V into S, then we call L(V,S) the set dual space of V. Infact L(V, S) is also a set vector space over S.*

The study of relation between set dimension of V and that of L(V, S) is an interesting problem.

***Example 2.2.5:*** Let V = {0, 1, 2, …, ∞} be a set vector space over the set S = {0, 1, 2, ..., ∞}. Clearly the set dimension of V over S is 1 and this has a unique generating set B = {1}. No other proper subset of B can ever generate V.

What is the dimension of L (V, S)?

We cannot define the notion of set hyperspace in case of set vector space using set linear functionals. We can define the concept of set annihilator if and only if the set vector space V and the set over which it is defined contains the zero element, otherwise we will not be in a position to define the set annihilator.

**DEFINITION 2.2.5:** *Let V be a set vector space over the set S, both V and S has zero in them. Let A be a proper subset of the*

*set V, the set annihilator of A is the set $A^o$ of all set linear functionals f on V such that $f(\alpha) = 0$ for every $\alpha$ in A.*

## 2.3 Set Linear Algebra and its Properties

In this section we define the new notion of set linear algebra and enumerate a few of its properties.

**DEFINITION 2.3.1:** *Let V be a set on which is defined a binary operation '+' which is closed on V i.e., for every pair a, b in V a + b $\in$ V. S any set we call V to be set linear algebra over S if for all s $\in$ S and v $\in$ V, s v $\in$ V and s(a+b) = sa + sb for all a, b $\in$ V and s $\in$ S.*

It is interesting and important to note that all set linear algebras are set vector spaces but a set vector space is not in general a set linear algebra. For we see in case of set vector space V we may not have any closed binary operation on the set V. In some cases we may not be in a position to define a closed binary operation on V. We see on the other hand, a set linear algebra V over the set S is such that V is endowed with a closed binary operation on it.

***Example 2.3.1:*** Let V = $Z^o$ = {0, 1, 2, …, ∞}, Let '+' be the operation on V. Suppose S = {0, 2, 4, …, ∞}, V is a set linear algebra over S. We see s (a + b) = sa + sb, ∀ s ∈ S and for all a, b ∈ V.

***Example 2.3.2:*** Let

$$V = M_{2\times 3} = \left\{ \begin{pmatrix} a & b & c \\ d & e & f \end{pmatrix} \middle| \text{ a, b, c, d, e and f} \in Z^+ \cup \{0\} \right\}.$$

$M_{2\times 3}$ is a set with addition of matrices as a closed binary operation. Take S = {0, 2, 4, …, ∞}. V is a set linear algebra over S.

***Example 2.3.3:*** Let $V = Z^+ \times Z^+ \times Z^+ = \{(a, b, c) | a, b c \in Z^+\}$; $(a, b, c) + (d, e, f) \in V$. Take $S = Z^+$; V is a set linear algebra over S.

Now we define the notion of set linear subalgebra of a set linear algebra V.

**DEFINITION 2.3.2:** *Let V be a set linear algebra over the set S. Suppose W is a proper non empty subset of V. If W itself is a set linear algebra over S then we call W to be a set linear subalgebra of V over the set S.*

We illustrate this by the following examples.

***Example 2.3.4:*** Let $V = Z^+ \times Z^+ \times Z^+ \times Z^+ \times Z^+$ be a set linear algebra over the set $S = Z^+$. Take $W = 2Z^+ \times Z^+ \times \phi \times 2Z^+ \times Z^+ \subseteq V$. Clearly W is also a set linear algebra over the set $S = Z^+$. Thus W is set linear subalgebra of V over the set S.

***Example 2.3.5:*** Let $V = Z^+$ be a set linear algebra over the set $S = \{2, 4, \ldots, \infty\}$. Take $W = \{2, 4, 6, 8, \ldots\} \subseteq V$. W is a set linear subalgebra of V.

***Example 2.3.6:*** Let

$$V = \left\{ \begin{pmatrix} a & b & c \\ d & e & f \end{pmatrix} \middle| \; a, b, c, d, e, f \in Z^+ \right\}$$

be a set linear algebra over the set $S = Z^+$. Take

$$W = \left\{ \begin{pmatrix} x & y & z \\ u & v & w \end{pmatrix} \middle| \; x, y, z, u, v \text{ and } w \in 2Z^+ \right\} \subseteq V.$$

W is also a set linear algebra over the set $S = Z^+$ i.e., W is a set linear subalgebra of V over the set S.

Just like a set vector space a set linear algebra may or may not contain the zero. If a set linear algebra contains zero then we call it a set linear algebra with zero. Examples 2.3.1 and 2.3.2 are set linear algebras with zero. It is still important to note that a set linear subalgebra may or may not contain the zero element even if the set linear algebra is one with zero.

It is interesting to note that a set linear subalgebra will contain zero if and only if the set S over which the set linear algebra is defined contains zero.

***Example 2.3.7:*** Let $V = Z^+ \cup \{0\} \times Z^+ \cup \{0\} \times Z^+ \cup \{0\}$ be a set linear algebra with zero defined over the set $S = \{0, 2, \ldots, \infty\}$. Clearly every set linear subalgebra contains zero as S the set over which it is defined contains zero.

***Example 2.3.8:*** Let $V = Z \times Z$ be a set vector space over $S = Z^+ = \{1, 2, \ldots, \infty\}$. V is a set linear algebra over S. Take $W = Z^+ \times Z^+ \subset V$, W is a set linear subalgebra of V over S. V is a set linear algebra with zero but W is not a set linear subalgebra with zero but W is only a set linear subalgebra of V over S.

Now we proceed on to define the notion of subset linear subalgebra of a set linear algebra over a set S.

**DEFINITION 2.3.3:** *Let V be a set linear algebra over the set S, suppose W is a subset of V and W is a set linear algebra over the subset $P \subseteq S$. Then we call W to be a subset linear subalgebra of V over P.*

We can have many subset linear subalgebras over a subset $P \subset S$. Likewise we can have the same subset W of V to be a subset linear subalgebra over different subsets of S.

We illustrate a few of these situations by examples.

***Example 2.3.9:*** Let $V = Z^+ \times Z^+ \times Z^+ \times Z^+$ be a set linear algebra over the set $S = Z^+ = \{1, 2, \ldots, \infty\}$. Take $W = Z^+ \times \{\phi\} \times Z^+ \times$ $\{\phi\} \subset V$ W is a subset linear subalgebra over P = {2, 4, 6, …} $\subseteq S = Z^+$.

Take $W_1 = 2Z^+ \times 2Z^+ \times Z^+ \times Z^+ \subseteq V$, $W_1$ is also a subset linear subalgebra over P. If $W_2 = \{\phi\} \times 2Z^+ \times 2Z^+ \times \{0\} \subseteq V$. $W_2$ is also a subset linear subalgebra over P. However if we take $W_3 = \{3, 3^2, \ldots\} \times \{1, 3, 3^2 \ldots\} \times \{\phi\} \times \{\phi\} \subseteq V$. $W_3$ is not a subset linear subalgebra over $P \subseteq S$. However $W_3$ is a subset sublinear algebra over the set $\{1, 3, 3^2, \ldots\}$.

Thus we see in general a subset of V may be a subset linear subalgebra over different subsets of S. So we can speak of intersection of subset linear subalgebras or set linear subalgebra. We may have the intersection to be at times empty. We have to overcome all these problems which happen to be challenging.

Since in case of set linear algebra we have no means to associate matrices with set linear operators or set linear transformation.

Now having defined set linear algebra and set linear subalgebra we now proceed to give some of its properties.

The set linear transformation of a set linear algebra would be a little different from that of set vector spaces, we define this in the following.

**DEFINITION 2.3.4:** *Let V and W be two set linear algebras defined over the same set S. A map T from V to W is called a set linear transformation if the following condition holds.*

$$T(c\alpha + \beta) = cT(\alpha) + T(\beta)$$

*for all $\alpha, \beta \in V$ and $c \in S$.*

*Note:*
1. The set linear transformation of set linear algebras are defined if and only if both the set linear algebras are defined over the same set S.
2. The set linear transformation of set linear algebras are different from the set linear transformation of set vector spaces.

Now we also define the notion of set basis.

Set linearly independent set etc are different for set linear algebras; as in case of set vector spaces, they are not the same. Before we define these notions we just define set linear operator of a set linear algebra.

**DEFINITION 2.3.5:** *Let V be a set linear algebra over the set S. T be a map from V to V which is a set linear transformation, then T is defined to be a set linear operator on V.*

***Example 2.3.10:*** Let V = $Z^+ \times Z^+ \times Z^+$ be a set linear algebra over the set $Z^+$ = S. Define T: V $\to$ V by T (x, y, z) = (x + y, y, x – z). T is set linear operator on V.

***Example 2.3.11:*** Let

$$V = \left\{ \begin{pmatrix} a & b & c \\ d & e & f \end{pmatrix} \middle| a,b,c,d,e,f \in Z^+ \cup \{0\} \right\},$$

V is a closed set under the binary operation matrix addition. Let S = $Z^+ \cup \{0\}$. V is a set linear algebra over the set S.
Define

$$T\left(\begin{pmatrix} a & b & c \\ d & e & f \end{pmatrix}\right) = \begin{pmatrix} a & 0 & b \\ d & 0 & 0 \end{pmatrix}$$

for all

$$\begin{pmatrix} a & b & c \\ d & e & f \end{pmatrix} \in V;$$

T is a set linear operator on V.

We see the concept of set vector spaces and set linear algebras will find its application in coding theory. They can be used to give best type of cryptosystems.

Now having defined set linear operators of a set linear algebra V over the set S we ask the reader to find out whether the set of all set linear operators of V, a set linear algebra over S. What is its dimension if V is a set linear algebra of dimension n?

Presently we proceed on to define yet a particular form of vector spaces which are different from set vector spaces, vector spaces and semivector spaces in the following section.

## 2.4 Semigroup Vector Spaces and their Generalizations

In this section for the first time we proceed on to define a new notion of semigroup vector spaces. Clearly from the definition the reader will be able to follow all semigroup vector spaces will be set vector spaces and not vice versa. Further all vector spaces and semivector spaces will be semigroup vector spaces and not conversely. We in this section define these new notions and illustrate them by examples.

**DEFINITION 2.4.1:** *Let V be a set, S any additive semigroup with 0. We call V to be a semigroup vector space over S if the following conditions hold good.*

1. *$sv \in V$ for all $s \in S$ and $v \in V$.*
2. *$0. v = 0 \in V$ for all $v \in V$ and $0 \in S$; 0 is the zero vector.*
3. *$(s_1 + s_2) v = s_1 v + s_2 v$*

*for all $s_1, s_2 \in S$ and $v \in V$.*

We illustrate this by the following examples.

***Examples 2.4.1:*** Let $V = (Z^+ \cup \{0\}) \times 2Z^+ \cup \{0\} \times (3Z^+ \cup \{0\})$ be a set and $S = Z^+ \cup \{0\}$ be a semigroup under addition. V is a semigroup vector space over S.

***Example 2.4.2:*** Let

$$V = \left\{ \begin{pmatrix} a_1 & a_2 & a_3 & a_4 \\ a_5 & a_6 & a_7 & a_8 \end{pmatrix} \middle| a_i \in Z^+ \cup \{0\}; 1 \le i \le 8 \right\}$$

be a set. Suppose $S = 2Z^+ \cup \{0\}$ be a semigroup under addition. V is a semigroup vector space over S.

***Example 2.4.3:*** Let $V = 3Z^+ \cup \{0\}$ be a set and $S = Z^+ \cup \{0\}$ be a semigroup under addition. V is a semigroup vector space over S.

***Example 2.4.4:*** Let V = {(1 0 0 0 1), (1 0 1 1 1), (0 1 1 1 1), (1 1 1 1 1), (0 0 0 0 0), (0 1 0 1 0)} be a non empty set. S = {0, 1} a semigroup under addition 1 + 1 = 1. V is a semigroup vector space over S.

***Example 2.4.5:*** Let V = {(1 1 1 1 1 1), (1 0 1 0 1 0), (1 1 1 0 0 0), (0 0 1 1 0 0), (0 0 0 0 0 0), (0 1 0 1 0 1)} and S = {0, 1} a semigroup under '+' with 1 + 1 = 0 then V is a semigroup vector space over S.

***Example 2.4.6:*** Let $V = Z_{12} = \{0, 1, 2, \ldots, 11\}$ be a set. V is a semigroup vector space over the semigroup S = {0,1} under addition modulo 2.

***Example 2.4.7:*** Let $V = 5Z^+ \cup \{0\}$ be a semigroup vector space over the semigroup $S = Z^+ \cup \{0\}$.

***Example 2.4.8:*** Let $\{Z^+ \cup \{0\}\}$ [x] be a semigroup vector space over $S = Z^+ \cup \{0\}$. S, a semigroup under addition.

***Example 2.4.9:*** Let

$$V = \left\{ \begin{pmatrix} a & b \\ a & b \\ a & b \end{pmatrix} \middle| a, b \in S = \{0\} \cup Z^+ \right\}$$

be a semigroup vector space over the semigroup $Z^+ \cup \{0\}$.

***Example 2.4.10:*** Let $V = Z_{12} \times Z_{12} \times Z_{12}$, the set of triples modulo 12. Suppose S = {0, 3, 6, 9} is a semigroup under addition modulo 12. V is a semigroup vector space over the semigroup S, sv ≡ $v_1$(mod 12) i.e., if 6 ∈ S and v = (5, 2, 7); sv = (6, 0, 6) = $v_1$ ∈ V.

Now we proceed on to define the notion of independent subset and generating subset of a semigroup vector space over the semigroup S.

**DEFINITION 2.4.2:** *Let V be semigroup vector space over the semigroup S. A set of vectors {$v_1$, …, $v_n$} in V is said to be a semigroup linearly independent set if*

$$\text{(i) } v_i \neq sv_j$$

*for any s ∈ S for i ≠ j; 1 ≤ i, j ≤ n.*

We first illustrate this by some examples before we proceed to define the generating set of a semigroup vector space.

***Example 2.4.11:*** Let V = $Z_{12}$ = {0, 1, 2, …, 11} be the set of integers modulo 12. S = {0, 1} be the semigroup under addition modulo 2. V is a semigroup vector space over the semigroup S. T = {1, 2, …, 11} is a linearly independent subset of V.

It is interesting and important to note that always all subsets of an independent set is also an independent set of V.

***Example 2.4.12:*** Let V = $Z^+ \cup \{0\}$ be a set. S = {0, 2, 4, …} be a semigroup under addition. Take T = {2, $2^5$, $2^9$, …, ∞} ⊆ V. T is not an independent subset of V. Now take P = {3, 9, 27, …} ⊆ V. P is an independent subset of V. Infact every subset of P is also an independent subset of V.

We now proceed on to define the notion of generating subset of a semigroup vector space V over a semigroup S.

**DEFINITION 2.4.3:** *Let V be a semigroup vector space over the semigroup S under addition. Let T = {$v_1$, …, $v_n$} ⊆ V be a subset of V we say T generates the semigroup vector space V over S if every element v of V can be got as v = $sv_i$, $v_i$ ∈ T; s ∈ S.*

We now illustrate this situation by the following example.

***Example 2.4.13:*** Let $V = 3Z^+ \cup \{0\}$ be a semigroup vector space over the semigroup $S = Z^+ \cup \{0\}$. Take $T = \{3\} \subseteq V$; T generates V over S.

***Example 2.4.14:*** Let $V = Z_{20} = \{0, 1, 2, \ldots, 19\}$ integers modulo 20. $S = \{0, 5, 10\}$ be the semigroup under + modulo 20. V is a semigroup vector space over S. $T = \{1, 2, 3, 4, 6, 7, 8, 9, 11, 12, 13, 14, 16, 17, 18, 19\}$ is a generating set of V.

***Example 2.4.15:*** Let $V = Z_{20} = \{0, 2, 4, 6, \ldots, 18\}$ integers modulo 20. $S = \{0, 1\}$ be the semigroup under addition modulo 2. V is a semigroup vector space over the semigroup S. Take the set $T = \{2, 4, 6, 8\ldots 18\} \subseteq V$; T is a generating subset of V.

***Example 2.4.16:*** Let $V = Z_{20} = \{0, 1, \ldots, 19\}$ integers modulo 20. Let $S = \{0, 10\}$ be a semigroup under addition modulo 20. V is a semigroup vector space over the semigroup S. $T = \{1, 2, 3, 4, 5, 6, 7, 8, 9, 11, 12, 13, 14, 15, 16, 17, 18, 19\} \subseteq V$ is a generating set of V over S.

Thus we see from the above examples the generating set of V is highly dependent on the semigroup S over which it is defined. Further nothing can be said about the dimension of V.

In the above cases when $|V| = 20$, we see the dimension of V is $\leq 20$ and nothing more so we can say, that is if $|V| = n < \infty$ then dimension of V is less than n and nothing more i.e., V is finite dimensional.

Now it is left as an exercise to find out that the generating set of semigroup vector space is also an independent subset of V.

**DEFINITION 2.4.4:** *Let V be a semigroup vector space over the semigroup S. Suppose P is a proper subset of V and P is also a semigroup vector space over the semigroup S, then we call P to be semigroup subvector space of V.*

We now illustrate this situation by the following examples.

***Example 2.4.17:*** Let V = $(Z^+ \cup \{0\}) \times (Z^+ \cup \{0\}) \times (2Z^+ \cup \{0\})$ be a semigroup vector space over the semigroup S = $Z^+ \cup \{0\}$ , under addition. Take W = $Z^+ \cup \{0\} \times \{0\} \times \{0\} \subseteq V$; W is also a semigroup vector space over S and we see W is a semigroup subvector space of V over S.

***Example 2.4.18:*** Let

$$V = \left\{ \begin{pmatrix} a_1 & a_2 & a_3 \\ a_4 & a_5 & a_6 \end{pmatrix} \middle| \begin{array}{c} a_i \in Z^+ \cup \{0\} \\ 1 \le i \le 6 \end{array} \right\}$$

be a semigroup vector space over the semigroup S = $\{Z^+ \cup \{0\}\}$. Take

$$W = \left\{ \begin{pmatrix} a_1 & a_2 & a_3 \\ 0 & 0 & 0 \end{pmatrix} \middle| \begin{array}{c} a_i \in Z^+ \cup \{0\} \\ 1 \le i \le 3 \end{array} \right\} \subset V;$$

W is a semigroup vector subspace of V over S.

***Example 2.4.19:*** Let V = $\{Q^+ \cup \{0\}\}$ be a semigroup vector space over the semigroup S = $Z^+ \cup \{0\}$. Take W = $2Z^+ \cup \{0\}$, W is a semigroup subvector space of V over S.

***Example 2.4.20:*** Take V = {(1 1 0 1 0), (0 0 0 0 0), (1 1 1 1 1), (1 0 1 0 1), (0 1 1 1 0), (0 0 0), (0 1 1), (0 1 0), (1 0 1)} be a set. V is a semigroup vector space over the semigroup S = {0, 1} under addition 1 + 1 = 1.

Take W = {(1 0 1), (0 1 1), (0 0 0), (0 1 0)} $\subseteq$ V; W is a semigroup subvector space over the semigroup S = {0, 1}. In fact every subset of V with the two elements (0 0 0) and or (0 0 0 0 0) is a semigroup subvector space over the semigroup S = {0, 1}.

***Example 2.4.21:*** Let V = $Z^+ \cup \{0\}$ be a semigroup vector space over the semigroup S = $Z^+ \cup \{0\}$ under '+'. The dimension of V is one and the generating set of V is {1}. Take W = $\{2Z^+$ $\cup \{0\}\} \subseteq V$ is a semigroup subvector space over S. The dimension of W is also one given by the set {2}. This is only possible when we have semigroup vector spaces V over the semigroup S or set vector spaces over the set S. In fact V has infinite number of semigroup vector subspaces defined in the semigroup $S = Z^+ \cup \{0\}$.

Now we proceed on to define the new notion of semigroup linear algebra over the semigroup.

**DEFINITION 2.4.5:** *Let V be a semigroup vector space over the semigroup S. If V is itself a semigroup under '+' then we call V to be a semigroup linear algebra over the semigroup S, if* $s(v_1 + v_2) = sv_1 + sv_2$ *;* $v_1, v_2 \in V$ *and* $s \in S$*.*

We illustrate this definition by examples before we proceed to study their properties.

***Example 2.4.22:*** Let $V = (Z^+ \cup \{0\})$ be a semigroup under addition $S = \{2Z^+ \cup \{0\}\}$ also a semigroup under addition. V is a semigroup linear algebra over S. Also S is a semigroup linear algebra over V.

***Example 2.4.23:*** Let $V = (Z^+ \cup \{0\}) \times (Z^+ \cup \{0\}) \times (Z^+ \cup \{0\})$ be a semigroup under addition of elements component wise. $S = Z^+ \cup \{0\}$ be a semigroup under addition. V is a semigroup linear algebra over S. Clearly S is not a semigroup linear algebra over V. Only in example 2.4.22 it so happened that V was a semigroup linear algebra over S and S happened to be a semigroup linear algebra over V.

Now we proceed on to give more examples of semigroup linear algebra.

***Example 2.4.24:*** Let

$$V = \left\{ \begin{pmatrix} a_1 & a_2 & a_3 \\ a_4 & a_5 & a_6 \end{pmatrix} \middle| \, a_i \in Z^+ \cup \{0\}; 1 \le i \le 6 \right\}$$

be a semigroup under matrix addition. $S = Z^+ \cup \{0\}$ a semigroup under '+'. Clearly V is a semigroup linear algebra over S.

***Example 2.4.25:*** Let

$$V = \left\{ \begin{pmatrix} a_1 & a_2 \\ a_3 & a_4 \end{pmatrix} \middle| a_i \in Z^+ \cup \{0\} \right\}$$

be a semigroup under addition. V is a semigroup linear algebra over the semigroup $S = Z^+ \cup \{0\}$.

***Example 2.4.26:*** Let $V = Q^+ \cup \{0\}$ be a semigroup under '+'. V is a semigroup linear algebra over the semigroup $S = Z^+ \cup \{0\}$ under usual addition.

We observe the following important properties.

**THEOREM 2.4.1:** *Every semigroup linear algebra is a semigroup vector space and a semigroup vector space in general is not a semigroup linear algebra.*

*Proof:* We see every semigroup linear algebra over a semigroup is obviously a semigroup vector space over a semigroup. Hence the claim. However we prove the converse using only a counter example.

Consider the semigroup vector space V = {(0 1 0 0), (0 0 0 0), (0 0 1 1), (0 0 0 1), (0 0 0), (1 1 1), (1 0 1)} over the semigroup S = {0, 1} under addition with 1 + 1 = 1. We see V is not a semigroup. But V is a semigroup vector space over S and not a semigroup linear algebra over S. Hence the theorem.

Another property about semigroup linear algebra V over the semigroup S is that all linear algebras and semilinear algebras are semigroup linear algebras. Thus the notion of semigroup linear algebras are a more generalized concept. However all semigroup linear algebras are set linear algebras and not conversely. Thus we have examples to show that set linear algebras are not semigroup linear algebras.

For instance take the following example.

***Example 2.4.27:*** Let V = {0, 1, 2, …, ∞} be a set. Take S = {0, 1, 3, 5, …}, S is not a semigroup under addition it is only a set. Clearly V is a set linear algebra over S for V is an additively closed set. We see V is not a semigroup linear algebra as S can never be a semigroup even if V is a semigroup under addition.

Thus the class of all semigroup linear algebras is a subset of the class of all set linear algebras.

Now we proceed on to define the substructure properties and transformations on semigroup linear algebras.

**DEFINITION 2.4.6:** *Let V be a semigroup linear algebra over the semigroup S. Suppose P is a proper subset of V and P is a subsemigroup of V. Further if P is a semigroup linear algebra over the same semigroup S then we call P a semigroup linear subalgebra of V over S.*

We now illustrate this by some examples before we enumerate some of its properties.

***Example 2.4.28:*** Let

$$\begin{aligned} V &= (Z^+ \cup \{0\}) \times (Z^+ \cup \{0\}) \times (3Z^+ \cup \{0\}) \\ &= \{(x, y, z) / x, y \in Z^+ \cup \{0\} \text{ and } z \in 3Z^+ \cup \{0\}\} \end{aligned}$$

be a semigroup under component wise addition. Take $S = Z^+ \cup \{0\}$, the semigroup under addition. V is a semigroup linear algebra over S. Take $W = \{(3Z^+ \cup \{0\}) \times \{0\} \times (3Z^+ \cup \{0\})\} \subseteq V$. W is a semigroup linear subalgebra of V over S. Take $P = \{Z^+ \cup \{0\}\} \times \{0\} \times \{0\}$, P is also a semigroup linear subalgebra of V over S.

However $\{0, 2, 2^2, 2^3, \ldots\} \times \{0\} \times \{0\}$ is not a semigroup linear subalgebra over S.

***Example 2.4.29:*** Let

$$V = \left\{ \begin{pmatrix} a_1 & \cdots & a_n \\ a_{n+1} & \cdots & a_{2n} \end{pmatrix} \middle| a_i \in Q^+ \cup \{0\};\ 1 \le i \le 2n \right\}$$

be a semigroup under matrix addition. S = $Z^+ \cup \{0\}$ be the semigroup under addition. V is a semigroup linear algebra over S. Let

$$W = \left\{ \begin{pmatrix} 0 & \cdots & 0 \\ a_{n+1} & \cdots & a_{2n} \end{pmatrix} \middle| a_{n+1}, \cdots, a_{2n} \in Z^+ \cup \{0\} \right\} \subseteq V,$$

W is a semigroup linear subalgebra over S.

***Example 2.4.30:*** Let V = $2Z^+ \cup \{0\}$ be a semigroup under addition. S = $\{2Z^+ \cup \{0\}\}$ is a semigroup V = S. V is a semigroup linear algebra over S.

***Example 2.4.31:*** Let V = $\{(Z^+ \cup \{0\}) \times (3Z^+ \cup \{0\}) \times (5Z^+ \cup \{0\}) \times (7Z^+ \cup \{0\})\}$ be a semigroup linear algebra over the semigroup S = $Z^+ \cup \{0\}$. Take W = $\{0\} \cup (3Z^+ \cup \{0\}) \times \{0\} \times \{0\} \subseteq V$, W is a semigroup linear subalgebra of V over S. V has infinite number of semigroup linear subalgebras.

***Example 2.4.32:*** Let V = $(Z^+ \cup \{0\}) \times (4Z^+ \cup \{0\}) \times (Q^+ \cup \{0\})$ be a semigroup linear algebra over S = $Z^+ \cup \{0\}$. This V too has infinite number of semigroup linear subalgebras.

Now as in case of semigroup vector spaces we can define the basis and the generating set of a semigroup linear algebra with a small difference. B = $\{v_1, \ldots, v_n\} \subset V$ is an independent set if $v_i \neq sv_j$; $i \neq j$ for some $s \in S$ and $v_k \neq \sum_{i=1}^{m} s_i v_i$ ; $i \le k \le n$; $m < n$, $s_i \in S$.

We say B is a generating subset of V if B is a linearly independent set and every element v ∈ V can be represented as $v = \sum_{i=1}^{n} \alpha_i v_i$ ; $a_i \in S$; $1 \leq i \leq n$.

***Example 2.4.33:*** Let $V = 2Z^+ \cup \{0\}$ be a semigroup linear algebra over the semigroup $S = Z^+ \cup \{0\}$ . The set $B = \{2\}$ is a generating set of V and the number of elements in B is one thus V is set dimension 1 over S.

***Example 2.4.34:*** Let $V = \{Z^+ \cup \{0\}\} \times \{Z^+ \cup \{0\}\} \times \{Z^+ \cup \{0\}\}$ be a semigroup linear algebra over the semigroup $S = Z^+ \cup \{0\}$. V is a generated by the set B = {(100), (010), (001)} and V is of dimension three over S.

***Example 2.4.35:*** Let $V = \{Q^+ \cup \{0\}\}$ be the semigroup linear algebra over the semigroup $Z^+ \cup \{0\}= S$. V is an infinite dimensional semigroup linear algebra over S.

***Example 2.4.36:*** Let $V = \{Q^+ \cup \{0\}\}$ be a semigroup linear algebra over the semigroup $S = \{Q^+ \cup \{0\}\}$. V is generated by the set B = {1} and dimension of V is 1 over $S = Q^+ \cup \{0\}$.

***Example 2.4.37:*** Let

$$V = \left\{ \begin{pmatrix} a & b \\ c & d \end{pmatrix} \middle| a, b, c, d \in Z^+ \cup \{0\} \right\}$$

be a semigroup linear algebra over the semigroup $Z^+ \cup \{0\} = S$. Then V is generated by

$$B = \left\{ \begin{pmatrix} 1 & 0 \\ 0 & 0 \end{pmatrix}, \begin{pmatrix} 0 & 1 \\ 0 & 0 \end{pmatrix}, \begin{pmatrix} 0 & 0 \\ 1 & 0 \end{pmatrix}, \begin{pmatrix} 0 & 0 \\ 0 & 1 \end{pmatrix} \right\}$$

and dimension of V over S is four. This is the unique subset of V which generates V over S and no other subset of V can generate V.

***Example 2.4.38:*** Let $V = Z^+ \cup \{0\}$ be a semigroup linear algebra over the semigroup $S = 2Z^+ \cup \{0\}$. B = {1} is a generating subset of V as $1 \in V$, $2.1 \in V$ and $2 + 1 \in V$ as V is a semigroup.

Thus we see dimension of V over $S = 2Z^+ \cup \{0\}$ is also 1. If we consider S as a semigroup vector space over V then dimension of V is not 1 and the generating subset of V is not {2}.

***Example 2.4.39:*** Let V[x] = {all polynomials in the variable x with coefficients from the set $Z^+ \cup \{0\}$}, V[x] is a semigroup under addition. V[x] is a semigroup linear algebra over the semigroup $S = \{Z^+ \cup \{0\}\}$. Clearly $B = \{1, x, x^2, x^3, \ldots, \infty\}$ is the generating subset of V[x] and the dimension of V[x] is infinite over S.

***Example 2.4.40:*** Let P[x] = {all polynomials of degree less than or equal 7 in the variable x with coefficient from $Z^+ \cup \{0\}$}. P[x] is a semigroup under polynomial addition. Let $S = Z^+ \cup \{0\}$ be a semigroup under addition. P[x] is a semigroup linear algebra over the semigroup S. Let $B = \{1, x, \ldots, x^7\} \subset P[x]$; B generates P[x] and dimension of P[x] is 8.

Now having seen examples of semigroup linear algebra and their dimension we now proceed on to define the notion of subsemigroup subvector spaces of semigroup vector spaces and subsemigroup linear subalgebra of semigroup linear algebra.

**DEFINITION 2.4.8:** *Let V be a semigroup vector space over the semigroup S. Let $P \subset V$ be a proper subset of V and T a subsemigroup of S. If P is a semigroup vector space over T then we call P to be a subsemigroup subvector space over T.*

We now illustrate this situation by the following example.

***Example 2.4.41:*** Let V = $Z^+ \cup \{0\} \times Z^+ \cup \{0\}$ be a semigroup vector space over the semigroup S = $Z^+ \cup \{0\}$. Let P = $\{2Z^+ \cup \{0\} \times \{0\}\} \subseteq V$ and T = $2Z^+ \cup \{0\}$ subsemigroup under addition. P is a subsemigroup subvector space over T = $2Z^+ \cup \{0\} \subseteq S$

***Example 2.4.42:*** Let V = $(3Z^+ \cup \{0\} \times 5Z^+ \cup \{0\} \times 7Z^+ \cup \{0\} \times 2Z^+ \cup \{0\})$ be a semigroup vector space over the semigroup S = $Z^+ \cup \{0\}$.

Let P = $\{\{0\} \times 5Z^+ \cup \{0\} \times \{0\} \times \{0\}\} \subseteq V$. P is a subsemigroup subvector space over the subsemigroup T = $5Z^+ \cup \{0\} \subseteq S$.

***Example 2.4.43:*** Let V[x] = {set of all polynomials in the variable x with coefficients from $Z^+ \cup \{0\}$}. Let P [x] = {set of all polynomials of even degree with coefficients from $Z^+ \cup \{0\}\} \subseteq$ V [x]. P [x] is a subsemigroup subvector space over the subsemigroup T = $2Z^+ \cup \{0\}$.

***Example 2.4.44:*** Let V = {(0 0 0), (1 1 1), (0 0 1), (1 0 0), (1 1 0)} be a set which is a semigroup vector space on the semigroup S = {0, 1}, under addition such that 1 + 1 = 1. Let P = {(0 0 0), (1 0 0)}. P is a subsemigroup subvector space over the set T = {1}.

In fact every subset of V with (0, 0, 0) is a subsemigroup subvector space over the set T = {1} which is a subsemigroup of S under addition.

Now we proceed on to define the notion of semigroup linear subalgebra of a semigroup linear algebra.

**DEFINITION 2.4.9:** *Let V be a semigroup linear algebra over the semigroup S. Let P $\subset$ V be a proper subset of V which is a subsemigroup under '+'. Let T be a subsemigroup of S. If P is a semigroup linear algebra over the semigroup T then we call P to be a subsemigroup linear subalgebra over the subsemigroup T.*

*Note:* A subsemigroup linear subalgebra is different from the notion of semigroup linear subalgebra. Also a semigroup linear subalgebra is not a subsemigroup linear subalgebra.

We now proceed onto describe them with examples.

***Example 2.4.45:*** Let

$$V = \left\{ \begin{pmatrix} a_1 & a_2 & \cdots & a_6 \\ a_7 & a_8 & \cdots & a_{12} \end{pmatrix} \middle| \; a_i \in Z^+ \cup \{0\}; 1 \le i \le 12 \right\}$$

be a semigroup linear algebra over the semigroup $S = Z^+ \cup \{0\}$. Let

$$P = \left\{ \begin{pmatrix} a_1 & \cdots & a_6 \\ a_7 & \cdots & a_{12} \end{pmatrix} \middle| \; a_i \in 2Z^+ \cup \{0\}; 1 \le i \le 12 \right\}$$

be a semigroup under matrix addition, P a proper subset of V. T = $\{2Z^+ \cup \{0\}\}$ the proper subsemigroup of S. P is a subsemigroup linear subalgebra of V over T.

***Example 2.4.46:*** Let {(100), (000), (010), (001), (111), (110), (011), (101)} = V be a semigroup under addition modulo 2. Let $Z_2 = \{0,1\}$ be the semigroup under addition modulo 2. V is a semigroup linear algebra over $Z_2 = \{0,1\}$ we see V has no subsemigroup linear subalgebras.

Now we are interested in studying and defining such algebras.

**DEFINITION 2.4.10:** *Let V be a semigroup linear algebra over the semigroup S. If V has no subsemigroup linear subalgebras over any subsemigroup of S then we call V to be a pseudo simple semigroup linear algebra.*

In fact we have a large class of such semigroup linear algebras. Further these classes of semigroup linear algebras find their applications in coding theory.

***Example 2.4.47:*** Let V = {(00), (10), (01), (11)} be a semigroup S = $Z_2$ = {0,1} modulo addition 2. We see V is a pseudo simple semigroup linear algebra.

***Example 2.4.48:*** Let V = {(0000), (1111), (1100), (0011)} be a semigroup linear algebra over $Z_2$ = {0,1} = S. V is also a pseudo simple semigroup linear algebra. (V is a semigroup under addition modulo 2).

***Example 2.4.49:*** Let V the $2^n$ elements of the form {(000…0) (10…0), …, (011 111 … 1) (111 … 1)}. V is a semigroup linear algebra over $Z_2$ = {0,1} the semigroup under addition, ($n \geq 2$). V is a pseudo simple semigroup linear algebra. (V is a semigroup under addition modulo 2).

Now having seen a very large class of pseudo simple linear algebras we now proceed on to define the notion of pseudo simple group subvector spaces in a semigroup linear algebra which are not semigroup linear subalgebras.

**DEFINITION 2.4.11:** *Let V be a semigroup linear algebra over a semigroup S. Suppose V has a proper subset P which is only a semigroup vector space over the semigroup S and not a semigroup linear algebra then we call P to be the pseudo semigroup subvector space over S.*

We illustrate this situation by the following examples.

***Example 2.4.50:*** Let V = {all polynomials of degree n with co-efficient from $Z^+ \cup \{0\}$}. V is a semigroup linear algebra over the semigroup S = $Z^+ \cup \{0\}$. V is also a semigroup under addition of polynomials.

P = {$ax + b$, $px^3 + dx + e$, $qx^9 + rx^2 + s$ / a, b, p, d, e, q, r, s $\in Z^+ \cup \{0\}$}; P is just a set, P is a semigroup vector space over S. However P is not a semigroup linear algebra over S. So P is a pseudo semigroup subvector space of V.

***Example 2.4.51:*** Let

$$M_{2\times 2} = \left\{ \begin{pmatrix} a & b \\ c & d \end{pmatrix} \middle| \; a, b, c, d \in Z^+ \cup \{0\} \right\}$$

be a semigroup under matrix addition. $M_{2\times2}$ is a semigroup linear algebra over the semigroup $S = Z^+ \cup \{0\}$.
Take

$$P = \left\{ \begin{pmatrix} a & b \\ 0 & 0 \end{pmatrix}, \begin{pmatrix} 0 & 0 \\ d & e \end{pmatrix} \middle| \; a, b, c, d \in P \right\};$$

P is only a pseudo semigroup subvector space over S. We see P is not closed under matrix addition.

***Example 2.4.52:*** Let

$$V = \left\{ \begin{pmatrix} a_1 & a_2 & a_3 \\ a_4 & a_5 & a_6 \end{pmatrix} | \; a_i \in Z^+ \cup \{0\} \right\};$$

V is a semigroup under matrix addition. V is a semigroup linear algebra over the semigroup $S = Z^+ \cup \{0\}$. Take P

$$= \left\{ \begin{pmatrix} a_1 & 0 & a_2 \\ a_3 & a_4 & 0 \end{pmatrix}, \begin{pmatrix} 0 & b & 0 \\ 0 & c & d \end{pmatrix} \middle| b, c, d, a_1, a_2, a_3, a_4 \in 2Z^+ \cup \{0\} \right\}.$$

$P \subset V$ but P is not a semigroup under matrix addition. Thus P is only a pseudo semigroup subvector space over S.

Now we proceed on to define the new notion of pseudo subsemigroup vector subspace of a semigroup linear algebra.

**DEFINITION 2.4.12:** *Let V be a semigroup linear algebra over the semigroup S. Let P be a proper subset of V and P is not a*

*semigroup under the operations of V. Suppose T $\subseteq$ S, a proper subset of S and T is also a semigroup under the same operations of S; i.e., T a subsemigroup of S, then we call P to be a pseudo subsemigroup subvector space over T if P is a semigroup vector space over T.*

We illustrate this situation by the following examples.

***Example 2.4.53:*** Let

$$V = \left\{ \begin{pmatrix} a_1 & a_2 & a_3 & a_4 \\ a_5 & a_6 & a_7 & a_8 \end{pmatrix} \middle| a_i \in Z^+ \cup \{0\}; 1 \le i \le 8 \right\}$$

be the semigroup under matrix addition. Let S = $Z^+ \cup \{0\}$ a semigroup under addition. V is a semigroup linear algebra over S. Take P

$$\left\{ \begin{pmatrix} a_1 & 0 & 0 & 0 \\ a_2 & 0 & 0 & 0 \end{pmatrix}, \begin{pmatrix} 0 & 0 & c_1 & 0 \\ 0 & 0 & c_2 & 0 \end{pmatrix} \middle| a_1, a_2, c_1, c_2 \in Z^+ \cup \{0\} \right\} \subseteq V .$$

P is just a subset of V and P is not closed under matrix addition. Take T = $2Z^+ \cup \{0\}$, T is a subsemigroup of S. Clearly P is a semigroup vector space over T, hence P is a subsemigroup pseudo subvector space over T.

***Example 2.4.54:*** Let $V_5$ = {$(Z^+ \cup \{0\})$ [x], i.e., set of all polynomials of degree less than or equal to 5 with coefficients from the semigroup S = $Z^+ \cup \{0\}$}. $V_5$ is a semigroup linear algebra over the semigroup S. Take $P_5$ = {$ax^2 + bx + c$, $px^3 + d$, $qx^4 + e$/ a, b, c, p, d, q and e $\in Z^+ \cup \{0\}$}. Clearly $P_5$ is only a proper subset of $V_5$. $P_5$ is not closed under the polynomial addition, so $P_5$ is not a semigroup. Take T = $3Z^+ \cup \{0\} \subseteq$ S = {$Z^+ \cup \{0\}$}. T is a semigroup under addition. Thus $P_5$ is a semigroup vector space over the semigroup T. $P_5$ is the pseudo subsemigroup vector subspace of $V_5$.

An important natural question would be that will every semigroup linear algebra have a pseudo subsemigroup vector subspace. The answer is no. We prove this by the following example.

***Example 2.4.55:*** Consider the semigroup $Z_2$ [x] of all polynomials with coefficients from the field $Z_2$ under polynomial addition. $Z_2$ [x] is a semigroup linear algebra over the semigroup $Z_2$.

Take P = {all polynomials $x^3 + 1$, $x^5 + 1$,…, $x^n + 1$, $n \in Z^+$}, P is only a proper subset of $Z_2$[x]. P is not a closed set under polynomial addition. P is in fact a pseudo semigroup subvector space of $Z_2$ [x]. Now $Z_2$ has no proper subsemigroups other than the trivial {0} semigroup. So P is not a pseudo subsemigroup subvector space of $Z_2$[x].

Thus we see every semigroup linear algebra need not contain a pseudo subsemigroup vector subspace.

In fact we have a class of such semigroup linear algebras which we state in the form of theorem.

**THEOREM 2.4.1:** *Let $Z_p$ [x] be the collection of all polynomials with coefficient from the prime field $Z_p$ of characteristic p. $Z_p$[x] is a semigroup under polynomial addition. Further $Z_p$ is also a semigroup under addition modulo p. $Z_p$ [x] is a semigroup linear algebra over $Z_p$. In fact $Z_p$ [x] has no subsemigroup linear subalgebras and $Z_p$ [x] has no pseudo subsemigroup subvector spaces.*

*Proof:* Given $Z_p$ [x] is a semigroup linear algebra over the semigroup $Z_p$ = {0, 1, …, p – 1}. Clearly $Z_p$ has no subsemigroups other than {0} and itself. So $Z_p$ [x] cannot have any non trivial subsemigroup linear subalgebras or pseudo subsemigroup vector subspaces. It can have only the {0} to be both these structure over {0}.

Now in view of this we define two new algebraic structures.

**DEFINITION 2.4.13:** *Let V be a semigroup linear algebra over the semigroup S. If V has no subsemigroup linear algebras over subsemigroups of S then we call V to be a simple semigroup linear algebra.*

We have non trivial classes of simple semigroup linear algebras given by the example.

***Example 2.4.56:*** Let

$$M_{n \times m} = \left\{ \left. \begin{pmatrix} a_{11} & \cdots & a_{1m} \\ \vdots & & \vdots \\ a_{n1} & \cdots & a_{nm} \end{pmatrix} \right| a_{ij} \in Z_p ; 1 \leq i \leq n ; i \leq j \leq m \right\} .$$

This is a simple semigroup linear algebra. (m=n) can also occur. We see $M_{n \times m}$ is taken only as a semigroup under matrix addition.

***Example 2.4.57:*** Let $V = Z_p \times \ldots \times Z_p = \{(x_1, \ldots, x_n) / x_i \in Z_p / 1 \leq i \leq n\}$, V is a semigroup under addition. V is also a semigroup linear algebra which is a simple semigroup linear algebra.

In fact both these semigroup linear algebras do not contain any proper pseudo subsemigroup subvector spaces. In view of all these we can have the following theorem before which, we just recall the definition of a simple semigroup. A semigroup S is S-simple if S has no proper subsemigroups. The only trivial subsemigroups of S being {0} or $\phi$ and S itself.

**THEOREM 2.4.2:** *Let V be a semigroup. S a semigroup such that it is S-simple. If V is a semigroup linear algebra over S then V is a simple semigroup linear algebra over S.*

*Proof:* Given V is a semigroup linear algebra over the semigroup S. Also it is given the semigroup S has no proper subsemigroups i.e., {0} and S are the only subsemigroups of S which are trivial. So if $W \subseteq V$; W cannot be a subsemigroup linear subalgebra over any $T \subseteq S$, T a subsemigroup of S. Hence the claim.

Now we proceed on to define the new notion of pseudo simple semigroup linear algebra.

**DEFINITION 2.4.14:** *Let V be a semigroup under addition and S a semigroup such that V is a semigroup linear algebra over the semigroup S. If V has no proper subset P ($\subseteq$V) such that V is a pseudo subsemigroup vector subspace over a subsemigroup, T of S then we call V to be a pseudo simple semigroup linear algebra.*

We illustrate this situation by the following examples.

***Example 2.4.58:*** Let V = $Z_5$ [x] be the collection of all polynomials with coefficients from $Z_5$, $Z_5$ a semigroup under addition modulo 5. V is semigroup linear algebra over the semigroup $Z_5$. $Z_5$ has no proper subsemigroups. Hence for any subset P of V; P cannot be a pseudo subsemigroup vector subspace. Hence V is a pseudo simple semigroup linear algebra.

***Example 2.4.59:*** Let

$$V = M_{3\times5} = \left\{\left(a_{ij}\right) \middle| a_{ij} \in Z_7;\ 1 \le i \le 3; 1 \le j \le 5\right\}$$

be a semigroup under matrix addition modulo 7, with entries from $Z_7$. S = $Z_7$ be the semigroup under addition modulo 7. V is a semigroup linear algebra over $Z_7$. $Z_7$ has no proper subsemigroup. So for any subset P of V, P is not a pseudo subsemigroup vector subspace of V. So V is a pseudo simple semigroup linear algebra over $Z_7$.

We prove the following interesting theorem.

**THEOREM 2.4.3:** *Let V be a semigroup, S a S-simple semigroup i.e. S has no subsemigroups other than {0} or $\phi$ or S. V be the semigroup linear algebra over S. V is a pseudo simple semigroup linear algebra over S.*

*Proof:* Given V is a semigroup linear algebra over the semigroup S, where S is a S-simple semigroup, i.e., S has no proper subsemigroups. So for any subset P of V, P cannot be a subsemigroup algebraic structure. In particular P cannot be subsemigroup subvector space of V. So V is a pseudo simple semigroup linear algebra over S.

Now we proceed onto define the notion of linear transformation of semigroup linear algebras defined over the same semigroup S. As in case of linear algebra transformation where both the linear algebras must be defined over the same field we see in case of semigroup linear algebras to have a linear transformation both of them must be defined over the same semigroup S.

**DEFINITION 2.4.15:** *Let V and W be any two semigroup linear algebras defined over the same semigroup, S we say T from V to W is a semigroup linear transformation if* $T(c\alpha + \beta) = cT(\alpha) + T(\beta)$ *for all* $c \in S$ *and* $\alpha, \beta \in V$.

It is left, as an exercise to the reader to prove the set of all semigroup linear transformations from V to W is a semigroup linear algebra over S with composition of maps as the operation.

Now we give few examples of semigroup linear algebras defined over the same semigroup S.

***Example 2.4.60:*** Let

$$V = \left\{ \begin{pmatrix} a_1 & a_2 \\ a_3 & a_4 \end{pmatrix} \middle| a_i \in Z^+ \cup \{0\}; 1 \leq i \leq 4 \right\}$$

be the semigroup under addition of matrices. $S = Z^+ \cup \{0\}$ is a semigroup under addition. V is a semigroup linear algebra over $S = Z^+ \cup \{0\}$. Let $W = \{P \times P \,/\, P = Z^+ \cup \{0\}\}$, W is a semigroup under component-wise addition. W is a semigroup linear algebra over S. Define T from V into W by

$$T\left(\begin{pmatrix} a_1 & a_2 \\ a_3 & a_4 \end{pmatrix}\right) = \left(a_1 + a_2, a_3 + a_4\right).$$

T is a semigroup linear transformation of V to W.

***Example 2.4.61:*** Let V and W be as in example 2.4.60. Define $T_1$ from W into V by

$$T_1 (x, y) = \begin{pmatrix} x & y \\ y & x \end{pmatrix}$$

for all $(x, y) \in P \times P$.

Prove T is a semigroup linear transformation from W into V.

Now we proceed onto define a new notion of semigroup linear operators.

**DEFINITIONS 2.4.16:** *Let V be a semigroup linear algebra over the semigroup S. A map T from V to V is said to be a semigroup linear operator on V if T (cu + v) = cT (u) + T (v) for every c $\in$ S and u, v $\in$ V.*

The reader is left with the task of proving the collection of all semigroup linear operators on V is again a semigroup linear algebra over S.

We now illustrate this situation by the following examples.

***Example 2.4.62:*** Let V = {set of all 2 × 3 matrices with entries from S}, be a semigroup under matrix addition. V is a semigroup linear algebra on the semigroup $S = Z^+ \cup \{0\}$.

Define T: V → V by

$$T\left[\begin{pmatrix} a_1 & a_2 & a_3 \\ a_4 & a_5 & a_6 \end{pmatrix}\right]$$

$$= \left[ \begin{pmatrix} 0 & a_1 + a_2 & 0 \\ a_4 & 0 & a_5 + a_6 \end{pmatrix} \middle| a_1, a_2, a_4, a_5, a_6 \in Z^+ \cup \{0\} \right].$$

T is a semigroup linear operator on V.

***Example 2.4.63:*** Let V = P × P × P × P where P = $Z^+ \cup \{0\}$ be a semigroup linear algebra over the semigroup S = $Z^+ \cup \{0\}$. Define T (x, y, z, w) = (x + y, y + z, z – w, x – w) for every v = (x, y, z, w) ∈ V. T is a semigroup linear operator on V.

***Example 2.4.64:*** Let V = P × P × P be a semigroup linear algebra over the semigroup S = $2Z^+ \cup \{0\}$, where P = $Z^+ \cup \{0\}$. Define T: V → V by T (x, y, z) = (y, z, x). Prove T is a semigroup linear operator which is one to one and invertible.

***Example 2.4.65:*** Let V = {all polynomials of degree less than or equal to 7 with coefficients from the semigroup S = $Z^+ \cup \{0\}$}. V is a semigroup under polynomial addition. V is a semigroup linear algebra over S. Define T : V → V by T (x) = $x^2$, T ($x^2$) = $x^3$ , …, T ($x^6$) = x ; i.e., T($x^n$) = $x^{n+1}$ if $1 \leq n \leq 5$ and T($x^6$) = x. Is T a 1 – 1 invertible semigroup linear operator on V?.

Now we define yet another new type of semigroup linear operator on a semigroup linear algebra V over the semigroups.

**DEFINITION 2.4.17:** *Let V be a semigroup linear algebra over the semigroup S. Let W ⊆ V be a subsemigroup linear algebra over the semigroup P, P a proper subsemigroup of S. Let T : V → W be a map such that T (αv + u) = T(α) T(v) + T(u) for all u, v ∈ V and T (α) ∈ P. We call T a pseudo semigroup linear operator on V.*

We first illustrate this situation by the following example.

***Example 2.4.66:*** Let V = P × P × P where P = $Z^+ \cup \{0\}$ be a semigroup linear algebra over the semigroup S = $Z^+ \cup \{0\}$. Let W = $2Z^+ \cup \{0\} \times \{2Z^+ \cup \{0\}\} \times \{2Z^+ \cup \{0\}\}$ be a subset of V and W be a subsemigroup linear subalgebra over the subsemigroup L = $2Z^+ \cup \{0\}$. Let T: V → W be defined by T(αu + v) = 2α (2u) + 2v, T is a pseudo linear operator on V.

We call this map T to be a pseudo projection.

We just give the definition of semigroup projection of a linear algebra.

**DEFINITION 2.4.18:** *Let V be a semigroup linear algebra over the semigroup S. Let W be a semigroup linear subalgebra of V over S. Let T be a linear operator on V. T is said to be a semigroup linear projection on W if*

$$T(v) = w, w \in W$$

*and*

$$T(\alpha u + v) = \alpha T(u) + T(v)$$
$$T(v) \text{ and } T(u) \in W$$

*for all* $\alpha \in S$ *and* $u, v \in V$.

We illustrate this situation by the following example.

***Example 2.4.67:*** Let V = P × P × P × P where P = $Z^+ \cup \{0\}$, V a semigroup linear algebra over the semigroup P = $Z^+ \cup \{0\}$. Let W = $2Z^+ \cup \{0\} \times \{2Z^+ \cup \{0\}\} \times \{0\} \times \{0\} \subseteq V$ be a semigroup linear subalgebra of V over P. Define T : V → V by T (x, y, z, w) = (2x, 2y, 0, 0). Clearly T is a semigroup linear projection of V onto W.

***Example 2.4.68:*** Let

$$V = \left\{ \begin{pmatrix} a_1 & a_2 & a_3 & a_4 \\ a_5 & a_6 & a_7 & a_8 \end{pmatrix} \middle| a_i \in Z^+ \cup \{0\}; 1 \le i \le 8 \right\}$$

be a semigroup linear algebra over the semigroup $S = Z^+ \cup \{0\}$.

Let

$$W = \left\{ \begin{pmatrix} a_1 & a_2 & a_3 & a_4 \\ 0 & 0 & 0 & 0 \end{pmatrix} \middle| a_i \in 2Z^+ \cup \{0\}; 1 \le i \le 4 \right\} \subseteq V;$$

W is a semigroup linear subalgebra of V. Define $T : V \to V$ by

$$T\left(\begin{pmatrix} a_1 & a_2 & a_3 & a_4 \\ a_5 & a_6 & a_7 & a_8 \end{pmatrix}\right) = \begin{pmatrix} 2a_1 & 2a_2 & 2a_3 & 2a_4 \\ 0 & 0 & 0 & 0 \end{pmatrix}.$$

Clearly T is a semigroup linear operator on V which is a semigroup linear projection of V into W.

***Example 2.4.69:*** Let V = P[x] where $P = Z^+ \cup \{0\}$, i.e., all polynomials in the variable x with coefficients from P. V is a semigroup linear algebra over P.

Let W = {all polynomials of even degree with coefficients from P} $\subseteq$ V; W is a semigroup linear subalgebra over P. Define a map

$$T : V \to W$$

by

$$T(\alpha_1 x) = \alpha_1(x^2),$$
$$T(\alpha_2 x^2) = \alpha_2 x^4, \ldots,$$
$$T(\alpha_n x^n) = \alpha_n x^{2n}$$

; $1 \le n \le \infty$.

T is clearly a semigroup linear operator which is a semigroup linear projection of V into W.

Now having defined the notion of semigroup linear projection we proceed on to define semigroup projection of semigroup vector spaces.

**DEFINITION 2.4.19:** *Let V be a semigroup vector space over the semigroup S. Let W $\subseteq$ V be a semigroup vector subspace of V. A linear operator on V is said to be a semigroup projection of V into W if T: V $\to$ W i.e., T (v) = w for every v $\in$ V and w $\in$ W.*

We illustrate this situation by the following example.

***Example 2.4.70:*** Let

$$V = \left\{ \begin{pmatrix} a_1 & a_2 & \cdots & a_5 \\ 0 & 0 & \cdots & 0 \end{pmatrix}, \begin{pmatrix} 0 & 0 & \cdots & 0 \\ b_1 & b_2 & \cdots & b_5 \end{pmatrix} \middle| \begin{matrix} a_i, b_j \in Z^+ \cup \{0\}; \\ 1 \leq i, j \leq 5 \end{matrix} \right\}$$

be a semigroup vector space over the semigroup $S = Z^+ \cup \{0\}$. Let

$$W = \left\{ \begin{pmatrix} a_1 & a_2 & \cdots & a_5 \\ 0 & 0 & \cdots & 0 \end{pmatrix} \middle| a_i \in 2Z^+ \cup \{0\}; 1 \leq i \leq 5 \right\} \subseteq V$$

be a semigroup subvector space of V. Let T: V → V defined by

$$T = \begin{pmatrix} a_1 & \cdots & a_5 \\ 0 & \cdots & 0 \end{pmatrix} = \begin{pmatrix} a_1 & \cdots & a_5 \\ 0 & \cdots & 0 \end{pmatrix}$$

and

$$T \begin{pmatrix} 0 & 0 & \cdots & 0 \\ b_1 & b_2 & \cdots & b_5 \end{pmatrix} = \begin{pmatrix} 0 & 0 & \cdots & 0 \\ 0 & 0 & \cdots & 0 \end{pmatrix};$$

then T is a semigroup projection of V on W.

We give yet another example of the semigroup projection of the semigroup vector spaces.

***Example 2.4.71:*** Let V = {0, 1, 3, 5, 7, …, (2n + 1)} be a semigroup vector space over the semigroup S = {0, 1} where 1 + 1 = 1. Let W = {0, 3, $3^2$, …} ⊆ V, W is a semigroup vector subspace of V. Let T be a semigroup linear operator on V defined by T(x) = x if x is of the form $3^n$. T(x) = 0 otherwise. T is a semigroup linear operator on V, which is a semigroup linear projection of V on W.

Now we proceed onto define direct union of semigroup vector subspaces of a semigroup vector space.

**DEFINITION 2.4.20:** *Let V be a semigroup vector space over the semigroup S. Let $W_1, \ldots, W_n$ be semigroup vector subspaces of V if $V = \cup W_i$ and $W_i \cap W_j = \phi$ or {0}, if $i \neq j$ then we say V is the direct union of the semigroup vector subspaces of the semigroup vector space V over S.*

We illustrate this situation by some examples.

***Example 2.4.72:*** Let

$$V = \left\{ \begin{pmatrix} a_1 & a_2 & a_3 \\ 0 & 0 & 0 \end{pmatrix}, \begin{pmatrix} 0 & 0 & 0 \\ b_1 & b_2 & b_3 \end{pmatrix} \middle| \; a_i, b_j \in Z^+ \cup \{0\}; 1 \leq i \leq 3 \right\}$$

be a semigroup vector space over the semigroup $S = Z^+ \cup \{0\}$. Take

$$W_1 = \left\{ \begin{pmatrix} a_1 & a_2 & a_3 \\ 0 & 0 & 0 \end{pmatrix} \middle| \; a_i \in Z^+ \cup \{0\}; 1 \leq i \leq 3 \right\}$$

and

$$W_2 = \left\{ \begin{pmatrix} 0 & 0 & 0 \\ b_1 & b_2 & b_3 \end{pmatrix} \middle| \; b_i \in Z^+ \cup \{0\}; 1 \leq i \leq 3 \right\}$$

be semigroup vector subspaces of V over the semigroup S. Clearly $V = W_1 \cup W_2$ and

$$W_1 \cap W_2 = \begin{pmatrix} 0 & 0 & 0 \\ 0 & 0 & 0 \end{pmatrix}.$$

Thus V is the direct union of vector subspaces over the semigroups.

***Example 2.4.73:*** Let

$$V = \left\{ \begin{pmatrix} a & b \\ c & d \end{pmatrix}, \begin{pmatrix} a_1 & a_2 & a_3 \\ a_4 & a_5 & a_6 \end{pmatrix}, \begin{pmatrix} a_1 & a_4 \\ a_2 & a_5 \\ a_3 & a_6 \end{pmatrix} \middle| \begin{array}{c} a, b, c, d, a_i \in Z^+ \cup \{0\}; \\ 1 \leq i \leq 3 \end{array} \right\};$$

V is a semigroup vector space over the semigroup S = $\{Z^+ \cup \{0\}\}$. Take

$$W_1 = \left\{ \begin{pmatrix} a & b \\ c & d \end{pmatrix} \middle| a, b, c, d \in Z^+ \cup \{0\} \right\},$$

$W_1$ is a semigroup vector subspace of V.

$$W_2 = \left\{ \begin{pmatrix} a_1 & a_2 & a_3 \\ a_4 & a_5 & a_6 \end{pmatrix} \middle| a_i \in Z^+ \cup \{0\}; 1 \leq i \leq 6 \right\}$$

is a semigroup vector subspace of V over S.

$$W_3 = \left\{ \begin{pmatrix} a_1 & a_4 \\ a_2 & a_5 \\ a_3 & a_6 \end{pmatrix} \middle| a_i \in Z^+ \cup \{0\}; 1 \leq i \leq 6 \right\}$$

is a semigroup vector subspace of V over the semigroup S. Thus $V = W_1 \cup W_2 \cup W_3$ with $W_i \cap W_j = \phi$ if $i \neq j$, $1 \leq i, j \leq 3$. Hence V is a direct union of vector subspaces of the semigroup vector space V.

***Example 2.4.74:*** Let V = $\{3Z^+ \cup \{0\}, 2Z^+ \cup \{0\}, 5Z^+ \cup \{0\}, \ldots, nZ^+ \cup \{0\} / 2 \leq n \leq \infty\}$ be a semigroup vector space over the semigroup S = $Z^+ \cup \{0\}$ . Let $W_1 = 2Z^+ \cup \{0\}$, $W_2 = (3Z^+ \cup \{0\}\}$, …, $W_n = (n + 1) Z^+ \cup \{0\}$. $2 \leq n \leq \infty$ be a semigroup vector subspaces of V or S.

Clearly $V = \bigcup_{i=1}^{\infty} W_i$ but $W_i \cap W_j \neq \phi$ or {0} so V is not a direct union of semigroup vector subspaces of V.

In view of this we define yet another new notion called pseudo direct union of semigroup vector subspaces of a semigroup vector space.

**DEFINITION 2.4.21:** *Let V be a semigroup vector space over the semigroup S. Let $W_1, \ldots, W_n$ be a semigroup subvector spaces of V over the semigroup S. If $V = \bigcup_{i=1}^{n} W_i$ but $W_i \cap W_j \neq \phi$ or {0} if i $\neq$ j then we call V to be the pseudo direct union of semigroup vector spaces of V over the semigroup S.*

We illustrate this situation by the following example.

***Example 2.4.75:*** Let

$$V = \left\{ \begin{pmatrix} a_1 & a_2 \\ a_3 & a_4 \end{pmatrix}, \begin{pmatrix} b_1 & b_2 \\ b_3 & b_4 \end{pmatrix}, \begin{pmatrix} c_1 & c_2 \\ c_3 & c_4 \end{pmatrix}, \right.$$

$$\begin{pmatrix} d_1 & d_2 \\ d_3 & d_4 \end{pmatrix}, \begin{pmatrix} y_1 & y_2 \\ y_3 & y_4 \end{pmatrix} \text{and} \begin{pmatrix} x_1 & x_2 \\ x_3 & x_4 \end{pmatrix}$$

where $a_1, a_2, a_3, a_4 \in 2Z^+ \cup \{0\}$, $b_1, b_2, b_3, b_4 \in 3Z^+ \cup \{0\}$, $c_1, c_2, c_3, c_4 \in 5Z^+ \cup \{0\}$, $d_1, d_2, d_3, d_4 \in 7Z^+ \cup \{0\}$, $y_1, y_2, y_3, y_4 \in 11Z^+ \cup \{0\}$ and $x_1, x_2, x_3, x_4 \in 19Z^+ \cup \{0\}\}$ be the semigroup vector space over the semigroup $S = Z^+ \cup \{0\}$.

Let

$$W_1 = \left\{ \begin{pmatrix} a_1 & a_2 \\ a_3 & a_4 \end{pmatrix} \middle| a_i \in 2Z^+ \cup \{0\}; 1 \leq i \leq 4 \right\}$$

$$W_2 = \left\{ \begin{pmatrix} b_1 & b_2 \\ b_3 & b_4 \end{pmatrix} \middle| b_i \in 3Z^+ \cup \{0\}; 1 \le i \le 4 \right\}$$

$$W_3 = \left\{ \begin{pmatrix} c_1 & c_2 \\ c_3 & c_4 \end{pmatrix} \middle| c_j \in 5Z^+ \cup \{0\}; 1 \le j \le 4 \right\}$$

$$W_4 = \left\{ \begin{pmatrix} d_1 & d_2 \\ d_3 & d_4 \end{pmatrix} \middle| d_j \in 7Z^+ \cup \{0\}; 1 \le j \le 4 \right\}$$

$$W_5 = \left\{ \begin{pmatrix} y_1 & y_2 \\ y_3 & y_4 \end{pmatrix} \middle| y_k \in 11Z^+ \cup \{0\} \right\}$$

and

$$W_6 = \left\{ \begin{pmatrix} x_1 & x_2 \\ x_3 & x_4 \end{pmatrix} \middle| x_i \in 19Z^+ \cup \{0\}; 1 \le i \le 4 \right\}$$

be semigroup vector subspaces of V over the semigroup $S = Z^+ \cup \{0\}$.

Clearly $V = \bigcup_{i=1}^{6} W_i$ we see $W_i \cap W_j \neq \phi$ or $\{0\}$, $i \neq j$, $1 \le i, j \le 6$. So V is the pseudo direct union of semigroup subvector spaces over S.

Now we proceed onto define the new notion of direct sum of semigroup linear subalgebras of a semigroup linear algebra over a semigroup S.

**DEFINITION 2.4.22:** *Let V be a semigroup linear algebra over the semigroup S. We say V is a direct sum of semigroup linear subalgebras $W_1, \ldots, W_n$ of V if*

1. $V = W_1 + \ldots + W_n$
2. $Wi \cap W_j = \{0\}$ *or* $\phi$ *if* $i \neq j$ $(1 \le i, j \le n)$.

We first illustrate this situation by the following example.

***Example 2.4.76:*** Let V = P × P × P × P be a semigroup linear algebra over P = {0} ∪ $Z^+$. Let $W_1$ = P × {0} × {0} × {0}, $W_2$ = {0} × P × {0} × {0}, $W_3$ = {0} × {0} × P × {0} and $W_4$ = {0} × {0} × {0} × P be the semigroup linear subalgebras of V. Clearly V = $W_1 + W_2 + W_3 + W_4$ and $W_i \cap W_j = \{0\}$ if i ≠ j.

We see this way of representation in general is not unique. For if we take $W_1'$ = {0} × P × {0} × {P} and $W_2'$ = {P} × {0} × P × {0} we get V = $W_1' + W_2'$ and $W_1' \cap W_2' = \{0\}$ thus V is also a direct sum of $W_1'$ and $W_2'$. Thus the direct sum in general is not unique.

We give yet another example.

***Example 2.4.77:*** Let

$$V = \left\{ \begin{pmatrix} a_1 & a_2 & a_3 \\ a_4 & a_5 & a_6 \\ a_7 & a_8 & a_9 \end{pmatrix} \middle| a_i \in Z^+ \cup \{0\} \text{ and } 1 \leq i \leq 9 \right\}$$

be a semigroup linear algebra over the semigroup S = $Z^+ \cup \{0\}$. Let

$$W_1 = \left\{ \begin{pmatrix} a_1 & 0 & 0 \\ 0 & a_5 & 0 \\ 0 & 0 & a_9 \end{pmatrix} \middle| a_i, a_5 \text{ and } a_9 \in Z^+ \cup \{0\} \right\} \subseteq V$$

be a semigroup linear subalgebra of V. Take

$$W_2 = \left\{ \begin{pmatrix} 0 & a_2 & a_3 \\ 0 & 0 & a_6 \\ 0 & 0 & 0 \end{pmatrix} \middle| a_2, a_3, a_6 \in Z^+ \cup \{0\} \right\} \subseteq V,$$

$W_2$ is also a semigroup linear subalgebra of V. Suppose

$$W_3 = \left\{ \begin{pmatrix} 0 & 0 & 0 \\ a_4 & 0 & 0 \\ a_7 & a_8 & 0 \end{pmatrix} \middle| a_4, a_7, a_8 \in Z^+ \cup \{0\} \right\} \subseteq V;$$

$W_3$ a semigroup linear subalgebra of V. Then we see $V = W_1 + W_2 + W_3$ with

$$W_i \cap W_j = \begin{pmatrix} 0 & 0 & 0 \\ 0 & 0 & 0 \\ 0 & 0 & 0 \end{pmatrix} \text{ if } i \neq j.$$

Thus V is a direct sum of $W_1$, $W_2$, $W_3$ of V. We see this is not the only way of representing V. For take

$$P_1 = \left\{ \begin{pmatrix} a_1 & 0 & 0 \\ 0 & 0 & 0 \\ 0 & 0 & a_9 \end{pmatrix} \middle| a_1, a_9 \in Z^+ \cup \{0\} \subseteq V \right\}$$

is a semigroup linear subalgebra of V.

$$P_2 = \left\{ \begin{pmatrix} 0 & a_2 & 0 \\ a_4 & 0 & a_6 \\ 0 & 0 & 0 \end{pmatrix} \middle| a_2, a_4, a_6 \in Z^+ \cup \{0\} \right\}$$

is a semigroup linear subalgebra of V over $S = Z^+ \cup \{0\}$.
Take

$$P_3 = \left\{ \begin{pmatrix} 0 & 0 & a_3 \\ 0 & a_5 & 0 \\ 0 & 0 & 0 \end{pmatrix} \middle| a_3, a_5 \in Z^+ \cup \{0\} \right\} \subseteq V$$

is a semigroup linear subalgebra of V.

$$P_4 = \left\{ \begin{pmatrix} 0 & 0 & 0 \\ 0 & 0 & 0 \\ a_7 & a_8 & 0 \end{pmatrix} \middle| \, a_7, a_8 \in Z^+ \cup \{0\} \right\} \subseteq V$$

is a semigroup linear subalgebra of V. We see $V = P_1 + P_2 + P_3 + P_4$ with

$$P_i \cap P_j = \begin{pmatrix} 0 & 0 & 0 \\ 0 & 0 & 0 \\ 0 & 0 & 0 \end{pmatrix} \; i \neq j \; 1 \leq i, j \leq 4 .$$

Thus V is a direct sum of semigroup linear subalgebras over the semigroup $Z^+ \cup \{0\}$.

Thus we see there exists more than one way of writing the semigroup linear algebra as the direct sum of semigroup linear subalgebras.

A semigroup linear algebra is said to be strongly simple if it cannot be written as a direct sum of semigroup linear subalgebras and has no proper semigroup linear subalgebra. Clearly the class of semigroup linear algebras $V = Z_p = \{0, 1, \ldots, p - 1\}$; p a prime over $S = Z_p = V$ are strong simple for in the first place they do have any semigroup linear subalgebras and it cannot be written as direct sum . All simple semigroup linear algebras are strongly simple however it is left as an open problem for the reader to find whether strongly simple implies simple.

## 2.5 Group Linear Algebras

Next we proceed onto define yet another new special class of linear algebras called group linear algebras and their generalizations group vector spaces. In this section we also enumerate a few of its properties.

**DEFINITION 2.5.1:** *Let V be a set with zero, which is non empty. Let G be a group under addition. We call V to be a group vector space over G if the following condition are true.*

1. *For every $v \in V$ and $g \in G$ gv and vg $\in V$.*
2. *$0.v = 0$ for every $v \in V$, 0 the additive identify of G.*

We illustrate this by the following examples.

***Example 2.5.1:*** Let V = {0, 1, 2, …, 15} integers modulo 15. G = {0, 5, 10} group under addition modulo 15. Clearly V is a group vector space over G, for $gv \equiv v_1 \pmod{15}$, for $g \in G$ and $v, v_1 \in V$.

***Example 2.5.2:*** Let V = {0, 2, 4, …, 10} integers 12. Take G = {0, 6}, G is a group under addition modulo 12. V is a group vector space over G, for $gv \equiv v_1 \pmod{12}$ for $g \in G$ and $v, v_1 \in V$.

***Example 2.5.3:*** Let

$$M_{2\times 3} = \left\{ \begin{pmatrix} a_1 & a_2 & a_3 \\ a_4 & a_5 & a_6 \end{pmatrix} \middle| a_i \in \{-\infty, \ldots, -4, -2, 0, 2, 4, \ldots, \infty\} \right\}.$$

Take G = Z be the group under addition. $M_{2\times 3}$ is a group vector space over G = Z.

***Example 2.5.4:*** Let $V = Z \times Z \times Z = \{(a, b, c) / a, b, c \in Z\}$. V is a group vector space over Z.

***Example 2.5.5:*** Let V = {0, 1} be the set. Take G = {0, 1} the group under addition modulo two. V is a group vector space over G.

***Example 2.5.6:*** Let

$$V = \left\{ \begin{pmatrix} 0 & 1 \\ 0 & 0 \end{pmatrix}, \begin{pmatrix} 1 & 1 \\ 0 & 0 \end{pmatrix}, \begin{pmatrix} 1 & 0 \\ 1 & 0 \end{pmatrix}, \begin{pmatrix} 0 & 1 \\ 1 & 0 \end{pmatrix}, \begin{pmatrix} 1 & 0 \\ 0 & 1 \end{pmatrix}, \begin{pmatrix} 0 & 0 \\ 1 & 0 \end{pmatrix}, \right.$$
$$\left. \begin{pmatrix} 0 & 0 \\ 0 & 1 \end{pmatrix}, \begin{pmatrix} 0 & 0 \\ 0 & 0 \end{pmatrix} \right\}$$

be set. Take G = {0, 1} group under addition modulo 2. V is a group vector space over G.

***Example 2.5.7:*** Let

$$V = \left\{ \begin{pmatrix} a_1 & a_2 & \ldots & a_n \\ 0 & 0 & \ldots & 0 \end{pmatrix}, \begin{pmatrix} 0 & 0 & \ldots & 0 \\ 0 & 0 & \ldots & 0 \end{pmatrix}, \right.$$
$$\left. \begin{pmatrix} 0 & 0 & \ldots & 0 \\ b_1 & b_2 & \ldots & b_n \end{pmatrix} \middle| \; a_i, b_i \in Z; 1 \leq i \leq n \right\}$$

be the non empty set. Take G = Z the group of integers under addition. V is the group vector space over Z.

***Example 2.5.8:*** Let

$$V = \left\{ \begin{pmatrix} a_1 & 0 & \ldots & 0 \\ a_2 & 0 & \ldots & 0 \end{pmatrix}, \begin{pmatrix} 0 & 0 & \ldots & 0 \\ 0 & 0 & \ldots & 0 \end{pmatrix}, \begin{pmatrix} 0 & b_1 & 0 & \ldots & 0 \\ 0 & b_2 & 0 & \ldots & 0 \end{pmatrix}, \right.$$
$$\left. \ldots, \begin{pmatrix} 0 & 0 & \ldots & t_1 \\ 0 & 0 & \ldots & t_2 \end{pmatrix} \middle| \; a_i, b_i, \ldots, t_i \in Z; \; 1 \leq i \leq 2 \right\}$$

be the set of 2 × n matrices of this special form. Let G = Z be the group of integers under addition. V is a group vector space over Z.

***Example 2.5.9:*** Let

$$V = \left\{ \begin{pmatrix} a_1 & 0 \\ 0 & 0 \end{pmatrix}, \begin{pmatrix} 0 & a_2 \\ 0 & 0 \end{pmatrix}, \begin{pmatrix} 0 & 0 \\ a_3 & 0 \end{pmatrix}, \right.$$

$$\begin{pmatrix} 0 & 0 \\ 0 & 0 \end{pmatrix}, \begin{pmatrix} 0 & 0 \\ 0 & a_4 \end{pmatrix} \Bigg| a_1, a_2, a_3, a_4 \in Z \Bigg\}$$

be the set. Z = G the group of integers V is a group vector space over Z.

Now having seen examples of group vector spaces which are only set defined over an additive group.

***Example 2.5.10:*** Let V = {(0 1 0 0), (1 1 1), (0 0 0), (0 0 0 0), (1 1 0 0), (0 0 0 0 0), (1 1 0 0 1), (1 0 1 1 0)} be the set. Take $Z_2$ = G = {0, 1} group under addition modulo 2. V is a group vector space over $Z_2$.

***Example 2.5.11:*** Let

$$V = \left\{ \begin{pmatrix} a_1 & a_2 & a_2 \\ 0 & 0 & 0 \\ 0 & 0 & 0 \end{pmatrix}, \begin{pmatrix} b_1 & 0 & 0 \\ b_2 & 0 & 0 \\ b_3 & 0 & 0 \end{pmatrix}, \begin{pmatrix} 0 & c_1 & 0 \\ 0 & c_2 & 0 \\ 0 & c_3 & 0 \end{pmatrix}, \begin{pmatrix} 0 & 0 & a'_1 \\ 0 & 0 & a'_2 \end{pmatrix}, \right.$$

$$\left. \begin{pmatrix} 0 & 0 & 0 \\ 0 & 0 & 0 \\ 0 & 0 & 0 \end{pmatrix}, \begin{pmatrix} 0 & 0 & 0 \\ 0 & 0 & 0 \end{pmatrix} \right| a_i\, b_i\, c_i \in Z;\ a'_1, a'_2 \in Z; 1 \leq i \leq 3 \Bigg\}$$

be the set, Z = G the group under addition. V is just a set but V is a group vector space over Z.

It is important and interesting to note that this group vector spaces will be finding their applications in coding theory.

Now we proceed onto define the notion of substructures of group vector spaces.

**DEFINITION 2.5.2:** *Let V be the set which is a group vector space over the group G. Let P $\subseteq$ V be a proper subset of V. We say P is a group vector subspace of V if P is itself a group vector space over G.*

***Example 2.5.12:*** Let

$$V = \left\{ \begin{pmatrix} 0 & 0 \\ 0 & 0 \end{pmatrix}, \begin{pmatrix} a_1 & a_2 \\ 0 & 0 \end{pmatrix}, \begin{pmatrix} a_1 & 0 \\ 0 & 0 \end{pmatrix}, \begin{pmatrix} 0 & a_2 \\ 0 & 0 \end{pmatrix}, \begin{pmatrix} 0 & 0 \\ a_1 & a_2 \end{pmatrix} \middle| a_1, a_2 \in Z \right\}$$

be the set. V is a group vector space over the group G = Z the group of integers under addition. Take

$$P = \left\{ \begin{pmatrix} 0 & 0 \\ 0 & 0 \end{pmatrix}, \begin{pmatrix} a_1 & a_2 \\ 0 & 0 \end{pmatrix}, \begin{pmatrix} 0 & 0 \\ a_1 & a_2 \end{pmatrix} \middle| a_1, a_2 \in Z \right\} \subseteq V \ .$$

P is a group vector subspace of V over Z. It is important and interesting to note that every proper subset of V need not be a group vector subspace of V. Take

$$T = \left( \begin{pmatrix} a_1 & a_2 \\ 0 & 0 \end{pmatrix}, \begin{pmatrix} 0 & 0 \\ a_1 & 0 \end{pmatrix} \right) \subseteq V \, .$$

T is not a group vector subspace of V it is only a set and has no additional properties.

***Example 2.5.13:*** Let V = {(1 1 0 0 1), (0 0 0 0 0), (1 0 0 1 0), (0 0 0), (1 1 1), (1 1 1 1), (0 0 0 0), (1 1 0 0), (1 0 0 1)} be a proper set. Take G = {0, 1} be a group under addition modulo 2. V is a group vector space over G. P = {(1 1 0 01), (0 0 0 0 0)} $\subseteq$ V; P is a group vector subspace of V over G.

$P_1$ = {(0 0 0), (1 1 1)} $\subset$ V,

$P_1$ is also a group vector subspace of V over G.

T = {(1 1 1), (1 1 0 0)} $\subset$ V,

T is not a group vector subspace of V over G.

***Example 2.5.14:*** Let

$$V = \left\{ \begin{pmatrix} a_1 & a_2 \\ a_3 & a_4 \end{pmatrix} \middle| a_i \in Z; 1 \le i \le 4 \right\}$$

be a group vector space over the group G = Z. Let

$$P = \left\{ \begin{pmatrix} b_1 & b_2 \\ b_3 & b_4 \end{pmatrix} \middle| b_i \in 2Z;\ 1 \leq i \leq 4 \right\} \subseteq V .$$

P is a group vector subspace over G.

We now define the notion of linearly independent subset of a group vector space.

**DEFINITION 2.5.3:** *Let V be a group vector space over the group G. We say a proper subset P of V to be a linearly dependent subset of V if for any $p_1, p_2 \in P$, $(p_1 \neq p_2)$ $p_1 = ap_2$ or $p_2 = a'p_1$ for some $a, a' \in G$. If for no distinct pair of elements $p_1, p_2 \in P$ we have $a, a_1 \in G$ such that $p_1 = ap_2$ or $p_2 = a_1p_1$ then we say the set P is a linearly independent set.*

We now illustrate this situation by some examples.

***Example 2.5.15:*** Let

$$V = \left\{ \begin{pmatrix} a_1 & a_2 \\ 0 & 0 \end{pmatrix}, \begin{pmatrix} 0 & 0 \\ a_1 & a_2 \end{pmatrix}, \begin{pmatrix} 0 & 0 \\ 0 & 0 \end{pmatrix} \middle| a_1, a_2 \in Z \right\}$$

be the group vector space over the group of integers Z. Take

$$P = \left\{ \begin{pmatrix} 2 & 4 \\ 0 & 0 \end{pmatrix}, \begin{pmatrix} 1 & 2 \\ 0 & 0 \end{pmatrix}, \begin{pmatrix} 6 & 12 \\ 0 & 0 \end{pmatrix}, \begin{pmatrix} 5 & 10 \\ 0 & 0 \end{pmatrix} \right\} \subseteq V ;$$

P is a linearly dependent subset in V over Z. Take

$$T = \left\{ \begin{pmatrix} 1 & 1 \\ 0 & 0 \end{pmatrix}, \begin{pmatrix} 0 & 0 \\ 1 & 1 \end{pmatrix} \right\}.$$

T is a linearly independent subset of V over Z.

$$T_1 = \left\{ \begin{pmatrix} 1 & 1 \\ 0 & 0 \end{pmatrix}, \begin{pmatrix} 4 & 6 \\ 0 & 0 \end{pmatrix}, \begin{pmatrix} 2 & 3 \\ 0 & 0 \end{pmatrix} \right\} \subseteq V .$$

$T_1$ is a linearly dependent set over Z.

$$T_2 = \left\{ \begin{pmatrix} 1 & 1 \\ 0 & 0 \end{pmatrix}, \begin{pmatrix} 4 & 6 \\ 0 & 0 \end{pmatrix} \right\} \subseteq V ,$$

T is a linearly independent set over Z.

An observation which is important and interesting is that both T and $T_2$ are linearly independent subsets of V but both of them are distinctly different in their behaviour. To this end we proceed onto define the notion of a generating subset of a group vector space V over the group G.

**DEFINITION 2.5.4:** *Let V be a group vector space over the group G. Suppose T is a subset of V which is linearly independent and if T generates V i.e., using t $\in$ T and g $\in$ V we get every v $\in$ V as v = gt for some g $\in$ G then we call T to be the generating subset of V over G. The number of elements in V gives the dimension of V. If T is of finite cardinality V is said to be finite dimensional otherwise V is said to be of infinite dimension.*

We illustrate this situation by the following example.

***Example 2.5.16:*** Let V = P = {(1 1 0 0), (0 0 0 0), (0 0 0 1), (1 1 1), (0 1 1), (0 1 0), (0 0 0)} be the given set. V is a group vector space over the group G = $Z_2$ = {0, 1} addition modulo 2. Take T = {(1 1 0 0), (0 0 0 1), (1 1 1), (0 1 1), (0 1 0) $\subseteq$ V; V is linearly independent set and dimension of V is 5 as T generates V.

It is important to note that no proper subset of T will generate V. Thus T is the only generating set of V and dimension of V is 5.

***Example 2.5.17:*** Let

$$V = \left\{ \begin{pmatrix} 0 & 0 \\ a_1 & a_1 \end{pmatrix}, \begin{pmatrix} a_1 & a_1 \\ 0 & 0 \end{pmatrix} \middle| a_1 \in Z \right\}$$

V is a group vector space over the group Z.

$$T = \left\{ \begin{pmatrix} 1 & 1 \\ 0 & 0 \end{pmatrix}, \begin{pmatrix} 0 & 0 \\ 1 & 1 \end{pmatrix} \right\},$$

is the generating set of V and no other set can generate V. Thus the group vector space V is of dimension two over Z. Clearly T is a linearly independent set.

***Example 2.5.18:*** Let

$$V = \left\{ \begin{pmatrix} a_1 & a_2 \\ 0 & 0 \end{pmatrix}, \begin{pmatrix} 0 & 0 \\ a_1 & a_2 \end{pmatrix} \middle| a_1, a_2 \in Z \right\} .$$

V is a group vector space over the group Z. Take

$$T = \left\{ \begin{pmatrix} 1 & 0 \\ 0 & 0 \end{pmatrix}, \begin{pmatrix} 0 & 1 \\ 0 & 0 \end{pmatrix}, \begin{pmatrix} 0 & 0 \\ 1 & 0 \end{pmatrix}, \begin{pmatrix} 0 & 0 \\ 0 & 1 \end{pmatrix} \right\} \subseteq V;$$

T is a linearly independent subset of V but T is not a generating subset of V. Take

$$T_1 = \left\{ \begin{pmatrix} 1 & 1 \\ 0 & 0 \end{pmatrix}, \begin{pmatrix} 0 & 0 \\ 1 & 1 \end{pmatrix}, \begin{pmatrix} 1 & 0 \\ 0 & 0 \end{pmatrix}, \begin{pmatrix} 0 & 1 \\ 0 & 0 \end{pmatrix}, \begin{pmatrix} 0 & 0 \\ 1 & 0 \end{pmatrix}, \begin{pmatrix} 0 & 0 \\ 0 & 1 \end{pmatrix} \right\} \subseteq V$$

$T_1$ is a linearly independent subset of V but T is not the generating subset of V over Z. In fact V cannot be generated over Z by any finite subset of V. Thus dimension of V over Z is infinite.

Take

$$P = \left\{ \begin{pmatrix} 0 & 0 \\ 1 & 2 \end{pmatrix}, \begin{pmatrix} 0 & 0 \\ 5 & 7 \end{pmatrix}, \begin{pmatrix} 1 & 7 \\ 0 & 0 \end{pmatrix}, \begin{pmatrix} 2 & 5 \\ 0 & 0 \end{pmatrix}, \right.$$

$$\left. \begin{pmatrix} 2 & 2 \\ 0 & 0 \end{pmatrix}, \begin{pmatrix} 7 & 0 \\ 0 & 0 \end{pmatrix}, \begin{pmatrix} 5 & 0 \\ 0 & 0 \end{pmatrix} \right\} \subseteq V.$$

P is a linearly independent subset of V but not a generating subset of V over Z.

Take

$$S = \left\{ \begin{pmatrix} 0 & 0 \\ 1 & 2 \end{pmatrix}, \begin{pmatrix} 0 & 0 \\ 3 & 6 \end{pmatrix}, \begin{pmatrix} 7 & 0 \\ 0 & 0 \end{pmatrix}, \begin{pmatrix} 1 & 0 \\ 0 & 0 \end{pmatrix} \right\} \subseteq V$$

S is not a linearly independent subset of V over Z.

***Example 2.5.19:*** Let V = Z × Z × Z = $\{(x_1\ x_2\ x_3) \mid x_i \in Z; 1 \le i \le 3\}$, V is a group vector space over Z. V is of infinite dimension over Z. Take T = {(1 1 0), (1 1 0), (0 0 1), (1 0 0), (0 1 0), (0 1 1)} ⊆ V is a linearly independent subset of V but T cannot generate V over Z. Take $T_1$ = {(1 1 1), (5 7 8), (7 8 1), (0 0 1)} ⊆ V. $T_1$ is again a linearly independent subset of V but not a generating subset of V over Z. Take W = Z × {0} × {0} ⊂ V to be the group vector subspace of V over Z. If T = {(1 0 0)} ⊂ W is the generating subset of W over Z and dimension of W over Z is 1. Suppose U = Z × Z × {0} ⊆ V; U is a group vector subspace of V over Z.

$T_1$ = {(1 1 0), (0 1 0), (1 0 0)} ⊂ V is a linearly independent subset of U but $T_1$ cannot generate U. In fact no finite subset of U can generate U. Thus the group vector subspace U of V is of infinite dimension over Z. Thus the group vector space V over Z has both group vector subspaces of finite and infinite dimension over Z.

***Example 2.5.20:*** Let V = {(a a a a) | a ∈ Z} be a group vector space over the group Z. Take T = {(1 1 1 1)} ⊆ V. T is the generating subset of V. In fact dimension of V over Z is one.

Further if we take W = $\{\langle(5\ 5\ 5\ 5)\rangle\} \subseteq$ V. W is a proper subset of V and W is a proper group vector subspace of V generated by the set $\langle(5\ 5\ 5\ 5)\rangle$ and dimension of W is also one. Thus V is a group vector space over Z of dimension one. W = $\{(x, x, x, x) \mid x \in 5Z\} \subseteq$ V is also of dimension one over Z but W is a proper group vector subspace of V.

It is still interesting to note that V has infinite number of proper group vector subspaces of dimension one. Take S = $\{(x\ x\ x\ x), (y\ y\ y\ y) \mid x \in 2Z \text{ and } y \in 3Z\} \subseteq$ V. S is a subset of V, V is of dimension one over Z. But S is a proper group vector subspace of V over Z and dimension of S over Z is two. The generating proper subset of S which generates S is given by $T_1 = \{(2\ 2\ 2\ 2), (3\ 3\ 3\ 3)\} \subseteq$ S. $T_1$ is a linearly independent subset of S and generates S over Z.

Thus it is still interesting and important to note that a one dimensional group vector space over the group has proper group vector subspaces of dimension greater than one. This sort of situations can occur only in case of group vector spaces.

This looks as if one cannot algebraically comprehend but concrete examples confirm the statement and establish it. In fact this one dimensional group vector space has proper group vector subspaces of infinite dimension also.

For take $S_1 = \{(x_n, x_n, x_n, x_n) \mid x_n \text{ a prime}\} \subseteq$ V. Thus

$$S_1 = \{\langle(2\ 2\ 2\ 2)\rangle, \langle(3\ 3\ 3\ 3)\rangle, \langle(5\ 5\ 5\ 5)\rangle, \ldots\} \subseteq V.$$

Thus $S_1$ is generated by

T = $\{(2\ 2\ 2\ 2), (3\ 3\ 3\ 3), (5\ 5\ 5\ 5), \ldots, (p\ p\ p\ p), \ldots\ p, \text{ a prime}\}$.

Clearly, cardinality of T is infinite. Thus V = $\{(x\ x\ x\ x) \mid x \in Z\}$ is of dimension one as it is generated by $\{(1\ 1\ 1\ 1)\}$ but it has a proper group vector subspace $S_1$ which is of infinite dimension as number of primes in Z is infinite.

Now having seen such types of group vector spaces we proceed onto give more examples of infinite dimensional group vector spaces.

***Example 2.5.21:*** Let V = {R × R × R, R reals} be a group vector space over the group Z. V is infinite dimensional. Let W = {(x x x) | x ∈ Z} ⊆ V be a proper subset of V, W is a finite dimensional group vector subspace of V. In fact T = {(1 1 1)} is the generator of W and W is of dimension one over Z.

Thus we see an infinite dimension group vector space can have as group vector subspaces of dimension one. This is a case just opposite to the case given in the earlier example where a one dimensional group vector space can have infinite dimensional group vector subspaces.

Now we proceed onto define yet another new type of substructures in a group vector space called the subgroup vector subspaces and illustrate them with examples.

**DEFINITION 2.5.5:** *Let V be a group vector space over the group G. Let W ⊆ V be a proper subset of V. H ⊂ G be a proper subgroup of G. If W is a group vector space over H and not over G then we call W to be a subgroup vector subspace of V.*

***Example 2.5.22:*** Let V = $Z_6 \times Z_6 \times Z_6$ be a group vector space over $Z_6$. W = {(2 2 2), (0 0 0), (1 1 1), (4 4 4)} ⊆ V. W is a subgroup vector subspace over the subgroup {0, 2, 4} = H ⊆ $Z_6$. Clearly W is not a vector subspace over $Z_6$ as 3 (1 1 1) = (3 3 3) ∉ W.

We give yet another example of subgroup vector subspace over a subgroup of the group over which it is defined.

***Example 2.5.23:*** Let V = $Z_{12} \times Z_{12} \times Z_{12}$ be a group vector space over the group G = $Z_{12}$. Take H = {0, 6}. Let W = {(1 1 1), (2 2 2), (6 6 6), (0 0 0), (3 3 3), (4 4 4)} ⊆ V. W is a subgroup vector subspace over the subgroup H. Clearly W is not a group vector subspace over $Z_{12}$.

Now it may so happen that a subset W may be a group vector subspace as well as subgroup vector subspace. We call in this situation W to be a duo subgroup vector subspace.

**DEFINITION 2.5.6:** *Let V be a group vector space over the group G. Let W $\subseteq$ V. If W is a subgroup vector subspace over a proper subgroup, H of G as well as W is a group vector subspace of V over G then we call W to be the duo subgroup vector subspace of V.*

We illustrate this with some examples before we proceed onto describe a few of its properties.

***Example 2.5.24:*** Let V = $Z_{12} \times Z_{12} \times Z_{12}$ be the group vector space over the group G = {0, 2, 4, 6, 8, 10}. Let W = $\{0\} \times Z_{12} \times \{0\} \subseteq$ V. W is a group vector subspace of V. Clearly W is also a subgroup vector subspace over H = {0, 6} a subgroup of G. Suppose S = {(0 0 0), (1 1 1), (6 6 6)} $\subset$ V. S is a subgroup vector subspace over H = {0, 6}. Clearly S is not a group vector subspace over G.

In view of this we prove the following theorem.

**THEOREM 2.5.1:** *Let V be a group vector space over the group G. If W is a group vector subspace of V then W need not be a subgroup vector subspace of V for some subgroup H of G.*

*Proof:* We illustrate this situation by examples. Let V = $Z_{12} \times Z_{12} \times Z_{12} \times Z_{12}$ be a group vector space over the group G = $Z_{12}$. W = $Z_{12} \times \{0\} \times \{0\} \times Z_{12} \subset$ V, W is a group vector subspace of V. W is also a subgroup vector subspace of V for every subgroup H of G. V = $Z_{11} \times Z_{11}$ is a group vector space over the group $Z_{11}$. W = $Z_{11} \times \{0\}$, W is only a group vector subspace of V and not a subgroup vector subspace of V as $Z_{11}$ has no proper subgroups. Thus V has no subgroup vector subspaces.

Conversely we have the following theorem.

**THEOREM 2.5.2:** *Let V be a group vector space over a group G. Suppose S $\subset$ V is a subgroup vector subspace of V then S need not in general be a group vector subspace of V.*

*Proof:* We prove this theorem only by a counter example. Let V = $Z_{10} \times Z_{10} \times Z_{10} \times Z_{10}$ be a group vector space over the group $Z_{10}$. Take S = {(1 1 1 1), (0 0 0 0), (5 5 5 5)} $\subseteq$ V. S is a subgroup vector subspace of V over the subgroup H = {0, 5} $\subset$ $Z_{10}$. Clearly S is not a group vector subspace of V over $Z_{10}$ as for 3 $\in$ $Z_{10}$, 3 (1 1 1 1) = (3 3 3 3) $\notin$ S. Hence the claim.

Thus a subgroup vector subspace of a group vector space V in general need not be a group vector subspace of V.

**THEOREM 2.5.3:** *Let V be a group vector space over the group G if W $\subset$ V is a duo subgroup vector subspace of V then W is both a group vector subspace of V as well as W is a subgroup vector subspace of V.*

*Proof:* The proof follows from the very definition of duo subgroup vector subspaces.

It may so happen we may find group vector spaces which has no subgroup vector subspaces over a proper subgroup. We define them in the following.

**DEFINITION 2.5.7:** *Let V be a group vector space over the group G. Suppose V has no subgroup vector subspaces then we call V to be a simple group vector space.*

We first illustrate this situation by the following examples.

***Example 2.5.25:*** Let V = $Z_7 \times Z_7 \times Z_7$ be a group vector space over the group $Z_7$. Since $Z_7$ has no proper subgroups under addition; V cannot have any subgroup vector subspaces. Thus V is a simple group vector space over $Z_7$.

***Example 2.5.26:*** Let V = $Z_5 \times Z_5$ be a group vector space over the group G = $Z_5$. V is a simple group vector space over $Z_5$.

***Example 2.5.27:*** Let V = {(1 1 1 1), (0 0 0 0), (1 0 1 1 0), (0 0 0 0 0), (1 1 0 0 1), (1 1 1), (0 0 0), (1 0 0), (0 1 1)}; V is a group vector space over the group $Z_2$ = {0, 1}. Clearly V has no proper subset which can be subgroup vector subspace of V i.e., V is a simple group vector space.

In view of the above examples we have the following theorem.

**THEOREM 2.5.4:** *Let V be a group vector space over a group G which has no proper subgroups then V is a simple group vector space over G.*

*Proof:* Obvious from the fact that the group G has no proper subgroup for a proper subset W to be a subgroup vector subspace; we need a proper subgroup in G over which W is a group vector space.

If G has no proper subgroup the existence of subgroup vector subspace is impossible.

Now we show we have a large class of simple group vector spaces.

**THEOREM 2.5.5:** *Let*

$$V = \underbrace{Z_p \times ... \times Z_p}_{n-times}$$

*be a group vector space over the group $Z_p$ where p is a prime i.e., $Z_p$ is a group under addition modulo p. V is a simple group vector space.*

*Proof:* Clear from the fact that $Z_p$ has no proper subgroups. Hence the claim.

Next we proceed on to define the notion of semigroup vector subspace of the group vector space V over G.

**DEFINITION 2.5.8:** *Let V be a group vector space over the group G. Let $W \subset V$ and $S \subset G$ where S is a semigroup under '+'. If W is a semigroup vector subspace over S then we call W to be pseudo semigroup vector subspace of V.*

We illustrate this by some examples.

***Example 2.5.28:*** Let V = Z × Z × Z be a group vector space over Z. $W = Z^+ \cup \{0\} \times Z^+ \cup \{0\} \times \{0\} \subset V$. W is a semigroup vector space over the semigroup $S = Z^+ \cup \{0\} \subseteq Z$. Thus W is a pseudo semigroup vector subspace of V.

***Example 2.5.29:*** Let

$$V = \left\{ \begin{pmatrix} a & b \\ c & d \end{pmatrix} \middle| \, a, b, c, d \in Z \right\}$$

be a group vector space over the group Z.

$$W = \left\{ \begin{pmatrix} a & b \\ c & d \end{pmatrix} \middle| \, a, b, c, d \in Z^+ \cup \{0\} \right\}$$

is a subset of V and W is a semigroup vector subspace over S = $Z^+ \cup \{0\}$. Thus W is a pseudo semigroup vector subspace of V.

***Example 2.5.30:*** Let V = {(1 1 0), (0 0 0), (0 1 1), (1 0 1), (0 0 1)} be a group vector space over the group $Z_2$ = {0, 1}. Clearly V has no pseudo semigroup vector subspace.

***Example 2.5.31:*** Let V = $Z_5 \times Z_5 \times Z_5 \times Z_5$ be a group vector space over $Z_5$; V has no pseudo semigroup vector subspace.

Now we define yet another type of subspace viz. pseudo set vector subspace of a group vector space V.

**DEFINITION 2.5.9:** *Let V be a group vector space over the group G. Suppose W $\subset$ V is a subset of V. Let S be a subset of G. If W is a set vector space over S then we call W to be a pseudo set vector subspace of the group vector space.*

We now give some illustrations.

***Example 2.5.32:*** Let V = {(1 1 1 0), (0 0 0 0), (1 1 0 0), (0 0 1 0)} be a group vector space over the group $Z_2$ = {0, 1}. Take W = {(1 1 1 0), (0 0 1 0)} $\subset$ V be a subset of V W is a pseudo set vector space over the set S = {1} $\subseteq Z_2$.

***Example 2.5.33:*** Let P = Z × Z × Z × Z be a group vector space over the group Z.

Take W = {(1 1 1 1), (0 0 0 0), (1 2 0 1), (3 3 5 1), (7 1 2 3), (1 0 0 1)} a proper subset of P. W is a pseudo set vector subspace over the set S = {0, 1} $\subseteq$ Z.

***Example 2.5.34:*** Let V = $Z_3 \times Z_3 \times Z_3 \times Z_3$ be a group vector space over $Z_3$. Take W = {(1 1 1), (2 2 2), (0 0 0), (1 0 1), (2 0 2)} $\subset$ V. Let S = {1, 2} a proper subset of $Z_3$ W is a set vector space over S. Thus W is a pseudo vector subspace of V.

It is an open problem whether there exists a group vector space, which has no pseudo, set vector subspaces.

Now we proceed onto define the notion of transformations of group vector spaces, which will be known as group linear transformations.

**DEFINITION 2.5.10:** *Let V and W be two group vector spaces defined over the same group G. A map T from V to W will be called as the group linear transformation if*

$$T(\alpha v) = \alpha T(v)$$

*for all $\alpha \in G$ and for all $v \in V$.*

We illustrate this by the following examples.

***Example 2.5.35:*** Let V = Z × Z × Z and W = Q × Q × Q × Q be two group vector spaces over the group Z. Let T : V $\rightarrow$ W be defined by T (x y z) = (z y x y). Clearly T is a group linear transformation of V into W.

***Example 2.5.36:*** Let V = {(0 0 0), (1 1 1), (0 1 0), (1 1 1 1), (0 0 0 0), (1 1 0 1), (0 1 1 1)} and

$$W = \left\{\begin{pmatrix}1 & 1\\0 & 1\end{pmatrix}, \begin{pmatrix}0 & 0\\1 & 1\end{pmatrix}, \begin{pmatrix}0 & 0\\0 & 0\end{pmatrix}, \begin{pmatrix}1 & 1\\1 & 1\end{pmatrix},\right.$$

$$\left.\begin{pmatrix}0 & 1\\0 & 0\end{pmatrix}, \begin{pmatrix}1 & 1\\1 & 0\end{pmatrix}, \begin{pmatrix}0 & 0\\1 & 1\end{pmatrix}, \begin{pmatrix}1 & 0\\0 & 1\end{pmatrix}, \begin{pmatrix}1 & 0\\0 & 0\end{pmatrix}\right\}$$

be the group vector space over the group $G = Z_2 = (0, 1)$.

$$T(x\ y\ z) = \begin{pmatrix} x & y \\ z & 0 \end{pmatrix}$$

and

$$T(x\ y\ z\ w) = \begin{pmatrix} x & y \\ z & w \end{pmatrix}$$

for (x y, z) and (x y z w) $\in$ V

$$T(0\ 0\ 0) = \begin{pmatrix} 0 & 0 \\ 0 & 0 \end{pmatrix}$$

$$T(1\ 1\ 1) = \begin{pmatrix} 1 & 1 \\ 1 & 0 \end{pmatrix}$$

$$T(0\ 1\ 0) = \begin{pmatrix} 0 & 1 \\ 0 & 0 \end{pmatrix}$$

$$T(1\ 1\ 1\ 1) = \begin{pmatrix} 1 & 1 \\ 1 & 1 \end{pmatrix}$$

$$T(0\ 0\ 0\ 0) = \begin{pmatrix} 0 & 0 \\ 0 & 0 \end{pmatrix}$$

$$T(1\ 1\ 0\ 1) = \begin{pmatrix} 1 & 1 \\ 0 & 1 \end{pmatrix}$$

$$T(0\ 1\ 1\ 1) = \begin{pmatrix} 0 & 1 \\ 1 & 1 \end{pmatrix}.$$

T is a group linear transformation of V to W.

We have for a group linear transformation T, $T^{-1}$ to exist provided the inverse mapping from W to V exists, otherwise we may not have $T^{-1}$ to exist for the T. Thus for a given T, $T^{-1}$ may or may not exist.

***Example 2.5.37:*** Let V = {(a a a a)| a ∈ Z} and

$$W = \left\{ \begin{pmatrix} a & a \\ a & a \end{pmatrix} \middle| a \in Z \right\}$$

be two group vector spaces over the group Z. A map

$$T\{(a\ a\ a\ a)\} = \begin{pmatrix} a & a \\ a & a \end{pmatrix}$$

for every (a a a a) ∈ V is both one to one and on to for define

$$T^{-1}\begin{pmatrix} a & a \\ a & a \end{pmatrix} = (a\ a\ a\ a).$$

$T^{-1}$ exists.

$$\begin{aligned} T^{-1} \text{ o } T\,(a\ a\ a\ a) &= T^{-1}\begin{pmatrix} a & a \\ a & a \end{pmatrix} \\ &= (a\ a\ a\ a) \end{aligned}$$

and

$$\begin{aligned} T \text{ o } T^{-1}\begin{pmatrix} a & a \\ a & a \end{pmatrix} &= T(a\ a\ a\ a) \\ &= \begin{pmatrix} a & a \\ a & a \end{pmatrix}. \end{aligned}$$

Thus T o $T^{-1}$ is identity map on W and $T^{-1}$ o T is the identity map on V. We call the group linear transformation T to be an invertible one.

***Example 2.5.38:*** Let $V = Z_{12} \times Z_{12} \times Z_{12} \times Z_{12}$ and

$$W = \left\{ \begin{pmatrix} a & b \\ c & d \end{pmatrix} \middle| a, b, c, d \in Z_{12} \right\}$$

be two group vector spaces over the group $Z_{12}$. T be a map such that

$$T(a, b, c, d) = \begin{pmatrix} a & b \\ c & d \end{pmatrix}$$

for every a, b, c, d $\in Z_{12}$.

T is a group linear transformation of V into W, in fact T is one to one and onto $T^{-1}$ exists.

Next we proceed onto define the notion of group linear operations on V, V a group vector space over the group G.

**DEFINITION 2.5.11:** *Let V be a group vector space over the group G. Let T from V to V be a group linear transformation then we call T to be a group linear operator on V.*

We now illustrate group linear operator on V by some examples.

***Example 2.5.39:*** Let V = {(a b c d) | a, b, c, d $\in$ Z} be a group vector space over Z. Define T from V to V by T {(a b c d)} = (d c b a) for every (a, b, c, d) $\in$ V. Clearly T is a group linear operator on V.

In fact it can further be verified $T^{-1}$ exists and $T^{-1}$ o T = T o $T^{-1}$ = identity group linear operator on V for

$$T^{-1} \text{ o } T \{(a\ b\ c\ d)\} = T^{-1} \{(d\ c\ b\ a)\} = (a\ b\ c\ d);$$

i.e., $T^{-1}$ o T is identity on V. Now

$$T \text{ o } T^{-1} \{(a\ b\ c\ d)\} = T \{(d\ c\ b\ a)\} = (a\ b\ c\ d).$$

$T \circ T^{-1}$ is also the group linear operator which is the identity map in this case.

All identity maps on V are identity group linear operator on V.

***Example 2.5.40:*** Let

$$V = \left\{ \begin{pmatrix} a & b \\ c & d \end{pmatrix} \middle| \; a, b, c, d \in Z_{10} \right\}$$

be the group vector space over the group $Z_{10}$. Define a map T from V to V by

$$T\left(\begin{pmatrix} a & b \\ c & d \end{pmatrix}\right) = \begin{pmatrix} a & b \\ 0 & 0 \end{pmatrix}$$

for every a, b, c, d ∈ V. T is a group linear operator on V, but T is clearly not an invertible group linear operator on V.

The reader is left with the task of finding

1. The algebraic structure given by the set of all group linear operators from the group vector space V to the group vector space W both V and W defined over the same group.
2. The algebraic structure of the set of all group linear operators from V into V, V the group vector space defined over the group G.

We denote the set of all group linear transformations from V to W defined over the group G by $M_G$ (V,W) and that the set of all group linear operators of V by $M_G$ (V, V).

1. What is the algebraic structure of $M_G$ (V,W) ?
2. What is the algebraic structure of $M_G$ (V,V)?

Now we proceed on to define the notion of group linear algebra over a group.

**DEFINITION 2.5.12:** *Let V be a group vector space over the group G. If V is again a group under addition then we call V to be a group linear algebra over G.*

It is clear from the very definition every group linear algebra defined over a group G is a group vector space over the group G but a group vector space is never a group linear algebra i.e., $\alpha_1 (v_1 + v_2) = \alpha_1 v_1 + \alpha_1 v_2$ for all $\alpha_1 \in G$ and $v_1, v_2 \in V$ may not be always true in V.

***Example 2.5.41:*** Let

$$V = \left\{ \begin{pmatrix} a & a \\ a & a \end{pmatrix} \middle| a \in Z \right\}.$$

V is a group linear algebra over Z with matrix addition on V.

***Example 2.5.42:*** Let V = {(0 0 0), (1 1 0), (0 0 1), (1 1 1), (0 0 0 0), (0 0 1 0), (0 0 1 1)} be a group vector space over $Z_2$ = {0, 1}. V is not closed under any additive operation so V is not a group linear algebra over G = $Z_2$.

Thus we see in general all group vector spaces are not group linear algebras.

***Example 2.5.43:*** Let $V = \{(a_{ij})_{m \times n} \mid a_{ij} \in Z_{12}\}$ be the collection of all m × n matrices with entries from $Z_{12}$. V is a group under matrix addition. V is a group linear algebra over $Z_{12}$.

***Example 2.5.44:*** Let V = {(0 0 0 0), (1 0 0 0), (0 1 0 0), (0 0 1 0), (0 0 0 1), (0 0 1 1), (0 1 0 1), (1 1 0 0), (0 1 1 0), (1 1 1 0), (0 1 1 1), (1 1 0 1), (1 0 1 1), (1 1 1 1), (1 0 0 1), (1 0 1 0)} with entries from $Z_2$ = {0, 1}. V is a group linear algebra over the group $Z_2$ = {0, 1}.

***Example 2.5.45:*** Let V = {(a b c) | a, b, c ∈ Z}, V under component wise addition is a group; V is a group linear algebra over Z.

***Example 2.5.46:*** Let V = {$Z_{10} \times Z_{10} \times Z_{10} \times Z_{10} \times Z_{10}$ = ($a_1$, $a_2$, …, $a_5$) | $a_i \in Z_{10}$; $1 \leq i \leq 5$} be a group under component wise addition, V is a group linear algebra over $Z_{10}$.

In case of group linear algebras the linear independence and the dimension are little different from that of the group vector spaces.

**DEFINITION 2.5.13:** *Let V be a group linear algebra over the group G. Let X ⊂ V be a proper subset of V, we say X is a linearly independent subset of V if X = {$x_1$, …, $x_n$} and for some $\alpha_i \in G$, $1 \leq i \leq n$, $\alpha_1 x_1 + \ldots + \alpha_n x_n = 0$ if and only if each $\alpha_i$ = 0.*

*A linearly independent subset X of V is said to generate V if every element of v ∈ V can be represented as v = $\sum_{i=1}^{n} \alpha_i x_i$, $\alpha_i \in$ G ($1 \leq i \leq n$).*

We illustrate this situation by the following examples.

***Example 2.5.47:*** Let V = {(0 0 0), (1 0 0), (0 0 1), (0 1 0), (1 1 0), (1 0 1), (0 1 1), (1 1 1)} be the group linear algebra over the group $Z_2$ = {0, 1}. V is generated by the set X = {(1 0 0), (0 1 0), (0 0 1)}. Clearly X is also a linearly independent subset of V over $Z_2$ = {0, 1}.

Thus dimension of V is 3. Hence as in case of usual vector spaces we say the dimension of a group linear algebra is also the cardinality of the linearly independent subset X of V which generates V.

***Example 2.5.48:*** Let

$$V = \left\{ \begin{pmatrix} a_1 & a_2 & a_3 \\ a_4 & a_5 & a_6 \end{pmatrix} \middle| \, a_i \in Z \right\}$$

be the group linear algebra over the group Z. Let

$$X = \left\{ \begin{pmatrix} 1 & 0 & 0 \\ 0 & 0 & 0 \end{pmatrix}, \begin{pmatrix} 0 & 1 & 0 \\ 0 & 0 & 0 \end{pmatrix}, \begin{pmatrix} 0 & 0 & 1 \\ 0 & 0 & 0 \end{pmatrix}, \right.$$

$$\left. \begin{pmatrix} 0 & 0 & 0 \\ 1 & 0 & 0 \end{pmatrix}, \begin{pmatrix} 0 & 0 & 0 \\ 0 & 1 & 0 \end{pmatrix}, \begin{pmatrix} 0 & 0 & 0 \\ 0 & 0 & 1 \end{pmatrix} \right\} \subseteq V$$

be the generating subset of V which is also linearly independent. Thus the dimension of the group linear algebra V is 6.

***Example 2.5.49:*** Let

$$V = \left\{ \begin{pmatrix} a & b \\ c & d \end{pmatrix} \middle| a, b, c, d \in Z_6 \right\}.$$

V is a group linear algebra over the group $Z_6$. The set

$$X = \left\{ \begin{pmatrix} 0 & 1 \\ 0 & 0 \end{pmatrix} \begin{pmatrix} 1 & 0 \\ 0 & 0 \end{pmatrix} \begin{pmatrix} 0 & 0 \\ 1 & 0 \end{pmatrix} \begin{pmatrix} 0 & 0 \\ 0 & 1 \end{pmatrix} \right\} \subseteq V$$

is the linearly independent subset of V which generates V. Clearly dimension of V is four.

***Example 2.5.50:*** Let $V = Z_6 \times Z_6 \times Z_6 \times Z_6$ be the group linear algebra over the group $Z_6$. Now X = {(1 0 0 0), (0 0 0 1), (0 1 0 0), (0 0 1 0)} is the generating set of V. The dimension of V is four over $Z_6$.

***Example 2.5.51:*** Let $V = Z_6 \times Z_6 \times Z_6 \times Z_6 \times Z_6$ be the group linear algebra over the group $Z_6$. For this group linear algebra also X = { (1 0 0 0 0), (0 1 0 0 0), (0 0 1 0 0), (0 0 0 1 0), (0 0 0 0 1)} is the linearly independent subset which generates V and dimension is 5.

Now we proceed onto define the notion of group linear subalgebra of a group linear algebra.

**DEFINITION 2.5.14:** *Let V be a group linear algebra over the group G. Let W $\subseteq$ V be a proper subset of V. We say W is a group linear subalgebra of V over G if W is itself a group linear algebra over G.*

We illustrate this situation by the following examples.

***Example 2.5.52:*** Let

$$V = \left\{ \begin{pmatrix} a & b \\ c & d \end{pmatrix} \middle| a, b, c, d \in Z \right\},$$

V is a group linear algebra over the group Z. Take

$$W = \left\{ \begin{pmatrix} a & a \\ a & a \end{pmatrix} \middle| a \in Z \right\};$$

W $\subseteq$ V and W is a group linear subalgebra of V over Z.

***Example 2.5.53:*** Let V = {(1 0 0), (0 0 0), (0 1 0), (0 0 1), (1 1 1), (0 1 1), (1 1 0), (1 0 1)} be the group linear algebra over the group $Z_2$ = {0, 1}. Take W = {(0 0 0), (1 1 1)} $\subseteq$ V, W is a group linear subalgebra of V.

***Example 2.5.54:*** Let V = $Z_9 \times Z_9 \times Z_9$ be the group linear algebra over $Z_9$. Let W = $Z_9 \times \{0\} \times Z_9 \subseteq$ V; W is the group linear subalgebra of V over $Z_9$.

***Example 2.5.55:*** Let V = $Z_8 \times Z_8$ be the group linear algebra over the group G = {0, 2, 4, 6} addition modulo 8. Let W = {0, 2, 4, 6} × {0, 2, 4, 6} $\subseteq$ V; W is a group linear subalgebra of V.

Now having seen several examples of group linear subalgebras over the group linear algebra we proceed onto define the notion of direct sum of group linear subalgebras of a group linear algebra.

**DEFINITION 2.5.15:** *Let V be a group linear algebra over the group G. Let $W_1, \ldots, W_n$ be group linear subalgebras of V over G.*

*We say V is a direct sum of the group linear subalgebras $W_1, W_2, \ldots, W_n$ if*

1. $V = W_1 + \ldots + W_n$
2. $W_i \cap W_j = \{0\}$ *if* $i \neq j$; $1 \leq i, j \leq n$.

Now we illustrate this situation by the following examples.

***Example 2.5.56:*** Let V be $Z_{14} \times Z_{14} \times Z_{14}$ be the group linear algebra over the group $Z_{14}$, the group under addition modulo 14. Let $W_1 = Z_{14} \times \{0\} \times \{0\}$, $W_2 = \{0\} \times Z_{14} \times \{0\}$ and $W_3 = \{0\} \times \{0\} \times Z_{14}$ be the group linear subalgebras of V. We see $V = W_1 + W_2 + W_3$ and $W_i \cap W_j = \{0\}$ if $i \neq j$; $1 \leq i, j \leq 3$.

***Example 2.5.57:*** Let

$$V = \left\{ \begin{pmatrix} a & b \\ c & d \end{pmatrix} \middle| a, b, c, d \in Z \right\}$$

be the group linear algebra over Z. Let

$$W_1 = \left\{ \begin{pmatrix} a & 0 \\ 0 & 0 \end{pmatrix} \middle| a \in Z \right\},$$

$$W_2 = \left\{ \begin{pmatrix} 0 & b \\ 0 & 0 \end{pmatrix} \middle| b \in Z \right\},$$

$$W_3 = \left\{ \begin{pmatrix} 0 & 0 \\ c & 0 \end{pmatrix} \middle| c \in Z \right\}$$

and

$$W_4 = \left\{ \begin{pmatrix} 0 & 0 \\ 0 & d \end{pmatrix} \middle| \, d \in Z \right\}$$

be the group linear subalgebras of V. Clearly $W_1 + W_2 + W_3 + W_4 = V$ and

$$W_i \cap W_j = \begin{pmatrix} 0 & 0 \\ 0 & 0 \end{pmatrix} \text{ if } i \neq j \text{ ; } 1 \leq i, j \leq 4.$$

Thus V is the direct sum of group linear subalgebras over Z. Now take

$$S_1 = \left\{ \begin{pmatrix} a & 0 \\ b & 0 \end{pmatrix} \middle| \, a, b \in Z \right\}$$

$$S_2 = \left\{ \begin{pmatrix} 0 & a \\ 0 & 0 \end{pmatrix} \middle| \, a \in Z \right\}$$

and

$$S_3 = \left\{ \begin{pmatrix} 0 & 0 \\ 0 & b \end{pmatrix} \middle| \, b \in Z \right\}$$

be group linear subalgebras of V. We see $V = S_1 + S_2 + S_3$ and

$$S_i \cap S_j = \begin{pmatrix} 0 & 0 \\ 0 & 0 \end{pmatrix}, \text{ if } i \neq j \ (1 \leq i, j \leq 3).$$

Suppose

$$T_1 = \left\{ \begin{pmatrix} a & 0 \\ 0 & b \end{pmatrix} \middle| \, a, b \in Z \right\}$$

and

$$T_2 = \left\{ \begin{pmatrix} 0 & c \\ d & 0 \end{pmatrix} \middle| \, c, d \in Z \right\}$$

be group linear subalgebras of V then we see $V = T_1 + T_2$ and

$$T_1 \cap T_2 = \begin{pmatrix} 0 & 0 \\ 0 & 0 \end{pmatrix}.$$

Thus V is the direct sum of group linear subalgebras. From this example it is evident that we have many ways of writing V as a direct sum of group linear subalgebras of V.
Further suppose

$$R_1 = \left\{ \begin{pmatrix} a & b \\ 0 & 0 \end{pmatrix} \middle| \, a, b \in Z \right\},$$

$$R_2 = \left\{ \begin{pmatrix} a & 0 \\ c & d \end{pmatrix} \middle| \, a, c, d \in Z \right\}$$

and

$$R_3 = \left\{ \begin{pmatrix} a & 0 \\ 0 & b \end{pmatrix} \middle| \, a, b \in Z \right\}$$

be three group linear subalgebras of V over Z we see

$$V \neq R_1 + R_2 + R_3 \text{ but } R_i \cap R_j \neq \begin{pmatrix} 0 & 0 \\ 0 & 0 \end{pmatrix}, \; 1 \leq i, j \leq 3.$$

Thus we see any set of group subalgebra need not lead to the direct sum. Also if

$$V_1 = \left\{ \begin{pmatrix} 0 & 0 \\ a & 0 \end{pmatrix} \middle| \, a \in Z \right\}$$

and

$$V_2 = \left\{ \begin{pmatrix} a & b \\ 0 & 0 \end{pmatrix} \middle| \, a, b \in Z \right\}$$

be two group subalgebras of V still $V \neq V_1 + V_2$ though

$$V_1 \cap V_2 \neq \begin{pmatrix} 0 & 0 \\ 0 & 0 \end{pmatrix}.$$

Thus we see in this case $\begin{pmatrix} a & b \\ c & d \end{pmatrix}$ with $d \neq 0$ cannot find its place in $V_1 + V_2$. Hence $V = V_1 + V_2$ is not a direct sum of V.

Now we proceed onto give yet another example of direct sum of group linear subalgebras.

***Example 2.5.58:*** Let V = {(1 0 0 0), (0 1 0 0), (0 0 1 0), (1 1 1 1), (0 0 0 1), (0 0 1 1), (1 1 0 0), (1 0 0 1), (0 1 1 0), (1 0 1 0), (0 0 0 0), (0 1 0 1), (1 1 1 0), (0 1 1 1), (1 0 1 1), (1 1 0 1)} be the group linear algebra over the group $Z_2 = \{0, 1\}$. Write V as a direct sum of group linear subalgebras. Can we represent V in more than one way as a direct sum?

Now we proceed onto define the notion of pseudo direct sum of a group linear algebra as a sum of group linear subalgebras.

**DEFINITION 2.5.16:** *Let V be a group linear algebra over the group G. Suppose $W_1, W_2, \ldots, W_n$ are distinct group linear subalgebras of V. We say V is a pseudo direct sum if*

1. $W_1 + \ldots + W_n = V$
2. $W_i \cap W_j \neq \{0\}$, *even if* $i \neq j$
3. *We need $W_i$'s to be distinct i.e., $W_i \cap W_j \neq W_i$ or $W_j$ if $i \neq j$.*

We now illustrate this situation by the following example.

***Example 2.5.59:*** Let

$$V = \left\{ \begin{pmatrix} a_1 & a_2 & a_3 \\ a_4 & a_5 & a_6 \end{pmatrix} \middle| a_i \in Z, 1 \leq i \leq 6 \right\}$$

be a group linear algebra over Z. Take

$$W_1 = \left\{ \begin{pmatrix} a_1 & a_2 & 0 \\ a_3 & 0 & 0 \end{pmatrix} \middle| \; a_1, a_3, a_2 \in Z \right\},$$

$$W_2 = \left\{ \begin{pmatrix} 0 & a_2 & a_3 \\ 0 & 0 & 0 \end{pmatrix} \middle| \; a_1, a_3 \in Z \right\},$$

$$W_3 = \left\{ \begin{pmatrix} a_1 & 0 & a_3 \\ a_4 & 0 & a_6 \end{pmatrix} \middle| \; a_1, a_3, a_4, a_6 \in Z \right\}$$

and

$$W_4 = \left\{ \begin{pmatrix} a_1 & a_4 & 0 \\ 0 & a_5 & a_6 \end{pmatrix} \middle| \; a_1, a_4, a_5, a_6 \in Z \right\}.$$

Clearly $W_1$, $W_2$, $W_3$ and $W_4$ are group linear subalgebras of V. We see

$$W_i \cap W_j \neq \begin{pmatrix} 0 & 0 & 0 \\ 0 & 0 & 0 \end{pmatrix}, \text{ for } 1 \leq i, j \leq 4 \; i \neq j.$$

Further $W_i \cap W_j \neq W_i$ or $W_i \cap W_j \neq W_j$ if $i \neq j$. Finally $V = W_1 + W_2 + W_3 + W_4$, so we say V is the pseudo direct sum of group linear subalgebras of V.

We give yet another example before we proceed on to describe further properties about group linear algebras.

***Example 2.5.60:*** Let $V = \{Z_{18} \times Z_{18} \times Z_{18}\}$ be a group linear algebra over $Z_{18}$.

Take $W_1 = Z_{18} \times Z_{18} \times \{0\}$, $W_2 = Z_{18} \times \{0\} \times Z_{18}$, $W_3 = \{0\ 2\ 4\ 6\ 8\ 10\ 12\ 14\ 16\} \times Z_{18} \times Z_{18}$ be three group linear subalgebras of V. Then $V = W_1 + W_2 + W_3$ and $W_i \cap W_j \neq \{0\}$; $i \neq j$; $1 \leq i, j \leq 3$ so V is the pseudo direct sum of group linear subalgebras.

It is important to note that a group linear algebra can both have a pseudo direct sum as well as direct sum. We see we do not have any relation among them. It can so happen a group linear algebra can have both way of decomposition.

Now we proceed onto define yet another new algebraic structure of a group linear algebra.

**DEFINITION 2.5.17:** *Let V be a group linear algebra over the group G. Let $W \subseteq V$ be a proper subgroup of V. Suppose $H \subset G$ be a semigroup in G. If W is a semigroup linear algebra over H then we call W to be a pseudo semigroup linear subalgebra of the group linear algebra V.*

We illustrate this situation by the following example.

***Example 2.5.61:*** Let

$$V = \left\{ \begin{pmatrix} a & b \\ c & d \end{pmatrix} \middle| a, b, c, d \in Z \right\}$$

be a group linear algebra over the group Z. Let $Z^+ \cup \{0\} = H \subset Z$ be the proper semigroup of Z under addition.

$$W = \left\{ \begin{pmatrix} a & a \\ a & a \end{pmatrix} \middle| a \in Z \right\} \subseteq V.$$

W is a semigroup linear algebra over H. We call W to be the pseudo semigroup linear subalgebra of the group linear algebra V.

We give yet another example.

***Example 2.5.62:*** Let $V = \{Z \times Z \times Z\}$ be the group linear algebra over the group Z. Let $P = Z^+ \cup \{0\}$ be the semigroup contained in Z. Let $W = 2Z \times 2Z \times 2Z \subset V$; W is a pseudo semigroup linear subalgebra over P.

***Example 2.5.63:*** Let $V = Z_2 \times Z_2 \times Z_2$ be the group linear algebra over $Z_2$. $Z_2$ has no proper subset which is a semigroup, so V cannot have pseudo semigroup linear subalgebra. In fact we have a class of group linear algebras which has no pseudo semigroup linear subalgebras.

Let $V = Z_p \times \ldots \times Z_p$, $Z_p$ the set of primes $\{0, 1, \ldots, p - 1\}$ under addition modulo p. V is a group linear algebra over $Z_p$. But $Z_p$ has no proper subset P which is a semigroup. Thus V has no pseudo semigroup linear subalgebras. Thus we have a class of group linear algebras which has no pseudo semigroup linear subalgebras.

Suppose we consider

$$\left\{ V = \left(a_{ij}\right)_{m\times n} \middle| a_{ij} \in Z_p;\ \text{p a prime}, 1 \leq i \leq m \text{ and } 1 \leq i \leq n \right\}.$$

V is a group linear algebra over $Z_p$. This has no pseudo semigroup linear subalgebras.

For varying primes p we get different classes of group linear algebras which has no pseudo semigroup linear subalgebras. We have yet another class of group linear algebras which has no pseudo semigroup linear subalgebras.

Consider $Z_p$ [x] = {all polynomials in the variable x with coefficients from $Z_p$; p a prime}; $Z_p$ [x] is a group linear algebra over $Z_p$. Clearly $Z_p$ has no proper subset which is a semigroup under addition modulo p. So $Z_p$ [x] has no pseudo semigroup linear subalgebras.

In fact we can have yet another substructure in group linear algebras which will be known as group vector subspaces of the group linear algebras.

**DEFINITION 2.5.18:** *Let V be a group linear algebra over the group G. Let P be a proper subset of V. P is just a set and it is not a closed structure. If P is a group vector space over G we call P to be the pseudo group vector subspace of V.*

We illustrate this by the following examples.

***Example 2.5.64:*** Let

V = {(1 0 0), (0 0 1), (0 1 0), (0 0 0), (1 1 0), (0 1 1), (1 0 1), (1 1 1)}

be the group linear algebra over the group $Z_2 = \{0, 1\}$.
Take

P = {(0 0 0), (1 1 0), (0 0 1), (0 1 0)} $\subseteq$ V.

P is a group vector subspace of V. Thus P is a pseudo group vector subspace of V.

***Example 2.5.65:*** Let

$$V = \left\{ \begin{pmatrix} a & b \\ c & d \end{pmatrix} \middle| a, b, c, d \in Z_3 \right\}$$

be the group linear algebra over $Z_3$.
Take

$$P = \left\{ \begin{pmatrix} 1 & 0 \\ 0 & 1 \end{pmatrix}, \begin{pmatrix} 0 & 2 \\ 0 & 1 \end{pmatrix}, \begin{pmatrix} 0 & 0 \\ 0 & 0 \end{pmatrix}, \begin{pmatrix} 2 & 0 \\ 0 & 2 \end{pmatrix}, \begin{pmatrix} 0 & 1 \\ 0 & 2 \end{pmatrix} \right\} \subset V;$$

P is the pseudo group vector subspace of V over $Z_3$.
Take

$$X = \left\{ \begin{pmatrix} 0 & 0 \\ 0 & 0 \end{pmatrix}, \begin{pmatrix} 0 & 0 \\ 1 & 1 \end{pmatrix}, \begin{pmatrix} 0 & 0 \\ 2 & 2 \end{pmatrix}, \begin{pmatrix} 1 & 1 \\ 0 & 0 \end{pmatrix}, \begin{pmatrix} 2 & 2 \\ 0 & 0 \end{pmatrix} \right\} \subset V.$$

X is also a pseudo group vector subspace of V over $Z_3$. However every proper subset of V is not a pseudo group vector subspace of V.

For instance

$$T = \left\{ \begin{pmatrix} 1 & 0 \\ 0 & 0 \end{pmatrix}, \begin{pmatrix} 1 & 1 \\ 1 & 1 \end{pmatrix}, \begin{pmatrix} 0 & 0 \\ 0 & 0 \end{pmatrix} \right\} \subset V$$

is not a pseudo group vector subspace of V as

$$2\begin{pmatrix} 1 & 0 \\ 0 & 0 \end{pmatrix} = \begin{pmatrix} 2 & 0 \\ 0 & 0 \end{pmatrix} \notin T$$

or

$$2\begin{pmatrix} 1 & 1 \\ 1 & 1 \end{pmatrix} = \begin{pmatrix} 2 & 2 \\ 2 & 2 \end{pmatrix} \notin T.$$

**Chapter Three**

# SET FUZZY LINEAR ALGEBRAS AND THEIR PROPERTIES

In this chapter we define the new notion of set fuzzy linear algebra analogous to set vector space; for these algebraic set up will be of immense use in application to fuzzy models or in any other models for in these set vector spaces and set fuzzy vector spaces we can induct any wanted elements without affecting the system and the structure. We now just recall the definition of fuzzy vector spaces before we proceed on to define set fuzzy vector spaces.

**DEFINITION 3.1:** *A fuzzy vector space (V, η) or ηV is an ordinary vector space V with a map $\eta : V \to [0, 1]$ satisfying the following conditions;*

1. $\eta(a + b) \geq \min\{\eta(a), \eta(b)\}$
2. $\eta(-a) = \eta(a)$
3. $\eta(0) = 1$
4. $\eta(ra) \geq \eta(a)$

*for all a, b, ∈ V and r ∈ F where F is the field.*

We now define the notion of set fuzzy vector space or $V_\eta$ or $V\eta$ or $\eta V$.

**DEFINITION 3.2:** *Let V be a set vector space over the set S. We say V with the map $\eta$ is a fuzzy set vector space or set fuzzy vector space if $\eta$: V $\to$ [0, 1] and $\eta$ (ra) $\geq$ $\eta$(a) for all a $\in$V and r $\in$ S. We call $V_\eta$ or V$\eta$ or $\eta$V to be the fuzzy set vector space over the set S.*

We now illustrate this situation by the following example.

***Example 3.1:*** Let V = {(1 3 5), (1 1 1), (5 5 5), (7 7 7), (3 3 3), (5 15 25), (1 2 3)} be set which is a set vector space over the set S = {0, 1}.

Define a map $\eta$: V $\to$ [0, 1] by

$$\eta(x, y, z) = \left(\frac{x+y+z}{50}\right) \in [0,1]$$

for (x, y, z) $\in$ V. $V_\eta$ is a fuzzy set vector space.

***Example 3.2:*** Let V = $Z^+$ the set of integers. S = $2Z^+$ be the set. V is a set vector space over S. Define $\eta$: V $\to$ [0, 1] by, for every v $\in$ V; $\eta(v) = \frac{1}{v}$ . $\eta$V is a set fuzzy vector space or fuzzy set vector space.

***Example 3.3:*** Let V ={($a_{ij}$) | $a_{ij}$ $\in$ $Z^+$; 1 $\leq$ i, j $\leq$ n} be the set of all n $\times$ n matrices with entries from $Z^+$.

Take S = $3Z^+$ to be the set. V is a fuzzy set vector space where $\eta$: V $\to$ [0, 1] is defined by

$$\eta(A = (a_{ij})) = \begin{cases} \dfrac{1}{5|A|} & \text{if } |A| \neq 0 \\ 1 & \text{if } |A| = 0. \end{cases}$$

V$\eta$ is the fuzzy set vector space.

The main advantage of defining set vector spaces and fuzzy set vector spaces is that we can include elements x in the set vector spaces V provided for all $s \in S$, $sx \in V$ this cannot be easily done in usual vector spaces. Thus we can work with the minimum number of elements as per our need and work with them by saving both time and money.

We give yet some more examples.

***Example 3.4:*** Let $V = 2Z^+ \times 5Z^+ \times 7Z^+$ be a set vector space over the set $Z^+$; with $\eta: V \to [0, 1]$ defined by

$$\eta((x, y, z)) = \frac{1}{x+y+z}$$

makes, $\eta V$ a fuzzy set vector space.

Now we define the notion of set fuzzy linear algebra.

**DEFINITION 3.3:** *A set fuzzy linear algebra (or fuzzy set linear algebra) (V, η) or ηV is an ordinary set linear algebra V with a map such η: V → [0, 1] such that η(a + b) ≥ min (η(a), η(b)) for a, b ∈ V.*

*Since we know in the set vector space V we merely take V to be a set but in case of the set linear algebra V we assume V is closed with respect to some operation usually denoted as '+' so the additional condition η(a + b) ≥ min (η(a), η(b)) is essential for every a, b ∈ V.*

We illustrate this situation by the following examples.

***Example 3.5:*** Let $V = Z^+[x]$ be a set linear algebra over the set $S = Z^+$; $\eta: V \to [0, 1]$.

$$\eta(p(x)) = \begin{cases} \dfrac{1}{\deg(p(x))} & \\ 1 & \text{if } p(x) \text{ is a constant.} \end{cases}$$

Clearly Vη is a set fuzzy linear algebra.

***Example 3.6:*** Let

$$V = \left\{ \begin{pmatrix} a & b \\ c & d \end{pmatrix} \middle| a, b, d, c \in Z^{+} \right\}$$

be set linear algebra over $2Z^+ = S$. Define

$$\eta\left(\begin{pmatrix} a & b \\ c & d \end{pmatrix}\right) = \begin{cases} \dfrac{1}{|ad - bc|} & \text{if } ad \neq bc \\ 0 & \text{if } ad = bc \end{cases}$$

for every a, b, c, d $\in Z^+$. Clearly Vη is a fuzzy set linear algebra.

***Example 3.7:*** Let $V = Z^+$ be a set linear algebra over $Z^+$. Define $\eta : V \to [0, 1]$ as $\eta(a) = \frac{1}{a}$. Vη is a fuzzy set linear algebra.

Now we proceed onto define the notion of fuzzy set vector subspace and fuzzy set linear subalgebra.

**DEFINITION 3.4:** *Let V be a set vector space over the set S. Let $W \subset V$ be the set vector subspace of V defined over S. If $\eta$: $W \to [0, 1]$ then $W_\eta$ is called the fuzzy set vector subspace of V.*

We illustrate this by the following example.

***Example 3.8:*** Let V = {(1 1 1), (1 0 1), (0 1 1), (0 0 0), (1 0 0)} be a set vector space defined over the set {0, 1}. Define $\eta : V \to [0, 1]$ by

$$\eta(x\ y\ z) = \frac{(x + y + z)}{9} \pmod{2}.$$

So that

$$\eta\ (0\ 0\ 0) = 0$$
$$\eta\ (1\ 1\ 1) = \frac{1}{9}$$
$$\eta\ (1\ 0\ 1) = 0$$
$$\eta\ (1\ 0\ 0) = \frac{1}{9}$$
$$\eta\ (0\ 1\ 1) = 0$$

$V\eta$ is a set fuzzy vector space. Take W = {(1 1 1), (0 0 0), (0 1 1)} $\subset$ V. W is a set vector subspace of V. $\eta$: W $\to$ [0, 1].

$$\eta\ (0\ 0\ 0) = 0$$
$$\eta\ (111) = \frac{1}{9}$$
$$\eta\ (011) = 0.$$

$W_\eta$ is the fuzzy set vector subspace of V.

***Example 3.9:*** Let V = {(111), (1011), (11110), (101), (000), (0000), (0000000), (00000), (1111111), (11101), (01010), (1101101)} be a set vector space over the set S = {0,1}.

Let W = {(1111111), (0000000), (000), (00000), (11101), (01010) (101)} $\subset$ V. Define $\eta$: W $\to$ [0, 1] by

$$\eta(x_1, x_2, \ldots, x_r) = \frac{1}{12}.$$

$\eta$W is a fuzzy set vector subspace.

We now proceed on to define the notion of fuzzy set linear subalgebra.

**DEFINITION 3.5:** *Let V be a set linear algebra over the set S. Suppose W is a set linear subalgebra of V over S. Let $\eta : W \to [0, 1]$, $\eta W$ is called the fuzzy set linear subalgebra if $\eta\,(a + b) \geq \min \{\eta\,(a), \eta\,(b)\}$ for a, b, $\in$ W.*

We give some examples before we define some more new concepts.

***Example 3.10:*** Let V = $Z^+ \times Z^+ \times Z^+$ be a set linear algebra over the set S = $2Z^+$. W = $Z^+ \times 2Z^+ \times 4Z^+$ is a set linear subalgebra over the set S = $2Z^+$. Define η: W → [0, 1]

$$\eta\,(a\ b\ c) = 1 - \frac{1}{a+b+c}.$$

Clearly η (x, y) ≥ min {η (x), η (y)} where x = $(x_1, x_2, x_3)$ and y = $(y_1, y_2, y_3)$; x, y ∈ W. Wη is a fuzzy set linear subalgebra.

***Example 3.11:*** Let

$$V = \left\{ \begin{pmatrix} a & b \\ c & d \end{pmatrix} \middle| a, b, c, d, \in Z^+ \right\}.$$

V is a set linear algebra over the set S = {1, 3, 5, 7} ⊆ $Z^+$. Let

$$W = \left\{ \begin{pmatrix} a & a \\ a & a \end{pmatrix} \middle| a \in Z^+ \right\}$$

be the set linear subalgebra of V. Define η: W → [0, 1] by

$$\eta \begin{pmatrix} a & a \\ a & a \end{pmatrix} = 1 - \frac{1}{a}.$$

$W_\eta$ or Wη is a set fuzzy linear subalgebra.

Now we proceed on to define fuzzy semigroup vector spaces.

**DEFINITION 3.6:** *A semigroup fuzzy vector space or a fuzzy semigroup vector space (V, η) or Vη where V is an ordinary semigroup vector space over the semigroup S; with a map η : V → [0, 1] satisfying the following condition;*

*η (ra) ≥ η (a) for all a ∈ V and r ∈ S.*

Let us illustrate this structure by some examples.

***Example 3.12:*** Let V = {(1000), (1011), (1110), (0111), (0100), (0000), (0001)} be the semigroup vector space over the semigroup S = $Z_2$ = {0, 1}. Define $\eta : V \to [0, 1]$ as

$$\eta\ (a\ b\ c\ d) = \frac{(a+b+c+d)}{6} \pmod 2 .$$

Clearly $V\eta$ is the semigroup fuzzy vector space.

We give yet another example.

***Example 3.13:*** Let V = $Z_3 \times Z_3 \times Z_3$ be the semigroup vector space over the semigroup $Z_3$. Define $\eta : V \to [0, 1]$ as

$$\eta\ (x\ y\ z) = \frac{(x+y+z)}{7} \pmod 3 .$$

$V\eta$ is a semigroup fuzzy vector space.

In fact given a semigroup vector space V, we can get many semigroup fuzzy vector spaces.

For define $\eta_1 : V \to [0, 1]$ as

$$\eta_1\ (x, y, z) = \begin{cases} \dfrac{1}{(x+y+z) \pmod 3} & \text{if } x+y+z \not\equiv 0 \pmod 3 \\ 0 & \text{if } x+y+z \equiv 0 \pmod 3 \end{cases} .$$

$V\eta_1$ is a semigroup fuzzy vector space different from $V\eta$. Define $\eta_3 : V \to [0, 1]$ as

$$\eta_3\ (x, y, z) = \begin{cases} \frac{1}{2} \text{ if } x+y+z \equiv 0 \pmod 3 \\ \frac{1}{4} \text{ if } x+y+z \equiv 1 \pmod 3 \\ \frac{1}{6} \text{ if } x+y+z \equiv 2 \pmod 3 \end{cases}$$

$V\eta_3$ is a semigroup fuzzy vector space.

Thus study of semigroup fuzzy vector spaces gives us more freedom for it solely depends on η which is defined from V to [0, 1].

Next we define semigroup fuzzy vector subspaces of a semigroup vector space V.

**DEFINITION 3.7:** *Let V be a semigroup vector space over the semigroup S. Let $W \subset V$ be a semigroup vector subspace of V over S. We say $W\eta$ is a semigroup fuzzy vector subspace if $\eta$ : $W \to [0, 1]$, such that*

*(i) $\eta(x, y) \geq \min(\eta(x), \eta(y))$*
*(ii) $\eta(rx) \geq \eta(x)$ for all $r \in S$ and $x, y \in W$.*

We illustrate this by the following example.

***Example 3.14:*** Let $V = Z_7 \times Z_7 \times Z_7 \times Z_7$ be a semigroup vector space over the semigroup $S = Z_7$. Let $W = Z_7 \times \{0\} \times Z_7 \times \{0\}$ be the semigroup vector subspace of V. Define the map $\eta : W \to [0, 1]$ by

$$\eta\,(x\ 0\ y\ 0) = \begin{cases} \frac{1}{4} \text{ if } x + y \equiv 1 \pmod 7 \\ \frac{1}{3} \text{ if } x + y \equiv 2 \text{ or } 4 \pmod 7 \\ \frac{1}{8} \text{ if } x + y \equiv 3 \text{ or } 5 \pmod 7 \\ \frac{1}{6} \text{ if } x + y \equiv 6 \text{ or } 0 \pmod 7. \end{cases}$$

$\eta W$ is the semigroup fuzzy vector subspace of V.

Here also using one semigroup vector subspace W we can define several semigroup fuzzy vector subspaces.

Define $\eta_1 : W \to [0, 1]$

$$\text{by } \eta_1 \,(x\ 0\ y\ 0) = \begin{cases} 0 & \text{if } x + y \equiv 0 \pmod 7 \\ 1 & \text{if } x + y \not\equiv 0 \pmod 7). \end{cases}$$

$W\eta_1$ is a semigroup fuzzy vector subspace.

We see the definition of fuzzy semigroup vector space is not in any way different from the fuzzy set vector space. So this will enable one to go from one type of space to another using fuzzy concepts defined on them.

So one can easily transfer a study from semigroup vector space to set vector space by defining the corresponding fuzzy set vector space and the semigroup fuzzy vector space as the map $\eta$ and hence $\eta V$ does not give different structures but same type of structures.

Now we see even in case of semigroup linear algebra and set linear algebra the fuzzy structures are identical.

**DEFINITION 3.8:** *Let V be a semigroup linear algebra defined over the semigroup S. We say $\eta V$ is a semigroup fuzzy linear algebra if $\eta : V \to [0, 1]$ such that $\eta\,(x + y) \geq \min\,(\eta\,(x), \eta(y))$; $\eta\,(rx) \geq \eta\,(x)$ for every $r \in S$ and $x \in V$.*

Now we illustrate this situation by the following examples.

***Example 3.15:*** Let $V = Z_7 \times Z_7 \times Z_7 \times Z_7$ be the semigroup linear algebra defined over the semigroup $S = Z_7$. Define $\eta : V \to [0, 1]$ by

$$\eta\,(x\ y\ z\ \omega) = \begin{cases} 1 & \text{if } x + y + z \equiv 0 \pmod 7) \\ 0 & \text{otherwise.} \end{cases}$$

$V\eta$ is a semigroup fuzzy linear algebra.

***Example 3.16:*** Let

$$V = \left\{ \begin{pmatrix} a_1 & a_2 & a_3 \\ a_4 & a_5 & a_6 \end{pmatrix} \middle| a_i \in Z^+; 1 \le i \le 6 \right\}$$

be the semigroup linear algebra over the semigroup S = $Z^+$. Define η: V → [0, 1] as

$$\eta \begin{pmatrix} a_1 & a_2 & a_3 \\ a_4 & a_5 & a_6 \end{pmatrix} = \left\{ 1 - \frac{1}{(a_1 + a_4)} \right.$$

Vη is a semigroup fuzzy linear algebra.

***Example 3.17:*** Let V = $Z^+$ [x] be the polynomials with coefficient from $Z^+$ in the variable x; V under addition is a semigroup. V is a semigroup linear algebra over $Z^+$.

Define η : V → [0, 1] as

$$\eta(p(x)) = \begin{cases} \dfrac{1}{\deg p(x)} \\ 1 \text{ if } \deg p(x) = 0 \end{cases} .$$

η is the semigroup fuzzy linear algebra.

Define $\eta_1$ : V → [0, 1] as

$$\eta_1(p(x)) = \begin{cases} \dfrac{1}{\deg p(x)} \\ 0 \text{ if } \deg p(x) = 0 \end{cases}$$

then also $V\eta_1$ is a semigroup fuzzy linear algebra.

We see as in case of fuzzy set vector spaces and semigroup fuzzy vector spaces the notion of fuzzy set linear algebra and semigroup fuzzy linear algebra are also identical. This sort of making them identical using fuzzy tool will find its use in certain applications. We shall define such structures as fuzzy equivalent structures.

We see fuzzy set vector spaces and semigroup fuzzy vector spaces are fuzzy equivalent structures though set vector spaces are distinctly different from semigroup vector spaces. Like wise set linear algebras and semigroup linear algebras are fuzzy equivalent although set linear algebras are different from semigroup linear algebras.

Now we proceed on to define group fuzzy vector spaces and group fuzzy linear algebras.

**DEFINITION 3.9:** *Let V be a group linear algebra over the group G. Let $\eta : V \to [0, 1]$ such that*

$$\begin{aligned} \eta(a + b) &\geq \min(\eta(a), \eta(b)) \\ \eta(-a) &= \eta(a) \\ \eta(0) &= 1 \\ \eta(ra) &\geq \eta(a) \text{ for all } a, b \in V \text{ and } r \in G. \end{aligned}$$

*We call $V\eta$ the group fuzzy linear algebra.*

We illustrate this by an example.

***Example 3.18:*** Let $V = Z \times Z \times Z$ be the group linear algebra.

Define $\eta : V \to [0, 1]$ by

$$\eta(a) = 1 - \frac{1}{|a|}$$

for every $a \in Z$

$$\eta(0) = 1.$$

$\eta V$ is the group fuzzy linear algebra.

It is pertinent to mention here that we have not so far defined group fuzzy vector spaces. We first mention group vector spaces are fuzzy equivalent with set vector spaces and semigroup vector spaces. However for the sake of completeness we just define group fuzzy vector spaces.

**DEFINITION 3.10:** *Let V be a group vector space over the group G. $\eta : V \to [0, 1]$ is such that $\eta(ra) \geq \eta(a)$ for all $r \in G$ and $a \in V$. We call $V\eta$ or $\eta V$ to be the group fuzzy vector space.*

We see from the very definition the group vector spaces are fuzzy equivalent with set vector spaces.

Now we give an example of this concept.

***Example 3.19:*** Let V = Z [x] be a group vector space over the group G. Define $\eta : V \to [0, 1]$ by

$$\eta(p(x)) = \frac{1}{\deg p(x)}$$

and

$$\eta\ (\text{constant}) = 0.$$

$V\eta$ is a group fuzzy vector space.

In the same example if we view Z [x] to be a group under addition. Clearly Z [x] = V can be viewed as a group linear algebra over the group Z. $\eta$ defined above is such that $\eta V$ is a group fuzzy linear algebra.

**DEFINITION 3.11:** *Let V be a group linear algebra over the group G. Let $W \subset V$, where W is a subgroup of V and W is a group linear subalgebra over the group G. $\eta : W \to [0, 1]$ such that*

$$\begin{aligned} \eta\,(a + b) &\geq \min(\eta\,(a), \eta\,(b)) \\ \eta\,(a) &= \eta\,(-a) \\ (0) &= 1 \\ \eta\,(ra) &\geq r\eta\,(a) \end{aligned}$$

*for all $a, b \in W$ and $r \in G$; we call $W\eta$ or $\eta W$ to be the group fuzzy group linear subalgebra.*

We illustrate this by the following.

***Example 3.20:*** Let V = Z × Z × Z × Z × Z be a group linear algebra over the group G = Z. W = 3Z × {0} × 5Z × {0} ×{0} be the group linear subalgebra over Z.
Define η: W → [0, 1] by

$$\eta(x, y, z, \omega, t) = \begin{cases} 0 & \text{if} \quad x+y+z+\omega+t=0 \\ \dfrac{1}{x+y+z+\omega+t} & \text{if} \quad x+y+z+\omega+t \neq 0 \end{cases};$$

ηW is a fuzzy group linear subalgebra.

The importance of this structure is that we do not demand for a field or any other perfect nice structure to work with. Even a set will do the work for we ultimately see when we define fuzzy vector spaces the field does not play any prominent role. Also we see the group linear algebra is the same as ordinary vector space, when they are made into respective fuzzy structures. In fact these two structures are basically fuzzy equivalent. Any one will like to work with least algebraic operations only. So as we have already mentioned set vector spaces happens to be the most generalized concept of ordinary vector spaces and it is easy to work with them.

Another advantage of working with these special vector spaces is we see most of them happen to be fuzzy equivalent with some other special space or the ordinary vector space. In certain models or study we may have meaning for the solution only when they are positive. In such circumstances we need not define the vector spaces over a field instead we can define it over the set S which is a subset of $Z^+ \cup \{0\}$ or over the semigroup $Z^+ \cup \{0\}$ or even just Z.

Further as our transformation to a fuzzy set up always demands only values from the positive unit interval [0, 1] these semigroup vector spaces or set vector spaces would be more appropriate than the ordinary vector spaces.

Further we see even in case of Markov process or Markov chain the transition probability matrix is a square matrix with entries which are non negative and the column sum adding up to one. So in such cases one can use set vector space where

$$V = \left\{ (a_{ij})_{nxn} \big/ a_{ij} \in [0, 1] \quad \text{with} \quad \sum_{i=1}^{n} a_{ik} = 1 \text{ for } 1 \leq k \leq n \right\}$$

is a set vector space over the set [0, 1]. So these new notions not only comes handy but involve lesser complication and lesser algebraic operations.

## Chapter Four

# SET BIVECTOR SPACES AND THEIR GENERALIZATION

In this chapter for the first time we define the notion of set bivector spaces and generalize them to set n vector spaces. We enumerate some of the properties. In fact these set n-vector spaces happens to be the most generalized form of n-vector spaces. They are useful in mathematical models which do not seek much abstract algebraic concepts.

**DEFINITION 4.1:** *Let $V = V_1 \cup V_2$ where $V_1$ and $V_2$ are two distinct set vector spaces defined over the same set S. That is $V_1 \not\subseteq V_2$ and $V_2 \not\subseteq V_1$ we may have $V_1 \cap V_2 = \phi$ or non empty. Then we call V to be a set bivector space over S.*

We illustrate this by the following examples.

***Example 4.1:*** Let $V = V_1 \cup V_2$ where $V_1 = Z_5 \times Z_5$ and

$$V_2 = \left\{ \begin{pmatrix} a & b \\ c & d \end{pmatrix} \middle| a, b, c, d, \in Z_5 \right\}$$

are set vector spaces over the set S = {0, 1}. V is a set bivector space over the set {0, 1} = S.

***Example 4.2:*** Let V = $V_1 \cup V_2$ where $V_1$ = {(1 1 1), (0 0 0), (1 1 0), (1 1 1 1 1), (0 0 0 0 0), (1 1 0 1 1), (1 1 0 0 0), (1 0 0)} and $V_2$ = {(0 1), (1 0), (0 0), (1 1 1 1), (0 0 0 0), (0 1 1 1), (1 1 1 1 1 1), (0 0 0 0 0 0), (1 0 0 0), (0 0 0 1)} be set vector spaces over the set S = {0, 1}. V = $V_1 \cup V_2$ is a set bivector space over the set S.

***Example 4.3:*** Let V = $V_1 \cup V_2$ where

$$V_1 = \left\{ \begin{pmatrix} a & a & a \\ a & a & a \end{pmatrix} \middle| \, a \in Z_{12} \right\}$$

and

$$V_2 = \left\{ \begin{pmatrix} a & a \\ a & a \\ a & a \end{pmatrix} \middle| \, a \in Z_{12} \right\}$$

be two set vector spaces over the set S = {0, 1}. V is a set bivector space over the set S.

Now we have seen that how a set bivector space is constructed from these examples.

***Example 4.4:*** Let

$$V_1 = \left\{ \begin{pmatrix} a & b \\ c & d \end{pmatrix} \middle| \, a, b, c, d \in Z_{12} \right\}$$

and

$$V_2 = \left\{ \begin{pmatrix} a & a \\ a & a \end{pmatrix} \middle| \, a \in Z_{12} \right\}$$

be set vector spaces over the set S = {0, 1}. Clearly V = $V_1 \cup V_2$ is not a set bivector space over S as $V_2 \subseteq V_1$. Thus we cannot say the union of two set vector spaces defined over the same set gives a set bivector space.

Now we proceed on to define the notion of set bivector bisubspaces of a set bivector space.

**DEFINITION 4.2:** *Let V = $V_1 \cup V_2$ be a set bivector space defined over the set S. A proper biset W = $W_1 \cup W_2$ ($W_1 \subset V_1$ and $W_2 \subset V_2$) such that $W_1$ and $W_2$ are distinct and contained in V is said to be a set bivector bisubspace of V (or set bivector subspace) if W is a set bivector space defined over S.*

We now illustrate situation by the following examples.

***Example 4.5:*** Let V = $V_1 \cup V_2$

$$= \left\{ \begin{pmatrix} a & b \\ c & d \end{pmatrix} \middle| a, b, c, d \in Z_{12} \right\} \cup \{(a\ b\ c\ d) \mid a, b, c, d \in Z_{12}\}$$

be a set bivector space over the set S = {0, 1}.
Let W = $W_1 \cup W_2$

$$= \left\{ \begin{pmatrix} a & a \\ a & a \end{pmatrix} \middle| a \in Z_{12} \right\} \cup \{(a\ a\ a\ a) \mid a \in Z_{12}\}$$

$\subseteq V_1 \cup V_2$ = V. W is a set bivector space over S = {0,1}.

Thus W is a set bivector subspace of V over S.

***Example 4.6:*** Let

$$V_1 = \left\{ \begin{pmatrix} a_1 & a_2 & a_3 \\ a_4 & a_5 & a_6 \end{pmatrix} \middle| a_i \in \{0, 2, 6, 8, 10, 12\} \subseteq Z_{14} \right\}$$

and

$$V_2 = \left\{ \begin{pmatrix} a & a & a \\ a & a & a \end{pmatrix} \middle| a_a \in Z_{14} \right\}.$$

$V = V_1 \cup V_2$ is a set bivector space over $S = \{0,1\}$. Clearly $V_1 \cap V_2 \neq \phi$ but $V_1$ and $V_2$ are distinct. Take

$$W_1 = \left\{ \begin{pmatrix} a & a & 0 \\ a & a & 0 \end{pmatrix} \middle| \, a \in (0,2,6,8,10,12) \right\} \subseteq V_1$$

and

$$W_2 = \left\{ \begin{pmatrix} a & a & a \\ 0 & 0 & 0 \end{pmatrix} \middle| \, a \in Z_{14} \right\} \subset V_2 .$$

$W = W_1 \cup W_2$ is a set bivector subspace (or set bivector subbispace) of V.

Suppose

$$P_1 = \left\{ \begin{pmatrix} a & b & c \\ a & a & 0 \end{pmatrix} \middle| a, b, c \in (0, 2, ..., 12) \subseteq Z_{14} \right\} \subseteq V_1$$

and

$$P_2 = \left\{ \begin{pmatrix} a & a & a \\ 0 & 0 & 0 \end{pmatrix} \middle| \, a \in Z_{14} \right\} \subseteq V_2$$

then also $P = P_1 \cup P_2$ is a set bivector bisubspace of V. Thus we can have several such set bivector subspaces of a given set bivector space V.

Now we define the bidimension and the generating biset of a set bivector space $V = V_1 \cup V_2$.

**DEFINITION 4.3:** *Let $V = V_1 \cup V_2$ be a set bivector space defined over the set S. Let $X = X_1 \cup X_2 \subset V_1 \cup V_2$, we say X is a bigenerating subset of V if $X_1$ is the generating set of the set vector space $V_1$ over S and $X_2$ is the generating set of the set vector space $V_2$ over S.*

The number of elements in $X = X_1 \cup X_2$ is the bidimension of V and is denoted by $(|X_1|; |X_2|)$ or $|X_1| \cup |X_2|$.

We shall illustrate this definition by some examples.

***Example 4.7:*** Let $V = V_1 \cup V_2$ = {(111), (000), (100), (010), (001)} ∪ {(1111), (0000), (1110), (1000)} be the set bivector space over the set S = {0,1}. Take X = {(111), (100), (010), (001)} ∪ {(1111), (1110), (1000)} $\subseteq V_1 \cup V_2$ is the generating bisubset of V over the set S. Clearly dim V = (4, 3) or (4 ∪ 3).

***Example 4.8:*** Let

$$V = V_1 \cup V_2 = \left\{ \begin{pmatrix} a & a \\ a & a \end{pmatrix} \middle| a \in Z \right\} \cup \{Z \times Z \times Z\}$$

be the set bivector space defined over the set S = Z.

$$X = \left\{ \begin{pmatrix} 1 & 1 \\ 1 & 1 \end{pmatrix} \right\} \cup \{(111), (100), (010), (001), (110), (a\ b\ c) \ldots\}$$
$$= X_1 \cup X_2$$

where $X_2$ is an infinite subset of $V_2$, this alone can generate V so bidimension of V is infinite i.e., bidimension of V = $\{1 \cup \infty\}$ or $(1, \infty)$.

***Example 4.9:*** Let $V = V_1 \cup V_2 = \{Z_{10} \times Z_{10} \times Z_{10}\} \cup$ {(a a a a a) | a ∈ {0, 2, 4, 6, 8} a proper subset of $Z_{10}$}. V is a set bivector space over $Z_{10}$. Prove bidimension of V is finite over $Z_{10}$.

Now we proceed onto define the notion of set bilinear algebra or equivalently we can call it as set linear bialgebra.

**DEFINITION 4.4:** *Let $V = V_1 \cup V_2$ be such that $V_1$ is a set linear algebra over the set S and $V_2$ is also a set linear algebra over S. Further $V_1 \neq V_2$, $V_1 \not\subseteq V_2$ or $V_2 \not\subseteq V_1$. Then $V = V_1 \cup V_2$ is defined to be the set linear bialgebra over the set S.*

We now illustrate this situation by some examples.

***Example 4.10:*** Let

$$V = V_1 \cup V_2 = \left\{ \begin{pmatrix} a & a & a \\ a & a & a \end{pmatrix} \middle| a \in Z_{16} \right\} \cup \{Z_{16} \times Z_{16}\}$$

be the set linear bialgebra over the set $Z_{16}$.

***Example 4.11:*** Let

$$\begin{aligned} V &= \{(000), (010), (100), (001), (110), (011), (101), (111)\} \\ &\quad \cup \left\{ \begin{pmatrix} a & b \\ c & d \end{pmatrix} \middle| a, b, c, d \in Z_2 = \{0, 1\} \right\} \\ &= V_1 \cup V_2\,. \end{aligned}$$

V is a set linear bialgebra over the set S = {0 1}.

***Example 4.12:*** Let

$$\begin{aligned} V &= V_1 \cup V_2 \\ &= \{Z\} \cup \{Z^+ \times Z^+ \times Z^+\} \end{aligned}$$

be the set linear bialgebra over the set $S = 2Z^+$, here $V_1$ and $V_2$ are set linear algebras over $S = 2Z^+$.

***Example 4.13:*** Let

$$\begin{aligned} V &= V_1 \cup V_2 \\ &= \{(a\ b) \,/\, a, b \in \{0\ 1\}\} \cup \{1110\}, \{0000\}, (0011)\}; \end{aligned}$$

V is not a set linear bialgebra over the set S = {0, 1}. V is only a set bivector space over the set S because $V_2$ is not closed under the operation '+'.

In view of this we have the following result which is left for the reader to prove.

*Result:* Every set linear bialgebra is a set bivector space but all set bivector spaces need not in general be a set linear bialgebras.

The example 4.13 is one such algebraic structure.

Now we proceed on to define the notion of set linear subbialgebra or equivalently the notion of set bilinear subalgebra.

**DEFINITION 4.5:** *Let $V = V_1 \cup V_2$ be a set linear bialgebra over the set S. If $W \subset V$ i.e., $W = W_1 \cup W_2 \subset V_1 \cup V_2$ ($W_i \subset V_i$, i = 1, 2) is a set linear bialgebra over the set S then we call W to be the set linear subbialgebra of V.*

We now illustrate this by some examples.

***Example 4.14:*** Let

$$\begin{aligned} V &= V_1 \cup V_2 \\ &= \left\{ \begin{pmatrix} a & a \\ a & a \end{pmatrix} \middle| a \in Z \right\} \cup \left\{ \begin{pmatrix} a & a & a & a \\ a & a & a & a \end{pmatrix} \middle| a \in Z \right\}; \end{aligned}$$

V is a set linear bialgebra over Z. Take

$$\begin{aligned} W &= W_1 \cup W_2 \\ &= \left\{ \begin{pmatrix} a & a \\ 0 & 0 \end{pmatrix} \middle| a \in Z \right\} \cup \left\{ \begin{pmatrix} a & 0 & a & 0 \\ a & 0 & a & 0 \end{pmatrix} \middle| a \in Z \right\} \end{aligned}$$

$\subseteq V_1 \cup V_2 = V$, W is a set linear subalgebra of V over S.

***Example 4.15:*** Let $V = V_1 \cup V_2 = Z[x] \cup Q$ be a set linear bialgebra over the set S = 2Z. Suppose $W = W_1 \cup W_2$ = {all polynomial of even degree with coefficient from 2Z} $\cup$ {5Z} $\subseteq V_1 \cup V_2 = V$. Clearly W is a set linear bisubalgebra over the set S = 2Z.

***Example 4.16:*** Let $V = V_1 \cup V_2$

$$= \left\{ \begin{pmatrix} x & 0 \\ 0 & x \end{pmatrix} \middle| x \in Z_{12} \right\} \cup \left\{ \begin{pmatrix} x & y & z \\ 0 & y & 0 \end{pmatrix} \middle| x, y, z \in Z_{12} \right\}.$$

V is a set linear bialgebra over $S = Z_{12}$. Take

$$\begin{aligned} W &= W_1 \cup W_2 \\ &= \left\{ \begin{pmatrix} x & 0 \\ 0 & 0 \end{pmatrix} \middle| x \in Z_{12} \right\} \cup \left\{ \begin{pmatrix} x & x & x \\ 0 & 0 & 0 \end{pmatrix} \middle| x \in Z_{12} \right\} \end{aligned}$$

$\subseteq V_1 \cup V_2 = V$.

W is a set linear bisubalgebra of V over the set $S = Z_{12}$.

**DEFINITION 4.6:** *Let $V = V_1 \cup V_2$ be a set linear bialgebra over the set S. Let $X = X_1 \cup X_2 \subset V_1 \cup V_2 = V$, if $X_1$ is a generating set of $V_1$ and $X_2$ is a generating set of $V_2$ then $X = X_1 \cup X_2$ is the generating subset of V. The bidimension of V is the cardinality of $(|X_1|, |X_2|)$.*

We illustrate this situation by some examples.

***Example 4.17:*** Let $V = V_1 \cup V_2 = \{Z\} \cup \{(a\ a\ a) \mid a \in Z\}$ be a set linear bialgebra over the set $S = Z$. Let $X = X_1 \cup X_2 = \{1\} \cup \{1\ \ 1\ \ 1\} \subseteq V_1 \cup V_2$. X is the generating biset of V. The bidimension of V is (1, 1).

***Example 4.18:*** Let

$$\begin{aligned} V &= V_1 \cup V_2 \\ &= \left\{ \begin{pmatrix} a & a \\ a & a \end{pmatrix} \middle| a \in Z \right\} \cup \{(a\ b\ c) \mid a, b, c \in Z\} \end{aligned}$$

be the set linear bialgebra over the set S = Z. Let

$$\begin{aligned} X &= \begin{pmatrix} 1 & 1 \\ 1 & 1 \end{pmatrix} \cup \{(100), \{010), (001)\} \\ &= X_1 \cup X_2 \subset V_1 \cup V_2 \\ &= V. \end{aligned}$$

The bidimension of V is $\{1\} \cup \{3\}$ or (1, 3).

***Example 4.19:*** Let $V = V_1 \cup V_2 = Z[x] \cup Z \times Z \times Z$ be a set linear bialgebra over the set Z. Let $X = \{1, x, x^2, \ldots, x^n, \ldots\} \cup \{(100), (010), (001)\}$ generates V and the bidimension of V is $\{\infty\} \cup \{3\} = (\infty, 3)$.

***Example 4.20:*** $V = V_1 \cup V_2 = Z[x] \cup Z \times Z \times Z$ as a set bivector space over the set Z is of bidimension $(\infty,\infty)$.

This is the marked difference between the set linear bialgebras and set bivector spaces.

Now we proceed onto define the notion of semigroup bivector spaces, biset bivector spaces and bisemigroup bivector spaces and illustrate them by examples.

**DEFINITION 4.7:** *Let $V = V_1 \cup V_2$ be such that $V_1$ is a set vector space over the set $S_1$ and $V_2$ be a set vector space over the set $S_2$. $S_1 \neq S_2$; $S_1 \not\subseteq S_2$ and $S_2 \not\subseteq S_1$ we define $V = V_1 \cup V_2$ to be the biset bivector space over the biset $S_1 \cup S_2$.*

Now we will illustrate this definition by some examples.

***Example 4.21:*** Let $V = V_1 \cup V_2 = Z \times Z \times Z \cup Z_{12} \times Z_{12} \times Z_{12} \times Z_{12}$ be a biset bivector space over the biset $S = Z \cup Z_{12}$.

***Example 4.22:*** Suppose

$$\begin{aligned} V &= V_1 \cup V_2 \\ &= Z_{12}[x] \cup \left\{ \left. \begin{pmatrix} a & a & a \\ a & a & a \end{pmatrix} \right| a \in Z_{10} \right\}; \end{aligned}$$

take the biset $S = Z_{12} \cup Z_{10}$; then V is the biset bivector space over the biset $S = Z_{12} \cup Z_{10}$. i.e., $V_1$ is a set vector space over the set $Z_{12}$ and $V_2$ is a set vector space over the set $Z_{10}$.

It is interesting as well as important to observe that set bivector spaces and biset bivector spaces are two different and distinct notions. They will find their applications in different sets of mathematical models.

***Example 4.23:*** Let

$$V = V_1 \cup V_2$$

$$= \{(111), (000), (11111), (00000), (110), (100), (001), (10100), (0000), (1100), (1010)\}$$

$$\cup \left\{ \begin{pmatrix} a & b \\ c & d \end{pmatrix} \middle| a, b, c, d \in Z_5 \right\}.$$

$V_1$ is a set vector space over the set $S_1 = Z_2 = \{0, 1\}$ and $V_2$ is a set vector space over the set $S_2 = Z_5$. Thus $V = V_1 \cup V_2$ is a biset bivector space over the biset $S = S_1 \cup S_2 = Z_2 \cup Z_5$.

***Example 4.24:*** Let

$$V = V_1 \cup V_2$$
$$= \{(1111), (0011), (0000), (1000), (11), (01), (00)\}$$

$$\cup \left\{ \begin{pmatrix} a & a & a \\ b & b & b \end{pmatrix} \middle| a, b \in Z_3 \right\}$$

be a biset bivector space over the biset $S = S_1 \cup S_2 = Z_2 \cup Z_3$.

Now we define the notion of biset bivector subspaces.

**DEFINITION 4.8:** *Let $V = V_1 \cup V_2$ be a biset bivector space over the biset $S = S_1 \cup S_2$. Let $W = W_1 \cup W_2 \subseteq V_1 \cup V_2$; if W is a biset bivector space over S then we call W to be the biset bivector subspace of V over the biset $S = S_1 \cup S_2$.*

***Example 4.25:*** Let $V = V_1 \cup V_2$

$$= \{(000000), (111000), (110110), (111), (000), (100), (011)\} \cup \left\{ \begin{pmatrix} a & a & a \\ a & a & a \end{pmatrix} \middle| a \in Z_4 \right\}$$

be the biset bivector space over the biset $S = S_1 \cup S_2 = \{0, 1\} \cup Z_4$.
Let

$$\begin{aligned} W &= W_1 \cup W_2 \\ &= \{(000000)\ (111000)\} \cup \left\{ \begin{pmatrix} a & a & a \\ 0 & 0 & 0 \end{pmatrix} \middle| a \in Z_4 \right\} \\ &\subseteq V_1 \cup V_2, \end{aligned}$$

W is a biset bivector subspace of V over the biset $S = S_1 \cup S_2$.

***Example 4.26:*** Let

$$V = V_1 \cup V_2 = \left\{ \begin{pmatrix} a & a \\ a & a \end{pmatrix} \middle| a \in Z_4 \right\} \cup \{Z_5[x]\};$$

$\{Z_5[x]$ i.e., all polynomial in the variable x with coefficients from $Z_5\}$ be the biset bivector space over the biset $S = S_1 \cup S_2 = Z_4 \cup Z_5$ . Let

$$\begin{aligned} W &= W_1 \cup W_2 \\ &= \left\{ \begin{pmatrix} a & 0 \\ 0 & a \end{pmatrix} \middle| a \in Z_4 \right\} \cup \\ &\quad \{\text{all polynomials of degree 1 with coefficient from } Z_5\} \\ &\subseteq V_1 \cup V_2 \\ &= V; \end{aligned}$$

W is a biset bivector subspace of V over the biset $S = Z_4 \cup Z_5$ .

***Example 4.27:*** Let

$$V = V_1 \cup V_2 = \{Z_6 \times Z_6\} \cup \left\{ \begin{pmatrix} a & b & c \\ d & e & f \end{pmatrix} \middle| a, b, c, d, e, f \in Z_7 \right\}$$

be the biset bivector space over the biset $S = S_1 \cup S_2 = \{0\ 2\ 4\} \cup Z_7$. Take

$$W = W_1 \cup W_2 = \{\{0, 3\} \times Z_6\} \cup \left\{ \begin{pmatrix} a & a & a \\ 0 & 0 & 0 \end{pmatrix} \middle| a \in Z_7 \right\} \subset V_1 \cup V_2 \,;$$

W is the biset bivector subspace over the biset $S = S_1 \cup S_2$ .

**DEFINITION 4.9:** *Let $V = V_1 \cup V_2$ be a biset bivector space over the biset $S = S_1 \cup S_2$. If $X = X_1 \cup X_2 \subset V_1 \cup V_2$ is such that $X_1$ generates $V_1$ as a set vector space over the set $S_1$ and $V_2$ is generated by the set $X_2$ over the set $S_2$ then we say the biset $X_1 \cup X_2$ is the bigenerator of the biset bivector space $V = V_1 \cup V_2$ over the biset $S = S_1 \cup S_2$. The bicardinality of $X = X_1 \cup X_2$ denoted by $(|X_1|, |X_2|)$ gives the bidimension of V over S.*

We illustrate this by the following examples.

***Example 4.28:*** Let $V = V_1 \cup V_2$ be a biset bivector space over the biset $S = S_1 \cup S_2$ where $V_1 = Z_5 \times Z_5 \times Z_5$ with $S_1 = Z_5$ and

$$V_2 = \left\{ \begin{pmatrix} a & b & c \\ d & e & f \end{pmatrix} \middle| a, b, c, d, e, f \in Z_2 = \{0, 1\} \right\}$$

and $S_2 = Z_2$. Take

$$X = X_1 \cup X_2$$

$$= \{(a, b) / a, b \in Z_5)$$

$$\cup \left\{ \begin{pmatrix} 1 & 1 & 1 \\ 1 & 1 & 1 \end{pmatrix}, \begin{pmatrix} 1 & 1 & 1 \\ 1 & 1 & 0 \end{pmatrix}, \begin{pmatrix} 1 & 1 & 1 \\ 0 & 0 & 1 \end{pmatrix} \right.$$

$$\begin{pmatrix} 1 & 1 & 1 \\ 1 & 0 & 1 \end{pmatrix}, \begin{pmatrix} 1 & 1 & 1 \\ 0 & 1 & 1 \end{pmatrix}, \begin{pmatrix} 1 & 1 & 1 \\ 0 & 1 & 0 \end{pmatrix},$$

$$\left. \begin{pmatrix} 1 & 1 & 1 \\ 1 & 0 & 0 \end{pmatrix}, \begin{pmatrix} 1 & 1 & 1 \\ 0 & 0 & 0 \end{pmatrix}, \begin{pmatrix} 0 & 0 & 0 \\ 1 & 1 & 1 \end{pmatrix}, \ldots \right\}$$

X is finite biset and bidimension of V = $\{24 \cup 63\} = (24, 63)$.

***Example 4.29:*** Let

$$V = V_1 \cup V_2 = \left\{ \begin{pmatrix} a & a & a \\ a & a & a \end{pmatrix} \middle| a \in Z \right\} \cup \{Z_3 \times Z_3\}$$

be a biset bivector space over the biset $S = S_1 \cup S_2 = Z \cup Z_3$ . Take

$$X = \left\{ \begin{pmatrix} 1 & 1 & 1 \\ 1 & 1 & 1 \end{pmatrix} \right\} \cup \{(11), (12), (10), (01)\} = X_1 \cup X_2 ,$$

X is the bigenerating biset of the biset bivector space V. The bidimension of V is $\{1\} \cup \{4\}$.

***Example 4.30:*** Let = $V_1 \cup V_2 = \{Z_5 \times Z_5\} \cup (Z \times Z)$ be the biset bivector space over the biset $S = S_1 \cup S_2 = Z_5 \cup Z$. Take $X = X_1 \cup X_2 = \{(11), (10), (01), (12), (13), (14)\} \cup \{(a, b) / a\ b \in Z\}$ bigenerates V. Clearly $|X_1| = 6$ and $|X_2| = \infty$ so the bidimension of V is $(6, \infty)$.

Thus even if one of the set vector spaces in the biset bivector space V is of infinite dimension we say the biset bivector space to be of bidimension infinity. Only if both $V_1$ and $V_2$ are of finite dimension we say V, the biset bivector space V is of finite bidimension.

Now we proceed on to define the notion of biset bilinear algebra or equivalently biset linear bialgebra.

**DEFINITION 4.10:** *Let $V = V_1 \cup V_2$, if $V_1$ is a set linear algebra over the set $S_1$ and $V_2$ a different set linear algebra on the set $S_2$ ($S_1 \neq S_2$, $S_1 \not\subseteq S_2$, $S_2 \not\subseteq S_1$) ($V_1 \neq V_2$, $V_1 \not\subseteq V_2$ or $V_2 \not\subseteq V_1$) then we call $V = V_1 \cup V_2$ to be the biset bilinear algebra over the biset S $= S_1 \cup S_2$ .*

We illustrate this definition by examples.

***Examples 4.31:*** Let

$$V = V_1 \cup V_2 = Z[x] \cup \left\{ \begin{pmatrix} a & a \\ a & a \end{pmatrix} \middle| a \in Z_5 \right\}$$

be the biset bilinear algebra over the biset $S = Z \cup Z_5$. Clearly Z[x] is the set linear algebra over Z and

$$\left\{ \begin{pmatrix} a & a \\ a & a \end{pmatrix} \middle| a \in Z_5 \right\}$$

is the set linear algebra over the set $Z_5$.

***Example 4.32:*** Let $V = V_1 \cup V_2 =\{Z_7 \times Z_7\} \cup \{(a\ a\ a\ a) / a \in Z_6\}$; $V_1$ is a set linear algebra over the set $Z_7$ and $V_2$ is a set linear algebra over the set $Z_6$. Thus $V = V_1 \cup V_2$ is the biset bilinear algebra over the biset $S = Z_7 \cup Z_6$.

***Example 4.33:*** Let $V = V_1 \cup V_2 = \{(a, b) / a, b \in Z\} \cup \{Z_9[x] =$ all polynomials in the variable x with coefficients from the set $Z_9\}$. Take $S = Z \cup Z_9 = S_1 \cup S_2$. Now V is the biset bilinear algebra over the biset $S = S_1 \cup S_2$.

***Example 4.34:*** Let

$$\begin{aligned} V &= V_1 \cup V_2 \\ &= \{Z_5 \times Z_5 \times Z_5\} \cup \left\{\begin{pmatrix} a & a & a \\ b & b & b \end{pmatrix} \middle| a, b \in Z_2 \right\}. \end{aligned}$$

Take the biset S = = $S_1 \cup S_2 = Z_5 \cup Z_2$. Clearly V is a biset bilinear algebra over the biset $S = Z_5 \cup Z_2$.

Now we proceed on to define the bidimension of a biset bilinear algebra over the biset $S = S_1 \cup S_2$.

**DEFINITION 4.11:** *Let $V = V_1 \cup V_2$ be a biset linear bialgebra over the biset $S = S_1 \cup S_2$. Let $W = W_1 \cup W_2 \subset V_1 \cup V_2$, if W is a biset bilinear algebra over the biset $S = S_1 \cup S_2$ then we call W to be the biset bilinear subalgebra of V over the biset $S = S_1 \cup S_2$.*

We illustrate this by some simple examples.

***Example 4.35:*** Let

$$\begin{aligned} V &= V_1 \cup V_2 \\ &= \left\{\begin{pmatrix} a & a & a \\ a & a & a \end{pmatrix} \middle| a \in Z_5 \right\} \cup \{Z_4 \times Z_4\} \end{aligned}$$

be the biset bilinear algebra over the biset $S = S_1 \cup S_2 = Z_5 \cup Z_4$. Let

$$\begin{aligned} W &= W_1 \cup W_2 \\ &= \left\{\begin{pmatrix} a & a & a \\ 0 & 0 & 0 \end{pmatrix} \middle| a \in Z_5 \right\} \cup \{\{0,2\} \times \{0,2\}\} \\ &\subseteq V_1 \cup V_2; \end{aligned}$$

W is a biset bilinear subalgebra of V over the biset $S = S_1 \cup S_2$ .

***Example 4.36:*** Let $V = V_1 \cup V_2 = Z_7[x] \cup \{(000), (111), (100), (001), (010), (110), (101), (011)\}$ be the biset bilinear algebra over the biset $S = Z_7 \cup Z_2$ . Take $W = W_1 \cup W_2 = \{$all polynomials of degree 2 with coefficient from $Z_7\} \cup \{000), (111)\} \subset V_1 \cup V_2$; W is a biset bilinear subalgebra of V over the biset $S = Z_7 \cup Z_2$.

***Example 4.37:*** Let

$$\begin{aligned} V &= V_1 \cup V_2 \\ &= \{Z_2 \times Z_2 \times Z_2 \times Z_2\} \cup \left\{ \begin{pmatrix} a & b \\ c & d \end{pmatrix} \middle| a, b, c, d \in Z_3 \right\}. \end{aligned}$$

V is a biset bilinear algebra over the biset $Z_2 \cup Z_3$.
Take

$$\begin{aligned} W &= W_1 \cup W_2 \\ &\subseteq V_1 \cup V_2 \end{aligned}$$

where

$$W_1 = Z_2 \times \{0\} \times Z_2 \times \{0\}$$

and

$$W_2 = \left\{ \begin{pmatrix} a & 0 \\ 0 & d \end{pmatrix} \middle| a, d \in Z_3 \right\}.$$

W is the biset bilinear subalgebra of V over the set $Z_2 \cup Z_3$.

Now having defined the substructure we now proceed on to define the notion of the bigenerating set and bidimension of the biset bilinear algebra.

**DEFINITION 4.12:** *Let $V = V_1 \cup V_2$ be a biset bilinear algebra defined over the biset $S = S_1 \cup S_2$. Let $X = X_1 \cup X_2 \subset V_1 \cup V_2$ where $X_1$ generates $V_1$ as a set linear algebra over $S_1$ and $X_2$ generates $V_2$ as a set linear algebra over $S_2$. Clearly $X = X_1 \cup X_2$ bigenerates V and the bidimension of V is ($|X_1|$; $|X_2|$).*

We give some examples to illustrate this concept.

***Example 4.38:*** Let

$$\begin{aligned} V &= V_1 \cup V_2 \\ &= \left\{ \begin{pmatrix} a & a & a \\ a & a & a \end{pmatrix} \middle| a \in Z \right\} \cup Z_2 \times Z_2 \times Z_2 \end{aligned}$$

be the biset bilinear algebra over the biset $S = Z \cup Z_2 = S_1 \cup S_2$.

Let

$$\begin{aligned} X &= X_1 \cup X_2 \\ &= \left\{ \begin{pmatrix} 1 & 1 & 1 \\ 1 & 1 & 1 \end{pmatrix} \right\} \cup \{(001), (000), (100), (010)\} \\ &\subseteq V_1 \cup V_2, \end{aligned}$$

we see X bigenerates V and bidimension of V is $\{1\} \cup \{3\}$ or $\{1, 3\}$.

***Example 4.39:*** Let

$$\begin{aligned} V &= V_1 \cup V_2 \\ &= \left\{ \begin{pmatrix} a & a & a \\ a & a & a \end{pmatrix} \middle| a \in Z_{12} \right\} \cup \{Z_6 \times Z_6\} \end{aligned}$$

be a biset bilinear algebra over the biset $S = Z_{12} \cup Z_6$ . Take

$$\begin{aligned} X &= X_1 \cup X_2 \\ &= \left\{ \begin{pmatrix} 1 & 1 & 1 \\ 1 & 1 & 1 \end{pmatrix} \right\} \cup \{(11), (10)\} \end{aligned}$$

is the bigenerator of V. The bidimension of V is $\{1\} \cup \{2\}$ or $\{1, 2\}$.

Now having seen the bidimension we wish to mention that for the same set V treated as a biset bilinear algebra and as a biset bivector space over the same biset $S = S_1 \cup S_2$, their bidimensions are distinct and not the same.

Thus in certain cases it is advantageous to work with biset bilinear algebra for it will make the cardinality of the bidimension relatively small when compared with the biset bivector space.

In spaces where it is possible we can make use of biset bilinear algebra instead of biset bivector spaces.

Now we proceed on to define the notion of semigroup bivector spaces.

**DEFINITION 4.13:** *Let $V = V_1 \cup V_2$ where $V_1$ is a semigroup vector space over the semigroup S.*

*If $V_2$ is also a semigroup vector space over the same S and if $V_1$ and $V_2$ are distinct ($V_1 \neq V_2$, $V_1 \not\subseteq V_2$ and $V_2 \not\subseteq V_1$) then we say $V = V_1 \cup V_2$ to be the semigroup bivector space over the semigroup S.*

We illustrate this by the following examples.

***Example 4.40:*** Let

$$\begin{aligned} V &= V_1 \cup V_2 \\ &= \{Z^+ \times Z^+ \times 2Z^+\} \cup \left\{ \begin{pmatrix} a & b & c \\ d & e & f \end{pmatrix} \middle| a, b, c, d, e, f \in 2Z^+ \right\} \end{aligned}$$

be a semigroup bivector space over the semigroup $S = 2Z^+$, we see $V_1$ is a semigroup vector space over the semigroup $2Z^+ = S$ and $V_2$ is also a semigroup vector space over the same semigroup $S = 2Z^+$.

***Example 4.41:*** Let $V = \{(111), (000), (100), (001)\} \cup \{(a, b) / a, b \in Z_2 = \{0, 1\}\} = V_1 \cup V_2$. V is a semigroup bivector space over the semigroup $Z_2 = \{0, 1\} = S$.

***Example 4.42:*** Let

$$V = \{(a, b, c) \mid a, b, c \in Z^+\} \cup \left\{ \begin{pmatrix} a & c \\ 0 & a \end{pmatrix} \middle| a, c \in 2Z^+ \right\}$$

$$= V_1 \cup V_2 .$$

$V = V_1 \cup V_2$ is a semigroup bivector space over the semigroup $S = Z^+$.

***Example 4.43:*** Let

$$V = \{(000), (00), (01), (111), (001), (011)\} \cup \{(1111), (0000), (0101), (1101), (000), (111), (11111), (00000), (10101)\}$$

$$= V_1 \cup V_2 .$$

V is a semigroup bivector space over the semigroup $Z_2 = \{0,1\}$ under addition modulo 2. $V_1 \cap V_2 = \{(000), (111)\} \neq \phi$. But $V_1 \not\subseteq V_2$ and $V_2 \not\subseteq V_1$.

Now we proceed on to define the new notion of semigroup bivector subspace of a semigroup bivector space V.

**DEFINITION 4.14:** *Let $V = V_1 \cup V_2$ be a semigroup bivector space over the semigroup S. Let $W = W_1 \cup W_2 \subset V_1 \cup V_2 = V$ be a proper biset of V, if W is a semigroup bivector space over S then we call W to be the semigroup bivector subspace of V over the semigroup S. Clearly $W_1 \neq W_2$ and $W_1 \not\subseteq W_2$ and $W_2 \not\subseteq W_1$ with $W_1 \subseteq V_1$ and $W_2 \subseteq V_2$.*

We now illustrate this definition by some examples.

***Example 4.44:*** Let $V = V_1 \cup V_2$ be a semigroup bivector space over the semigroup $S = Z^+$. Let

$$V_1 = \left\{ \begin{pmatrix} a & a & a \\ a & a & a \end{pmatrix} \middle| a \in Z^+ \cup \{0\} \right\}$$

and

$$V_2 = \{3Z^+ \times 5Z^+\}.$$

$V_1$ is a semigroup vector space over the semigroup $S = Z^+$ and $V_2$ is a semigroup vector space over the semigroup $S = Z^+$.

Thus V is a semigroup bivector space over the semigroup $S = Z^+$.

Take

$$W_1 = \left\{ \begin{pmatrix} a & a & a \\ a & a & a \end{pmatrix} \middle| a \in 3Z^+ \cup \{0\} \right\} \subseteq V_1$$

and

$$W_2 = 3Z^+ \times \{0\} \subseteq V_2.$$

$W = W_1 \cup W_2$ is a semigroup bivector space over the semigroup $S = Z^+$. W is a semigroup bivector subspace of $V = V_1 \cup V_2$ over the semigroup $S = Z^+$.

***Example 4.45:*** Let V = {(1110), (0000), (1010), (1000), (00), (11), (10)} ∪ {(11111), (00000), (000), (111), (11011), (101)} = $V_1 \cup V_2$ be the semigroup bivector space over the semigroup $S = Z_2$. W = {(0000), (1111)} ∪ {(000), (101)} $\subseteq V_1 \cup V_2$ is a semigroup bivector subspace of V over $S = Z_2$.

Now we proceed on to define the bidimension and bigenerator of the semigroup bivector space.

**DEFINITION 4.15:** *Let $V = V_1 \cup V_2$ be a semigroup bivector space over the semigroup S.*

*Let $X = X_1 \cup X_2 \subseteq V_1 \cup V_2$, if $X_1$ generates the semigroup vector space $V_1$ over the semigroup S and $X_2$ generates the semigroup vector space $V_2$ over the semigroup S then, $X = X_1 \cup X_2$ is the bigenerator of the semigroup bivector space V over the semigroup S.*

*The bidimension of V is $|X_1| \cup |X_2|$ or $(|X_1|, |X_2|)$ over the semigroup S. If even one of $|X_1|$ or $|X_2|$ is infinite we say the bidimension of V is infinite.*

***Example 4.46:*** Let $V = V_1 \cup V_2$ = {(000), (111), (010), (001), (00), (01)} $\cup$ {(a a a a) / a $\in Z_2$ = {0, 1}} be the semigroup bivector space over the semigroup $S = Z_2$ = {0,1}. X = {(111), (010), (001), (01)} $\cup$ ((1111)}; bigenerates V over $S = Z_2$ = {0,1} so the bidimension of V is (4,1) or {4} $\cup$ {1) over S.

***Example 4.47:*** Let V = {(a b c) / a, b, c $\in Z^+$} $\cup$ {(a a a a a) / a $\in Z^+$} be the semigroup bivector space over the semigroup S = $Z^+$. Take X = {(a b c) / a, b, c $\in Z^+$} $\cup$ {(11111)} $\subseteq V_1 \cup V_2$ , X is a bigenerator of V and the bidimension of V over S is {$\infty$} $\cup$ {1} = {$\infty$, 1}. Thus V is an infinite bidimensional semigroup bivector space over S = $Z^+$.

Now we proceed on to define the notion of semigroup bilinear algebra over the semigroup.

**DEFINITION 4.16:** *Let $V = V_1 \cup V_2$ be such that $V_1$ is a semigroup linear algebra over the semigroup S and $V_2$ is a semigroup linear algebra over the semigroup S. with $V_1 \neq V_2$ , $V_1 \not\subseteq V_2$ and $V_2 \not\subseteq V_1$.*

*Then we call V to be the semigroup bilinear algebra over the semigroup S.*

We illustrate this by few examples.

***Example 4.48:*** Let

$$\begin{aligned} V &= V_1 \cup V_2 \\ &= \{(111), (000), (110), (101), (100), (010), (001), \\ &\quad (011)\} \cup \left\{ \begin{pmatrix} a & a \\ a & a \end{pmatrix} \middle| a \in Z_2 = \{0,1\} \right\}. \end{aligned}$$

V is a semigroup bilinear algebra over the semigroup $S = Z_2$ = {0,1}.

***Example 4.49:*** Let $V = V_1 \cup V_2 = \{Z_5 [x]\} \cup \{Z_5 \times Z_5 \times Z_5\}$. V is a semigroup bilinear algebra over the semigroup $Z_5$.

***Example 4.50:*** Let

$$V = V_1 \cup V_2 = \{Z^+ \times Z^+ \times Z^+\} \cup \left\{ \left. \begin{pmatrix} a & a & a \\ a & a & a \end{pmatrix} \right| a \in 2Z^+ \right\}.$$

V is a semigroup bilinear algebra over the semigroup $S = Z^+$.

We see all semigroup bilinear algebras defined over the semigroup are semigroup bivector spaces but a semigroup bivector space in general is not a semigroup bilinear algebra.

To this end we give an example.

***Example 4.51:*** Let $V = V_1 \cup V_2 = \{(111), (000), (11), (00)\} \cup \{(0000), (1111), (1101), (0110)\}$. V is a semigroup bivector space over the semigroup $S = Z_2 = \{0,1\}$.

Clearly V is not a semigroup bilinear algebra over $Z_2 = \{0,1\}$ as $V_1$ is not a semigroup under addition and $V_2$ is also not a semigroup under addition. Hence the claim.

It may so happen in $V = V_1 \cup V_2$ we may have $V_1$ to be a semigroup linear algebra over the semigroup S and $V_2$ is only a semigroup vector space, in such cases we define a new algebraic structure.

**DEFINITION 4.17:** *Let $V = V_1 \cup V_2$ be such that $V_1$ is a semigroup linear algebra over the semigroup S and $V_2$ is only a semigroup vector space over S with $V_1 \neq V_2$, $V_1 \not\subseteq V_2$ and $V_2 \not\subseteq V_1$. Then we call $V = V_1 \cup V_2$ to be a quasi semigroup bilinear algebra over S.*

We illustrate this by the following examples.

***Example 4.52:*** Let

$$\begin{aligned} V &= V_1 \cup V_2 \\ &= \{(000), (111), (01), (10), (11), (00), (100)\} \\ &\quad \cup \left\{ \left. \begin{pmatrix} a & a \\ a & a \end{pmatrix} \right| a \in Z_2 \right\} \end{aligned}$$

$V_1$ is only a semigroup vector space over the semigroup S = $Z_2$ = {0,1}. $V_2$ is a semigroup linear algebra over the semigroup S = $Z_2$ = {0, 1}; $V_1 \neq V_2$. So V is a quasi semigroup bilinear algebra over the semigroup S = $Z_2$ = {0,1}.

It is interesting to note that all semigroup bilinear algebras are quasi semigroup bilinear algebras but converse is never true. Also all quasi semigroup bilinear algebras are semigroup bivector spaces but the converse is not true.

***Example 4.53:*** Let

$$\begin{aligned} V &= \left\{ \begin{pmatrix} a & a \\ a & a \end{pmatrix}, \begin{pmatrix} b & b & b \\ b & b & b \end{pmatrix} \middle| a, b \in Z_2 \right\} \\ &\cup \{(110), (0000), (00000), (11111), (1101), (000)\} \\ &= V_1 \cup V_2\,. \end{aligned}$$

V is a semigroup bivector space over the semigroup S = $Z_2$ = {0, 1}. $V_1$ is only a semigroup vector space also $V_2$ is only a semigroup vector space over $Z_2$. So V = $V_1 \cup V_2$ is only a semigroup bivector space over $Z_2$ and never a quasi semigroup bilinear algebra over S = $Z_2$ = {0,1}.

***Example 4.54:*** Let

$$\begin{aligned} V &= V_1 \cup V_2 \\ &= \{(000), (111), (110), (111111), (000000), (111000), (101010)\} \\ &\cup \left\{ \begin{pmatrix} a & b \\ c & d \end{pmatrix} \middle| a, b, c, d \in Z_2 \right\}. \end{aligned}$$

$V_1$ is only a semigroup vector space over the semigroup S = $Z_2$ = {0, 1}. $V_2$ is a semigroup linear algebra over the semigroup S = $Z_2$ = {0, 1}. Thus V = = $V_1 \cup V_2$ is only a quasi semigroup bilinear algebra over $Z_2$ = {0, 1}.

Thus a quasi semigroup bilinear algebra can have quasi semigroup bilinear subalgebra as well as quasi semigroup bivector subspaces.

**DEFINITION 4.17:** *Let $V = V_1 \cup V_2$ be a quasi semigroup bilinear algebra over the semigroup S. Here $V_1$ is a semigroup linear algebra over S and $V_2$ is a semigroup vector space over S.*

*Let $W = W_1 \cup W_2 \subset V_1 \cup V_2$ where $W_1$ is a semigroup linear subalgebra of $V_1$ and $W_2$ is only a semigroup vector subspace of $V_2$. Then $W = W_1 \cup W_2$ is the quasi semigroup bilinear subalgebra of V.*

If $P = P_1 \cup P_2 \subset V_1 \cup V_2$ is such that $P_1$ is only a semigroup vector subspace of the semigroup linear algebra $V_1$ over S and $P_2$ is a semigroup vector subspace of the semigroup vector space $V_2$ then we call $P = P_1 \cup P_2$ to be the quasi semigroup bivector subspace of $V = V_1 \cup V_2$ over the semigroup S.

***Example 4.55:*** Let

$$\begin{aligned} V &= V_1 \cup V_2 \\ &= \{(000), (111), (100), (010), (001), (110), (011), \\ &\quad (101)\} \cup \left\{ \begin{pmatrix} a & a \\ a & a \end{pmatrix}, \begin{pmatrix} b & b & b \\ b & b & b \end{pmatrix}, \begin{pmatrix} a & 0 \\ a & 0 \end{pmatrix} \middle| a, b \in Z_2 \right\}. \end{aligned}$$

V is a quasi semigroup bilinear algebra over the semigroup $S = Z_2 = \{0,1\}$. Take

$$\begin{aligned} W &= W_1 \cup W_2 \\ &= \{(000), (100), (010)\} \cup \left\{ \begin{pmatrix} a & a \\ a & a \end{pmatrix}, \begin{pmatrix} a & 0 \\ a & 0 \end{pmatrix} \middle| a \in Z_2 \right\} \\ &\subseteq V_1 \cup V_2. \end{aligned}$$

W is only a semigroup bivector space over the semigroup $S = Z_2 = \{0, 1\}$. So W is a quasi semigroup bivector subspace of V over $Z_2 = \{0, 1\}$. Let

$$\begin{aligned} P &= P_1 \cup P_2 \\ &= \{(000), (111)\} \cup \left\{ \begin{pmatrix} a & a & a \\ a & a & a \end{pmatrix}, \begin{pmatrix} a & a \\ a & a \end{pmatrix} \middle| a \in Z_2 \right\} \\ &\subseteq V = V_1 \cup V_2. \end{aligned}$$

P is a quasi semigroup bilinear subalgebra over the semigroup $S = Z_2 = \{0,1\}$.

Now have seen the substructures of a quasi semigroup bilinear algebra we proceed on to define bidimension and bigenerating subset of the quasi semigroup bilinear algebra.

**DEFINITION 4.18:** *Let $V = V_1 \cup V_2$ be a quasi semigroup bilinear algebra over the semigroup S. Let $X = X_1 \cup X_2 \subseteq V_1 \cup V_2$ where $X_1$ generates the semigroup linear algebra $V_1$ and $X_2$ generates the semigroup vector space $V_2$ over S.*

*Then $X = X_1 \cup X_2$ is called the bigenerator of V and the bidimension of V is $(|X_1|, |X_2|)$ or $(|X_1| \cup |X_2|$.*

We illustrate this situation by some examples.

***Example 4.56:*** Let $V = V_1 \cup V_2 = \{(a, a, a) \mid a \in Z_2\} \cup \{(1\ 1\ 1), (0\ 0\ 0), (1\ 1\ 0), (1\ 1\ 1\ 0), (0\ 0\ 0\ 0), (1\ 1\ 0\ 0), (1\ 1\ 0\ 1), (1\ 1\ 0\ 0\ 1), (0\ 0\ 0\ 0\ 0), (1\ 1\ 1\ 0\ 1)\}$ be a quasi semigroup linear algebra over $Z_2$. Let $X = \{(1\ 1\ 1)\} \cup \{(1\ 1\ 1), (1\ 1\ 0), (1\ 1\ 1\ 0), (1\ 1\ 0\ 0), (1\ 1\ 0\ 1), (1\ 1\ 0\ 0\ 1), (1\ 1\ 1\ 0\ 1)\} \subset V_1 \cup V_2$, X is a bisubset of V which bigenerates V. The bidimension of X is $\{1\} \cup \{7\} = (1, 7)$.

***Example 4.57:*** Let $V = V_1 \cup V_2 = \{(1\ 1), (1\ 0), (0\ 0), (1\ 1\ 1), (0\ 0\ 0), (1\ 1\ 1\ 1\ 1), (0\ 0\ 0\ 0\ 0), (1\ 1\ 0\ 0\ 0), (0\ 1\ 1\ ), (1\ 0\ 1\ 0\ 1)\}$

$$\cup \left\{ \begin{pmatrix} a & a \\ a & a \end{pmatrix} \middle| a \in Z_2 \right\}$$

be a quasi semigroup linear algebra over the semigroup $Z_2 = \{0, 1\}$.

$$X = \{(1\ 1), (1\ 0), (1\ 1\ 1), (1\ 1\ 1\ 1\ 1), (1\ 1\ 0\ 0\ 0), (0\ 1\ 1), (1\ 0\ 1\ 01\ )\} \cup \left\{\begin{pmatrix} 1 & 1 \\ 1 & 1 \end{pmatrix}\right\}$$

is the bigenerator of V. The bidimension of V is $\{7\} \cup \{1\} = (7, 1)$.

We see from these examples the dimension of the semigroup linear algebra V is less than the dimension of semigroup vector space V. Same V for which '+' is taken and for the other '+' operation is not taken.

We illustrate this situation by an example.

***Example 4.58:*** Let V = {(1 1 1), (1 0 0), (0 1 0), (0 0 1), (1 1 0), (1 0 1), (0 1 1), (0 0 0)} be a semigroup linear algebra over the semigroup $Z_2 = \{0, 1\}$. Suppose V = {(1 1 1), (1 0 0), (0 1 0), (0 0 1), (1 1 0), (1 0 1), (0 1 1), (0 0 0)} be a semigroup vector space over $Z_2 = \{0,1\}$, the semigroup under addition. Dimension of V as a semigroup linear algebra is three given by the generating set X = {(1 0 0), (0 1 0), (0 0 1)}. The dimension of V as a semigroup vector space is 7 given by the generating set X = {(1 1 1), (1 0 0), (0 0 1) (0 1 0), (1 1 0), (0 1 1), (1 0 1)}. Thus we see dimension varies or the dimension is small when the structure is a semigroup linear algebra V and the dimension is large for the same V when it is a semigroup vector space.

Now we proceed onto define the new notion of group bivector spaces and group bilinear algebras.

**DEFINITION 4.19:** *Let $V = V_1 \cup V_2$ be such that $V_1 \neq V_2$, $V_1 \not\subseteq V_2$ and $V_2 \not\subseteq V_1$, $V_1$ and $V_2$ group vector spaces over the same group G, then we call V to be a group bivector space defined over the group G.*

We illustrate this by the following examples.

***Example 4.59:*** Let $V = V_1 \cup V_2$ where

$$V_1 = \left\{\begin{pmatrix} a & a \\ a & a \end{pmatrix} \middle| a \in Z_3 \right\}$$

and

$$V_2 = \{(a\ a\ a\ a\ a) \mid a \in Z_3\}.$$

V is a group bivector space over the group $Z_3 = \{0, 1, 2\}$ addition modulo 3.

***Example 4.60:*** Let

$$\begin{aligned} V &= V_1 \cup V_2 \\ &= \{(0\ 0), (1\ 1), (1\ 1\ 1\ 1\ 1), (0\ 1), (0\ 0\ 0\ 0\ 0)\} \cup \\ &\quad \left\{ \begin{pmatrix} a & a \\ a & a \end{pmatrix}, \begin{pmatrix} a & a & a \\ a & a & a \end{pmatrix} \middle| a \in Z_2 = \{0,1\} \right\}. \end{aligned}$$

V is a group bivector space over the group $G = Z_2 = \{0, 1\}$ addition modulo 2.

***Example 4.61:*** Let $V = V_1 \cup V_2 = \{Z[x]\} \cup \{(Z \times Z \times Z)\}$ be a group bivector space over the group G = Z, group under addition.

We now define some interesting substructures of group bivector spaces.

**DEFINITION 4.20:** *Let $V = V_1 \cup V_2$ be a group bivector space over the group G. $W = W_1 \cup W_2 \subseteq V_1 \cup V_2$ is said to be a group bivector subspace of V over G if W itself is a group bivector space over G.*

We illustrate this by some examples.

***Example 4.62:*** Let

$$\begin{aligned} V &= V_1 \cup V_2 \\ &= \{(0\ 0\ 0), (1\ 1\ 1), (0\ 0\ 0\ 0), (1\ 1\ 1\ 1), (1\ 1\ 0\ 0), (0\ 0\ 1\ 1)\} \\ &\quad \cup \left\{ \begin{pmatrix} a & a \\ a & a \end{pmatrix}, \begin{pmatrix} a & a & a \\ a & a & a \end{pmatrix} \middle| a \in Z_2 \right\} \end{aligned}$$

be the group bivector space over the group G = $Z_2$ = {0, 1} under addition. Take

$$\begin{aligned} W &= W_1 \cup W_2 \\ &= \{(0\ 0\ 0), (1\ 1\ 1), (0\ 0\ 0\ 0), (1\ 1\ 1\ 1)\} \cup \\ &\quad \left\{ \begin{pmatrix} a & a & a \\ a & a & a \end{pmatrix} \middle| a \in Z_2 \right\} \subseteq V_1 \cup V_2 ; \end{aligned}$$

Clearly W is a group bivector subspace of V over G = $Z_2$.

***Example 4.63:*** Let

$$\begin{aligned} V &= V_1 \cup V_2 \\ &= \left\{ \begin{pmatrix} a & b \\ c & d \end{pmatrix}, \begin{pmatrix} a & a & a \\ a & a & a \end{pmatrix} \middle| a, b, c, d \in Z_2 \right\} \cup \{(1\ 1\ 1\ 0\ 1\ 1), \\ &\quad (0\ 0\ 0\ 0\ 0\ 0), (1\ 1\ 1\ 0\ 0\ 0), (1\ 1), (0\ 0), (1\ 0)\} \end{aligned}$$

be a group bivector space over $Z_2$ = {0, 1}. Take

$$\begin{aligned} W &= \left\{ \begin{pmatrix} a & b \\ c & d \end{pmatrix} \middle| a, b, c, d \in Z_2 \right\} \cup \{(0\ 0\ 0\ 0\ 0\ 0), (1\ 1\ 1\ 0\ 0\ 0)\} \\ &\subset V_1 \cup V_2, \end{aligned}$$

W is a group bivector subspace of V over $Z_2$ = {0, 1}.

***Example 4.64:*** Let

$$\begin{aligned} V &= V_1 \cup V_2 \\ &= \{Z[x]\} \cup \left\{ \begin{pmatrix} a & b \\ c & d \end{pmatrix}, \begin{pmatrix} a & b & c \\ d & e & f \end{pmatrix} \middle| a, b, c, d, e, f \in Z \right\} \end{aligned}$$

be a group bivector space of V over Z. Take
W = $W_1 \cup W_2$

$$= \quad \{2Z[x]\} \cup \left\{\begin{pmatrix} a & b \\ c & d \end{pmatrix} \middle| a, b, c, d \in Z\right\} \subset V_1 \cup V_2;$$

W is a group bivector subspace of V over Z.

Now we define the notion of pseudo semigroup bivector subspace.

**DEFINITION 4.21:** *Let $V = V_1 \cup V_2$ be a group bivector space over the group G. Let $W = W_1 \cup W_2 \subseteq V_1 \cup V_2$ and $H \subseteq G$ be a semigroup of the group G. If W is a semigroup bivector space over the semigroup H then we call W to be a pseudo semigroup bivector subspace of V.*

***Example 4.65:*** Let $V = V_1 \cup V_2 = \{Z[x]\} \cup \{ (a, b, c) \mid a, b, c, \in Z\}$ be a group bivector space over the group Z. Take $W = W_1 \cup W_2 = \{Z^+ [x]\} \cup \{a, b, c) \mid a, b, c \in Z^+\} \subset V_1 \cup V_2$; W is a semigroup bivector space. W is a pseudo semigroup bivector subspace of V over the semigroup $Z^+ \subseteq Z$.

***Example 4.66:*** Let

$$V = V_1 \cup V_2$$
$$= \left\{\begin{pmatrix} a & b & c \\ c & d & b \end{pmatrix} \middle| a, b, c, d, e, f \in Z\right\} \cup$$
$$\{(a, b\ c\ d) \mid a, b, c, d \in Z\}$$

be a group bivector space over the group Z. Take

$$W = W_1 \cup W_2$$
$$= \left\{\begin{pmatrix} a & b & c \\ d & e & f \end{pmatrix} \middle| a, b, c, d, e, f \in Z^+ \cup \{0\}\right\} \cup$$
$$\{(a, b, c, d) \mid a, b, c, d \in 2Z^+ \cup \{0\}\} \subseteq V_1 \cup V_2;$$

W is a pseudo semigroup bivector subspace of V over the semigroup $2Z^+ \cup \{0\}$.

***Example 4.67:*** Let

$$V = V_1 \cup V_2$$
$$= \{(a\ b\ c\ d) \mid a, b, c, d \in Z\} \cup$$
$$\left\{ \begin{pmatrix} a & b \\ c & d \end{pmatrix}, \begin{pmatrix} a & a & a \\ a & a & a \end{pmatrix} \middle| a,b,c,d \in Z \right\}$$

be a group bivector space over the group Z. Take

$$W = W_1 \cup W_2$$
$$= \{(a\ 0\ b\ 0) \mid a, b \in Z^+ \cup \{0\}\} \cup$$
$$\left\{ \begin{pmatrix} a & b \\ c & d \end{pmatrix} \middle| a,b,c,d \in Z^+ \cup \{0\} \right\} \subseteq V_1 \cup V_2 .$$

W is a pseudo semigroup bivector subspace of V over the semigroup $Z^+ \cup \{0\}$.

Now we proceed onto define the substructure pseudo set bivector subspace of a group bivector space.

**DEFINITION 4.22:** *Let $V = V_1 \cup V_2$ be a group bivector space over the group G. Let $W = W_1 \cup W_2 \subseteq V_1 \cup V_2 = V$, take S a proper subset of G. If W is a set vector space over the set S then we call W to be the pseudo set bivector subspace of V over the set S.*

We now illustrate this by the following example.

***Example 4.68:*** Let $V = V_1 \cup V_2$ be a group bivector space over the group G = Z; where

$$V_1 = \{(a, b, c, d) \mid a, b, c, d \in 2Z\}$$

and

$$V_2 = \left\{ \begin{pmatrix} a & b \\ c & d \end{pmatrix} \middle| a,b,c,d \in Z \right\}.$$

Take

$$
\begin{aligned}
W &= W_1 \cup W_2 \\
&= \{(a\ b\ c\ d) \mid a, b, c, d \in 2Z^+ \cup \{0\}\} \cup \\
&\quad \left\{\begin{pmatrix} a & b \\ c & d \end{pmatrix} \middle| a, b, c, d \in 2Z^+ \cup \{0\}\right\} \\
&\subseteq V_1 \cup V_2 = V,
\end{aligned}
$$

W is a pseudo set bivector subspace over the set $S = \{0, 2, 2^2, 2^4, 2^8, \ldots, 2^{2n} \mid n \in N\}$.

Next we proceed onto define the notion of bisemigroup bivector group space.

**DEFINITION 4.23:** *Let $V = V_1 \cup V_2$, where $V_1$ is a semigroup vector space over the semigroup $S_1$ and $V_2$ is a semigroup vector space over the semigroup $S_2$, ($S_1 \neq S_2$, $S_1 \not\subseteq S_2$ and $S_2 \not\subseteq S_1$). Also $V_1 \neq V_2$, $V_1 \not\subseteq V_2$ and $V_2 \not\subseteq V_1$. We call V to be the bisemigroup bivector space over the bisemigroup $S = S_1 \cup S_2$.*

We illustrate this by the following examples.

***Example 4.69:*** Let $V = V_1 \cup V_2$ where

$$V_1 = \{(a, a, a) \mid (a \in Z_6\}$$

and

$$V_2 = \left\{\begin{pmatrix} a & a \\ a & a \end{pmatrix} \middle| a \in Z^+ \right\};$$

V is a bisemigroup bivector space over the bisemigroup $S = Z_6 \cup Z^+$.

***Example 4.70:*** Let $V = V_1 \cup V_2 = \{Z^+ [x]\} \cup \{Q^+ \times Q^+\}$, V is a bisemigroup bivector space over the bisemigroup $3Z^+ \cup 5Z^+$.

***Example 4.71:*** Let $V = V_1 \cup V_2$ where

$$V_1 = \left\{\begin{pmatrix} a & a \\ a & a \end{pmatrix} \middle| a \in 2Z^+ \right\}$$

is a semigroup vector space over the semigroup $2Z^+$ and $V_2 = \{a, a, a\} \mid a \in 3Z^+\}$ is the semigroup vector space over the semigroup $3Z^+$. V is a bisemigroup bivector space over the bisemigroup $S = 2Z^+ \cup 3Z^+$.

***Example 4.72:*** Let $V = V_1 \cup V_2$ where

$$V_1 = \{Z^+ \times Z^+ \times Z^+\}$$

and

$$V_2 = \left\{ \begin{pmatrix} a & b \\ c & d \end{pmatrix} \middle| \, a, b, c, d \in 3Z^+ \right\}.$$

V is bisemigroup bivector space over the bisemigroup $S = 2Z^+ \cup 3Z^+$.

The notion of bisemigroup bivector subspace can be defined as in case of other bivector spaces.

Next we proceed onto define the notion of bisemigroup bilinear algebra defined over the bisemigroup $S = S_1 \cup S_2$.

**DEFINITION 4.24:** *Let $V = V_1 \cup V_2$ where $V_1$ is a semigroup linear algebra over the semigroup $S_1$ and $V_2$ is a semigroup linear algebra over the semigroup $S_2$ ($V_1 \neq V_2$, $V_1 \not\subseteq V_2$; $V_2 \not\subseteq V_1$) ($S_1 \neq S_2$; $S_1 \not\subseteq S_2$ and $S_2 \not\subseteq S_1$). Then we call V to be the bisemigroup bilinear algebra over the bisemigroup $S = S_1 \cup S_2$.*

We illustrate this by the following examples.

***Example 4.73:*** Let $V = V_1 \cup V_2 =$

$$V_1 = \left\{ \begin{pmatrix} a & a \\ a & a \end{pmatrix} \middle| \, a \in 2Z^+ \right\}$$

and

$$V_2 = (a\ b\ c\ d) \mid a, b, c, d \in 3Z^+\}$$

be the bisemigroup bilinear algebra over the bisemigroup $S = 2Z^+ \cup 3Z^+$.

***Example 4.74:*** Let $V = V_1 \cup V_2 = \{(1\ 1), (0\ 0), (0\ 0\ 0), (1\ 1\ 1)\} \cup \{(a\ b\ c\ d) \mid a, b, c, d \in Z^+$ be the bisemigroup bivector space over the bisemigroup $S = Z_2 \cup Z^+$. Clearly V is not a bisemigroup bilinear algebra over S.

In view of this we have the following.

Every bisemigroup bilinear algebra is a bisemigroup bivector space but in general a bisemigroup bivector space is not a bisemigroup bilinear algebra. The above example is a semigroup bivector space which is not a bisemigroup bilinear algebra.

We now proceed onto define the notion of bisemigroup bilinear subalgebra.

**DEFINITION 4.25:** *Let $V = V_1 \cup V_2$ be a semigroup bilinear algebra over the bisemigroup $S = S_1 \cup S_2$. Let $W = W_1 \cup W_2 \subset V_1 \cup V_2$ if ($W_1 \neq W_2$ $W_1 \not\subseteq W_2$ and $W_2 \not\subseteq W_1$)and W is itself a bisemigroup bilinear algebra over S then we call W to be a bisemigroup bilinear subalgebra over the bisemigroup $S = S_1 \cup S_2$.*

We illustrate this by the following example.

***Example 4.75:*** Let $V = V_1 \cup V_2 = \{Z^+ [x]\} \cup Z^+ \times Z^+ \times Z^+\}$ be the bisemigroup bilinear algebra over the bisemigroup $S = 3Z^+ \cup 5Z^+$. Take $W = W_1 \cup W_2 = \{3Z^+ [x]\} \cup \{Z^+ \times \{0\} \times Z^+\} \subseteq V_1 \cup V_2$, W is a bisemigroup bilinear subalgebra of V over the bisemigroup S.

***Example 4.76:*** Let

$$\begin{aligned} V &= V_1 \cup V_2 \\ &= \{(a, b) \text{ such that } a, b \in Z_{10}\} \cup \\ &\quad \left\{ \begin{pmatrix} a & b \\ c & d \end{pmatrix} \middle| \, a, b, c, d \in Z^+ \right\} \end{aligned}$$

be a bisemigroup bilinear algebra over the bisemigroup $S = Z_{10} \cup Z^+$. Take

$$\begin{aligned} W &= W_1 \cup W_2 \\ &= \{(a, a) \mid a \in Z_{10}\} \cup = \left\{ \begin{pmatrix} a & a \\ a & a \end{pmatrix} \middle| a \in Z^+ \right\}, \\ &\subseteq V_1 \cup V_2 \end{aligned}$$

W is a bisemigroup bilinear subalgebra of V over the bisemigroup $S = Z_{10} \cup Z^+$.

***Example 4.77:*** Let

$$\begin{aligned} V &= V_1 \cup V_2 \\ &= \{(a, a) \mid a \in 5Z^+\} \cup \left\{ \begin{pmatrix} a & a \\ a & a \\ a & a \end{pmatrix} \middle| a \in 7Z^+ \right\} \end{aligned}$$

be a bisemigroup bilinear algebra over the bisemigroup $S = 5Z^+ \cup 7Z^+$. Take

$$\begin{aligned} W &= W_1 \cup W_2 \\ &= \{(a, a) \mid a \in 15\, Z^+\} \cup \left\{ \begin{pmatrix} a & a \\ a & a \\ a & a \end{pmatrix} \middle| a \in 14Z^+ \right\} \\ &\subseteq V_1 \cup V_2; \end{aligned}$$

W is a bisemigroup bilinear subalgebra of V over the bisemigroup $S = 5Z^+ \cup 7Z^+$.

***Example 4.78:*** Let $V = V_1 \cup V_2 = \{7Z^+ [x]\} \cup \{5Z^+ \times 5Z^+ \times 5Z^+\}$ be a bisemigroup bilinear algebra over the bisemigroup S $= 7Z^+ \cup 5Z^+$. Take $W = W_1 \cup W_2 = \{14Z^+ [x]\} \cup \{5Z^+ \times \phi \times$ $5Z^+\} \subseteq V_1 \cup V_2$, W is a bisemigroup bilinear subalgebra over the bisemigroup $S = 7Z^+ \cup 5Z^+$.

Now having seen the definition of semigroup bilinear subalgebra we now proceed on to define the notion of bidimension of the bisemigroup bilinear algebra over the bisemigroup.

**DEFINITION 4.26:** *Let $V = V_1 \cup V_2$ be the bisemigroup bilinear algebra over the bisemigroup $S = S_1 \cup S_2$. Take $X = X_1 \cup X_2 \subseteq V_1 \cup V_2$; if $X_1$ generates $V_1$ and $X_2$ generates $V_2$ then we say X bigenerates V over the bisemigroup $S = S_1 \cup S_2$.*

*The cardinality of $X_1 \cup X_2$ is given by $|X_1| \cup |X_2|$ or $(|X_1|, |X_2|)$, called the bidimension of the bisemigroup bilinear algebra $V = V_1 \cup V_2$. If even one of $X_1$ or $X_2$ is of infinite dimension then we say the bidimension of V is infinite, only when both $X_1$ and $X_2$ are of finite cardinality we say V is of finite bidimension over the bisemigroup S.*

***Example 4.79:*** Let $V = V_1 \cup V_2$ = {(1 0), (0 1), (0 0), (1 1 1), (0 0 0), (1 0 1)} ∪ {(1 1 1), (3 3 3), (2 2 2), (0 0 0), (1 0 0), (2 0 0), (3 0 0)} be the bisemigroup bivector space over the bisemigroup $S = Z_2 \cup Z_4$. Let X = {(1 0), (0 1), (1 1 1), (1 0 1)} ∪ {(1 1 1), (1 0 0)} $= X_1 \cup X_2$ be the bigenerator of V. The bidimension of V is (4, 2).

***Example 4.80:*** Let $V = V_1 \cup V_2 = \{Z_2 [x\} \cup \{(a\ a\ a) \mid a \in Z_5\}$ be the bisemigroup bilinear bialgebra over the bisemigroup $S = Z_2 \cup Z_5$. Let $X = \{1, x, \ldots, x^\infty \} \cup \{(1\ 1\ 1)\} = X_1 \cup X_2 \subseteq V_1 \cup V_2$ be the bigenerator of V. Clearly bidimension of V is $(\infty, 1)$ or $\infty \cup \{1\}$ so V is of infinite bidimension over S.

***Example 4.81:*** Suppose $V = V_1 \cup V_2$ = {(1 1 1), (1 1 1 1), (0 0 1), (0 0 0), (1 0 0), (1 1 0 0), (0 0 0 0)} ∪ $\{Z_3 \times Z_3 \times Z_3\}$ be the bisemigroup bivector space over the bisemigroup $S = Z_2 \cup Z_3$. Let X = {(1 1 1), (1 1 1 1), (0 0 1), (1 0 0), (1 1 0 0)} ∪ {(1 0 0), (0 1 0) (0 0 1), (1 1 0), (0 1 1), (1 0 1), (1 2 0), (1 1 1), (1 0 2), (0 1 2), (1 2 2), (2 1 2), (2 2 1) etc} $= X_1 \cup X_2$ which is

the bigenerator of the bisemigroup bivector space over the bisemigroup. The bidimension of V over S is $(5 \cup 26)$.

***Example 4.82:*** Let

$$V = V_1 \cup V_2 = \{2Z^+ [x]\} \cup \left\{ \begin{pmatrix} a & b \\ c & d \end{pmatrix} \middle| \; a, b, c, d \in 3Z^+ \right\},$$

be a bisemigroup bivector space over the bisemigroup $S = 2Z^+ \cup 3Z^+$. Clearly the bidimension of V is infinite.

Now we proceed onto define the notion of bigroup bivector space over the bigroup.

**DEFINITION 4.27:** *Let $V = V_1 \cup V_2$, such that $V_1 \neq V_2$, $V_1 \not\subseteq V_2$, $V_2 \not\subseteq V_1$. If $V_1$ is a group vector space over the group $G_1$ and $V_2$ is a group vector space over the group $G_2$ ($G_1 \neq G_2$ , $G_1 \not\subseteq G_2$ and $G_2 \not\subseteq G_1$) then we say $V = V_1 \cup V_2$ is a bigroup bivector space over the bigroup $G = G_1 \cup G_2$. Clearly if $V_1$ and $V_2$ are just set it is sufficient to define bigroup bivector space over a bigroup.*

***Example 4.83:*** Let $V = V_1 \cup V_2$ be a bigroup bivector space over the bigroup $G = Z_3 \cup Q$ where $V_1 = Z_3[x]$ and $V_2 = Q \times Q$.

***Example 4.84:*** Let $V = V_1 \cup V_2 = \{(1\ 1), (0\ 0\ ), (0\ 1), (1\ 1\ 1\ 0), (0\ 0\ 0\ 0), (0\ 1\ 0\ 0), (1\ 1\ 1\ 0\ 0), (0\ 0\ 0\ 0\ 0), (0\ 0\ 0\ 1\ 0)\} \cup$

$$\left\{ \begin{pmatrix} a & a \\ a & a \end{pmatrix} \text{and} \begin{pmatrix} a & a & a \\ a & a & a \end{pmatrix} \middle| \; a \in Z \right\}$$

be the bigroup bivector space over the bigroup $G = G_1 \cup G_2 = Z_2 \cup Z$. Clearly both $V_1$ and $V_2$ are just sets.

***Example 4.85:*** Let $V = V_1 \cup V_2 = \{Z_5 \times Z_5 \times Z_5\} \cup \{Z_7[x]\}$ be the bigroup bivector space over the bigroup $G = Z_5 \cup Z_7$.

***Example 4.86:*** Let $V = V_1 \cup V_2 = \{\langle(1\ 1\ 1), (2\ 0\ 0), (2\ 2\ 0), (0\ 0\ 0)\rangle\} \cup \{\langle(1\ 5\ 7\ 2), (0\ 1\ 0\ 0), (0\ 0\ 0\ 0), (3\ 5\ 4\ 1), (1\ 2\ 3\ 4), (5\ 6\ 7\ 0)\rangle\}$ be the bigroup bivector space over the bigroup $G = Z_3 \cup Z_8$. ($\langle, \rangle$; denotes generated over the related groups).

***Example 4.87:*** Let

$$\begin{aligned} V &= V_1 \cup V_2 \\ &= \left\{\begin{pmatrix} a & b & c \\ d & e & f \end{pmatrix} \middle| \ a, b, c, d, e, f \in Z_6\right\} \\ &\quad \cup \left\{\begin{pmatrix} x & x \\ y & y \end{pmatrix} \middle| \ x, y \in Z_{10}\right\}, \end{aligned}$$

be a bigroup bivector space over the bigroup bivector space over the bigroup space over the bigroup $G = G_1 \cup G_2 = Z_6 \cup Z_{10}$.

Now we proceed onto define the notion of substructures in bigroup bivector spaces.

**DEFINITION 4.28:** *Let $V = V_1 \cup V_2$ be a bivector bispace over the bigroup $G = G_1 \cup G_2$. Let $W = W_1 \cup W_2 \subset V_1 \cup V_2 = V$ be such that $W_1 \neq W_2$. $W_1 \not\subseteq W_2$ or $W_2 \not\subseteq W_1$ if W itself is a bigroup bivector space over the bigroup $G = G_1 \cup G_2$ then we say W is a bigroup bivector subspace of V over the bigroup $G = G_1 \cup G_2$.*

We now illustrate this situation by the following examples.

***Example 4.88:*** Let

$$\begin{aligned} V &= V_1 \cup V_2 \\ &= \{(a, a, a) \mid a \in Z_6\} \cup \left\{\begin{pmatrix} a & a & a \\ a & a & a \end{pmatrix}, \begin{pmatrix} a & a \\ a & a \end{pmatrix} \middle| \ a \in Z_8\right\} \end{aligned}$$

be the bigroup bivector space over the bigroup $G = Z_6 \cup Z_8$. Let

$$W = \{(a, a, a) \mid a \in \{0, 3\}\} \left\{ \begin{pmatrix} a & a \\ a & a \end{pmatrix} \middle| a \in Z_8 \right\}$$
$$= W_1 \cup W_2.$$

W is the bigroup bivector subspace of V over the bigroup $G = Z_6 \cup Z_8$.

***Example 4.89:*** Let

$$V = V_1 \cup V_2$$
$$= \{(a\ a\ a), (a\ a), (a\ a\ a\ a) \mid a \in Z\} \cup \left\{ \begin{pmatrix} a & b \\ c & d \end{pmatrix} \middle| a, b, c, d \in Z_2 \right\}$$

be the bigroup bivector space over the bigroup $G = Z \cup Z_2$.

$$W = \{(a, a, a) \mid a \in Z\} \cup \left\{ \begin{pmatrix} a & a \\ a & a \end{pmatrix} \middle| a \in Z_2 \right\}$$
$$= W_1 \cup W_2 \subseteq V_1 \cup V_2$$

is the bigroup bivector subspace of V over the bigroup $G = Z \cup Z_2$.

***Example 4.90:*** Let $V = V_1 \cup V_2 = \{(a\ a\ a), (a\ a\ a\ a\ a), (a\ a), (a\ a\ a\ a\ a\ a) \mid a \in Z_5\} \cup Z_2[x]$ be the bigroup bivector space over the bigroup $Z_5 \cup Z_2$.

Take $W = W_1 \cup W_2 = \{(a\ a\ a), (a, a, a, a, a, a, a) \mid a \in Z_5\} \cup$ {all polynomial of even degree with coefficients from $Z_2$} $\subseteq V_1 \cup V_2$ is a bigroup bivector subspace of V over the bigroup $G = Z_5 \cup Z_2$.

Next we define pseudo bisemigroup bivector subspace of a bigroup bivector space V.

**DEFINITION 4.29:** *Let $V = V_1 \cup V_2$ be a bigroup bivector space over the bigroup $G = G_1 \cup G_2$. Let $W = W_1 \cup W_2 \subseteq V_1 \cup V_2$ (such that $W_1 \neq W_2$, $W_1 \not\subseteq W_2$, $W_2 \not\subseteq W_1$) where $W_1$ is a semigroup vector space over the semigroup $H_1$ contained in $G_1$ and $W_2$ is a semigroup vector space over the semigroup $H_2$ contained in $G_2$ $H_1 \not\subseteq H_2$, $H_2 \not\subseteq H_1$, $H_1 \neq H_2$.*

*Then we call $W = W_1 \cup W_2$ to be the pseudo bisemigroup bivector subspace of the bigroup bivector space V over $H = H_1 \cup H_2 \subseteq G_1 \cup G_2$, H the bisemigroup contained in the bigroup.*

We illustrate thus by the following examples.

***Example 4.91:*** Let

$$V = V_1 \cup V_2 = \{(a\ a\ a\ a), (a\ a\ a) \mid a \in Z\} \cup \left\{\begin{pmatrix} a & b \\ c & d \end{pmatrix} \middle| a, b, c, d \in Z_{12}\right\}$$

be the bigroup bivector space over the bigroup $G = Z \cup Z_{12}$.

Take

$$W = W_1 \cup W_2 = \{(a\ a\ a\ a), \mid a \in Z\} \cup \left\{\begin{pmatrix} a & b \\ 0 & d \end{pmatrix} \middle| a, b, d \in \{0, 2, 4, 6, 8, 10\}\right\} \subseteq V_1 \cup V_2.$$

W is the pseudo bisemigroup bivector subspace over the bisemigroup

$$Z^+ \cup \{0, 6\} = H = H_1 \cup H_2.$$

***Example 4.92:*** Let

$$V = V_1 \cup V_2$$

$$= \left\{ \begin{pmatrix} a & a \\ a & a \end{pmatrix} \begin{pmatrix} a & a & a \\ a & a & a \\ a & a & a \end{pmatrix} \middle| \, a \in Z \right\} \cup$$

$$\{(a\ a\ a\ a), (a\ a), (a\ a\ a) \mid a \in 2Z\}$$

be the bigroup bivector space over the bigroup $3Z \cup 2Z$.
Take

$$W = W_1 \cup W_2$$

$$= \left\{ \begin{pmatrix} a & a \\ a & a \end{pmatrix} \middle| \, a \in 3Z \right\} \cup \{(a\ a\ a\ a) \mid a \in 2Z\},$$

W is the pseudo bisemigroup bivector subspace of V over the bisemigroup $H = 3Z^+ \cup 2Z^+ \subseteq 3Z \cup 2Z$.

***Example 4.93:*** Let

$$V = V_1 \cup V_2$$
$$= \{Z[x]\} \cup \{(a\ a\ a\ a), (a\ a) \mid a \in 2Z\}$$

be the bigroup bivector space over the bigroup $G = 3Z \cup 2Z$.
Let

$$W = \{Z^+[x]\} \cup \{(a\ a) \mid a \in 2Z\}$$

$$\subseteq V_1 \cup V_2$$

be the pseudo bisemigroup bivector subspace over the bisemigroup $H = 3Z^+ \cup 2Z^+$.

Now we define yet a new mixed structure which we call as quasi bigroup bivector space.

**DEFINITION 4.30:** *Let $V = V_1 \cup V_2$ where $V_1$ is a semigroup vector space over the semigroup $S_1$ and $V_2$ is the group vector space over the group $G_1$. Then we call $V = V_1 \cup V_2$ to be the pseudo bigroup bivector space over the pseudo bigroup $G = S_1 \cup G_1$.*

We illustrate this situation by the following example.

***Example 4.94:*** Let

$$\begin{aligned} V &= V_1 \cup V_2 \\ &= \{Z_5\,[x]\} \cup \{Z^+ \times Z^+ \times Z^+\} \end{aligned}$$

be the pseudo bigroup bivector space over the pseudo bigroup $G = Z_5 \cup Z^+$.

***Example 4.95:*** Let

$$\begin{aligned} V &= V_1 \cup V_2 \\ &= Z^+[x] \cup \left\{ \begin{pmatrix} a & a \\ a & a \end{pmatrix} \middle|\, a \in Z_{12} \right\} \end{aligned}$$

be the pseudo bigroup bivector space over the pseudo bigroup $G = Z^+ \cup Z_{12}$ where $Z^+$ is the semigroup and $Z_{12}$ is the group under addition modulo 12.

***Example 4.96:*** Let

$$\begin{aligned} V &= V_1 \cup V_2 \\ &= \left\{ \begin{pmatrix} a & a \\ b & b \end{pmatrix} \middle|\, a, b \in 3Z^+ \right\} \cup \{Z_{20}\,[x]\} \end{aligned}$$

be the pseudo bigroup bivector space over the pseudo bigroup G = $3Z^+ \cup Z_{20}$.

The author leaves it as an exercise for the reader to define various substructures of this structure and bigenerator and bidimension.

Now we proceed onto generalize this to n-set n-vector spaces set n-vector spaces $n \geq 3$, in the following chapter.

## Chapter Five

# SET n-VECTOR SPACES AND THEIR GENERALIZATIONS

In this chapter we for the first time introduce the notion of set n-vector spaces, semigroup n-vector spaces and group n-vector spaces ($n \geq 3$) when $n = 2$ we get set bivector spaces, semigroup bivector spaces and so on.

**DEFINITION 5.1:** *Let $V = V_1 \cup \ldots \cup V_n$, each $V_i$ is a distinct set with $V_i \not\subseteq V_j$ or $V_j \not\subseteq V_i$ if $i \neq j$; $1 \leq i, j \leq n$. Let each $V_i$ be a set vector space over the set S, i = 1, 2, …, n, then we call $V = V_1 \cup V_2 \cup \ldots \cup V_n$ to be the set n-vector space over the set S.*

We illustrate this by the following examples.

***Example 5.1:*** Let

$$\begin{aligned} V &= V_1 \cup V_2 \cup V_3 \cup V_4 \\ &= \{(1\ 1\ 1), (0\ 0\ 0), (1\ 0\ 0), (0\ 1\ 0), (1\ 1), (0\ 0), (1\ \ 1\ \ 1\ \ 1), \\ &\quad (1\ 0\ 0\ 0), (0\ 0\ 0)\} \cup \\ &\quad \left\{ \begin{pmatrix} a & a \\ a & a \end{pmatrix} \middle| a \in Z_2 \right\} \cup \left\{ \begin{pmatrix} a & a & a & a \\ a & a & a & a \end{pmatrix} \middle| a \in Z_2 = \{0, 1\} \right\} \\ &\quad \cup \{Z_2\,[x]\}. \end{aligned}$$

V is a set 4 vector space over the set S = {0, 1}.

***Example 5.2:*** Let

$$V = V_1 \cup V_2 \cup V_3 \cup V_4 \cup V_5 \cup V_6$$
$$= \left\{\begin{pmatrix} a & a \\ a & a \end{pmatrix} \middle| a \in Z^+\right\} \cup \{Z^+ \times Z^+ \times Z^+\} \cup \{(a, a, a),$$
$$(a, a, a, a, a) \mid a \in Z^+\} \cup \left\{\begin{pmatrix} a \\ a \\ a \\ a \end{pmatrix} \middle| a \in Z^+\right\} \cup \{Z^+[x]\}$$
$$\cup \left\{\begin{pmatrix} a_1 & a_2 & a_3 \\ a_4 & a_5 & a_6 \end{pmatrix} \middle| a_i \in Z^+; 1 \le i \le 6\right\}$$

be the set 6-vector space over the set $S = Z^+$.

***Example 5.3:*** Let

$$V = V_1 \cup V_2 \cup V_3$$
$$= \{Z_6 [x]\} \cup \{Z_6 \times Z_6 \times Z_6\} \cup \left\{\begin{pmatrix} a & a & a \\ a & a & a \end{pmatrix} \middle| a \in Z_6\right\}$$

be the set 3 vector space over the set S = {0, 2, 4}. We call this also as set trivector space over the set S. Thus when n = 3 we call the set n vector space as set trivector space.

We define set n-vector subspace of a set n-vector space V.

**DEFINITION 5.2:** *Let $V = V_1 \cup \ldots \cup V_n$ be a set n-vector space over the set S. If $W = W_1 \cup \ldots \cup W_n$ with $W_i \ne W_j$; $i \ne j$, $W_i \not\subseteq W_j$ and $W_j \not\subseteq W_i$, $1 \le i, j \le n$ and $W = W_1 \cup W_2 \cup \ldots \cup W_n \subseteq V_1 \cup V_2 \cup \ldots \cup V_n$ and W itself is a set n-vector space over the set S then we call W to be the set n vector subspace of V over the set S.*

We illustrate this by a simple example.

***Example 5.4:*** Let

$$
\begin{aligned}
V &= V_1 \cup V_2 \cup V_3 \cup V_4 \\
&= \{(a, a, a), (a, a) \mid a \in Z^+\} \cup \left\{ \begin{pmatrix} a & a \\ b & b \end{pmatrix} \middle| a, b \in Z^+ \right\} \cup \\
&\quad \{Z^+[x]\} \cup \left\{ \begin{pmatrix} a \\ a \\ a \end{pmatrix} \middle| a \in Z^+ \right\},
\end{aligned}
$$

V is a set 4-vector space over the set S = $Z^+$. Take

$$
\begin{aligned}
W &= W_1 \cup W_2 \cup W_3 \cup W_4 \\
&= \{(a, a, a) \mid a \in Z^+\} \cup \left\{ \begin{pmatrix} a & a \\ 0 & 0 \end{pmatrix} \middle| a \in Z^+ \right\} \cup \\
&\quad \{\text{all polynomial of even degree}\} \cup \left\{ \begin{pmatrix} a \\ a \\ a \end{pmatrix} \middle| a \in 2Z^+ \right\} \\
&\subseteq V_1 \cup V_2 \cup V_3 \cup V_4 = V,
\end{aligned}
$$

is a set 4-vector subspace of V over the set S = $Z^+$.
We can find several set 4-vector subspaces of V.

Now we proceed on to define the n-generating set of a set n-vector space over the set S.

**DEFINITION 5.3:** *Let $V = V_1 \cup \ldots \cup V_n$ be a set n-vector space over set S. Let $X = X_1 \cup \ldots \cup X_n \subset V_1 \cup V_2 \cup \ldots \cup V_n = V$. If each set $X_i$ generates $V_i$ over the set S, i = 1, 2, …, n then we say the set n vector space $V = V_1 \cup \ldots \cup V_n$ is generated by the n-set $X = X_1 \cup X_2 \cup \ldots \cup X_n$ and X is called the n-generator of V. If each of $X_i$ is of cardinality $n_i$, i = 1, 2, …, n then we say the n-cardinality of the set n vector space V is given by $|X_1| \cup \ldots \cup$*

*$|X_n| = \{|X_1|, |X_2|, \ldots, |X_n|\} = \{(n_1, n_2, \ldots, n_n)\}$. If even one of the $X_i$ is of infinite cardinality we say the n-cardinality of V is infinite. Thus if all the sets $X_1, \ldots, X_n$ have finite cardinality then we say the n-cardinality of V is finite.*

We now illustrate this by the following examples.

***Example 5.5:*** Let

$$\begin{aligned} V &= V_1 \cup V_2 \cup \ldots \cup V_5 \\ &= \{(a, a, a) \mid a \in \{0, 1\}\} \cup \left\{\begin{pmatrix} a & a \\ a & a \end{pmatrix} \middle| a \in \{0, 1\}\right\} \cup \\ &\quad \left\{\begin{pmatrix} a & a & a \\ a & a & a \end{pmatrix} \middle| a \in \{0, 1\}\right\} \cup \{(a, a, a, a)\} \; a \in \{0, 1\}\} \cup \\ &\quad \left\{\begin{pmatrix} a & b \\ a & b \end{pmatrix} \middle| a, b \in \{0, 1\}, a \neq b\right\} \end{aligned}$$

be a set 5-vector space over the set S = {0, 1}.
Choose

$$\begin{aligned} X &= \{(1\ 1\ 1)\} \cup \left\{\begin{pmatrix} 1 & 1 \\ 1 & 1 \end{pmatrix}\right\} \cup \left\{\begin{pmatrix} 1 & 1 & 1 \\ 1 & 1 & 1 \end{pmatrix}\right\} \cup \{(1\ 1\ 1\ 1)\} \cup \\ &\quad \left\{\begin{pmatrix} 1 & 0 \\ 1 & 0 \end{pmatrix}, \begin{pmatrix} 0 & 1 \\ 0 & 1 \end{pmatrix}\right\} \\ &\subseteq V_1 \cup V_2 \cup \ldots \cup V_5 = V. \end{aligned}$$

It is easily verified each $X_i$ generates $V_i$, i = 1, 2, …, 5. Thus X = $X_1 \cup X_2 \cup \ldots \cup X_5$ is the 5-generator set of the set 5-vector space over the set S. In fact each set $X_i$ is of finite cardinality, so V is a set 5-vector space of finite 5-dimension. In fact the 5 dimension of V is {1, 1, 1, 1, 2}.

One of the important and interesting factor to observe about these set n-vector spaces over the set S, is at times they can have one and only one generating set. The example 5.5 is one such case.

Now we proceed onto give yet another example of a set n-vector space V over a set S.

***Example 5.6:*** Let

$$
\begin{aligned}
V &= V_1 \cup V_2 \cup V_3 \cup V_4 \cup V_5 \\
&= \{(a, a, a) \mid a \in Z^+\} \cup \{Q^+ \times Q^+\} \cup \{Z^+ [x]\} \cup \\
&\quad \left\{\begin{pmatrix} a & b \\ c & d \end{pmatrix} \middle| a, b, c, d \in Z^+\right\} \cup \left\{\begin{pmatrix} a & a & a & a \\ a & a & a & a \end{pmatrix} \middle| a \in Q^+\right\}
\end{aligned}
$$

be a set 5-sector space over the set $S = Z^+$. Take

$$
\begin{aligned}
X &= X_1 \cup X_2 \cup X_3 \cup X_4 \cup X_5 \\
&= \{(1\ 1\ 1)\} \cup \{\text{an infinite set of pairs including } (1, 1)\} \cup \\
&\quad \{1, x, \ldots, x^n \text{ and an infinite set}\} \cup \left\{\begin{pmatrix} 1 & 1 \\ 1 & 1 \end{pmatrix} \text{ with an} \right. \\
&\quad \text{infinite set of } 2 \times 2 \text{ matrices}\} \cup \left\{\begin{pmatrix} 1 & 1 & 1 & 1 \\ 1 & 1 & 1 & 1 \end{pmatrix} \text{ together} \right. \\
&\quad \left. \text{with an infinite set } \begin{pmatrix} a & a & a & a \\ a & a & a & a \end{pmatrix} \middle| a \in Q^+ \right\}.
\end{aligned}
$$

X is a 5-generating set. In fact X is an infinite 5-generating subset of $V = V_1 \cup \ldots \cup V_5$.

***Example 5.7:*** Let

$$
\begin{aligned}
V &= V_1 \cup V_2 \cup V_3 \cup V_4 \cup V_5 \\
&= \{Z^+\} \cup \{(a\ a\ a) \mid a \in Z^+\} \cup \left\{\begin{pmatrix} a & a \\ a & a \end{pmatrix} \middle| a \in Z^+\right\} \cup
\end{aligned}
$$

$\{k_1 (1 + x^2 + x), k_2x, k_3x^3, k_4(x^7 + 1), k_5 (x^8+1), k_6(x^{13} + 1 + 2x); k_i \in Z^+; 1 \leq i \leq 6\} \cup \left\{ \begin{pmatrix} a & a & a & a & a \\ a & a & a & a & a \\ a & a & a & a & a \end{pmatrix} \middle| a \in Z^+ \right\}$

be a 5 set vector space over the set S = $Z^+$.
Take

$$X = \{1\} \cup \{(1\ 1\ 1)\} \cup \left\{ \begin{pmatrix} 1 & 1 \\ 1 & 1 \end{pmatrix} \right\} \cup \{1 + x^2 + x, x, x^3, x^7 + 1, x^8 + 1, x^{13} + 2x + 1\} \cup \left\{ \begin{pmatrix} 1 & 1 & 1 & 1 & 1 \\ 1 & 1 & 1 & 1 & 1 \\ 1 & 1 & 1 & 1 & 1 \end{pmatrix} \right\}$$

$$\subseteq V_1 \cup V_2 \cup V_3 \cup V_4 \cup V_5.$$

Clearly X, 5-generates V i.e., X is the unique 5-generator of V. Further V is 5-finitely generated and 5-dimension of V is $\{1 \cup 1 \cup 1 \cup 6 \cup 1\}$ = (1, 1, 1, 6, 1).

We now proceed onto define the n-set basis of the n-set vector space V over S.

**DEFINITION 5.4:** *Let $V = V_1 \cup V_2 \cup \ldots \cup V_n$ be a set n-vector space over the set S. Suppose $X = X_1 \cup \ldots \cup X_n$ is a n-set generating subset of $V = V_1 \cup \ldots \cup V_n$ then we call X to be the n-set basis of V. If $X = ( |X_1|, |X_2|, \ldots, |X_n| )$ in which each $|X_i| < \infty$, $1 \leq i \leq n$, then we say V is finitely n set generated by the n-set X and is of n-dimension $(|X_1|, |X_2|, \ldots, |X_n|)$. Even if one of $|X_i| = \infty$ we say V is infinitely n-set generated by X.*

Now we will give one or two examples of n-set basis before we proceed onto define other interesting notions about n-set vector spaces.

***Example 5.8:*** Let

$$V = \{(a, a, a) \mid a \in Z^+ \cup \{0\}\} \cup \left\{\begin{pmatrix} a & a \\ 0 & a \end{pmatrix} \middle| a \in Z^+ \cup \{0\}\right\} \cup$$
$$\{a_1(x^2 + 1 + 3x), a_2(x + 1), a_4x^7 a_3 (x^6 + x^2 + x^3 + x^4 + 1), a_5x, a_6x^2, a_7(x^8 + 5x^2 + 1) \mid a_i \in Z^+ \cup \{0\}; 1 \leq i \leq 7 \} \cup$$
$$\left\{\begin{pmatrix} a & a & a & a \\ a & a & a & a \end{pmatrix} \middle| a \in Z^+ \cup \{0\}\right\}$$
$$= V_1 \cup V_2 \cup V_3 \cup V_4$$

be a 4-set vector space over the set $S = Z^+ \cup \{0\}$. Now take

$$X = X_1 \cup X_2 \cup X_3 \cup X_4$$
$$= \{(1\ 1\ 1)\} \cup \left\{\begin{pmatrix} 1 & 1 \\ 0 & 1 \end{pmatrix}\right\} \cup \{(x^3 + 3x+1), x+1, x^7, x^6 + x^2 + x^3 + x^4 + 1, x, x^2, x^8 + 5x^2 + 1\} \cup \left\{\begin{pmatrix} 1 & 1 & 1 & 1 \\ 1 & 1 & 1 & 1 \end{pmatrix}\right\}$$
$$\subseteq V_1 \cup V_2 \cup V_3 \cup V_3.$$

Clearly V is 4 set generated by X and the basis 4-set of V is X and the n-dimension of the n-set X is (1, 1, 7, 1), (n = 4).

We give yet another example.

***Example 5.9:*** Let

$$V = V_1 \cup V_2 \cup \ldots \cup V_5$$
$$= \{(0\ 0\ 0\ 1\ 1\ 1), (0\ 0\ 1\ 1\ 1\ 1), (0\ 1\ 1\ 1\ 1\ 1), (0\ 0\ 0\ 0\ 0\ 0), (1\ 1\ 1), (0\ 1\ 0), (1\ 1\ 0), (0\ 0\ 0)\} \cup \{(1\ 1\ 1\ 1\ 1\ 1\ 1\ 1\ 1), (0\ 0\ 0\ 0\ 0\ 0\ 0\ 0\ 0)\} \cup \left\{\begin{pmatrix} a & a \\ a & a \end{pmatrix} \middle| a \in Z_2 = \{0,1\}\right\} \cup$$
$$\left\{\begin{pmatrix} a & a & a & a \\ a & a & a & a \end{pmatrix}, \begin{pmatrix} 1 & 1 & 1 & 1 & 1 \\ 1 & 1 & 0 & 1 & 1 \end{pmatrix}, \begin{pmatrix} a & a & a \\ a & a & a \end{pmatrix},\right.$$

$$\left.\begin{pmatrix}1&1&1&1\\0&0&0&0\end{pmatrix}, \begin{pmatrix}0&0&0&0&0\\0&0&0&0&0\end{pmatrix}, \begin{pmatrix}0&0&0\\0&0&0\end{pmatrix} \right| a \in Z_2 = \{0,1\} \Bigg\}$$

$\cup \{1 + x^2, 0, x + x^2 + 1, x^3 + 1\}$

be a 5-set vector space over the set $Z_2 = \{0, 1\}$.

$$\begin{aligned} X &= X_1 \cup \ldots \cup X_5 \\ &= \{(0\ 0\ 0\ 1\ 1\ 1), (0\ 0\ 1\ 1\ 1\ 1), (0\ 1\ 1\ 1\ 1\ 1), (1\ 1\ 1), \\ &\quad (0\ 1\ 0), (1\ 1\ 0)\} \cup \{(1\ 1\ 1\ 1\ 1\ 1\ 1\ 1\ 1)\} \cup \left\{\begin{pmatrix}1&1\\1&1\end{pmatrix}\right\} \cup \\ &\quad \left\{\begin{pmatrix}1&1&1&1\\1&1&1&1\end{pmatrix}, \begin{pmatrix}1&1&1\\1&1&1\end{pmatrix}, \begin{pmatrix}1&1&1&1\\0&0&0&0\end{pmatrix}, \begin{pmatrix}1&1&1&1&1\\1&1&0&1&1\end{pmatrix}\right\} \cup \\ &\quad \{1 + x^2, x^2 + x + 1, x^3 + 1\} \\ &\subseteq V_1 \cup V_2 \cup \ldots \cup V_5 \end{aligned}$$

be the 5 set which 5-generates V. Thus X is a 5-set- 5-basis of V of 5-dimension (6, 1, 1, 4, 3).

Now we proceed on to define n-set linear algebra over the S.

**DEFINITION 5.5:** *Let $V = V_1 \cup V_2 \cup \ldots \cup V_n$ be a n set vector space over the set S. If each $V_i$ is closed under addition then we call V to be a n-set linear algebra over S; $1 \leq i \leq n$.*

Now we illustrate this by the following examples.

***Example 5.10:*** Let

$$\begin{aligned} V &= V_1 \cup V_2 \cup V_3 \cup V_4 \\ &= \left\{\left.\begin{pmatrix}a&a&a&a\\a&a&a&a\end{pmatrix}\right| a \in Z^+ \cup \{0\}\right\} \cup \{Z^+ \cup \{0\} \times Z^+ \cup \\ &\quad \{0\} \times Z^+ \cup \{0\}\} \cup \{\text{all polynomials of degree less than or equal to 5 with coefficients from } Z^+ \cup \{0\}\} \cup \\ &\quad \left\{\left.\begin{pmatrix}a&b\\c&d\end{pmatrix}\right| a,b,c,d \in Z^+ \cup \{0\}\right\} \end{aligned}$$

be a 4-set vector space over the set $S = Z^+ \cup \{0\}$. Clearly V is a 4-set linear algebra over S. The basis of V as a 4 set linear algebra over S is only finite whereas V as a 4-set vector space over S is infinite.

A 4-set which generates V is given by

$$X = \left\{\begin{pmatrix} 1 & 1 & 1 & 1 \\ 1 & 1 & 1 & 1 \end{pmatrix}\right\} \cup \{(1\ 0\ 0), (0\ 1\ 0), (0\ 0\ 1)\} \cup \{1, x, x^2, x^3, x^4, x^5\} \cup \left\{\begin{pmatrix} 1 & 0 \\ 0 & 0 \end{pmatrix}, \begin{pmatrix} 0 & 1 \\ 0 & 0 \end{pmatrix}, \begin{pmatrix} 0 & 0 \\ 1 & 0 \end{pmatrix}, \begin{pmatrix} 0 & 0 \\ 0 & 1 \end{pmatrix}\right\}$$

$$\subseteq V_1 \cup V_2 \cup V_3 \cup V_4.$$

Clearly the 4-dimension of V is (1, 3, 6, 4).

*Note:* It is interesting to observe that V is a 4-set vector space over $S = Z^+ \cup \{0\}$ is not finitely generated. The given X in the example 5.10 does not 4-generate V as a 4-set vector space over $S = Z^+ \cup \{0\}$.

***Example 5.11:*** Let

$$V = V_1 \cup V_2 \cup V_3 \cup V_4 \cup V_5$$
$$= \{(1\ 1\ 1\ 1), (1\ 0\ 0\ 0), (1\ 1\ 0\ 0), (0\ 0\ 1\ 1), (0\ 0\ 0\ 1), (0\ 0\ 1\ 0), (0\ 1\ 0\ 0), (0\ 1\ 1\ 0), (1\ 0\ 0\ 1), (1\ 0\ 1\ 0), (0\ 1\ 0\ 1), (1\ 1\ 1\ 0), (1\ 0\ 1\ 1)\ (0\ 1\ 1\ 1), (1\ 1\ 0\ 1), (0\ 0\ 0\ 0)\} \cup \left\{\begin{pmatrix} a & b \\ c & d \end{pmatrix} \middle| a, b, c, d \in Z_2\right\} \cup$$
$$\left\{\begin{pmatrix} a_1 & a_2 & a_3 \\ a_4 & a_5 & a_5 \end{pmatrix} \middle| a_i \in Z_2 = \{0, 1\}, 1 \leq i \leq 5\right\} \cup \{Z_2[x] \mid x \text{ an indeterminate of degree less than } 4\} \cup \{(1\ 1\ 1\ 1\ 1), (0\ 0\ 0\ 0\ 0)\}$$

be a 5-set vector space over $Z_2 = \{0, 1\}$. Clearly V is a 5-set linear algebra over $Z_2 = \{0, 1\}$. Take

$$X = \{(1\ 0\ 0\ 0), (0\ 1\ 0\ 0), (0\ 0\ 1\ 0), (0\ 0\ 0\ 1) \cup \left\{\begin{pmatrix}1 & 0\\0 & 0\end{pmatrix}, \begin{pmatrix}0 & 1\\0 & 0\end{pmatrix}, \begin{pmatrix}0 & 0\\1 & 0\end{pmatrix}, \begin{pmatrix}0 & 0\\0 & 1\end{pmatrix}\right\} \cup \left\{\begin{pmatrix}1 & 0 & 0\\0 & 0 & 0\end{pmatrix},\right.$$

$$\begin{pmatrix}0 & 1 & 0\\0 & 0 & 0\end{pmatrix}, \begin{pmatrix}0 & 0 & 1\\0 & 0 & 0\end{pmatrix}, \begin{pmatrix}0 & 0 & 0\\1 & 0 & 0\end{pmatrix}, \begin{pmatrix}0 & 0 & 0\\0 & 1 & 0\end{pmatrix},$$

$$\left.\begin{pmatrix}0 & 0 & 0\\0 & 1 & 0\end{pmatrix}\right\} \cup \{1, x, x^2, x^3\} \cup \{(1\ 1\ 1\ 1\ 1)\}$$

$$= X_1 \cup X_2 \cup \ldots \cup X_5$$
$$\subseteq V_1 \cup V_2 \cup \ldots \cup V_5$$
$$= V.$$

Now the 5-dimension of V as a 5-linear algebra is (4, 4, 6, 4, 1). Whereas if we consider V as a 5-set vector space over the set $Z_2$, X is not the 5-generator of V then

$$Y = Y_1 \cup Y_2 \cup \ldots \cup Y_5$$
$= \{(1\ 1\ 1\ 1), (1\ 0\ 0\ 0), (0\ 1\ 0\ 0), (0\ 0\ 1\ 0), (0\ 0\ 0\ 1), (1\ 1\ 0\ 0), (1\ 0\ 1\ 0), (1\ 0\ 0\ 1), (0\ 1\ 1\ 0), (0\ 1\ 0\ 1), (0\ 0\ 1\ 1), (1\ 1\ 1\ 0), (1\ 0\ 1\ 1)\ (0\ 1\ 1\ 1), (1\ 1\ 0\ 1)\} \cup$ {A set with 15 elements) $\cup$ {A set with 63 elements} $\cup$ {A set with 15 elements} $\cup \{(1\ 1\ 1\ 1\ 1)\}$,

5 generates V as a 5 set vector space over {0, 1}. Now the 5-dimension of V is given by (15, 15, 63, 15, 1).

Thus we see by making when ever possible or feasible a n-set vector space into a n-set linear algebra, we can minimize the number of elements in the generating n-set.

Now we proceed on to define the notion of n-set linear transformation of n-set linear algebras.

**DEFINITION 5.6:** *Let $V = V_1 \cup V_2 \cup \ldots \cup V_n$ be a n-set linear algebra over the set S. Let $W = W_1 \cup W_2 \cup \ldots \cup W_n$ be a n-set linear algebra over the same set S. Suppose T is a map from V to W such that $T = T_1 \cup T_2 \cup \ldots \cup T_n : V_1 \cup V_2 \cup \ldots \cup V_n \rightarrow W_1 \cup W_2 \cup \ldots \cup W_n$ where $T_i : V_i \rightarrow W_i$ is a set linear algebra* *transformation i.e., $T_i(\alpha x + y) = T_i(\alpha x) + T_i(y)$; $1 \leq i \leq n$. Then we call T to be a n-set linear transformation from V to W.*

We illustrate this by the following example.

***Example 5.12:*** Let

$$
\begin{aligned}
V &= V_1 \cup V_2 \cup V_3 \cup V_4 \\
&= \{(a\ a\ a) \mid a \in Z_2 = \{0, 1\}\} \cup \left\{\begin{pmatrix} a & a \\ a & a \end{pmatrix} \middle| a \in Z^2 = \{0,1\}\right\} \cup \\
&\quad \{(0\ 0\ 0\ 0), (1\ 1\ 1\ 1), (1\ 1\ 0\ 0), (0\ 0\ 1\ 1)\} \cup \{Z_2 \times Z_2\}
\end{aligned}
$$

be a 4-set linear algebra over $Z_2 = \{0, 1\}$ and

$$
\begin{aligned}
W &= W_1 \cup W_2 \cup W_3 \cup W_4 \\
&= \{Z_2 \times Z_2 \times Z_2\} \cup \left\{\begin{pmatrix} a & b \\ c & d \end{pmatrix} \middle| a \in Z_2\right\} \cup \\
&\quad \left\{\begin{pmatrix} a & a & a & a \\ a & a & a & a \end{pmatrix} \middle| a \in Z_2\right\} \cup \{(1\ 1\ 1\ 1\ 1), (0\ 0\ 0\ 0\ 0), \\
&\quad (1\ 1\ 0\ 0\ 0), (0\ 0\ 1\ 1\ 1)\}
\end{aligned}
$$

is also a 4-set linear algebra over the set $Z_2$. Define $T: V \to W$ by $T = T_1 \cup T_2 \cup T_3 \cup T_4: V_1 \cup V_2 \cup V_3 \cup V_4 \to (W_1 \cup W_2 \cup W_3 \cup W_4)$ with $T_i : V_i \to W_i$ ; $1 \leq i \leq 4$. $T_1 : V_1 \to W_1$ is given by

$$T_1((a\ a\ a)) = (a\ a\ a).$$

$T_2: V_2 \to W_2$ is given by

$$T_2\left(\begin{pmatrix} a & a \\ a & a \end{pmatrix}\right) = \begin{pmatrix} a & a \\ a & a \end{pmatrix}.$$

$T_3 : V_3 \to W_3$ is defined by

$$T_3\,(0\ 0\ 0\ 0) = \begin{pmatrix} 0 & 0 & 0 & 0 \\ 0 & 0 & 0 & 0 \end{pmatrix};$$

$$T_3\,(1\ 1\ 1\ 1) = \begin{pmatrix} a & a & a & a \\ a & a & a & a \end{pmatrix};\ a = 0$$

$$T_3\,(1\ 1\ 0\ 0) = \begin{pmatrix} 1 & 1 & 1 & 1 \\ 1 & 1 & 1 & 1 \end{pmatrix}$$

and

$$T_3\,(0\ 0\ 1\ 1) = \begin{pmatrix} 1 & 1 & 1 & 1 \\ 1 & 1 & 1 & 1 \end{pmatrix}.$$

$T_4 : V_4 \to W_4$ is given by

$$\begin{aligned} T_4\,(x, y) &= (1\ 1\ 1\ 1\ 1) \quad x \neq 0 \\ T_4\,(0, 0) &= (0\ 0\ 0\ 0\ 0) \quad y \neq 0 \\ T_4\,(1, 0) &= (0\ 0\ 1\ 1\ 1) \end{aligned}$$

and

$$T_4\,(0, 1) = (1\ 1\ 0\ 0\ 0).$$

It is easily verified that $T = T_1 \cup T_2 \cup T_3 \cup T_4$ is a 4-set linear transformation from V to W.

Now we proceed onto define n-set linear operator on a n-set linear algebra.

**DEFINITION 5.7:** *Let $V = V_1 \cup V_2 \cup \ldots \cup V_n$ be a n-set linear algebra over the set S. A map $T = T_1 \cup T_2 \cup \ldots \cup T_n$ from V to V is said to be a n-set linear operator if $T_i : V_i \to V_i$ is such that $T_i\,(\alpha x + y) = \alpha T_i(x) + T_i(y)$ for $x, y \cup V_i$; $\alpha \in S$; $1 \leq i \leq n$.*

We will illustrate this situation by the following example.

***Example 5.13:*** Let

$$V = V_1 \cup V_2 \cup V_3 \cup V_4$$

$$= \{Z^+ \times Z^+\} \cup \left\{\begin{pmatrix} a & b \\ c & d \end{pmatrix} \middle| a, b, c, d \in Z^+\right\} \cup \{Z^+[x]\} \cup \left\{\begin{pmatrix} a & b & c & g \\ d & e & f & h \end{pmatrix} \middle| a,b,c,d,e,f,g,h \in Z^+\right\}$$

be a 4-set linear algebra over the set $Z^+$. A map $T = T_1 \cup T_2 \cup T_3 \cup T_4 : V \to V$; such that

$$T_1(x, y) = (y, x);$$

$$T_2\left(\begin{pmatrix} a & b \\ c & d \end{pmatrix}\right) = \begin{pmatrix} a & b \\ c & c \end{pmatrix}$$

$$T_3 (p(x) = p_0 + p_1x + \ldots + p_nx^n) = p_nx^n$$

and

$$T_4\left(\begin{pmatrix} a & b & c & g \\ d & e & f & h \end{pmatrix}\right) = \begin{pmatrix} a & a & a & a \\ d & d & d & d \end{pmatrix}.$$

It is easily verified that T is a 4-set linear operator on V.

Now in case of n-set linear algebra we can define a notion called n-set quasi linear operator on V.

**DEFINITION 5.8:** *Let $V = V_1 \cup V_2 \cup \ldots \cup V_n$ be a n-set linear algebra over the set S. An map $T = T_1 \cup T_2 \cup \ldots \cup T_n$ from $V = V_1 \cup V_2 \cup \ldots \cup V_n$ to $V = V_1 \cup V_2 \cup \ldots \cup V_n$ such that $T_i : V_i \to V_j$; $i \neq j$, $1 \leq i, j \leq n$ such that $T_i (\alpha u + v) = \alpha T_i(u) + T_i(v)$ for u, v $\in V_i$ and $\alpha \in S$ for each i, $1 \leq i \leq n$ is called the n-set quasi linear operator on V.*

*If on the other hand $T = T_1 \cup T_2 \cup \ldots \cup T_n$ is such that $T_i : V_i \to V_j$ for some $j \neq i$ and $T_k : V_k \to V_k$ for some $1 \leq k \leq n$ then we call T to be a n-set semiquasi linear operator on V.*

We cannot get any interrelation between n-set linear operator, n-set linear semiquasi operator and n-set linear quasi operator.

The reader can substantiate these definitions with examples. However problems based on these definitions are given in the last chapter of this book. Interested reader can refer to them.

Now we proceed on to define the notion of semigroup n-set vector space.

**DEFINITION 5.9:** *Let $V = V_1 \cup V_2 \cup \ldots \cup V_n$ be a n-set vector space over the set S. If S is an additive semigroup then we call V to be a semigroup n-set vector space over the semigroup S.*

We illustrate this situation by some examples.

***Example 5.14:*** Let $V = V_1 \cup V_2 \cup \ldots \cup V_n$ be a semigroup 4-set vector space over the semigroup $S = Z_2 = \{0, 1\}$ under addition modulo 2, where

$$V_1 = \left\{ \begin{pmatrix} a & b \\ c & d \end{pmatrix} \middle| \; a, b, c, d \in Z_2 \right\},$$

$$V_2 = \{Z_2 \times Z_2 \times Z_2\},$$

$$V_3 = \{(1\ 1\ 0\ 1), (0\ 1\ 0\ 0), (0\ 0\ 0\ 0), (1\ 1\ 1), (0\ 0\ 0), (1\ 0\ 0)\}$$

and

$$V_4 = \left\{ \begin{pmatrix} a & a \\ a & a \end{pmatrix}, \begin{pmatrix} a & a & a & a \\ a & a & a & a \end{pmatrix} \middle| \; a \in Z_2 \right\}.$$

V is a semigroup 4-set vector space over the semigroup $S = Z_2 = \{0,1\}$.

***Example 5.15:*** Let

$$\begin{aligned} V &= V_1 \cup V_2 \cup V_3 \cup V_4 \cup V_5 \cup V_6 \\ &= \{Z^+ \times Z^+ \times Z^+\} \cup \left\{ \begin{pmatrix} a & b & c \\ d & e & f \end{pmatrix} \middle| a, b, c, d, e, f \in Z^+ \right\} \cup \\ &\quad \{(a, b, c, d), (a, a, a, a, b, c) \mid a, b, c, d \in Z^+\} \cup \{(a, b), \end{aligned}$$

$(a, a, a, a, a) \mid a, b \in Z^+\} \cup \{Z^+[x] \mid x \text{ an indeterminate}\} \cup \{3Z^+ \times 2Z^+ \times 5Z^+ \times 7Z^+\}$

is a semigroup n-vector space over the semigroup $Z^+$.

***Example 5.16:*** Let

$$V = V_1 \cup V_2 \cup V_3 \cup V_4 \cup V_5$$

$$= \left\{\begin{pmatrix} a & b \\ c & d \end{pmatrix} \middle| a, b, c, d \in Z_3\right\} \cup \{Z_3 \times Z_3 \times Z_3\} \cup$$

$$\left\{\begin{pmatrix} a & 0 & 0 \\ b & c & d \end{pmatrix} \middle| a, b, c, d \in Z_3\right\} \cup Z_3[x] \cup$$

$$\left\{\begin{bmatrix} a & b & c & d \\ 0 & c & f & g \\ 0 & 0 & i & j \\ 0 & 0 & 0 & k \end{bmatrix} \middle| a, b, c, d, e, f, g, i, j, k \in Z_3\right\}$$

be a semigroup 5-set vector space over $Z_3$.

We now prove that all semigroup n-set vector spaces are n-set vector spaces but a n-set vector space in general is not a semigroup n-vector space.

**THEOREM 5.1:** *Let $V = V_1 \cup V_2 \cup \ldots \cup V_n$ be a semigroup n set vector space over the semigroup S, then $V = V_1 \cup V_2 \cup \ldots \cup V_n$ is the n-set vector space over S. Conversely if $V = V_1 \cup V_2 \cup \ldots \cup V_n$ is a n-set vector space over S; then V in general is not a semigroup n-set vector space over the set S.*

*Proof:* Let $V = V_1 \cup V_2 \cup \ldots \cup V_n$ be a semigroup n-set vector space over the semigroup S. Now every semigroup is a set S. So V is a n-set vector space over the set S. To prove a n-set vector space over the set in general is not a semigroup n-set vector space. Suppose

$$V = V_1 \cup V_2 \cup V_3 \cup V_4$$

$$= \{(0\ 0\ 0\ 0), (1\ 1\ 1\ 1), (2\ 2\ 2\ 2), \ldots, (2n\ 2n\ 2n\ 2n)\} \cup \left\{\begin{pmatrix} 2n & 2n \\ 2n & 2n \end{pmatrix} \middle| n \in Z^+\right\} \cup \left\{\begin{pmatrix} 2n \\ 2n \\ 2n \end{pmatrix} \middle| n \in Z^+\right\} \cup \{2Z^+[x]\}$$

is a 4-set vector space over the set S = {1, 2}. Clearly S is not a semigroup under addition. Hence the claim.

Now we proceed onto define some substructures on semigroup n-set vector spaces and illustrate them by examples.

**DEFINITION 5.10:** *Let $V = V_1 \cup V_2 \cup \ldots \cup V_n$ be a semigroup n-set vector space over the semigroup S. Suppose $W = W_1 \cup W_2 \cup \ldots \cup W_n \subseteq V = V_1 \cup V_2 \cup \ldots \cup V_n$ and each $W_i \subseteq V_i$ is a semigroup set vector space over the semigroup S then we call W to be semigroup n-set vector subspace of V over the semigroup S, $(1 \leq i \leq n)$.*

***Example 5.17:*** Let

$$V = V_1 \cup V_2 \cup V_3 \cup V_4 \cup V_5$$
$$= \left\{\begin{pmatrix} a & a & a \\ a & a & a \end{pmatrix} \middle| a \in Z_4\right\} \cup \{Z_4 \times Z_4\} \cup \left\{\begin{pmatrix} a & a \\ a & a \\ a & a \end{pmatrix} \middle| a \in Z_4\right\}$$
$$\cup \{(0\ 0\ 0\ 0), (2\ 2\ 2\ 2), (0\ 0\ 0), (2\ 2\ 2), (0\ 0\ 0\ 0\ 0), (2\ 2\ 2\ 2\ 2)\} \cup \left\{\begin{pmatrix} a & a \\ a & a \end{pmatrix} \middle| a \in Z_4\right\}$$

be a semigroup 5-set vector space over the semigroup H = {0, 2} $\subseteq Z_4$ = {0, 1, 2, 3} under addition modulo 4. Take

$$W = W_1 \cup W_2 \cup W_3 \cup W_4 \cup W_5$$
$$= \left\{\begin{pmatrix} 0 & 0 & 0 \\ 0 & 0 & 0 \end{pmatrix}, \begin{pmatrix} 2 & 2 & 2 \\ 2 & 2 & 2 \end{pmatrix}\right\} \cup \{(0\ 0), (2\ 2)\} \cup$$

$$\left\{\begin{pmatrix}0&0\\0&0\\0&0\end{pmatrix},\begin{pmatrix}2&2\\2&2\\2&2\end{pmatrix}\right\} \cup \{(0\ 0\ 0\ 0),\ (2\ 2\ 2\ 2)\} \cup$$

$$\left\{\begin{pmatrix}0&0\\0&0\end{pmatrix},\begin{pmatrix}2&2\\2&2\end{pmatrix}\right\}$$

$$\subseteq V_1 \cup V_2 \cup V_3 \cup V_4 \cup V_5$$
$$= V,$$

is a semigroup 5-set vector subspace over the semigroup $S = \{0, 2\}$ under addition modulo 4.

***Example 5.18:*** Let

$$V = V_1 \cup V_2 \cup V_3 \cup V_4$$
$$= \{Z_5 \times Z_5 \times Z_5\} \cup \left\{\begin{pmatrix}a&a&a\\a&a&a\end{pmatrix} \middle| a \in Z_5\right\} \cup$$
$$\left\{\begin{pmatrix}a&a&a\\a&a&a\\a&a&a\end{pmatrix} \middle| a \in Z_5\right\} \cup \{Z_5[x]\}$$

be a semigroup 4-set vector space over the semigroup $S = Z_5$, the semigroup under addition modulo 5.
Let

$$W = W_1 \cup W_2 \cup W_3 \cup W_4$$
$$= \{(0\ 0\ 0),\ (1\ 1\ 1)\} \cup \left\{\begin{pmatrix}1&1&1\\1&1&1\end{pmatrix},\begin{pmatrix}2&2&2\\2&2&2\end{pmatrix},\begin{pmatrix}0&0&0\\0&0&0\end{pmatrix}\right\}$$
$$\cup \left\{\begin{pmatrix}0&0&0\\0&0&0\\0&0&0\end{pmatrix},\begin{pmatrix}3&3&3\\3&3&3\\3&3&3\end{pmatrix}\right\} \cup \{a+ax+\ldots+ ax^n \mid a \in Z^5\}.$$

W is not a semigroup 4-set vector subspace over the semigroup $Z_5$.

***Example 5.19:*** Let

$$V = V_1 \cup V_2 \cup V_3 \cup V_4 \cup V_5$$

$$= \left\{\begin{pmatrix} a & a & a \\ a & a & a \end{pmatrix} \middle| a \in Z_6 \right\} \cup \{Z_6 \times Z_6 \times Z_6\} \cup \{Z_6[x]\} \cup$$

$$\left\{\begin{pmatrix} a & a \\ a & a \end{pmatrix}, \begin{pmatrix} a \\ a \\ a \end{pmatrix} \middle| a \in Z_6 \right\} \cup \left\{\begin{pmatrix} a & & & (0) \\ & b & & \\ & & c & \\ (0) & & & d \end{pmatrix} \middle| a, b, c, d \in Z_6 \right\}$$

be a semigroup 5-set vector space over the semigroup S = {0, 3} under addition modulo 6.

$$W = W_1 \cup W_2 \cup \ldots \cup W_5$$

$$= \left\{\begin{pmatrix} 1 & 1 & 1 \\ 1 & 1 & 1 \end{pmatrix}, \begin{pmatrix} 0 & 0 & 0 \\ 0 & 0 & 0 \end{pmatrix}, \begin{pmatrix} 2 & 2 & 2 \\ 2 & 2 & 2 \end{pmatrix}, \begin{pmatrix} 3 & 3 & 3 \\ 3 & 3 & 3 \end{pmatrix}\right\} \cup \{(1\ 1\ 1),$$

$$(0\ 0\ 0), (2\ 2\ 2), (3\ 3\ 3)\} \cup \{a + ax + ax^2 + \ldots + ax^n \mid a \in$$

$$Z_6\} \cup \left\{\begin{pmatrix} a & a \\ a & a \end{pmatrix} \middle| a \in Z_6 \right\} \cup \left\{\begin{pmatrix} a & & & (0) \\ & a & & \\ & & a & \\ (0) & & & a \end{pmatrix} \middle| a \in Z_6 \right\}$$

$$\subseteq V_1 \cup V_2 \cup \ldots \cup V_5$$

$$= V$$

is a semigroup 5-set subvector space over the semigroup S = {0, 3}.

***Example 5.20:*** Let

$$V = V_1 \cup V_2 \cup V_3 \cup V_4 \cup V_5 \cup V_6$$

$$= \{Z_2 \times Z_2\} \cup \{Z_2[x]\} \cup \left\{\begin{pmatrix} a & a \\ a & a \\ a & a \end{pmatrix} \middle| a \in Z_2 \right\}$$

$$\left\{\begin{pmatrix} a & a \\ a & a \end{pmatrix}, \begin{pmatrix} a & a & a \\ a & a & a \end{pmatrix} \middle| a \in Z_2 \right\} \cup \left\{ \begin{pmatrix} a & & & (0) \\ & a & & \\ & & a & \\ (0) & & & a \end{pmatrix} \middle| a \in Z_2 \right\}$$

$$\cup \{(a\ a\ a\ a\ a\ a), (a\ a\ a) \mid a \in Z_2\}$$

is a semigroup set 6-vector space over the semigroup S = $Z_2$ under addition modulo 2. Let

$$\begin{aligned} W &= W_1 \cup W_2 \cup \ldots \cup W_6 \\ &= \{(0\ 0), (1\ 1)\} \cup \{a + ax + \ldots + ax^n \mid a \in Z_2\} \cup \\ &\quad \left\{ \begin{pmatrix} 1 & 1 \\ 1 & 1 \\ 1 & 1 \end{pmatrix}, \begin{pmatrix} 0 & 0 \\ 0 & 0 \\ 0 & 0 \end{pmatrix} \right\} \cup \left\{ \begin{pmatrix} a & a & a \\ a & a & a \end{pmatrix} \middle| a \in Z_2 \right\} \\ &\quad \cup \left\{ \begin{pmatrix} 0 & 0 & 0 & 0 & 0 \\ 0 & 0 & 0 & 0 & 0 \\ 0 & 0 & 0 & 0 & 0 \\ 0 & 0 & 0 & 0 & 0 \\ 0 & 0 & 0 & 0 & 0 \end{pmatrix}, \begin{pmatrix} 1 & & (0) \\ & \ddots & \\ (0) & & 1 \end{pmatrix} \right\} \\ &\quad \cup \{(a\ a\ a) \mid a \in Z_2\} \end{aligned}$$

is a semigroup 6-set subvector space over S = {1, 0} ⊆ $Z_2$.

Now we proceed on to define the new type of substructure in semigroup n-set vector space.

**DEFINITION 5.11:** *Let V = $V_1 \cup V_2 \cup \ldots \cup V_n$ be a semigroup n-set vector space over the semigroup S. Let W = $W_1 \cup W_2 \cup \ldots \cup W_n \subseteq V_1 \cup V_2 \cup \ldots \cup V_n$, where W is a n-set vector space over the set P then W is a pseudo n-set vector space over the subset P ⊆ S.*

We illustrate them by some examples.

***Example 5.21:*** Let

$$
\begin{aligned}
V &= V_1 \cup V_2 \cup \ldots \cup V_5 \\
&= \left\{ \begin{pmatrix} a & a \\ a & a \end{pmatrix} \middle| a \in Z_{12} \right\} \cup \\
&\quad \left\{ \begin{pmatrix} a & a & a \\ a & a & a \end{pmatrix}, \begin{pmatrix} a & a & a & a & a \\ a & a & a & a & a \end{pmatrix} \middle| a \in Z_{12} \right\} \cup \\
&\quad \left\{ \begin{pmatrix} a_1 & & & & 0 \\ & a_2 & & & \\ & & a_3 & & \\ & & & a_4 & \\ 0 & & & & a_5 \end{pmatrix} \middle| a_1, a_2, \ldots, a_5 \in Z_{12} \right\} \cup \{Z_{12}[x]\} \\
&\quad \cup \{(a\ a\ a\ a\ a), (a\ a\ a), (a\ a\ a\ a) \mid a \in Z_{12}\}
\end{aligned}
$$

be a semigroup 5-set vector space over the semigroup $S = Z_{12}$, a semigroup under addition modulo 12. Take

$$
\begin{aligned}
W &= W_1 \cup W_2 \cup W_3 \cup W_4 \cup W_5 \\
&= \left\{ \left\{ \begin{pmatrix} 1 & 1 \\ 1 & 1 \end{pmatrix}, \begin{pmatrix} 0 & 0 \\ 0 & 0 \end{pmatrix}, \begin{pmatrix} 6 & 6 \\ 6 & 6 \end{pmatrix} \right\} \cup \left\{ \begin{pmatrix} a & a & a \\ a & a & a \end{pmatrix} \middle| a \in Z_{12} \right\} \cup \right. \\
&\quad \left\{ \begin{pmatrix} a & 0 & 0 & 0 & 0 \\ 0 & a & 0 & 0 & 0 \\ 0 & 0 & a & 0 & 0 \\ 0 & 0 & 0 & a & 0 \\ 0 & 0 & 0 & 0 & a \end{pmatrix} \middle| a \in Z_{12} \right\} \cup \{a + ax + \ldots + ax^n \mid a \\
&\quad \in Z_{12}\} \cup \{(a\ a\ a\ a\ a) \mid a \in Z_{12}\} \\
&\subseteq V_1 \cup V_2 \cup \ldots \cup V_5;
\end{aligned}
$$

W is a pseudo 5-set vector subspace of V over the set $S = \{0, 1, 6\} \subseteq Z_{12}$.

***Example 5.22:*** Let

$$V = V_1 \cup V_2 \cup V_3 \cup V_4$$
$$= \{Z^+ [x]\} \cup \{Z^+ \times Z^+ \times Z^+\} \cup \left\{ \begin{pmatrix} a & a \\ a & a \end{pmatrix} \middle| a \in Z^+ \right\} \cup$$
$$\left\{ \begin{pmatrix} a & a & a & a \\ a & a & a & a \end{pmatrix}, \begin{pmatrix} a & a & a \\ a & a & a \\ a & a & a \end{pmatrix} \middle| a \in Z^+ \right\};$$

be a semigroup 4-set vector space over the semigroup S = $Z^+$, semigroup under addition.

$$W = W_1 \cup W_2 \cup W_3 \cup W_4$$
$$= \{a + ax + \ldots + ax^n \mid a \in Z^+\} \cup \{(a\ a\ a) \mid a \in Z^+\} \cup$$
$$\left\{ \begin{pmatrix} 2n & 2n \\ 2n & 2n \end{pmatrix} \middle| n \in Z^+ \right\} \cup \left\{ \begin{pmatrix} a & a & a \\ a & a & a \\ a & a & a \end{pmatrix} \middle| a \in Z^+ \right\}$$
$$\subseteq V_1 \cup V_2 \cup V_3 \cup V_4,$$

is a pseudo 4-set vector subspace of V over the subset {0, 2, 4, 3, 7} $\subseteq Z^+$.

Now we proceed on to define the notion of semigroup n-set linear algebra over the semigroup.

**DEFINITION 5.12:** *Let $V = V_1 \cup \ldots \cup V_n$ be a semigroup n-set vector space over the semigroup S. If each $V_i$ is a semigroup under addition, $1 \leq i \leq n$; then we call V to be the semigroup n set linear algebra over the semigroup S.*

We now illustrate this definition by some examples.

***Example 5.23:*** Let

$$V = V_1 \cup V_2 \cup \ldots \cup V_5$$

$$= \left\{\begin{pmatrix} a & a \\ a & a \end{pmatrix} \middle| a \in Z^+ \right\} \cup \{(a\ a\ a\ a) \mid a \in Z^+\} \cup \{Z_2[x]\} \cup$$

$$\left\{\begin{pmatrix} a & b & c & d \\ e & f & g & h \end{pmatrix} \middle| a, b, c, d, e, f, g, h \in Z^+ \right\} \cup$$

$$\left\{\begin{pmatrix} a & b & c \\ d & e & f \\ g & h & i \end{pmatrix} \middle| a, b, c, d, e, f, g, h, i \in Z^+ \right\}$$

be the semigroup 5-set linear algebra over the semigroup $Z^+$.

***Example 5.24:*** Let

$$V = V_1 \cup V_2 \cup V_3 \cup V_4$$

$$= \left\{\begin{pmatrix} a & a & a \\ a & a & a \\ a & a & a \end{pmatrix} \middle| a \in Z_2 \right\} \cup \{Z_2\,[x]\} \cup \{Z_2 \times Z_2 \times Z_2 \times Z_2\}$$

$$\cup \left\{\begin{pmatrix} a & a & a & a & a \\ a & a & a & a & a \end{pmatrix} \middle| a \in Z_2 \right\}$$

be a semigroup 4-set linear algebra over the semigroup $Z_2$.

***Example 5.25:*** Let

$$V = V_1 \cup V_2 \cup V_3 \cup V_4 \cup V_5$$

$$= \{Z_3[x]\} \cup \{Z_3 \times Z_3 \times Z_3\} \cup \left\{\begin{pmatrix} a & a \\ a & a \end{pmatrix} \middle| a \in Z_3 \right\} \cup$$

$$\{(a\ a\ a\ a\ a\ a) \mid a \in Z_3\} \cup \left\{\begin{pmatrix} a & a & a & a \\ a & a & a & a \\ a & a & a & a \end{pmatrix} \middle|\ a \in Z_3\right\}$$

be a semigroup 5-set linear algebra over $Z_3$.

***Example 5.26:*** Let

$$\begin{aligned} V &= V_1 \cup V_2 \cup V_3 \cup V_4 \\ &= \left\{\begin{pmatrix} a & a & a \\ a & a & a \\ a & a & a \end{pmatrix} \middle|\ a \in Z_7\right\} \cup \{Z_7 \times Z_7 \times Z_7\} \cup \{Z_7[x]\} \cup \\ &\quad \left\{\begin{pmatrix} a & a & a & a \\ b & b & b & b \end{pmatrix} \middle|\ a, b \in Z_7\right\} \end{aligned}$$

be a semigroup 4-set linear algebra over the semigroup $Z_7$.

**DEFINITION 5.13:** *Let $V = V_1 \cup V_2 \cup V_3 \cup \ldots \cup V_n$ be a semigroup n-set linear algebra over the semigroup S. Let $W = W_1 \cup W_2 \cup \ldots \cup W_n \subseteq V_1 \cup V_2 \cup \ldots \cup V_n = V$ be a n-subset of V such that W is a semigroup n-set linear algebra over the semigroup S, then we call W to be the semigroup n-set linear subalgebra of V over S.*

We illustrate this by the following example.

***Example 5.27:*** Let

$$\begin{aligned} V &= V_1 \cup V_2 \cup V_3 \cup V_4 \cup V_5 \\ &= \{Z_6 \times Z_6 \times Z_6\} \cup \{Z_6[x]\} \cup \left\{\begin{pmatrix} a & a & a \\ a & a & a \end{pmatrix} \middle|\ a \in Z_6\right\} \cup \\ &\quad \left\{\begin{pmatrix} a & b \\ c & d \\ e & f \end{pmatrix} \middle|\ a, b, c, d, e, f \in Z_6\right\} \cup \end{aligned}$$

$$\left\{\begin{pmatrix} a & b & c \\ d & e & f \\ g & h & i \end{pmatrix} \middle| a,b,c,d,e,f,g,h,i \in Z_6 \right\}$$

be a semigroup 5-set linear algebra over the semigroup $Z_6$, $Z_6$ a semigroup under addition modulo 6. Take

$$\begin{aligned} W &= W_1 \cup W_2 \cup W_3 \cup W_4 \cup W_5 \\ &= \{Z_6 \times Z_6 \times \{0\}\} \cup \{a + ax + \ldots + ax^n \mid a \in Z\} \cup \\ &\quad \left\{\begin{pmatrix} a & a & a \\ a & a & a \end{pmatrix} \middle| a \in \{0,2,4\}\right\} \cup \left\{\begin{pmatrix} a & a \\ a & a \\ a & a \end{pmatrix} \middle| a \in Z_6 \right\} \cup \\ &\quad \left\{\begin{pmatrix} a & a & a \\ a & a & a \\ a & a & a \end{pmatrix} \middle| a \in Z_6 \right\} \end{aligned}$$

is a semigroup 5-set linear subalgebra over the semigroup $Z_6$.

***Example 5.28:*** Let

$$\begin{aligned} V &= V_1 \cup V_2 \cup \ldots \cup V_6 \\ &= \{Z^+ \times Z^+ \times Z^+\} \cup \left\{\begin{pmatrix} a & a & a \\ a & a & a \end{pmatrix} \middle| a \in Z^+ \right\} \cup \{Z^+[x]\} \cup \\ &\quad \left\{\begin{pmatrix} a & a & a & a \\ a & a & a & a \end{pmatrix} \middle| a \in Z^+ \right\} \cup \left\{\begin{pmatrix} a & a \\ a & a \end{pmatrix} \middle| a \in Z^+ \right\} \cup \\ &\quad \left\{\begin{pmatrix} a & a & a \\ a & b & b \\ c & c & c \end{pmatrix} \middle| a, b, c \in Z^+ \right\} \end{aligned}$$

be a semigroup 6 set linear algebra over the semigroup $Z^+$ under addition.

Take

$$W = \{3Z^+ \times 3Z^+ \times 3Z^+\} \cup \left\{\begin{pmatrix} a & a & a \\ a & a & a \end{pmatrix} \middle| a \in 2Z^+ \right\} \cup \{a + ax + \ldots + ax^6 \mid a \in 5Z^+\} \cup \left\{\begin{pmatrix} a & a & a & a \\ a & a & a & a \end{pmatrix} \middle| a \in 3Z^+ \right\} \cup$$

$$\left\{\begin{pmatrix} a & a \\ a & a \end{pmatrix} \middle| a \in 2Z^+ \right\} \cup \left\{\begin{pmatrix} a & a & a \\ a & a & a \\ a & a & a \end{pmatrix} \middle| a \in 3Z^+ \right\}$$

$$= W_1 \cup \ldots \cup W_6$$
$$\subseteq V_1 \cup V_2 \cup V_3 \cup V_4 \cup V_5 \cup V_6$$
$$= V.$$

W is a semigroup 6-set linear subalgebra of V over the semigroup $Z^+$. We call W to be a semigroup 6 – set linear subalgebra over the semigroup $Z^+$.

***Example 5.29:*** Let

$$V = V_1 \cup V_2 \cup V_3 \cup V_4$$
$$= \{Z_8 \times Z_8\} \cup \{Z_8 [x] \mid Z_8 [x] \text{ contains all polynomials of degree less than or equal to } 7\} \cup$$
$$\left\{\begin{pmatrix} a & b \\ c & d \end{pmatrix} \middle| a, b, c, d \in Z_8 \right\} \cup \left\{\begin{pmatrix} a & a & a \\ a & a & a \end{pmatrix}, \begin{pmatrix} a & a \\ a & a \\ a & a \end{pmatrix} \middle| a \in Z_8 \right\}$$

be a semigroup 4-set vector space over $Z_8$. This can never be made into a semigroup 4-set linear algebra over $Z_8$.

Thus all semigroup n-set vector spaces are not in general semigroup n-set linear algebras.

Now we define yet another new subalgebraic structure of semigroup n-set vector spaces and semigroup n-set linear algebras.

**DEFINITION 5.14:** *Let $V = V_1 \cup \ldots \cup V_n$ be a semigroup n – set vector space over the semigroup S. If $W = W_1 \cup W_2 \cup \ldots \cup W_n \subseteq V_1 \cup V_2 \cup \ldots \cup V_n = V$ is a proper subset of V and W is a semigroup n-set vector space over a proper subsemigroup P of S then we call W to be the subsemigroup n set vector subspace of V over the subsemigroup of the semigroup S.*

Suppose $V_1 \cup V_2 \cup \ldots \cup V_n$ is a semigroup n – linear algebra over the semigroup S and $W = W_1 \cup W_2 \cup \ldots \cup W_n$ be a proper subset of V such that $W_i \subseteq V_i$ and $W_i$ is a subsemigroup of $V_i$ for each i, $1 < i \leq n$ and if $P \subseteq S$ is a proper subsemigroup of the semigroup S and if W is a semigroup n-linear algebra over the semigroup $P \subset S$ then we call W to be the subsemigroup n-linear subalgebra of V over the subsemigroup P of the semigroup S.

Now we illustrate this situation by a few examples.

***Example 5.30:*** Let

$$\begin{aligned} V &= V_1 \cup V_2 \cup V_3 \cup V_4 \\ &= \{(1000), (0000), (0011), (1100)\} \cup \{Z_2 \times Z_2 \times Z_2\} \cup \\ &\quad \{(111), (001), (000), (11100), (00011), (100), \\ &\quad (00000)\} \cup \left\{ \begin{pmatrix} a & a & a \\ a & a & a \end{pmatrix} \middle| a \in Z_2 \right\} \end{aligned}$$

be a semigroup 4 set vector space over the semigroup $Z_2 = \{0,1\}$. Clearly V has no subsemigroup 4-set vector space as $Z_2$ has no proper subsemigroup.

***Example 5.31:*** Let

$$\begin{aligned} V &= V_1 \cup V_2 \cup V_3 \cup V_4 \\ &= \left\{ \begin{pmatrix} a & b \\ c & d \end{pmatrix} \middle| a, b, c, d \in Z_4 \right\} \cup \{(a, a, a, a) \,/\, a \in Z_4\} \cup \{Z_4 \\ &\quad \times Z_4 \times Z_4\} \cup \{(Z_4 [x]\} \end{aligned}$$

be a semigroup 4-set vector space over the semigroup $Z_4$. Take

$$W = \left\{\begin{pmatrix} a & a \\ a & a \end{pmatrix} \middle| a \in Z_4\right\} \cup \{(0000), (2222)\} \cup \{(a\ a\ a) |\ a \in Z_4\}$$
$$\cup \{a_1 + a_2x^2 + a_3x^4 + a_4x^6 + a_5x^8 \mid a_1, a_2, a_3, a_4, a_5 \in Z_4\}$$

is a subsemigroup 4-set vector subspace of V over the subsemigroup S = {0, 2} addition modulo 4.

On similar lines we can define the notion of subsemigroup n set linear subalgebra of V as follows.

**DEFINITION 5.15:** *Let $V = V_1 \cup V_2 \cup \ldots \cup V_n$ be a semigroup n-set linear algebra over the semigroup S. If $W = W_1 \cup W_2 \cup \ldots \cup W_n$ be a proper subset of V and if P is any proper subsemigroup of the semigroup S. We call W to be the subsemigroup n-set linear subalgebra of V over the subsemigroup P of the semigroup of S; if W is a semigroup n-set linear algebra over P.*

***Example 5.32:*** Let

$$V = V_1 \cup V_2 \cup \ldots \cup V_5$$
$$= \{Z_{16} \times Z_{16} \times Z_{16}\} \cup \{Z_{16}[x]; \text{ all polynomials of degree less than or equal to } 5\} \cup \left\{\begin{pmatrix} a & a & a \\ a & a & a \end{pmatrix} \middle| a \in Z_{16}\right\} \cup$$
$$\left\{\begin{pmatrix} a & b \\ c & d \end{pmatrix} \middle| a, b, c, d \in Z_{16}\right\} \cup \left\{\begin{pmatrix} a & a & a & a \\ 0 & a & a & a \\ 0 & 0 & a & a \\ 0 & 0 & 0 & a \end{pmatrix} \middle| a \in Z_{16}\right\}$$

be a semigroup 5-set linear algebra over the semigroup $S = Z_{16}$. Take

$$W = W_1 \cup W_2 \cup W_3 \cup W_4 \cup W_5$$
$$= \{S \times S \times S \mid S = \{0,4,8\}\} \cup \{a + ax + ax^2 + ax^3 + ax^4 +$$

$$ax^5 \mid a \in Z_{16}\} \cup \left\{\begin{pmatrix} a & a & a \\ a & a & a \end{pmatrix} \middle| a \in \{0,2,4,6,8,10,12,14\}\right\}$$

$$\cup \left\{\begin{pmatrix} a & a \\ a & a \end{pmatrix} \middle| a \in \{0,8\}\right\} \cup \left\{\begin{pmatrix} a & a & a & a \\ 0 & a & a & a \\ 0 & 0 & a & a \\ 0 & 0 & 0 & a \end{pmatrix} \middle| a \in \{0,4,8,12\}\right\}$$

is a subsemigroup 4-set linear subalgebra over the subsemigroup $S = \{0, 4, 8, 12\} \subseteq Z_{16}$.

Not all semigroup n-set linear algebras have subsemigroup n-set linear subalgebras. We give classes of semigroup n-set vector spaces and semigroup n-set linear algebras which do not contain this type of substructures.

**THEOREM 5.2:** *Let $V = V_1 \cup V_2 \cup \ldots \cup V_n$ be a semigroup n-set vector space over the semigroup $Z_p$ (p a prime) under addition modulo p. V does not contain any proper subsemigroup n-set semigroup vector subspaces.*

*Proof:* Given $V = V_1 \cup \ldots \cup V_n$ is a semigroup n-set vector space over the semigroup $Z_p$, p a prime. Clearly $Z_p$ has no proper subsemigroup. So even if $W = W_1 \cup \ldots \cup W_n \subseteq V_1 \cup V_2 \cup \ldots \cup V_n = V$ a proper subset of V then also W is not a subsemigroup n-set vector subspace of V as $Z_p$ has no proper subsemigroups under addition modulo p. Hence the claim.

Thus in view of this theorem we give the following interesting definition.

**DEFINITION 5.16:** *Let $V = V_1 \cup V_2 \cup \ldots \cup V_n$ be a n-set vector space over a semigroup S. If S has no proper subsemigroup n-set vector subspace then we call V to be a pseudo simple semigroup n-set vector space.*

Now we will also prove we have a class of semigroup n-set vector spaces which are pseudo simple.

**THEOREM 5.3:** *Let $V = V_1 \cup \ldots \cup V_n$ be a semigroup n-set linear algebra over the semigroup $Z_p$, p a prime. Then V has no proper subsemigroup n-set linear subalgebra.*

*Proof:* Given $V = V_1 \cup \ldots \cup V_n$ is a semigroup n-set linear algebra over the semigroup $Z_p$. Clearly $Z_p$ has no proper subsemigroup as p is a prime so even if W is a proper subset of V with each $W_i \subset V_i$ a subsemigroup of $V_i$ for each i; $1 \leq i \leq n$ still V has no proper subsemigroup n-set linear subalgebra as $Z_p$ has no proper subsemigroups. Thus V has no subsemigroup n-set linear subalgebra.

Hence we can define such V's described in this theorem as pseudo simple semigroup n-set linear algebras.

**THEOREM 5.4:** *Let $V = V_1 \cup \ldots \cup V_n$ be a semigroup n-set vector space over the semigroup $S= Z^+$ or $Z_n$, n a composite number for appropriate V's we can have in V subsemigroup n-set vector subspaces.*

*Proof:* When $V = V_1 \cup \ldots \cup V_n$ has a proper subset $W = W_1 \cup \ldots \cup W_n$ such that for some proper subsemigroups S of $Z^+$ or $Z_n$ (n-composite number) W is a semigroup n-set vector space of V over S then W is the subsemigroup n-set vector subspace of V. Hence the claim.

Now we give some examples of them.

***Example 5.33:*** Let

$$V = V_1 \cup V_2 \cup V_3 \cup V_4$$
$$= \{Z_{10} \times Z_{10}\} \cup \{Z_{10}[x]\} \cup \left\{\begin{pmatrix} a & b \\ c & d \end{pmatrix} \middle| a,b,c,d \in Z_{10}\right\} \cup$$
$$\left\{\begin{pmatrix} a & a & a \\ a & a & a \end{pmatrix}, \begin{pmatrix} a & a & a & a & a \\ a & a & a & a & a \end{pmatrix} \middle| a \in Z_{10}\right\}$$

be a semigroup 4-set vector space over the semigroup S = $Z_{10}$. Take

W = $\{(a, a) \mid a \in Z_{10}\} \cup$ {All polynomials of even degree in x with coefficients from $Z_{10}\} \cup \left\{\begin{pmatrix} a & a \\ a & a \end{pmatrix} \middle| a \in Z_{10}\right\} \cup$

$$\left\{\begin{pmatrix} a & a & a \\ a & a & a \end{pmatrix} \middle| a \in Z_{10}\right\}$$

$= W_1 \cup W_2 \cup W_3 \cup W_4 \subseteq V$;

W is a subsemigroup 4 vector subspace over the subsemigroup S = $\{0, 5\} \subseteq Z_{10}$. In fact W is also a subsemigroup 4-set vector subspace over the subsemigroup P = $\{0, 2, 4, 6, 8\} \subseteq Z_{10}$. Thus we see W is a subsemigroup 4 set vector subspace over all the subsemigroups of $Z_{10}$.

Now we proceed on to define yet another type of substructure of these semigroup n-set vector spaces and semigroup n-set linear algebras.

**DEFINITION 5.17:** *Let $V = V_1 \cup V_2 \cup \ldots \cup V_n$ be a semigroup n-set vector space over the semigroup S. If $P_1, \ldots, P_n$ is the complete set of subsemigroups of S (n can also be infinite). Suppose $W = W_1 \cup W_2 \cup \ldots \cup W_n \subseteq V_1 \cup V_2 \cup \ldots \cup V_n$ is a proper subset of V and W is a subsemigroup n – set vector subspace of V for every subsemigroup $P_i$ of S for i = 1, 2, …, n then we call W to be the strong subsemigroup n set vector subspace of V.*

(All subsemigroup n-set vector subspaces of V need not be a strong subsemigroup n-set vector subspaces of V). Similarly if $V = V_1 \cup V_2 \cup \ldots \cup V_n$ is a semigroup n set linear algebra over the semigroup S and if $P_1, P_2, \ldots P_n$ is the set of all subsemigroups of S and if $W = W_1 \cup W_2 \cup \ldots \cup W_n \subseteq V_1 \cup V_2 \cup \ldots \cup V_n$ is such that W is the subsemigroup n-set linear subalgebra of V over every subsemigroup $P_i$, for i = 1, 2, …, n then we call W to be the strong subsemigroup n-set linear subalgebra of V. As in case of semigroup n set vector spaces we see in case of semigroup n-set linear algebras all subsemigroup n-set linear algebras need not always be strong subsemigroup n-set linear subalgebras of V.

We now illustrate this situation by a simple example.

***Example 5.34:*** Let

$V = V_1 \cup V_2 \cup \ldots \cup V_5$

$= \{Z_6 \times Z_6 \times Z_6\} \cup \left\{\begin{pmatrix} a & a \\ a & a \end{pmatrix} \middle| a \in Z_6\right\} \cup \{Z_6[x]$; all polynomials of degree less than or equal to 5 with coefficients from $Z_6\} \cup \left\{\begin{pmatrix} a & b & c \\ d & e & f \end{pmatrix} \middle| a, b, c, d, e, f \in Z_6\right\}$
$\cup \left\{\begin{pmatrix} a & a & a \\ a & a & a \\ a & a & a \end{pmatrix} \middle| a \in Z_6\right\}$

be a semigroup 5-set linear algebra over the semigroup Z. The subsemigroups of $Z_6$ are $P_1 = \{0, 3\}$ and $P_2 = \{0, 2, 4\}$.
Take

$W = \{(a, a, a) / a \in Z_6\} \cup \left\{\begin{pmatrix} a & a \\ a & a \end{pmatrix} \middle| a \in \{0, 2, 4\}\right\} \cup \{$All polynomial of the form $a + ax + ax^2 + ax^3 + ax^4 + ax^5 \mid a \in Z_6\} \cup \left\{\begin{pmatrix} a & a & a \\ a & a & a \end{pmatrix} \middle| a \in Z_6\right\} \cup \left\{\begin{pmatrix} a & a & a \\ a & a & a \\ a & a & a \end{pmatrix} \middle| a \in \{0, 3\}\right\}$

$\subseteq V = V_1 \cup V_2 \cup V_3 \cup V_4 \cup V_5$

is a strong subsemigroup 5-set linear subalgebra as W is a subsemigroup 5 set linear subalgebra over both the subsemigroups $P_1$ and $P_2$.

Now we proceed on to define the notion of group n-set vector spaces and group n-set linear algebras.

**DEFINITION 5.18:** *Let $V = V_1 \cup V_2 \cup \ldots \cup V_n$, V is said to be a group n-set vector space over the group G where $V_i$ are sets such that $g\, v_i \in V_i$ for all $v_i \in V_i$ and $g \in G$, $1 \le i \le n$. Here G is just an additive abelian group.*

We now illustrate this definition by some examples.

***Example 5.35:*** Let

$$V = V_1 \cup V_2 \cup V_3 \cup V_4$$

$$= \left\{ \begin{pmatrix} a & a \\ a & a \end{pmatrix} \middle| a \in Z_7 \right\} \cup \{Z_7[x]\} \cup \{Z_7 \times Z_7 \times Z_7\} \cup \left\{ \begin{pmatrix} a & a & a \\ a & a & a \end{pmatrix}, \begin{pmatrix} a & a & a \\ a & a & a \\ a & a & a \end{pmatrix}, \begin{pmatrix} a & a & a & a \\ a & a & a & a \\ a & a & a & a \end{pmatrix} \middle| a \in Z_7 \right\},$$

V is a group 4-set vector space over the group $Z_7$ under addition modulo 7.

***Example 5.36:*** Let

$$V = V_1 \cup V_2 \cup V_3 \cup V_4 \cup V_5$$

$$= \{Z \times Z \times Z\} \cup \{Z\,[x]\} \cup \left\{ \begin{pmatrix} a & b \\ c & d \end{pmatrix} \middle| a, b, c, d \in Z \right\} \cup \left\{ \begin{pmatrix} a & a & a \\ a & a & a \end{pmatrix} \middle| a \in Z \right\} \cup \left\{ \begin{pmatrix} a & a & a & a \\ 0 & a & a & a \\ 0 & 0 & a & a \\ 0 & 0 & 0 & a \end{pmatrix} \middle| a \in Z \right\}$$

be a group n-set vector space (n = 5) over the group Z under addition.

***Example 5.37:*** Let

$$\begin{aligned} V &= V_1 \cup V_2 \cup V_3 \cup V_4 \\ &= \{(0000), (1100), (0011)\} \cup \{(000), (010), (00), (10)\} \cup \\ &\quad \left\{ \begin{pmatrix} a & a \\ a & a \end{pmatrix} \middle| a \in Z_2 \right\} \cup \{Z_2 [x]\} \end{aligned}$$

be a group n set vector space over the group $Z_2 = \{0, 1\}$, a group under addition modulo 2.

Now we proceed on to define the notion of group n-set linear algebra over a group and illustrate them by some examples.

**DEFINITION 5.18:** *Let $V = V_1 \cup V_2 \cup \ldots \cup V_n$ be a group n-set vector space over the additive group G. If each $V_i$ is an additive group then we call V to be a group n-set linear algebra over G for $1 \le i \le n$.*

***Example 5.38:*** Let

$$\begin{aligned} V &= V_1 \cup V_2 \cup V_3 \cup V_4 \cup V_5 \\ &= \{Z \times Z \times Z\} \cup \{Z [x]\} \cup \left\{ \begin{pmatrix} a & b \\ c & d \end{pmatrix} \middle| a, b, c, d \in Z \right\} \cup \\ &\quad \left\{ \begin{pmatrix} a & a & a \\ a & a & a \end{pmatrix} \middle| a \in 2Z \right\} \cup \{(a\ a\ a\ a\ a) / a \in Z\} \end{aligned}$$

be a group 5-set linear algebra over the additive group Z.

***Example 5.39:*** Let

$$\begin{aligned} V &= V_1 \cup V_2 \cup V_3 \cup V_4 \\ &= \{(000), (111), (100), (001), (010), (110), (011), (101)\} \\ &\quad \cup \{Z_2 [x]\} \cup \{(1111), (0000), (1100), (0011)\} \\ &\quad \cup \left\{ \begin{pmatrix} a & b \\ c & d \end{pmatrix} \middle| a, b, c, d \in Z_2 \right\} \end{aligned}$$

be a group 4 set linear algebra over the additive group $Z_2$ modulo 2.

***Example 5.40:*** Let

$$
\begin{aligned}
V &= V_1 \cup V_2 \cup V_3 \cup V_4 \\
&= \left\{\begin{pmatrix} a & b \\ c & d \end{pmatrix} \middle| a,b,c,d \in Z_{10}\right\} \cup \left\{\begin{pmatrix} a & a & a & a \\ a & a & a & a \end{pmatrix} \middle| a \in Z_{10}\right\} \cup \\
&\quad Z_{10}[x] \cup \{Z_{10} \times Z_{10} \times Z_{10} \times Z_{10}\};
\end{aligned}
$$

V is a group 4 set linear algebra over the group $Z_{10}$, group under addition modulo 10.

Now we define a few substructures of these two concepts.

**DEFINITION 5.19:** *Let $V = V_1 \cup ... \cup V_n$ be a group n set vector space over the group G. If $W = W_1 \cup W_2 \cup \ldots \cup W_n$ is a proper subset of V and W is a group n set vector space over the group G then we call W to be a group n-set vector subspace of V over G.*

We illustrate this by some examples.

***Example 5.41:*** Let

$$
\begin{aligned}
V &= V_1 \cup V_2 \cup V_3 \cup V_4 \\
&= \{Z_6 \times Z_6\} \cup \left\{\begin{pmatrix} a & b \\ c & d \end{pmatrix} \middle| a,b,c,d \in Z_6\right\} \\
&\quad \cup \left\{\begin{pmatrix} a & a & a & a \\ a & a & a & a \end{pmatrix}, \begin{pmatrix} a & a & a \\ a & a & a \end{pmatrix} \middle| a \in Z_6\right\} \cup Z_6[x]
\end{aligned}
$$

all polynomials of degree less than or equal to 4}, V is a group 4-set vector space over the group $Z_6$, group under addition modulo 6. Take

$W = W_1 \cup W_2 \cup W_3 \cup W_4$

$$= \{(a, a) / a \in Z_6\} \cup \left\{\begin{pmatrix} a & a \\ a & a \end{pmatrix} \middle| a \in Z_6\right\} \cup$$

$$\left\{\begin{pmatrix} a & a & a & a \\ a & a & a & a \end{pmatrix} \middle| a \in Z_6\right\} \cup \{\text{all polynomials of the form}$$

$$a + ax + ax^2 + ax^3 + ax^4 / a \in Z_6\}$$

$$\subseteq V_1 \cup V_2 \cup V_3 \cup V_4$$

$$= V,$$

W is a group 4-set vector subspace of V over the group $Z_6$.

***Example 5.42:*** Let

$$V = V_1 \cup V_2 \cup V_3 \cup V_4 \cup V_5$$

$$= \{Z \times Z\} \cup \{Z[x]\} \cup \left\{\begin{pmatrix} a & a & a \\ a & a & a \end{pmatrix} \middle| a \in Z\right\} \cup$$

$$\left\{\begin{pmatrix} a & b \\ c & d \end{pmatrix} \middle| a, b, c, d \in Z\right\} \cup$$

$$\left\{\begin{pmatrix} a & b & c & d \\ 0 & e & f & g \\ 0 & 0 & h & i \\ 0 & 0 & 0 & j \end{pmatrix} \middle| a, b, c, d, e, f, g, h, i, j \in Z\right\}$$

be a group 5-set vector space over the additive group Z. Take

$$W = W_1 \cup W_2 \cup W_3 \cup W_4 \cup W_5$$

$$= \{(a, a) \mid a \in Z\} \cup \{a + ax + \ldots + ax^n \mid a \in Z\}$$

$$\cup \left\{\begin{pmatrix} a & a & a \\ a & a & a \end{pmatrix} \middle| a \in 2Z\right\} \cup \left\{\begin{pmatrix} a & a \\ a & a \end{pmatrix} \middle| a \in Z\right\}$$

$$\cup \left\{ \left( \begin{matrix} a & a & a & a \\ 0 & a & a & a \\ 0 & 0 & a & a \\ 0 & 0 & 0 & a \end{matrix} \right) \middle| a \in 2Z \right\}$$

$$\subseteq V_1 \cup V_2 \cup V_3 \cup V_4 \cup V_5$$
$$= V.$$

W is clearly a group 5 set vector subspace of V over the group G.

We now prove the following interesting result.

**THEOREM 5.6:** *Every group n-set linear algebra $V = V_1 \cup ... \cup V_n$ over the group G is a group n-set vector space over the group G but a group n set vector space over the group G in general need not be a group n-set linear algebra over G.*

*Proof:* Let $V = V_1 \cup V_2 \cup \ldots \cup V_n$ be a group n set linear algebra over the additive group G. Clearly V is a group n-set vector space over the additive group G for every n-group is a n-set.

To prove a group n set vector space V over a group G is not in general a group n-set linear algebra over the group G. We give a counter example. Consider the group 5 set vector space

$$V = V_1 \cup V_2 \cup V_3 \cup V_4 \cup V_5$$
$$= \{(111), (000), (100), (00), (11), (01), (10)\} \cup$$
$$\left\{ \left( \begin{matrix} a & b \\ c & d \end{matrix} \right), \left( \begin{matrix} a & a & a \\ a & a & a \end{matrix} \right) \middle| a,b,c,d \in Z_2 = \{0,1\} \right\} \cup$$
$$\left\{ \left( \begin{matrix} a & a & a \\ 0 & a & a \\ 0 & 0 & a \end{matrix} \right), \left( \begin{matrix} a & a \\ a & a \\ a & a \\ a & a \end{matrix} \right) \middle| a \in Z_2 = \{0,1\} \right\} \cup \{Z_2 [x]\} \cup$$
$$\{Z_2 \times Z_2 \times Z_2 \times Z_2\}$$

over the additive group $Z_2 = \{0, 1\}$. Clearly $V_1$, $V_2$ and $V_3$ are not groups under addition so V is not a group 5-set linear algebra over the group $Z_2 = \{0, 1\}$ under addition modulo 2. Hence the claim.

Now we define group n set linear subalgebra over the group G.

**DEFINITION 5.21:** *Let $V = V_1 \cup V_2 \cup \ldots \cup V_n$ be the group n-set linear algebra over the additive group G, where each $V_i$ is a group under addition, i = 1, 2, …, n. Take $W = W_1 \cup W_2 \cup \ldots \cup W_n$, a subset of V such that W the group n set linear algebra over the same group G i.e., each $W_i \subseteq V_i$ is a proper subgroup of $V_i$, i = 1, 2, …, n. We call $W = W_1 \cup W_2 \cup \ldots \cup W_n$ to be a group n-set linear subalgebra of V over the additive group G.*

Now in case of group n set linear algebra V; we can define the notion of pseudo group n-set vector subspaces of V which is given below.

**DEFINITION 5.22:** *Let $V = V_1 \cup V_2 \cup \ldots \cup V_n$ be a group n set linear algebra over the group G. If $W = W_1 \cup W_2 \cup \ldots \cup W_n$ be a proper n-subset, where at least one of the $W_i$'s is not a subgroup of the group $V_i$ and if W is a group n-set vector space over G then we call $W = W_1 \cup W_2 \cup \ldots \cup W_n$ to be the pseudo group n-set vector subspace of V;$(1 \leq i \leq n)$.*

Now we illustrate these definitions by the following examples.

***Example 5.43:*** Let

$$V = V_1 \cup V_2 \cup V_3 \cup V_4$$

$$= \left\{ \begin{pmatrix} a & a \\ a & a \end{pmatrix} \middle| a \in Z_{12} \right\} \cup \{Z_{12}[x]\} \cup \{Z_{12} \times Z_{12} \times Z_{12}\} \cup$$

$$\left\{ \begin{pmatrix} a_1 & a_2 & a_3 \\ a_4 & a_5 & a_6 \\ a_7 & a_8 & a_9 \end{pmatrix} \middle| a_i \in Z_{12}; \quad 1 \leq i \leq 9 \right\}$$

be a group 4 set linear algebra over the group $Z_{12}$, addition modulo 12.

Let

$$W = \left\{\begin{pmatrix}6 & 6\\ 6 & 6\end{pmatrix}, \begin{pmatrix}0 & 0\\ 0 & 0\end{pmatrix}\right\} \cup \{(ax^2 + ax + a, ax^7 + ax^3 + a, a + ax^7 + ax^5 + ax^2 \mid a \in Z_{12}\} \cup \{(a\ a\ 0), (0\ a\ a), (a\ 0\ a) \mid a \in Z_{12}\} \cup \left\{\begin{pmatrix}a & a & a\\ 0 & 0 & 0\\ 0 & 0 & 0\end{pmatrix}, \begin{pmatrix}0 & 0 & 0\\ a & a & a\\ 0 & 0 & 0\end{pmatrix}, \begin{pmatrix}0 & 0 & 0\\ 0 & 0 & 0\\ a & a & a\end{pmatrix} \middle| a \in Z_{12}\right\}$$

$$= W_1 \cup W_2 \cup W_3 \cup W_4 \subseteq V_1 \cup V_2 \cup V_3 \cup V_4 = V.$$

Clearly $W_2$, $W_3$ and $W_4$ are not even closed under addition. So W is verified easily to be a pseudo group 4 set vector subspace of V over $Z_{12}$.

We further say all proper subsets of V need to be pseudo group 4 set vector subspaces of V over the group $Z_{12}$. For take

$$W = \left\{\begin{pmatrix}1 & 1\\ 1 & 1\end{pmatrix}, \begin{pmatrix}5 & 5\\ 5 & 5\end{pmatrix}\right\} \cup \{3x + 2x^2 + 1, 5x^3 + 1\} \cup \{(321), (123)\} \cup \left\{\begin{pmatrix}3 & 7 & 2\\ 1 & 2 & 0\\ 5 & 1 & 4\end{pmatrix}\right\}$$

$$= W_1 \cup W_2 \cup W_3 \cup W_4 \subseteq V_1 \cup V_2 \cup V_3 \cup V_4.$$

Clearly for $W_1$ we take

$$\begin{pmatrix}1 & 1\\ 1 & 1\end{pmatrix} \in W_1$$

and

$$6 \in Z_{12}\ ,\ 6\begin{pmatrix}1 & 1\\ 1 & 1\end{pmatrix} = \begin{pmatrix}6 & 6\\ 6 & 6\end{pmatrix} \notin W_1.$$

Likewise $3x + 2x^2 + 1 \in W_2$ and $7 \in Z_{12}$; $7\,(3x + 2x^2 + 1) = 9x + 2x^2 + 7 \notin W_2$, (3, 2, 1) is in $W_3$ and $5 \in Z_{12}$, $5(3, 2, 1) = (3, 10, 5) \notin W_3$.
Finally

$$\begin{pmatrix}3 & 7 & 2\\ 1 & 2 & 0\\ 5 & 1 & 4\end{pmatrix} \in W_4$$

and $0 \in Z_{12}$.

$$0.\begin{pmatrix}3 & 7 & 2\\ 1 & 2 & 0\\ 5 & 1 & 4\end{pmatrix} = \begin{pmatrix}0 & 0 & 0\\ 0 & 0 & 0\\ 0 & 0 & 0\end{pmatrix} \notin W_4\,.$$

So $W = W_1 \cup W_2 \cup W_3 \cup W_4 \subseteq V$ is just a 4-subset of V, but is not a pseudo group 4 set vector subspace of V over the group $Z_{12}$. Thus the 4-subsets of V are not pseudo group 4 set vector subspace of V.

Now we proceed on to define the notion of subgroup n-set linear subalgebra of a group n set linear subalgebra of a group n-set linear algebra over the group G.

**DEFINITION 5.23:** *Let $V = V_1 \cup V_2 \cup \ldots \cup V_n$ be a group n set linear algebra over the group G. Let $W = W_1 \cup \ldots \cup W_n$ be a proper subset of G ie. $W_i \subseteq V_i$ and $W_i$ is a subgroup of $V_i$, $1 \le i \le n$. Let H be a proper subgroup of G. If W is a group n-set linear algebra over the group H, then we call W to be the subgroup n-set linear subalgebra of V over the subgroup H of the group G.*

***Example 5.44:*** Let

$$\begin{aligned} V &= V_1 \cup V_2 \cup V_3 \cup V_4 \\ &= \{Z_6[x]\} \cup \{Z_6 \times Z_6 \times Z_6\} \cup \left\{\begin{pmatrix}a & b\\ c & d\end{pmatrix} \middle|\, a, b, c, d \in Z_6\right\} \end{aligned}$$

$$\cup \left\{ \begin{pmatrix} a & a & a & a \\ a & a & a & a \end{pmatrix} \middle| a \in Z_6 \right\}$$

be the group 4-set linear algebra over the group $Z_6$ under addition modulo 6. Let

$$W = W_1 \cup W_2 \cup W_3 \cup W_4$$
$$= \{a + ax + ax^2 + \ldots ax^n \,/\, a \in Z_6\} \cup \{(a, a, a) \,/\, a \in Z_6\}$$
$$\cup \left\{ \begin{pmatrix} a & a \\ a & a \end{pmatrix} \middle| a \in Z_6 \right\} \cup \left\{ \begin{pmatrix} a & a & a & a \\ a & a & a & a \end{pmatrix} \middle| a \in \{0,3\} \subset Z_6 \right\}$$

be the group 4 set linear subalgebra of V over the group $Z_6$.

***Example 5.45:*** Let

$$V = V_1 \cup V_2 \cup V_3 \cup V_4$$

$$= \{Z \times Z \times Z\} \cup \{Z[x]\} \cup$$

$$\left\{ \begin{pmatrix} a & a & a \\ 0 & a & a \\ 0 & 0 & a \end{pmatrix} \middle| a \in Z \right\} \cup$$
$$\left\{ \begin{pmatrix} a_1 & 0 & 0 & 0 \\ a_2 & a_3 & 0 & 0 \\ a_4 & a_5 & a_6 & 0 \\ a_7 & a_8 & a_9 & a_{10} \end{pmatrix} \middle| a_i \in Z; \quad 1 \le i \le 10 \right\}$$

be a group 4 set linear algebra over the group G = Z. Take

$$W = W_1 \cup W_2 \cup W_3 \cup W_4$$
$= \{(a\ a\ a) \mid a \in Z\} \cup \{$all polynomials of even degree with

coefficients from Z} $\cup$ $\left\{\begin{pmatrix} 2n & 2n & 2n \\ 0 & 2n & 2n \\ 0 & 0 & 2n \end{pmatrix} \middle| n \in Z\right\}$

$$\cup \left\{\begin{pmatrix} a & 0 & 0 & 0 \\ a & a & 0 & 0 \\ a & a & a & 0 \\ a & a & a & a \end{pmatrix} \middle| a \in Z\right\}$$

$\subseteq$ $V_1 \cup V_2 \cup V_3 \cup V_4 = V$.

It is easily verified W is a group 4 set linear subalgebra of V over the group G = Z.

***Example 5.46:*** Let

$$V = V_1 \cup V_2 \cup V_3 \cup V_4$$
$$= \{Z_{14} \times Z_{14} \times Z_{14} \times Z_{14}\} \cup \{Z_{14}[x]\} \cup \left\{\begin{pmatrix} a & a \\ a & a \\ a & a \end{pmatrix} \middle| a \in Z_{14}\right\}$$
$$\cup \left\{\begin{pmatrix} a & b \\ c & d \end{pmatrix} \middle| a, b, c, d \in Z_{14}\right\}$$

be a group 4 set linear algebra over the group $Z_{14}$. Take

$$W = W_1 \cup W_2 \cup W_3 \cup W_4$$
$$= \{(a, a, a, a) \mid a \cup Z_{14}\} \cup \{a + ax + \ldots + ax^n \mid a \in Z_{14}\}$$
$$\cup \left\{\begin{pmatrix} a & a \\ a & a \\ a & a \end{pmatrix} \middle| a \in \{0, 7\}\right\} \cup \left\{\begin{pmatrix} a & a \\ a & a \end{pmatrix} \middle| a \in Z_{14}\right\}$$
$$\subseteq V_1 \cup V_2 \cup V_3 \cup V_4 = V.$$

W is clearly a subgroup 4-set linear subalgebra over the subgroup $\{0,7\} \subseteq Z_{14}$.

***Example 5.47:*** Let

$$
\begin{aligned}
V &= V_1 \cup V_2 \cup V_3 \cup V_4 \cup V_5 \\
&= \{Z \times Z \times Z \times Z\} \cup \left\{\begin{pmatrix} a & a & a & a \\ a & a & a & a \end{pmatrix} \middle| a \in Z\right\} \cup \\
&\quad \left\{\begin{pmatrix} a & b & c & d \\ 0 & e & f & g \\ 0 & 0 & h & i \\ 0 & 0 & 0 & j \end{pmatrix} \middle| a, b, c, d, e, f, g, h, i, j \in Z\right\} \cup \{Z[x]\} \cup \\
&\quad \left\{\begin{pmatrix} a_1 & a_2 & a_3 \\ a_4 & a_5 & a_6 \\ a_7 & a_8 & a_9 \end{pmatrix} \middle| a_i \in Z;\ \ 1 \le i \le 9\right\}
\end{aligned}
$$

be a group 5-set linear algebra over the additive group Z. Take

$$
\begin{aligned}
W &= W_1 \cup W_2 \cup W_3 \cup W_4 \cup W_5 \\
&= \{(2Z \times 3Z \times 4Z \times 5Z)\} \cup \left\{\begin{pmatrix} a & a & a & a \\ a & a & a & a \end{pmatrix} \middle| a \in 2Z\right\} \cup \\
&\quad \left\{\begin{pmatrix} a & a & a & a \\ 0 & a & a & a \\ 0 & 0 & a & a \\ 0 & 0 & 0 & a \end{pmatrix} \middle| a \in Z\right\} \cup \{(a + ax + ax^2 + \ldots + ax^n \mid \\
&\quad a \in Z\} \cup \left\{\begin{pmatrix} a & a & a \\ a & a & a \\ a & a & a \end{pmatrix} \middle| a \in 2Z\right\} \\
&\subseteq V_1 \cup V_2 \cup V_3 \cup V_4 \cup V_5 \\
&= V.
\end{aligned}
$$

W is a subgroup 5-set linear subalgebra of V over the subgroup $2Z \subseteq Z$.

Let us define the n-generating subset of a group n-set vector space and group n-set linear algebra over the group G.

**DEFINITION 5.23:** *Let $V = V_1 \cup V_2 \cup \ldots \cup V_n$ be a group n-set vector space over the group G. Suppose $W = W_1 \cup W_2 \cup \ldots \cup W_n$ is a proper subset of V and if W is a semigroup n-set vector space over some proper subset H of G where H is a semigroup, then we call W to be the pseudo semigroup n-set vector subspace of V over the semigroup H of G.*

We illustrate it by the following examples.

***Example 5.48:*** Let

$$V = V_1 \cup V_2 \cup V_3 \cup V_4 \cup V_5$$

$$= \{Z \times Z \times Z\} \cup \{Z[x]\} \cup \left\{ \left. \begin{pmatrix} a & a & a \\ a & a & a \\ a & a & a \end{pmatrix} \right| a \in Z \right\} \cup$$

$$\left\{ \left. \begin{pmatrix} a & b \\ c & d \end{pmatrix} \right| a, b, c, d \in Z \right\} \cup$$

$$\left\{ \left. \begin{pmatrix} a_1 & 0 & 0 & 0 \\ a_2 & a_3 & 0 & 0 \\ a_4 & a_5 & a_6 & 0 \\ a_7 & a_8 & a_9 & a_{10} \end{pmatrix} \right| a_i \in Z; 1 \le i \le 10 \right\}$$

be a group 5-set vector space over the additive group Z.

$$W = W_1 \cup W_2 \cup W_3 \cup W_4 \cup W_5$$

$$= \{Z^+ \times Z^+ \times Z^+\} \cup \{Z^+[x]\} \cup \left\{ \left. \begin{pmatrix} a & a & a \\ a & a & a \\ a & a & a \end{pmatrix} \right| a \in Z^+ \right\} \cup$$

$$\left\{\begin{pmatrix} a & a \\ a & a \end{pmatrix} \middle| a \in Z^+ \right\} \cup \left\{\begin{pmatrix} a_1 & 0 & 0 & 0 \\ a_2 & a_3 & 0 & 0 \\ a_4 & a_5 & a_6 & 0 \\ a_7 & a_8 & a_9 & a_{10} \end{pmatrix} \middle| a_i \in Z^+; 1 \leq i \leq 10 \right\}$$

$$\subseteq \; V_1 \cup V_2 \cup V_3 \cup V_4 \cup V_5$$
$$= \; V$$

is a pseudo semigroup 5-set vector subspace of V over the semigroup $Z^+$ of Z.

***Example 5.49:*** Let

$$V = V_1 \cup V_2 \cup V_3 \cup V_4$$

$$= \left\{\begin{pmatrix} a & a & a \\ a & a & a \\ a & a & a \end{pmatrix} \middle| a \in 2Z \right\} \cup \left\{\begin{pmatrix} a & a & a & a \\ a & a & a & a \end{pmatrix} \middle| a \in 3Z \right\} \cup$$

$$\{Z \times Z \times Z \times 3Z\} \cup \left\{\begin{pmatrix} a & a \\ a & a \\ a & a \end{pmatrix} \middle| a \in 5Z \right\}$$

be a group n-set vector space over the group Z. Take

$$W = W_1 \cup W_2 \cup W_3 \cup W_4$$

$$= \left\{\begin{pmatrix} a & a & a \\ a & a & a \\ a & a & a \end{pmatrix} \middle| a \in 2Z^+ \right\} \cup \left\{\begin{pmatrix} a & a & a & a \\ a & a & a & a \end{pmatrix} \middle| a \in 3Z^+ \right\} \cup$$

$$\{Z^+ \times Z^+ \times Z^+ \times 3Z^+\} \cup \left\{\begin{pmatrix} a & a \\ a & a \\ a & a \end{pmatrix} \middle/ a \in 5Z^+ \right\}$$

$$\subseteq \; V_1 \cup V_2 \cup V_3 \cup V_4$$
$$= \; V$$

is easily verified to be a pseudo semigroup n-set vector subspace over the semigroup $Z^+ \subseteq Z$.

Now we proceed to define n-generating set and n-basis.

**DEFINITION 5.24:** *Let $V = V_1 \cup V_2 \cup \ldots \cup V_n$ be a group n-set vector space over the group G. If we have a n-set $X = X_1 \cup X_2 \cup \ldots \cup X_n \subseteq V_1 \cup V_2 \cup \ldots \cup V_n = V$ such that each $X_i$ generates $V_i$ over the group G, for i = 1, 2, …, n, then we call X to be the n-set generator of V. If the cardinality of each $X_i$ is finite we say V is generated finitely and the n-dimension of V is $(|X_1|, |X_2|, \ldots, |X_n|)$. If even one of the $X_i$ of X happens to have infinite cardinality then we say n-dimension of V is infinite. We call the n-generating n-subset of V to be the n-basis of V over the group G.*

We illustrate this by the following examples.

***Example 5.50:*** Let $V = V_1 \cup V_2 \cup V_3 \cup V_4$ be a group 4-set vector space over the group $Z_2 = \{0,1\}$ where $V_1 = \{(1100), (0011), (1111), (0000), (111), (101), (000)\}$, $V_2 = \{Z_2 [x]$ / every polynomial is of degree less than or equal to 2},

$$V_3 = \left\{ \begin{pmatrix} a & b \\ c & d \end{pmatrix} \middle| a, b, c, d \in Z_2 \right\}$$

and

$$V_4 = \left\{ \begin{pmatrix} a & a & a & a \\ a & a & a & a \end{pmatrix} \middle| a \in Z_2 \right\}.$$

Take

$$X = \{(1111), (1100), (0011), (111), (101)\} \cup \{1, x, x^2\ 1 + x, 1+x^2, x+x^2, 1+x+x^2\} \cup$$
$$\left\{ \begin{pmatrix} 0 & 1 \\ 0 & 0 \end{pmatrix}, \begin{pmatrix} 1 & 0 \\ 0 & 0 \end{pmatrix}, \begin{pmatrix} 0 & 0 \\ 1 & 0 \end{pmatrix}, \begin{pmatrix} 0 & 0 \\ 0 & 1 \end{pmatrix}, \begin{pmatrix} 1 & 1 \\ 0 & 0 \end{pmatrix}, \right.$$

$$\begin{pmatrix}1&0\\1&0\end{pmatrix},\begin{pmatrix}0&1\\0&1\end{pmatrix},\begin{pmatrix}0&0\\1&1\end{pmatrix},\begin{pmatrix}1&0\\0&1\end{pmatrix},\begin{pmatrix}0&1\\1&0\end{pmatrix}$$

$$\left.\begin{pmatrix}1&0\\1&1\end{pmatrix},\begin{pmatrix}1&1\\0&1\end{pmatrix},\begin{pmatrix}0&1\\1&1\end{pmatrix},\begin{pmatrix}1&1\\1&0\end{pmatrix},\begin{pmatrix}1&1\\1&1\end{pmatrix}\right\} \cup$$

$$\left\{\begin{pmatrix}1&1&1&1\\1&1&1&1\end{pmatrix}\right\}$$

$$= X_1 \cup X_2 \cup X_3 \cup X_4$$
$$\subseteq V_1 \cup V_2 \cup V_3 \cup V_4.$$

X n-generates V (n=4). Further V is finitely 4-generated and the number of 4-generators is given by (5,7,15,1).

We give yet another example.

***Example 5.51:*** Let

$$\begin{aligned} V &= V_1 \cup V_2 \cup V_3 \cup V_4 \cup V_5 \\ &= \{Z \times Z\} \cup Z[x] \cup \{(a, a, a) \mid a \in Z\} \cup \left\{\begin{pmatrix}a&a\\a&a\end{pmatrix} \middle| a \in Z\right\} \\ &\quad \cup \{(a, a, a, a, a), (a, a, a, a, a, a, a) \mid a \in Z\} \end{aligned}$$

be the group 5 set vector space over the group Z under addition. Now if we consider the 5-generating subset of V as $X = X_1 \cup X_2 \cup X_3 \cup X_4 \cup X_5 \subseteq V_1 \cup V_2 \cup V_3 \cup V_4 \cup V_5$ where $X_1$, $X_2$ are infinite sets which alone can generate $V_1$ and $V_2$ respectively, but $X_3 = \{(111)\}$,

$$X_4 = \left\{\begin{pmatrix}1&1\\1&1\end{pmatrix}\right\}$$

and $X_5 = \{(11111), (1111111)\}$. Thus the 5-dimension of V is infinite given by $(\infty,\infty, 1, 1, 2)$.

Now we find the n-basis and n-dimension of group n-set linear algebras over the group. We just define the n-basis in case of group n-set linear algebras.

**DEFINITION 5.25:** *Let $V = V_1 \cup V_2 \cup \ldots \cup V_n$ be a group n-set linear algebra over the additive group G. Let $X = X_1 \cup X_2 \cup \ldots \cup X_n \subseteq V_1 \cup V_2 \cup \ldots \cup V_n = V$ be the n-generator of V. If each $X_i$ generates $V_i$ as a group linear algebra over the group G then we call X to be the n-generating set of V; i = 1, 2, …, n. Here each $V_i$ is a group under addition; $1 \le i \le n$. If each $X_i$ of X is finite then we say V is n-finitely generated by X over G. If even one of the $X_i$'s is infinite we say V is infinitely generated over G.*

The n-cardinality of X is called as the n-dimension of V and X is called as the n-set basis of V over the group G.

Now we illustrate this situation by some examples.

***Example 5.52:*** Let

$$V = V_1 \cup V_2 \cup V_3 \cup V_4 = \left\{ \begin{pmatrix} a & a \\ a & a \end{pmatrix} \middle| a \in Z_{12} \right\} \cup \{Z_{12} \times Z_{12}\} \cup \{Z_{12}[x]; \text{ all polynomials of degree less than or equal to } 10\} \cup \left\{ \begin{pmatrix} a & b & c \\ 0 & d & e \end{pmatrix} \middle| a, b, c, d, e \in Z_{12} \right\}$$

be a group 4-set linear algebra over the group $G = Z_{12}$, group under addition modulo 12. Take

$$X = X_1 \cup X_2 \cup X_3 \cup X_4 = \left\{ \begin{pmatrix} 1 & 1 \\ 1 & 1 \end{pmatrix} \right\} \cup \{(1, 0), (0, 1)\} \cup \{1, x, x^2, \ldots, x^{10}\}$$
$$\cup \left\{ \begin{pmatrix} 1 & 0 & 0 \\ 0 & 0 & 0 \end{pmatrix}, \begin{pmatrix} 0 & 1 & 0 \\ 0 & 0 & 0 \end{pmatrix}, \begin{pmatrix} 0 & 0 & 1 \\ 0 & 0 & 0 \end{pmatrix}, \begin{pmatrix} 0 & 0 & 0 \\ 0 & 1 & 0 \end{pmatrix}, \begin{pmatrix} 0 & 0 & 0 \\ 0 & 0 & 1 \end{pmatrix} \right\}$$
$$\subseteq V_1 \cup V_2 \cup V_3 \cup V_4 .$$

X, 4-generates V over $Z_{12}$. Now the 4-dimension of V is (1, 2, 11, 5). Thus X is a 4-basis of V and V is 4-finite dimensional.

***Example 5.53:*** Let

$$V = V_1 \cup V_2 \cup V_3 \cup V_4 \cup V_5$$
$$= \{Z \times Z \times Z\} \cup \{Z[x]; \text{all polynomials of degree less than or equal to } 4\} \cup \left\{\begin{pmatrix} a & a \\ a & a \end{pmatrix} \middle| a \in Z\right\} \cup$$
$$\left\{\begin{pmatrix} a & b \\ 0 & d \\ 0 & c \end{pmatrix} \middle| a, b, c, d \in Z\right\} \cup \left\{\begin{pmatrix} 0 & 0 & 0 & 0 \\ a & 0 & 0 & 0 \\ a & a & 0 & 0 \\ a & a & a & 0 \end{pmatrix} \middle| a_i \in Z\right\}$$

be a group 5-set linear algebra over the additive group Z. Take

$$X = \{(100), (010), (001)\} \cup \{1, x, x^2\ x^3, x^4\} \cup \left\{\begin{pmatrix} 1 & 1 \\ 1 & 1 \end{pmatrix}\right\}$$
$$\left\{\begin{pmatrix} 1 & 0 \\ 0 & 0 \\ 0 & 0 \end{pmatrix}, \begin{pmatrix} 0 & 1 \\ 0 & 0 \\ 0 & 0 \end{pmatrix}, \begin{pmatrix} 0 & 0 \\ 0 & 1 \\ 0 & 0 \end{pmatrix}, \begin{pmatrix} 0 & 0 \\ 0 & 0 \\ 0 & 1 \end{pmatrix}\right\} \cup$$
$$\left\{\begin{pmatrix} 0 & 0 & 0 & 0 \\ 1 & 0 & 0 & 0 \\ 1 & 1 & 0 & 0 \\ 1 & 1 & 1 & 0 \end{pmatrix}\right\}$$

$$= X_1 \cup X_2 \cup X_3 \cup X_4 \cup X_5\,,$$

clearly X 5-generates V as a group 5-set linear algebra over the group Z. Now the 5-cardinality of V as a group 5-linear algebra over Z is finite and is given by (3, 5, 1, 4, 1).

It is important at this juncture to mention the essential fact that even if $V = V_1 \cup ... \cup V_n$ is a group n-linear algebra over the additive group G as well as if the same V is treated as a group n-vector space over the same group G, the dimension is distinctly different. One cannot say because V is the same the n-dimension is also the same over the group G. In view of this we give some concrete examples.

***Example 5.54:*** Let

$$
\begin{aligned}
V &= V_1 \cup V_2 \cup V_3 \cup V_4 \cup V_5 \\
&= \{Z \times Z \times Z\} \cup \{Z[x] \mid \text{all polynomials of degree less than or equal to 7}\} \cup \left\{\begin{pmatrix} a & b \\ c & d \end{pmatrix} \middle| a, b, c, d \in Z\right\} \cup \{(a, a, a, a, a) \mid a \in Z\} \cup \left\{\begin{pmatrix} a & 0 \\ 0 & b \\ d & 0 \\ 0 & c \end{pmatrix} \middle| a, b, c, d, \in Z\right\}
\end{aligned}
$$

be a group 5-set linear algebra over the group Z. Now take a 5-subset of V to be

$$
\begin{aligned}
X &= X_1 \cup X_2 \cup X_3 \cup X_4 \cup X_5 \\
&= \{(100), (010), (001)\} \cup \{1, x, x^2, \ldots, x^7\} \cup \left\{\begin{pmatrix} 1 & 0 \\ 0 & 0 \end{pmatrix}, \begin{pmatrix} 0 & 1 \\ 0 & 0 \end{pmatrix}, \begin{pmatrix} 0 & 0 \\ 0 & 1 \end{pmatrix}, \begin{pmatrix} 0 & 0 \\ 1 & 0 \end{pmatrix}\right\} \cup \{(11111)\} \cup \\
&\quad \left\{\begin{pmatrix} 1 & 0 \\ 0 & 0 \\ 0 & 0 \\ 0 & 0 \end{pmatrix}, \begin{pmatrix} 0 & 0 \\ 0 & 1 \\ 0 & 0 \\ 0 & 0 \end{pmatrix}, \begin{pmatrix} 0 & 0 \\ 0 & 0 \\ 1 & 0 \\ 0 & 0 \end{pmatrix}, \begin{pmatrix} 0 & 0 \\ 0 & 0 \\ 0 & 0 \\ 0 & 1 \end{pmatrix}\right\} \\
&\subseteq V_1 \cup V_2 \cup V_3 \cup V_4 \cup V_5 \\
&= V.
\end{aligned}
$$

Clearly X 5-generates V over Z and is a 5-basis of V. The 5-dimension of V over Z is {(3, 8, 4, 1, 4)}.

Now we can think of V as a group 5-set vector space over the group Z. Take $Y = = Y_1 \cup Y_2 \cup Y_3 \cup Y_4 \cup Y_5$ = {infinite generating subset of $V_1$ alone can generate the group vector space (as $V_1$ is just a set) over Z} ∪ {only an infinite subset can generate $V_2$} ∪ {only an infinite subset can generate $V_3$, as a group vector space over Z} ∪ {(11111)} ∪ {only an infinite subset can generate $V_5$ as a group vector space over the group G = Z}. Thus Y is an infinite 5-subset of V which 5-generates V as a group 5-set vector space over the group Z. So the 5-dimension of V is (∞,∞,∞, 1, ∞).

Thus we see from this example when we use the group n-set vector space over an additive group it can be infinite dimensional but at the same time if V is realized as a group n-set linear algebra it can be finite dimensional.

We provide yet another example

***Example 5.55:*** Let

$$V = V_1 \cup V_2 \cup V_3 \cup V_4 = \{Z_2 \times Z_2 \times Z_2\} \cup \{Z_2[x] \text{ all polynomials of degree less than or equal to } 4\} \cup \left\{ \begin{pmatrix} a & b \\ c & d \end{pmatrix} \middle| a,b,c,d \in Z_2 \right\} \cup \left\{ \begin{pmatrix} a & a & a & a \\ a & a & a & a \end{pmatrix} \middle| a \in Z_2 \right\}$$

be a group 4-vector space over the group $Z_2$ = {0,1}. V is also group liner algebra over the group $Z_2$ = {0,1}.

Now the 4-basis of V as a group 4-set vector space is given by the 4-set

X = $X_1 \cup X_2 \cup X_3 \cup X_4$
= {(100), (001), (010), (110), (011), (101), (111)} ∪ {1, x, $x^2$, $x^3$, $x^4$, 1+x, 1+$x^2$, x+$x^2$, 1+$x^3$, 1+$x^4$, x+$x^3$, x+$x^4$, $x^2$+$x^3$, $x^2$+$x^4$, $x^4$+$x^3$, 1+x+$x^2$, 1+x+$x^3$, 1+x+$x^4$, 1+$x^2$+$x^3$, 1+$x^2$+$x^4$, 1+$x^3$+$x^4$, $x^2$+$x^3$+$x^4$, x+$x^2$+$x^3$, x+$x^2$+$x^4$, x+$x^3$+$x^4$, 1+x+$x^2$+$x^3$, 1+x+$x^2$+$x^4$, 1+x+$x^3$+$x^4$, 1+$x^2$+$x^3$+$x^4$, x+$x^3$+$x^4$+$x^2$, 1+x + $x^2$ + $x^3$ + $x^4$} ∪

$$\left\{\begin{pmatrix}1&0\\0&0\end{pmatrix},\begin{pmatrix}0&1\\0&0\end{pmatrix},\begin{pmatrix}0&0\\1&0\end{pmatrix},\begin{pmatrix}0&0\\0&1\end{pmatrix},\begin{pmatrix}1&0\\1&0\end{pmatrix},\right.$$

$$\begin{pmatrix}1&1\\0&0\end{pmatrix},\begin{pmatrix}0&1\\0&1\end{pmatrix},\begin{pmatrix}0&0\\1&1\end{pmatrix},\begin{pmatrix}1&0\\0&1\end{pmatrix},\begin{pmatrix}0&1\\1&0\end{pmatrix},$$

$$\left.\begin{pmatrix}1&1\\0&1\end{pmatrix},\begin{pmatrix}1&1\\1&0\end{pmatrix},\begin{pmatrix}0&1\\1&1\end{pmatrix},\begin{pmatrix}1&0\\1&1\end{pmatrix},\begin{pmatrix}1&1\\1&1\end{pmatrix}\right\} \cup$$

$$\left\{\begin{pmatrix}1&1&1&1\\1&1&1&1\end{pmatrix}\right\}$$

$$\subseteq \quad V_1 \cup V_2 \cup V_3 \cup V_4$$
$$= \quad V$$

is a 4-basis of V and the 4-dimension of V is given by {7, 31, 15, 1)}. Now consider V as a group 4-set linear algebra over the group $Z_2 = \{0,1\}$. Take

$$X = X_1 \cup X_2 \cup X_3 \cup X_4$$
$$= \{(100), (001), (010)\} \cup \{1, x, x^2, x^3, x^4\} \cup$$
$$\left\{\begin{pmatrix}1&0\\0&0\end{pmatrix},\begin{pmatrix}0&1\\0&0\end{pmatrix},\begin{pmatrix}0&0\\1&0\end{pmatrix},\begin{pmatrix}0&0\\0&1\end{pmatrix}\right\} \cup \left\{\begin{pmatrix}1&1&1&1\\1&1&1&1\end{pmatrix}\right\}$$
$$\subseteq V_1 \cup V_2 \cup V_3 \cup V_4 ,$$

now X is a 4-basis and 4-generates V as a group 4-set linear algebra over the group $Z_2 = \{0,1\}$. The four dimension of V is given by {(3, 5, 4, 1)}. Thus we see when some additional structure is imposed the n-cardinality reduces drastically.

Now having defined the notion of n-generating set of group n-vector spaces and group n-linear algebra, we now proceed on to define the notion of n-group transformation.

**DEFINITION 5.26:** *Let $V = V_1 \cup V_2 \cup \ldots \cup V_n$ and $W = W_1 \cup W_2 \cup \ldots \cup W_n$ be two group n set vector space and group m-set vector space over the same group G; m > n. Let $T: V \to W$ be a n map from V into W defined by $T_i : V_i \to W_j$ $(i \neq j)$ (i = 1, 2, …, n) is a group set linear map where $T = T_1 \cup T_2 \cup \ldots \cup T_n$.*

*(Now when m>n we say we can have several ways of defining T's). We call $T = T_1 \cup T_2 \cup \ldots \cup T_n$ as the quasi group n-linear map of V into W over G or quasi n-linear transformation.*

We first illustrate this by an example.

***Example 5.56:*** Let $V = V_1 \cup V_2 \cup V_3 \cup V_4$ be a group 4-vector space over $Z_2 = \{0,1\}$ where

$$
\begin{aligned}
V &= V_1 \cup V_2 \cup V_3 \cup V_4 \\
&= \{Z_2 \times Z_2 \times Z_2\} \cup \left\{ \begin{pmatrix} a & b \\ c & d \end{pmatrix} \middle| a,b,c,d \in Z_2 \right\} \cup \\
&\quad \left\{ \begin{pmatrix} a & a & a & a \\ a & a & a & a \end{pmatrix} \middle| a \in Z_2 \right\} \cup \{Z_2[x] \text{ all polynomials of degree} \leq 4\}.
\end{aligned}
$$

Let

$$
\begin{aligned}
W &= W_1 \cup W_2 \cup W_3 \cup W_4 \cup W_5 \\
&= \left\{ \begin{pmatrix} a & b \\ c & d \end{pmatrix} \middle| a,b,c,d \in Z_2 \right\} \cup \{Z_2 \times Z_2 \times Z_2 \times Z_2\} \cup \\
&\quad \left\{ \begin{pmatrix} a & a & a \\ a & a & a \end{pmatrix}, \begin{pmatrix} a & a & a & a & a \\ a & a & a & a & a \end{pmatrix} \right\} \cup \\
&\quad \left\{ \begin{pmatrix} a & a \\ a & a \\ a & a \\ a & a \end{pmatrix}, \begin{pmatrix} a & a \\ a & a \end{pmatrix} \middle| a \in Z_2 \right\} \cup \{Z_2[x] \text{ all polynomials of degree less than or equal to } 8\}
\end{aligned}
$$

be a group 5 set vector space over the group $Z_2 = \{0,1\}$. Let T = $T = T_1 \cup T_2 \cup T_3 \cup T_4$ . $V \to W$, such that
$T_1 : V_1 \to W_2$ defined by

$$T_1 (x, y, z) = (x , y, z, x + y + z);$$

$T_2 : V_2 \to W_1$ defined by

$$T_2\left(\begin{pmatrix} a & b \\ c & d \end{pmatrix}\right) = \begin{pmatrix} a & b \\ c & d \end{pmatrix}$$

$T_3 : V_3 \to W_4$ defined by

$$T_3 \begin{pmatrix} a & a & a & a \\ a & a & a & a \end{pmatrix} = \begin{pmatrix} a & a \\ a & a \\ a & a \\ a & a \end{pmatrix}.$$

and $T_4 : V_4 \to W_5$

$$T_4\,(a) = a \text{ if } a \in Z_2,\; T_4\,(x^n) = x^{2n}\;;\; 1 \leq n \leq 4$$

i.e., $p\,(x) = x^2 + x + 1$ then $T_4\,(p(x)) = x^4 + x^2 + 1$. It is easily verified $T = T_1 \cup T_2 \cup T_3 \cup T_4$ is a quasi group n-linear map of V into W. We can have several T's for some group n-linear map $S = S_1 \cup S_2 \cup S_3 \cup S_4$; $V \to W$ which can be defined by

$S_1 : V_1 \to W_5$
$S_2 : V_2 \to W_1$
$S_3 : V_3 \to W_2$
$S_4 : V_4 \to W_3$

such that S is a quasi group n-linear map. So we have several such quasi group n-linear maps.

Next we define the notion of pseudo quasi group n linear transformations.

**DEFINITION 5.27:** *Let $V = V_1 \cup \ldots \cup V_n$ and $W = W_1 \cup \ldots \cup W_m$ (n>m) be group n set vector space and group m set vector space over the same group G.*

*A n-map $T = T_1 \cup \ldots \cup T_n$ from V into W defined by $T_i : V_i \to W_j$ where at least more than one $V_i$ will be mapped onto the same $W_j$ ($1 \leq i \leq n$ and $1 \leq j \leq m$). If each $T_i$ is a group set linear* *transformation then we call T to be a pseudo quasi n-set linear transformation of V into W.*

We give an example.

***Example 5.57:*** Let

$$V = V_1 \cup V_2 \cup V_3 \cup V_4 \cup V_5$$

$$= \{Z \times Z \times Z\} \cup \left\{\begin{pmatrix} a & b \\ c & d \end{pmatrix} \middle| a, b, c, d \in Z\right\} \cup$$

$$\left\{\begin{pmatrix} a & a & a \\ a & a & a \end{pmatrix}, \begin{pmatrix} a & a & a \\ a & a & a \\ a & a & a \\ a & a & a \end{pmatrix} \middle| a \in Z\right\} \cup \{Z[x] \text{ all polynomials of degree less than or equal to } 4\} \cup$$

$$\left\{\begin{pmatrix} a & 0 & 0 & 0 \\ b & c & 0 & 0 \\ d & e & f & 0 \\ g & h & i & j \end{pmatrix} \middle| a, b, c, d, e, f, g, h, i, j \in Z\right\}$$

is a group 5-set vector space over the group Z. Take

$$W = W_1 \cup W_2 \cup W_3 \cup W_4$$

$$= \{Z[x] \text{ all polynomials of degree less than or equal to } 2\}$$

$$\cup \left\{\begin{pmatrix} a & a & a & a \\ a & a & a & a \\ a & a & a & a \end{pmatrix}, \begin{pmatrix} a & a \\ a & a \\ a & a \end{pmatrix} \middle| a \in Z\right\} \cup \left\{\begin{pmatrix} a & a \\ a & a \end{pmatrix} \middle| a \in Z\right\}$$

$$\cup \left\{\begin{pmatrix} a & b & c & d \\ 0 & e & f & g \\ 0 & 0 & h & i \\ 0 & 0 & 0 & j \end{pmatrix} \middle| a, b, c, d, e, f, g, h, i, j \in Z\right\}$$

be a group 4-set vector space over Z.

Define $T = T_1 \cup T_2 \cup T_3 \cup \ldots \cup T_5$ from V to W by

$T_1 : V_1 \to W_1$ defined by

$$T_1 (a, b, c) = ax^2 + bx + c$$

$T_2 : V_2 \to W_3$ defined by

$$T_2 \left( \begin{pmatrix} a & b \\ c & d \end{pmatrix} \right) = \begin{pmatrix} a & a \\ a & a \end{pmatrix}$$

$T_3 : V_3 \to W_2$ defined by

$$T_3 \left( \begin{pmatrix} a & a & a \\ a & a & a \end{pmatrix} \right) = \begin{pmatrix} a & a \\ a & a \\ a & a \end{pmatrix}$$

$$T_3 \left( \begin{pmatrix} a & a & a \\ a & a & a \\ a & a & a \\ a & a & a \end{pmatrix} \right) = \begin{pmatrix} a & a & a & a \\ a & a & a & a \\ a & a & a & a \end{pmatrix}$$

$T_4 : V_4 \to W_1$

$$\begin{aligned} T_4 (c) &= c \quad \forall\, c \in Z \\ T_4 (x) &= x \\ T_4 (x^2) &= x^2 \\ T_4 (x^3) &= x \\ T_4 (x_4) &= x^2 \end{aligned}$$

and

$T_5 : V_5 \to W_4$ defined by

$$T_5 \left( \begin{pmatrix} a & 0 & 0 & 0 \\ b & c & 0 & 0 \\ d & e & f & 0 \\ g & h & i & j \end{pmatrix} \right) = \begin{pmatrix} a & b & c & d \\ 0 & e & f & g \\ 0 & 0 & h & i \\ 0 & 0 & 0 & j \end{pmatrix}$$

i.e.;

$$T_5 \begin{pmatrix} 4 & 0 & 0 & 0 \\ 7 & 8 & 0 & 0 \\ 9 & 6 & 1 & 0 \\ 2 & 1 & 3 & 5 \end{pmatrix} = \begin{pmatrix} 4 & 7 & 8 & 9 \\ 0 & 6 & 1 & 2 \\ 0 & 0 & 1 & 3 \\ 0 & 0 & 0 & 5 \end{pmatrix}.$$

Now $T = T_1 \cup T_2 \cup T_3 \cup T_4 \cup T_5 : V = V_1 \cup V_2 \cup V_3 \cup V_4 \cup V_5 \to W_1 \cup W_2 \cup W_3 \cup W_4$ is a pseudo quasi 5-set linear transformation V into W.

It is important to mention here that the study of what is the structure of $H_G^q$ (V, W) where $H_G^q$ (V, W) denotes the collection of all quasi group n-set linear transformations of V into W, is an interesting and innovative problem. Likewise the study of the algebraic structure of $\text{Hom}_G^{pq}$ (V, W) where $\text{Hom}_G^{pq}$ (V, W) is the set of all pseudo quasi group n-set linear transformations of V into W is left as an exercise for the readers.

Now we proceed onto discuss about the n-projection.

**DEFINITION 5.28:** *Let $V = V_1 \cup V_2 \cup \ldots \cup V_n$ be a $(m_1, \ldots, m_n)$ dimensional group n-set vector space over the group G; each $m_i$ is finite $1 \leq i \leq n$. Suppose $W = W_1 \cup W_2 \cup \ldots \cup W_n$ be a $(t_1, \ldots, t_n)$ dimensional group n-set vector subspace of V over the group G; $t_i < m_i$, $1 \leq i \leq n$.*

*Let $T = T_1 \cup T_2 \cup \ldots \cup T_n \to V$ into W be such that $T_i : V_i \to W_i$ is a projection of $V_i$ into $W_i$ ; $1 \leq i \leq n$. Then $T = T_1 \cup T_2 \cup \ldots \cup T_n$ is defined as the group n set projection of the group n set vector space V into the group n-set vector subspace W of V.*

We now illustrate this by simple examples.

**Example 5.58:** Let

$$V = V_1 \cup V_2 \cup V_3 \cup V_4$$

$$= \left\{ \begin{pmatrix} a & b \\ c & d \end{pmatrix} \middle| a, b, c, d \in Z \right\} \cup \{Z[x] \text{ all polynomials of degree less than or equal to } 5\} \cup \{Z \times Z \times Z \times Z\} \cup$$

$$\left\{\begin{pmatrix} a & b \\ c & d \\ e & f \\ g & h \end{pmatrix} \middle| a,b,c,d,e,f,g,h \in Z\right\}$$

be the group 4 set vector space over the group G.
Let

$W = W_1 \cup W_2 \cup W_3 \cup W_4$

$= \left\{\begin{pmatrix} a & a \\ a & a \end{pmatrix} \middle| a \in Z\right\} \cup \{Z[x]$ all polynomials of degree less than or equal to 3$\} \cup \{Z \times Z \times Z \times \{0\} \times \{0\}\}$

$$\cup \left\{\begin{pmatrix} a & a \\ a & a \\ a & a \\ a & a \end{pmatrix} \middle| a \in Z\right\}$$

be a proper group 4 set vector subspace of V over G. Define P = $P_1 \cup P_2 \cup P_3 \cup P_4$ from V into W by $P_i : V_i \to W_i$; i = 1, 2, 3, 4 by

$P_1 : V_1 \to W_1$ is such that

$$P_1\left(\begin{pmatrix} a & b \\ c & d \end{pmatrix}\right) = \begin{pmatrix} a & a \\ a & a \end{pmatrix}.$$

$P_2 : V_2 \to W_2$ is defined by

$$\begin{aligned} P_2(a) &= a, \text{ for all } a \in Z. \\ P_2(x) &= x \\ P_2(x^2) &= x^2 \\ P_2(x^3) &= x^3 \\ P_2(x^4) &= 0 \\ P_2(x^5) &= 0. \end{aligned}$$

$P_3 : V_3 \to W_3$ is such that by

$$P_3((x, y, z, w, t)) = (x, y, z, 0, 0)$$

and $P_4 : V_4 \to W_4$ is given by

$$P_4\left(\begin{pmatrix} a & b \\ c & d \\ e & f \\ g & h \end{pmatrix}\right) = \begin{pmatrix} a & a \\ a & a \\ a & a \\ a & a \end{pmatrix}.$$

It is easily verified $P = P_1 \cup P_2 \cup P_3 \cup P_4$ is a group 4 set projection of the group 4 set vector space $V = V_1 \cup V_2 \cup V_3 \cup V_4$ into the group 4-set vector subspace $W = W_1 \cup W_2 \cup W_3 \cup W_4$ of V.

$$P \text{ o } P = P_1 \text{ o } P_1 \cup P_2 \text{ o } P_2 \cup P_3 \text{ o } P_3 \cup P_4 \text{ o } P_4 \text{ is given by}$$

$$P \text{ o } P = P_1 \text{ o } P_1 \left(\begin{pmatrix} a & b \\ c & d \end{pmatrix}\right) \cup P_2 \text{ o } P_2 \left((a_o + a_1x + a_2x^2 + a_3x^3 + a_4x^4 + a_5x^5)\right) \cup P_3 \text{ o } P_3 ((x, y, z, \omega, t)) \cup P_4 \text{ o } P_4 \left(\begin{pmatrix} a & b \\ c & d \\ e & f \\ g & h \end{pmatrix}\right)$$

$$= P_1 \left(\begin{pmatrix} a & a \\ a & a \end{pmatrix}\right) \cup P_2 ((a_o + a_1x + a_2x^2 + a_3x^3)) \cup P_3 ((x\ y\ z\ 0\ 0)) \cup P_4 \begin{pmatrix} a & a \\ a & a \\ a & a \\ a & a \end{pmatrix}$$

$$= \begin{pmatrix} a & a \\ a & a \end{pmatrix} \cup (a_o + a_1x + a_2x^2 + a_3x^3) + (x\ y\ z\ 0\ 0) \cup \begin{pmatrix} a & a \\ a & a \\ a & a \\ a & a \end{pmatrix}$$

$$= P;$$

i.e.; the group 4-set projection of the group 4 set vector space is an idempotent 4-map.

It is left as an exercise for the reader to prove that the n-projection of a group n-set linear transformation from $V = V_1 \cup \ldots \cup V_n$ into its n-subspace is an idempotent n-map.

We illustrate this situation by another example.

***Example 5.59:*** Let

$$V = V_1 \cup V_2 \cup V_3 \cup V_4 \cup V_5$$

$$= \left\{ \begin{pmatrix} a & b \\ c & d \end{pmatrix} \middle| a, b, c, d \in Z_3 \right\} \cup \{Z_3 [x] \text{ all polynomials of degree less than or equal to } 5\} \cup \{ Z_3 \times Z_3 \times Z_3 \times Z_3 \times Z_3\} \cup \left\{ \begin{pmatrix} a & b & c & d & e \\ f & g & h & i & j \end{pmatrix} \middle| a, b, c, d, e, f, g, h, i, j \in Z_3 \right\} \cup \left\{ \begin{pmatrix} a & b \\ c & d \\ e & f \\ g & h \end{pmatrix} \middle| a, b, c, d, e, f, g, h \in Z_3 \right\}$$

be a group 5-set vector space over the group $Z_3$. Take a proper 5 subset W of V i.e.,

$$W = W_1 \cup W_2 \cup W_3 \cup W_4 \cup W_5$$

$$= \left\{ \begin{pmatrix} a & a \\ a & a \end{pmatrix} \middle| a \in Z_3 \right\} \cup \{\text{all polynomials of even degree only from } V_2\} \cup \{(a, a, a, a, a) / a \in Z_3\} \cup \left\{ \begin{pmatrix} a & a & a & a & a \\ a & a & a & a & a \end{pmatrix} \middle| a \in Z_3 \right\} \cup \left\{ \begin{pmatrix} a & a \\ a & a \\ a & a \\ a & a \end{pmatrix} \middle| a \in Z_3 \right\}$$

$$\subseteq V_1 \cup V_2 \cup V_3 \cup V_4 \cup V_5 = V.$$

Define a group 5 set linear transformation P: V → V by $P = P_1 \cup P_2 \cup P_3 \cup P_4 \cup P_5 : V_1 \cup V_2 \cup V_3 \cup V_4 \cup V_5 \to V_1 \cup V_2 \cup V_3 \cup V_4 \cup V_5$ by

$$P_1 = \left(\begin{pmatrix} a & b \\ c & d \end{pmatrix}\right) = \begin{pmatrix} a & a \\ a & a \end{pmatrix}$$

where $P_1 : V_1 \to V_1$

$P_2 : V_2 \to V_2$ given by

$$P_2 ((a_o + a_1x + a_2x^2 + a_3x^3 + a_4x^4 + a_5x^5)) = a_o + a_2x^2 + a_4x^4;$$

$P_3 : V_3 \to V_3$ defined by

$$P_3 ((x, y, z, \omega, t)) = (x, x, x, x, x);$$

$P_4 : V_4 \to V_4$ defined by

$$P_4 \left\{\begin{pmatrix} a & b & c & d & e \\ f & g & h & i & j \end{pmatrix}\right\} = \begin{pmatrix} a & a & a & a & a \\ a & a & a & a & a \end{pmatrix}$$

and $P_5 : V_5 \to V_5$

$$P_5 \left(\begin{pmatrix} a & b \\ c & d \\ e & f \\ g & h \end{pmatrix}\right) = \begin{pmatrix} a & a \\ a & a \\ a & a \\ a & a \end{pmatrix}.$$

We see $P = P_1 \cup P_2 \cup P_3 \cup P_4 \cup P_5$ can be treated as a 5-projection map of $V = V_1 \cup V_2 \cup V_3 \cup V_4 \cup V_5$ into its group 5 set vector subspace $W = W_1 \cup W_2 \cup W_3 \cup W_4 \cup W_5$ .

It is easily further verified that $P \circ P = P_1 \circ P_1 \cup P_2 \circ P_2 \cup P_3 \circ P_3 \cup P_4 \circ P_4 \cup P_5 \circ P_5 = P_1 \cup P_2 \cup \ldots \cup P_5$.

Now we can have for any group n-set linear algebra a new substructure a pseudo group semigroup n-set linear subalgebra.

**DEFINITION 5.29:** *Let $V = V_1 \cup V_2 \cup \ldots \cup V_n$ be a group n-set linear algebra over the group G. Let $W = W_1 \cup W_2 \cup \ldots \cup W_n$ be a proper n-subset of $V = V_1 \cup V_2 \cup \ldots \cup V_n$ ; each $W_i \subset V_i$; $W_i$ is an additive subsemigroup of the additive group $V_i$; i = 1, 2, …, n. If we have a proper subset H of G where H is a semigroup under the operations of G and $W = W_1 \cup W_2 \cup \ldots \cup W_n$ is a semigroup n set linear algebra over the semigroup H,* *then we call $W = W_1 \cup W_2 \cup \ldots \cup W_n$ to be the pseudo group semigroup n-set linear subalgebra of V over the subsemigroup H of G.*

We illustrate this situation by the following examples.

***Example 5.60:*** Let

$$\begin{aligned} V &= V_1 \cup V_2 \cup V_3 \cup V_4 \\ &= \{Z \times Z \times Z\} \cup \{Z[x] \text{ all polynomials of degree less than or equal to five}\} \cup \left\{ \begin{pmatrix} a & b \\ c & d \end{pmatrix} \middle| a,b,c,d \in Z \right\} \cup \\ &\quad \left\{ \begin{pmatrix} a & a & a & a & a \\ a & a & a & a & a \end{pmatrix} \middle| a \in Z \right\} \end{aligned}$$

be a group 4-set linear algebra over the additive group Z. Let

$$\begin{aligned} W &= W_1 \cup W_2 \cup W_3 \cup W_4 \\ &= \{(Z^+ \times Z^+ \times Z^+)\} \cup \{Z^+ [x], \text{ all polynomials of degree less than or equal to five}\} \cup \left\{ \begin{pmatrix} a & b \\ c & d \end{pmatrix} \middle| a,b,c,d \in Z^+ \right\} \cup \\ &\quad \left\{ \begin{pmatrix} a & a & a & a & a \\ a & a & a & a & a \end{pmatrix} \middle| a \in 2Z^+ \right\} \\ &\subseteq V_1 \cup V_2 \cup V_3 \cup V_4 \\ &= V. \end{aligned}$$

It is easily verified that W is a pseudo group semigroup 4 set linear subalgebra of V.

***Example 5.61:*** Let

$$\begin{aligned} V &= V_1 \cup V_2 \cup V_3 \cup V_4 \cup V_5 \\ &= \{Z \times Z \times Z \times Z \times Z\} \cup \{Q[x] \text{ all polynomials of degree less than or equal to 5}\} \cup \end{aligned}$$

$$\left\{\begin{pmatrix} a & b & d \\ c & f & g \\ h & i & j \end{pmatrix} \middle| a,b,c,d,e,f,g,h,i,j \in Q\right\} \cup$$

$$\left\{\begin{pmatrix} a & a & a & a & a \\ 0 & a & a & a & a \\ 0 & 0 & a & a & a \\ 0 & 0 & 0 & a & a \\ 0 & 0 & 0 & 0 & a \end{pmatrix} \middle| a \in Q\right\} \cup \left\{\begin{pmatrix} a & a & a & a & a \\ b & b & b & b & b \\ c & c & c & c & c \end{pmatrix} \middle| a,b,c \in Q\right\}.$$

V is a group 5-set linear algebra over the additive group Z. Take

$$\begin{aligned} W &= W_1 \cup W_2 \cup W_3 \cup W_4 \cup W_5 \\ &= \{3Z^+ \times 3Z^+ \times 3Z^+ \times 3Z^+ \times 3Z^+\} \cup \{2Z^+[x] \mid \text{all polynomials of degree less than or equal to 5 with coefficients for } 2Z^+\} \cup \\ &\quad \left\{\begin{pmatrix} a & b & d \\ c & f & g \\ h & i & j \end{pmatrix} \middle| a,b,c,d,e,f,g,h,i,j \in 5Z^+\right\} \cup \\ &\quad \left\{\begin{pmatrix} a & a & a & a & a \\ 0 & a & a & a & a \\ 0 & 0 & a & a & a \\ 0 & 0 & 0 & a & a \\ 0 & 0 & 0 & 0 & a \end{pmatrix} \middle| a \in 7Z^+\right\} \cup \\ &\quad \left\{\begin{pmatrix} a & a & a & a & a \\ a & a & a & a & a \\ a & a & a & a & a \end{pmatrix} \middle| a \in 3Z^+\right\} \\ &\subseteq V_1 \cup V_2 \cup \ldots \cup V_5. \end{aligned}$$

It is easily verified that W is a pseudo group semigroup 5 set linear subalgebra of V over the semigroup $Z^+$.

## Chapter Six

# SET FUZZY LINEAR BIALGEBRA AND ITS GENERALIZATION

In this chapter we define the notion of set fuzzy linear bialgebra and its generalized structure namely fuzzy linear n-algebra.

**DEFINITION 6.1:** *Let $V = V_1 \cup V_2$ be a set vector bispace over the set S. We say the set vector bispace $V = V_1 \cup V_2$ with the bimap $\eta = \eta_1 \cup \eta_2$ is a fuzzy set bivector space or set fuzzy bivector space if $\eta = \eta_1 \cup \eta_2 : V = V_1 \cup V_2 \to [0, 1]$ and $\eta (r_1a_1 \cup r_2a_2) = (\eta_1 \cup \eta_2) (r_1a_1 \cup r_2a_2) = \eta_1 (r_1a_1) \cup \eta_2 (r_2a_2) \geq \eta_1(a_1) \cup \eta_2(a_2)$ for all $a_i \in V_i$; i = 1, 2 and $r \in S$. We call $V_\eta = V\eta = (V_1 \cup V_2)_{\eta_1 \cup \eta_2}$ or $\eta V = (\eta_1 \cup \eta_2) (V_1 \cup V_2)$, to be the fuzzy set bivector space over the set S.*

We illustrate this by the following simple examples.

***Example 6.1:*** Let $V = V_1 \cup V_2$ where $V_1$ = {(4 0 0 1), (5 2 2 1) (1 5 5 1), (7 8 9 2), (1 1 1 1), (7 2 7 2), (1 0 5 3)} ∪ {(1 1 2), (2 3 4), (5 6 7), (2 3 1), (4 –1 2), (0 1 7), (0 8 0)} be the biset which is a set vector bispace over the set S = {0, 1}. Define a map $\eta = (\eta_1 \cup \eta_2): V = V_1 \cup V_2 \to [0, 1]$ by $\eta_i: V_i \to [0, 1]$ for i = 1, 2 defined by

$$\eta_1 (x, y, z, w) = \frac{x + y + z + w}{100} \in [0, 1]$$

$$\eta_2 (x, y, z, w) = \frac{x + y + z}{50} \in [0, 1]$$

$V_\eta = \left(V_1 \cup V_2\right)_{\eta_1 \cup \eta_2}$ is a fuzzy set vector bispace.

***Example 6.2:*** Let V = ($3Z^+ \cup 2Z^+$), $Z^+$ the set of positive integers. Take S = $Z^+$, V is a set vector bispace over S. Define $\eta = \eta_1 \cup \eta_2 : V_1 \cup V_2 \to [0, 1]$ by; for every $v_i \in V_i$ ,

$$\eta_i(v_i) = \frac{1}{v_i}, i = 1, 2.$$

Clearly $\eta V$ or $V_{\eta_1 \cup \eta_2} = \left(V_1 \cup V_2\right)_{\eta_1 \cup \eta_2}$ is a set fuzzy vector bispace or fuzzy vector bispace or fuzzy set vector bispace.

***Example 6.3:*** Let $V_1 = \{(a_{ij}) \mid a_{ij} \in Z^+, 1 \leq i, j \leq 7\} \cup \{(a_{ij}) \mid a_{ij} \in 2Z^+, 1 \leq i, j \leq 5\} = V_1 \cup V_2$ be a biset. $V = V_1 \cup V_2$ is a set vector bispace over the set S = $Z^+$.

Define $\eta = \eta_1 \cup \eta_2 : V_1 \cup V_2 \to [0, 1]$ by

$$\eta_1 (A_1 = (a_{ij}) ) = \frac{1 | A_1 |}{8 | A_1 |}$$

where $A_1 \in V_1$ and $|A_1|$ denotes the determinant value of A.

$$\eta_2 (A_1 = (a_{ij}) ) = \frac{1 | A_2 |}{7 | A_2 |}$$

where $A_2 \in V_2$ and if $|A_i| = 0$, $1 \leq i \leq 2$ then $\eta_i (A_i) = 0$. $V_\eta = (V_1 \cup V_2)_{\eta_1 \cup \eta_2}$ is a fuzzy set vector bispace.

One of the main advantages of using the concept of set vector bispaces and fuzzy set vector bispaces is that we can at liberty include any element in $V_1$ or $V_2$ provided $sx \in V_1$ or $V_2$. This cannot in general be done when we use usual vector bispaces. Another advantage of this system is that we can work with the minimum number of elements i.e., we can just delete any element in $V_1$ or $V_2$ ($V = V_1 \cup V_2$); if it is not going to be useful for us. So we can as per our need increase or decrease the number of elements in the set bivector spaces. This is one of the chief advantages, when the use of elements mean economy and time.

We give yet another example.

***Example 6.4:*** Let $V = V_1 \cup V_2$ where $V_1 = \{p(x) \in Z^+[x]$, all polynomials in x with coefficients from $Z^+\}$ and $V_2 = \{Z^+ \times 3Z^+ \times 5Z^+ \times 7Z^+\}$; V is a set vector bispace over the set $Z^+$. Define $\eta = \eta_1 \cup \eta_2$ from $V = V_1 \cup V_2$ to [0, 1] by

$$\eta_1(p(x)) = \begin{cases} \dfrac{1}{\deg p(x)} \\ 1 \quad \text{if } \deg(p(x)) = 0 \end{cases}$$

where $\eta_1 : V_1 \to [0, 1]$ and $\eta_2 : V_2 \to [0, 1]$ defined by

$$\eta_2 (x, y, z, w) = \frac{1}{x + y + z + w}.$$

Clearly $V_\eta = (V_1 \cup V_2)_{\eta_1 \cup \eta_2}$ is a set fuzzy vector bispace.

Now we proceed of to define notion of set fuzzy linear bialgebra.

**DEFINITION 6.2:** *A set fuzzy linear bialgebra $(V, \eta) = (V_1 \cup V_2, \eta_1 \cup \eta_2)$ or $(\eta_1 \cup \eta_2)$ $(V_1 \cup V_2)$ is an ordinary set linear bialgebra. V with a bimap $\eta = \eta_1 \cup \eta_2 : V_1 \cup V_2 \to [0, 1]$ such that $\eta_1(a_1 + b_1) \geq \min (\eta_1(a_1), \eta_1 (b_1))$; $\eta_2(a_2 + b_2) \geq \min (\eta_2(a_2), \eta_2(b_2))$ and $\eta_1(r_1 a_1) \cup \eta_2(r_2 a_2) \geq \eta_1(a_1) \cup \eta_2(a_2)$ for $a_1, b_1 \in V_1$ and $a_2, b_2 \in V_2$.*

*Since we know in the set vector bispace $V = V_1 \cup V_2$ we merely take V to be a biset, but in case of the set linear bialgebra $V = V_1 \cup V_2$ we assume V is closed with respect to some operation usually denoted by '+'; so the additional condition $\eta_1 (a_1+b_1) \geq \{\min (\eta_1 (a_1), \eta_1(b_1)\}$ and $\eta_2 (a_2 + b_2) \geq \min \{(\eta_2(a_2), \eta_2(b_2)\}$ is essential for every $a_i, b_i \in V_i$, i = 1,2.*

We denote this by some examples.

***Example 6.5:*** Let $V = V_1 \cup V_2$ where $V_1 = Z^+ \times Z^+ \times Z^+$ and

$$V_2 = \left\{ \begin{pmatrix} a & b \\ c & d \end{pmatrix} \middle| \; a, b, c, d \in Z^+ \right\}.$$

V is clearly a set linear bialgebra over the set $S = Z^+$. Define $\eta$ : $V \to [0, 1]$ where $\eta = \eta_1 \cup \eta_2 : V_1 \cup V_2 \to [0, 1]$ with $\eta : V_1 \to [0, 1]$ and $\eta_2 : V_2 \to [0, 1]$ as

$$\eta_1 ((x, y, z)) = \frac{1}{x+y+z}$$

and

$$\eta_2 \left( \begin{pmatrix} a & b \\ c & d \end{pmatrix} \right) = \begin{cases} \dfrac{ad-bc}{5(ad-bc)} & \text{if } ad \neq bc \\ 0 & \text{if } ad = bc \end{cases}.$$

Clearly $V_\eta = (V_1 \cup V_2)_{\eta_1 \cup \eta_2}$ is a fuzzy set linear bialgebra.

***Example 6.6:*** Let $V = V_1 \cup V_2$ where

$$V_1 = \left\{ \begin{pmatrix} a & b & c & d \\ e & f & g & h \end{pmatrix} \middle| \; a, b, c, d, e, f, g, h \in Z^+ \right\}$$

and

$$V_2 = \left\{ \begin{pmatrix} a \\ b \\ c \\ d \\ e \end{pmatrix} \middle| \; a, b, c, d, e \in 2Z^+ \right\}$$

be a set linear bialgebra over $Z^+$. Define $\eta = \eta_1 \cup \eta_2 : V = V_1 \cup V_2 \to [0, 1]$ by

$$\eta_1\begin{pmatrix} a & b & c & d \\ e & f & g & h \end{pmatrix} = \frac{1}{(a+b+c+d)\,(e+f+g+h)}$$

as both a + b + c + d ≠ 0 and e + f + g + h ≠ 0 which is true.

$$\eta_2\left(\begin{pmatrix} a \\ b \\ c \\ d \\ e \end{pmatrix}\right) = \frac{1}{a+c+e}.$$

Clearly $V_\eta = (V_1 \cup V_2)_{\eta 1 \cup \eta 2}$ is a set fuzzy linear bialgebra.

***Example 6.7:*** Let $V = V_1 \cup V_2$ where $V_1 = Z_{10}[x]$ and $V_2 = Z_{10} \times Z_{10} \times Z_{10} \times Z_{10}$. $V = V_1 \cup V_2$ is a set linear bialgebra. Define $\eta = \eta_1 \cup \eta_2 : V = V_1 \cup V_2 \to [0, 1]$ by

$$\eta_1(p(x)) = \begin{cases} \dfrac{1}{\deg p(x)} \\ 1 \text{ if } p(x) \text{ is constant or } 0 \end{cases}$$

and

$$\eta_2(a, b, c, d) = \begin{cases} \dfrac{1}{a+b+c+d} \text{ if } a+b+c+d \not\equiv (\text{mod}10) \\ 1 \text{ otherwise} \end{cases}$$

$V\eta$ is a set fuzzy linear bialgebra.

Now we proceed onto define the notion of fuzzy set vector subspace and fuzzy set linear subalgebra.

**DEFINITION 6.3:** *Let $V = V_1 \cup V_2$ be a set vector bispace over the set S. Let $W = W_1 \cup W_2 \subseteq V_1 \cup V_2$ be the set vector bisubspace of V defined over S. If $\eta = \eta_1 \cup \eta_2 :\ W = W_1 \cup W_2 \to [0, 1]$ then $W_\eta = (W_1 \cup W_2)_{\eta_1 \cup \eta_2}$ is called the fuzzy set vector bisubspace of V.*

We illustrate this by the following example.

***Example 6.8:*** Let V = {(1 1 1), (0 1 1), (1 0 1), (0 1 0), (1 1 0), (0 0 0)} $\cup$ {$Z_2[x]$} = $V_1 \cup V_2$ be a set vector bispace over $Z_2$ = {0, 1}. Define $\eta = \eta_1 \cup \eta_2$: $V = V_1 \cup V_2 \to [0, 1]$ by $\eta_1 : V_1 \to [0, 1]$ and $\eta_2 : V_2 \to [0, 1]$

$$\eta_1 (1\ 1\ 1) = 1 \qquad \eta_1 (1\ 0\ 1) = 0$$
$$\eta_1 (0\ 1\ 1) = 0 \qquad \eta_1 (0\ 1\ 0) = 0$$
$$\eta_1 (1\ 1\ 0) = 0 \qquad \eta (0\ 0\ 0) = 0$$

and

$$\eta_2 [p(x)] = \begin{cases} \dfrac{1}{\deg p(x)} \\ 0 \quad \text{if } p(x) \text{ is constant or } 0 \end{cases}$$

Take $W = W_1 \cup W_2$ = {(1 1 1), (0 0 0)} $\cup$ {All polynomials of degree less than or equal to 4}.

Now $\eta = \eta_1 \cup \eta_2 : W \to [0, 1]$ is well defined and $W\eta = (W_1 \cup W_2)_{\eta_1 \cup \eta_2}$ is a fuzzy set vector bisubspace of $V = V_1 \cup V_2$.

***Example 6.9:*** Let

$$\begin{aligned} V &= V_1 \cup V_2 \\ &= \{Z_{15} \times Z_{15} \times Z_{15}\} \cup \left\{ \begin{pmatrix} a & b \\ c & d \end{pmatrix} \middle|\ a, b. c, d \in Z_{15} \right\} \end{aligned}$$

be a set vector bispace over $Z_{15}$. Define $\eta = \eta_1 \cup \eta_2 : V = V_1 \cup V_2 \to [0, 1]$ by

$$\eta_1(a\ b\ c) = \begin{cases} 0 & \text{if } a+b+c \equiv 15 \text{ or } 0 (\text{mod } 15) \\ \dfrac{1}{a+b+c} & \text{if } a+b+c \neq 0 \end{cases}$$

$$\eta_2\left(\begin{pmatrix} a & b \\ c & d \end{pmatrix}\right) = \begin{cases} 0 & \text{if } a+b+c+d \equiv 0 (\text{mod } 15) \\ \dfrac{1}{a+b+c+d} & \text{if } a+b+c+d \not\equiv 0 (\text{mod} 15) \end{cases}$$

Now $V_\eta = (V_1 \cup V_2)_{\eta_1 \cup \eta_2}$ is a fuzzy set vector space. Take W = $W_1 \cup W_2 \subseteq V_1 \cup V_2$ where

$$W_1 = \{0\} \times Z_{15} \times Z_{15} \subseteq V_1$$

and

$$W_2 = \left\{ \begin{pmatrix} a & b \\ c & d \end{pmatrix} \middle| \; a, b, c, d \in Z_{15} \right\}$$

$W_\eta = (W_1 \cup W_2)_{\eta_1 \cup \eta_2}$ is a fuzzy set vector subspace of V.

$$\eta_1 (a\ b\ 0) = \begin{cases} 0 & \text{if } 0+a+b \equiv 0 (\text{mod } 15) \\ \dfrac{1}{0+a+b} & \text{if } 0+a+b \not\equiv 0 (\text{mod} 15) \end{cases}$$

$$\eta_2\left(\begin{pmatrix} a & b \\ 0 & d \end{pmatrix}\right) = \begin{cases} 0 & \text{if } a+b+0+d = 0 (\text{mod } 15) \\ \dfrac{1}{a+b+0+d} & \text{if } a+b+0+d \neq 0 (\text{mod} 15) \end{cases}$$

$W_\eta = (W_1 \cup W_2)_{\eta_1 \cup \eta_2}$ is a fuzzy set vector subspace of V.

***Example 6.10:*** Let $V = V_1 \cup V_2$ where

$$V_1 = \left\{ \begin{pmatrix} a & b & c \\ d & e & f \end{pmatrix} \middle| \; a, b, c, d, e, f \in Z^+ \right\}$$

and

$$V_2 = \{Q^+ \times Q^+ \times Z^+ \times Z^+\}.$$

$V = V_1 \cup V_2$ is a set vector bispace over the set $S = Z^+$.

Define $\eta = \eta_1 \cup \eta_2 : V_1 \cup V_2 = V \to [0, 1]$ by

$$\eta_1\left(\begin{pmatrix} a & b & c \\ d & e & f \end{pmatrix}\right) = \frac{1}{a+b+c+d+e+f}$$

$$\eta_2\,(a\ b\ c\ d) = \begin{cases} \dfrac{1}{a+b+c+d} & \text{if } a+b+c+d \geq 1 \text{ and an integer} \\ 0 & \text{if } a+b+c+d < 1 \text{ or a proper rational} \end{cases}$$

of the form p/q, q ≠ 0. $V_\eta = (V_1 \cup V_2)_{\eta_1 \cup \eta_2}$ is a fuzzy set vector bispace. Take $W = W_1 \cup W_2 \subseteq V_1 \cup V_2 = V$ where

$$W_1 = \left\{\begin{pmatrix} a & a & a \\ a & a & a \end{pmatrix} \middle|\ a \in Z^+\right\}$$

and

$$W_2 = 5Z^+ \times 3Z^+ \times 2Z^+ \times 7Z^+ \subseteq V_2.$$

Clearly $W_\eta = (W_1 \cup W_2)_{\eta_1 \cup \eta_2}$ is a fuzzy set vector bisubspace of V.

Now we proceed onto define the notion of fuzzy set bilinear algebra.

***Example 6.11:*** Let $V = V_1 \cup V_2$ where

$$V_1 = \left\{\begin{pmatrix} a & a \\ a & a \\ a & a \\ a & a \\ a & a \end{pmatrix} \middle|\ a \in Z^+\right\}$$

and
$V_2 = \{Z^+[x] \mid$ all polynomials in the variable x with coefficients from $Z^+\}$. $V = V_1 \cup V_2$ is a set linear bialgebra over $Z^+$. Define $\eta = \eta_1 \cup \eta_2 : V_1 \cup V_2 \to [0, 1]$ by

$$\eta_1 \begin{pmatrix} a & a \\ a & a \\ a & a \\ a & a \\ a & a \end{pmatrix} = \frac{1}{8a} \text{ and } \eta_2\,(p(x)) = \frac{1}{\deg\, p(x)}.$$

Clearly $V_\eta = (V_1 \cup V_2)_{\eta_1 \cup \eta_2}$ is a fuzzy set linear bialgebra. Take $W = W_1 \cup W_2 \subseteq V_1 \cup V_2 = V$ where

$$W_1 = \begin{pmatrix} 3a & 3a \\ 3a & 3a \\ 3a & 3a \\ 3a & 3a \\ 3a & 3a \end{pmatrix} \subseteq V_1$$

and

$W_2$ = {All polynomial in the variable x of even degree} $\subseteq V_2$;

$W_\eta = (W_1 \cup W_2)_{\eta_1 \cup \eta_2}$ is a fuzzy set linear subbialgebra of V.

***Example 6.12:*** Let $V = V_1 \cup V_2$ where

$$V_1 = \{Z_{12} \times Z_{12} \times Z_{12} \times Z_{12}\}$$

and

$$V_2 = \left\{ \begin{pmatrix} a & b & c \\ d & e & f \\ g & h & i \end{pmatrix} \middle|\, a, b, c, d, e, f, i \in Z_{12} \right\}.$$

We see V is a set linear subbialgebra of V over $Z_{12}$. Take

$W = W_1 \cup W_2 =$

$$\{Z_{12} \times \{0\} \times Z_{12} \times \{0\}\} \cup \left\{ \begin{pmatrix} a & b & c \\ 0 & e & f \\ 0 & 0 & i \end{pmatrix} \middle|\, a, b, c, d, e, f, i \in Z_{12} \right\}$$

Define $\eta = \eta_1 \cup \eta_2 : V_1 \cup V_2 \to [0, 1]$ by

$$\eta_1 (a\ b\ c\ d) = \begin{cases} 0 & \text{if } a + b + c + d \equiv 0 (\text{mod } 12) \\ \dfrac{1}{a+b+c+d} & \text{if } a + b + c + d \not\equiv 0 (\text{mod} 12) \end{cases}$$

$$\eta_2 \left\{ \begin{pmatrix} a & b & c \\ d & e & f \\ g & h & i \end{pmatrix} \right\} = \begin{cases} \dfrac{1}{a} & \text{if } a \neq 0 \\ 0 & \text{if } a = 0 \end{cases}.$$

Now $V_\eta = (V_1 \cup V_2)_{\eta_1 \cup \eta_2}$ is a fuzzy set linear bialgebra. $W_\eta = (W_1 \cup W_2)_{\eta_1 \cup \eta_2}$ is a fuzzy set linear bisubalgebra of $V_\eta$.

Now we proceed onto define fuzzy semigroup vector bispaces.

**DEFINITION 6.4:** *A semigroup fuzzy vector bispace or a fuzzy semigroup vector bispace $V_\eta$ (or $\eta V$) is an ordinary semigroup vector bispace V over the semigroup S; with a bimap $\eta = \eta_1 \cup \eta_2$; $V = V_1 \cup V_2 \to [0, 1]$ satisfying the following conditions $\eta(ra) = \eta_1(r_1a_1) \cup \eta_2(r_2a_2)$ is such that $\eta_1(r_1a_1) \geq \eta_1(a_1)$ and $\eta_2(r_2a_2) \geq \eta_2(a_2)$ for all $a_1 \in V_1$ and $a_2 \in V_2$ and $r_1, r_2 \in S$ i.e., $\eta(ra) = \eta_1(r_1a_1) \cup \eta_2(r_2a_2) \geq \eta_1 (a_1) \cup \eta_2(a_2)$.*

We illustrate this by some examples.

***Example 6.13:*** Let $V = V_1 \cup V_2$ where $V_1 = \{(1\ 1\ 1\ 1\ 1), (0\ 1\ 1\ 1\ 0), (1\ 0\ 0\ 0\ 1), (0\ 0\ 0\ 0\ 0)\}$ and $V_2 = \{Z_2 \times Z_2 \times Z_2 \times Z_2 \times Z_2\}$. V is a semigroup vector bispace over the semigroup $Z_2 = \{0, 1\}$. Define $\eta = \eta_1 \cup \eta_2 : V = V_1 \cup V_2 \to [0, 1]$ by

$$\eta_1 (a\ b\ c\ d) = \frac{(a+b+c+d)\ (\text{mod } 2)}{5}$$

and

$$\eta_2 (a\ b\ c\ d\ e) = \frac{(a+b+c+d+e)\ (\text{mod } 2)}{7}$$

Clearly $V_\eta = (V_1 \cup V_2)_{\eta_1 \cup \eta_2}$ is a fuzzy semigroup vector bispace.

***Example 6.14:*** Let V = $\{Z_5 \times Z_5 \times Z_5\} \cup \{Z_5[x]$; all polynomials in x of degree less than or equal to 10 with coefficients from $Z_5\} = V_1 \cup V_2$. V is a semigroup vector bispace over the semigroup $Z_5$.
Define $\eta = \eta_1 \cup \eta_2 : V \to [0, 1]$; $\eta_1 : V_1 \to [0, 1]$ with

$$\eta_1 (x\ y\ z) = \begin{cases} \dfrac{1}{(x+y+z)} & \text{if } x+y+z \not\equiv 0 (\text{mod } 5) \\ 0 & \text{if } x+y+z \equiv 0 (\text{mod } 5) \end{cases}.$$

$\eta_2 : V_2 \to [0, 1]$ given by

$$\eta_2 (p(x)) = \frac{1}{\deg (p(x))}; \ \eta_2(p(x)) = 0 \text{ if } p(x) \text{ is a constant.}$$

Clearly $V_\eta = (V_1 \cup V_2)_{\eta_1 \cup \eta_2}$ is a fuzzy semigroup vector bispace. Define another bimap $\delta : V \to [0, 1]$ as $\delta = \delta_1 \cup \delta_2 : V_1 \cup V_2 \to [0, 1]$ by $\delta_1 : V_1 \to [0, 1]$, $\delta_2 : V_2 \to [0, 1]$ where

$$\delta_1 (x\ y\ z) = \begin{cases} \dfrac{1}{5(x+y+z)} & \text{if } x+y \neq 0; \text{not modulo addition} \\ 0 & \text{if } x+y+z = 0; \text{i.e., } x = y = z = 0 \end{cases}$$

$$\delta_2 (p(x)) = \frac{1}{3 \deg p(x)}; \ \delta_2(p(x)) = 0 \text{ if } p(x) \text{ is a constant.}$$

$V_\delta = (V_1 \cup V_2)_{\delta_1 \cup \delta_2}$ is yet another fuzzy semigroup vector bispace. Thus given any semigroup vector bispace, we can obtain infinite number of fuzzy semigroup vector bispaces.

Thus we can have several semigroup fuzzy vector bispaces for a given semigroup vector bispace. This is the advantage of using the notion of fuzzy semigroup vector bispaces.

Now we proceed onto define the notion of fuzzy semigroup vector bispaces.

**DEFINITION 6.5:** *Let $V = V_1 \cup V_2$ be a semigroup vector bispace over the semigroup S. Suppose $V_\eta = (V_1 \cup V_2)_{\eta_1 \cup \eta_2}$ is a semigroup fuzzy vector bispace of V and if $W = W_1 \cup W_2 \subseteq V_1 \cup V_2$ is a semigroup vector subbispace of V, then $W_\eta = (W_1 \cup W_2)_{\eta_1 \cup \eta_2}$ is defined to be the semigroup fuzzy vector subbispace of V; where $\bar{\eta} : W \to [0, 1]$ i.e., $\bar{\eta}$ is the restriction of $\eta$ to W i.e., $\bar{\eta}_1 : W_1 \to [0, 1]$ and $\bar{\eta}_2 : W_2 \to [0, 1]$.*

We illustrate this by some examples.

***Example 6.15:*** Let

$$V = V_1 \cup V_2 = \left\{ \begin{pmatrix} a & b & c \\ d & e & f \\ g & h & i \end{pmatrix} \middle| \; a, b, c, d, e, f, g, h, i \in Z^+ \cup \{0\} \right\} \cup \{Z^0 \times Z^0 \times Z^0 \times Z^0 \mid Z^0 = Z^+ \cup \{0\}\}$$

be a semigroup vector bispace over the semigroup $Z^0 = Z^+ \cup \{0\}$. Define $\eta = \eta_1 \cup \eta_2 : V = V_1 \cup V_2 \to [0, 1]$ by

$$\eta_1 \begin{pmatrix} a & b & c \\ d & e & f \\ g & h & i \end{pmatrix} = \begin{cases} \dfrac{1}{a+e+i} & \text{if } a+e+i \neq 0 \\ 0 & \text{if } a+e+i = 0 \end{cases}.$$

$\eta_2 : V_2 \to [0, 1]$ is defined by

$$\eta_2 (a, b, c, d) = \begin{cases} \dfrac{1}{a+b+c+d} & \text{if } a+b+c+d \neq 0 \\ 0 & \text{if } a+b+c+d = 0 \end{cases}.$$

$V_\eta = (V_1 \cup V_2)_{\eta_1 \cup \eta_2}$ is a semigroup fuzzy vector bispace.
Take

$$\begin{aligned} W &= W_1 \cup W_2 \\ &= \left\{\left.\begin{pmatrix} a & a & a \\ a & a & a \\ a & a & a \end{pmatrix}\right| a \in Z^0 \right\} \cup \{(a\ a\ a\ a) \mid (a \in Z^0\} \\ &\subseteq V_1 \cup V_2 \end{aligned}$$

is a proper semigroup vector subbispace of V. Now $\overline{\eta}_1 : W_1 \to [0, 1]$ is given by

$$\overline{\eta}_1 \begin{pmatrix} a & a & a \\ a & a & a \\ a & a & a \end{pmatrix} = \begin{cases} \dfrac{1}{3a} & \text{if } 3a \neq 0 \\ 0 & \text{if } 3a = 0 \end{cases}$$

and $\overline{\eta}_2 : W_2 \to [0, 1]$ given by

$$\overline{\eta}_2\,(a\ a\ a\ a) = \begin{cases} \dfrac{1}{4a} & \text{if } a \neq 0 \\ 0 & \text{if } a = 0 \end{cases}$$

is such that $W_\eta = (W_1 \cup W_2)_{\eta_1 \cup \eta_2}$ is a semigroup fuzzy vector subbispace of V.

***Example 6.16:*** Let

$$\begin{aligned} V &= V_1 \cup V_2 \\ &= \{Z_{16} \times Z_{16} \times Z_{16} \times Z_{16}\} \cup \left\{\left.\begin{pmatrix} a & a & a & a & a \\ a & a & a & a & a \end{pmatrix}\right| a \in Z_{16} \right\} \end{aligned}$$

be a semigroup vector bispace over the semigroup $Z_{16}$. Let

$$\begin{aligned} W &= W_1 \cup W_2 \\ &= \{S \times S \times S \times S \mid S = \{0, 2, 4, 6, 8, 10, 12, 14\}\} \\ &\cup \left\{\left.\begin{pmatrix} a & a & a & a & a \\ a & a & a & a & a \end{pmatrix}\right| a \in S \right\} \end{aligned}$$

$\subseteq \quad V = V_1 \cup V_2$.

W is a semigroup vector bisubspace over the semigroup $Z_{16}$. Define $\eta = \eta_1 \cup \eta_2 : W_1 \cup W_2 \to [0, 1]$ as follows :
$\eta_1 : W_1 \to [0, 1]$ defined by

$$\eta_1 \text{ (a b c d)} = \begin{cases} \frac{1}{a} & \text{if } a \neq 0 \\ 0 & \text{if } a = 0 \end{cases}.$$

$\eta_2 : W_2 \to [0, 1]$ defined by

$$\eta_2 \begin{pmatrix} a & a & a & a & a \\ a & a & a & a & a \end{pmatrix} = \begin{cases} \frac{1}{10a} & \text{if } a \neq 0 \\ 0 & \text{if } a = 0 \end{cases}.$$

Clearly $W_\eta = (W_1 \cup W_2)_{\eta_1 \cup \eta_2}$ is a semigroup fuzzy vector bisubspace.

Now we proceed onto define semigroup fuzzy linear bialgebra.

**DEFINITION 6.6:** *Let $V = V_1 \cup V_2$ be a semigroup linear bialgebra defined over the semigroup S. We say $V_\eta$ or $\eta$ V or $(V_1 \cup V_2)_{\eta_1 \cup \eta_2}$ or $(\eta_1 \cup \eta_2)$ $(V_1 \cup V_2)$ is a semigroup fuzzy linear bialgebra if $\eta : V_1 \cup V_2 \to [0, 1]$ where $\eta = \eta_1 \cup \eta_2$ such that $\eta_1 : V_1 \to [0, 1]$ and $\eta_2 : V_2 \to [0, 1]$ satisfying the conditions $\eta_i(x_i + y_i) \geq \min(\eta_i(x), \eta_i(y))$; $\eta_i(rx_i) \geq \eta_i(x_i)$ ; i = 1, 2 for every $r \in S$; $y_i, x_i \in V_i$; i = 1,2.*

Now we illustrate this by a few examples.

***Example 6.17:*** Let

$$V = V_1 \cup V_2 = \{Z_8 \times Z_8 \times Z_8\} \cup \left\{ \begin{pmatrix} a_1 & a_2 & a_3 \\ a_4 & a_5 & a_6 \end{pmatrix} \middle| a_i \in Z_8; 1 \leq i \leq 6 \right\}$$

be a semigroup linear bialgebra over the semigroup $Z_8$. Define $\eta = \eta_1 \cup \eta_2 : V \to [0, 1]$ as follows :

$\eta_1 : V_1 \to [0, 1]$ is defined by

$$\eta_1 (x\ y\ z) = \begin{cases} \frac{1}{x} & \text{if } x \neq 0 \\ 0 & \text{if } x = 0 \end{cases}.$$

$\eta_2: V_2 \to [0, 1]$ defined by

$$\eta_2 \left( \begin{pmatrix} a_1 & a_2 & a_3 \\ a_4 & a_5 & a_6 \end{pmatrix} \right) = \begin{cases} \frac{1}{a_6} & \text{if } a_6 \neq 0 \\ 0 & \text{if } a_6 = 0 \end{cases}.$$

$V_\eta = (V_1 \cup V_2)_{\eta_1 \cup \eta_2}$ is a semigroup fuzzy linear bialgebra.

***Example 6.18:*** Let $V = V_1 \cup V_2 = \{p(x) \mid p(x) \in Z_9^+ [x]$ where $p(x)$ is of degree less than or equal to 9 with coefficients from $Z^+\} \cup$

$$\left\{ \begin{pmatrix} a & b & c & d \\ e & f & g & h \\ i & j & k & l \\ m & n & o & p \end{pmatrix} \middle| \; a, b, \ldots, p \in Z^+ \right\}$$

is a semigroup linear bialgebra over the semigroup $Z^+$. Define $\eta = \eta_1 \cup \eta_2 : V = V_1 \cup V_2 \to [0, 1]$ by

$\eta_1 : V_1 \to [0, 1]$ such that

$$\eta_1 (p(x)) = \frac{1}{\deg p(x)}; \; \eta_1(p(x)) = 0 \text{ if } p(x) \text{ is constant or } 0$$

$\eta_2 : V_2 \to [0, 1]$ is such that

$$\eta_2 \begin{pmatrix} a & b & c & d \\ e & f & g & h \\ i & j & k & l \\ m & n & o & p \end{pmatrix} = \frac{1}{a+b+k+p}.$$

$V_\eta = (V_1 \cup V_2)_{\eta_1 \cup \eta_2}$ is a semigroup fuzzy linear bialgebra.

We see as in case of fuzzy set vector spaces and semigroup fuzzy vector spaces the notion of fuzzy set vector bispaces and semigroup fuzzy vector bispaces are identical. Likewise the notion of fuzzy set linear bialgebra and fuzzy semigroup linear bialgebra are also identical. This sort of visualizing them as identical algebraic structures is done using the fuzzy tool which will be useful in applications. We as in case of fuzzy set linear algebras define these structures as fuzzy equivalent algebraic structures.

Thus we see fuzzy set linear bialgebra and semigroup fuzzy linear bialgebra are fuzzy equivalent though set linear bialgebra and semigroup linear bialgebra are different from each other.

Now we proceed onto define group fuzzy vector bispaces and group fuzzy linear bialgebras.

**DEFINITION 6.7:** *Let $V = V_1 \cup V_2$ be a group linear bialgebra over the group G. Let $\eta = \eta_1 \cup \eta_2 : V = V_1 \cup V_2 \to [0, 1]$ such that*

$$\eta_i (a_i + b_i) \geq \min (\eta_i (a_i), \eta_i(b_i))$$

$$\eta_i (-a_i) = \eta_i (a_i)$$

$$\eta_i (0) = 1$$

$$\eta_i (ra_i) \geq \eta_i (a_i) \text{ for all } a_i, b_i \in V_i \text{ and } r \in G, i = 1, 2.$$

*We call $V_\eta = (V_1 \cup V_2)_{\eta_1 \cup \eta_2}$ to be the group fuzzy linear bialgebra.*

Likewise the same definition holds good for group fuzzy vector bispaces as both group fuzzy vector bispaces and group fuzzy linear bialgebras are fuzzy equivalent.

Now we illustrate this by some examples.

***Example 6.19:*** Let

$$V = V_1 \cup V_2$$

$$= \{Z \times Z \times Z\} \cup \left\{ \begin{pmatrix} a & b & c \\ d & e & f \end{pmatrix} \middle| a, b, c, d, e, f \in Z \right\}$$

be a group vector bispace over the group Z. Define $\eta = \eta_1 \cup \eta_2$ : $V = V_1 \cup V_2 \to [0, 1]$ where

$\eta_1 : V_1 \to [0, 1]$ is such that

$$\eta_1 \text{ (a b c)} = \begin{cases} \dfrac{1}{5\,|(a+b+c)|} & \text{if } a+b+c \neq 0 \\ 1 & \text{if } a+b+c = 0 \end{cases}.$$

$\eta_2 : V_2 \to [0,1]$ such that

$$\eta_2 \begin{pmatrix} a & b & c \\ d & e & f \end{pmatrix} = \begin{cases} \dfrac{1}{|a+b+c+e|} & \text{if } a+b+c+e \neq 0 \\ 1 & \text{if } a+b+c+e = 0 \end{cases};$$

$\eta = \eta_1 \cup \eta_2 : V \to [0, 1]$ is well defined. $V_\eta = (V_1 \cup V_2)_{\eta_1 \cup \eta_2}$ is a group fuzzy linear bialgebra.

***Example 6.20:*** Let $V = V_1 \cup V_2$ be such that $V_1$ = {Z[x] all polynomials of degree less than or equal to 25} $\cup$ {Z × Z × Z × Z} be a group linear bialgebra over the group Z. Define

$\eta = \eta_1 \cup \eta_2 : V_1 \cup V_2 \to [0, 1]$ as follows

$\eta_1 : V_1 \to [0,1]$ by

$$\eta_1 \text{ (p(x))} = \begin{cases} \eta_1(p(x)) = \dfrac{1}{\deg p(x)} \\ \eta_1(a) = 1 \quad a \in Z \end{cases}$$

$\eta_2 : V_2 \to [0, 1]$ is defined by

$$\eta_2 \text{ (a b c d)} = \begin{cases} \dfrac{1}{|a+b+c+d|} & \text{if } a+b+c+d \neq 0 \\ 1 & \text{if } a+b+c+d = 0 \end{cases}.$$

$V_\eta = (V_1 \cup V_2)_{\eta_1 \cup \eta_2}$ is a group fuzzy linear bialgebra.

Now we proceed onto define the notion of fuzzy group vector bisubspace and fuzzy group vector bisubalgebra.

**DEFINITION 6.8:** *Let $V = V_1 \cup V_2$ be a group vector bispace over the group G. Let $W = W_1 \cup W_2 \subseteq V_1 \cup V_2 = V$, where W is a group vector subspace of V.*

*Define $\eta = \eta_1 \cup \eta_2$: $W = W_1 \cup W_2 \to [0, 1]$ such that $\eta_1 : W_1 \to [0, 1]$, $\eta_2 : W_2 \to [0, 1]$ where $\eta_i (a_i + b_i) \geq \min (\eta_i(a_i), \eta_i (b_i))$*

$$\eta_i (a_i) = \eta_i (-a_i)$$
$$\eta_i (0) = 1$$
$$\eta_i (ra_i) \geq r\eta_i (a_i) \text{ for all } a_i \in V_i \text{ and } r \in G;\ i = 1, 2.$$

*We call $W_\eta = (W_1 \cup W_2)_{\eta_1 \cup \eta_2}$ to be the group fuzzy bivector subspace.*

We illustrate this definition by an example.

***Example 6.21*:** Let

$$\begin{aligned} V &= V_1 \cup V_2 \\ &= \{Z \times Z \times Z \times Z\} \cup \left\{ \begin{pmatrix} a & a & a & a & a \\ a & a & a & a & a \end{pmatrix} \middle| \, a \in Z \right\}, \end{aligned}$$

be a group vector bispace over the group Z. Take

$$\begin{aligned} W &= W_1 \cup W_2 \\ &= \{2Z \times 2Z \times 2Z \times \{0\}\} \cup \left\{ \begin{pmatrix} a & a & a & a & a \\ a & a & a & a & a \end{pmatrix} \middle| \, a \in 2Z \right\} \\ &\subseteq V = V_1 \cup V_2 \end{aligned}$$

is a group vector subbispace of V. Define $\eta = \eta_1 \cup \eta_2 : W = W_1 \cup W_2 \to [0, 1]$; $\eta_1 : W_1 \to [0, 1]$ and $\eta_2 : W_2 \to [0, 1]$ with

$$\eta_1 (x\ y\ z\ 0) = \begin{cases} \dfrac{1}{x + y + z} & \text{if } x + y + z \neq 0 \\ 1 & \text{if } x + y + z = 0 \end{cases}.$$

$$\eta_2 \left( \begin{pmatrix} a & a & a & a & a \\ a & a & a & a & a \end{pmatrix} \right) = \begin{cases} \dfrac{1}{10a} & \text{if } a \neq 0 \\ 1 & \text{if } a = 0 \end{cases}.$$

$W_\eta = (W_1 \cup W_2)_{\eta_1 \cup \eta_2}$ is a group linear fuzzy subalgebra. Further as our transformation to a fuzzy set up always demands only values from the positive unit interval [0, 1] these semigroup vector bispaces or set vector spaces would be more appropriate than the ordinary vector bispaces.

Further we see even in case of Markov biprocess or Markov bichain, the transition probability bimatrix is a square bimatrix whose entries are non negative adding up to one. So in such cases one can use set vector bispaces where

$$V = V_1 \cup V_2 = \left\{ \left(a^1_{ij}\right)_{\eta_1 \times \eta_1} \middle| a^1_{ij} \in [0, 1] \right\} \cup \left\{ \left(a^2_{ij}\right)_{\eta_2 \times \eta_2} \middle| a^2_{ij} \in [0, 1] \right\}$$

with $\sum_{i_1=1}^{n_1} a^1_{i_k k_1} = 1$ for $1 \leq k_1 \leq n_1$

and $\sum_{i_2=1}^{n_2} a^2_{i_2 k_2} = 1$ for $1 \leq k_2 \leq n_2$

is a set vector bispace over the set [0, 1]. So these new notion not only comes handy but involve lesser complication and lesser algebraic structures.

**DEFINITION 6.9:** *Let $V = V_1 \cup V_2$ be such that $V_1$ is a set vector space over the set $S_1$ and $V_2$ is a set vector space over $S_2$ ; $S_1 \neq S_2$. We define $V = V_1 \cup V_2$ to be the biset bivector space or biset vector bispace over the biset $S = S_1 \cup S_2$.*
*Let $\eta = \eta_1 \cup \eta_2 : V = V_1 \cup V_2 \to [0, 1]$*
*where $\eta_1 : V_1 \to [0, 1]$ and*
*$\eta_2 : V_2 \to [0, 1]$*
*such that*

$$\eta_1(r_1 a_1) \geq \eta_1(a_1) \text{ for all } r_1 \in S, \text{ and } a_1 \in V_1.$$
$$\eta_2(r_2 a_2) \geq \eta_2(a_2) \text{ for all } r_2 \in S_2 \text{ and } a_2 \in V_2.$$

*$V_\eta = (V_1 \cup V_2)_{\eta_1 \cup \eta_2}$ is a biset fuzzy vector bispace.*

We illustrate this by some simple examples.

***Example 6.22:*** Let

$$V = V_1 \cup V_2$$
$$= \{Z^+ \times Z^+ \times Z^+) \cup \left\{ \begin{pmatrix} a & a & a \\ a & a & a \end{pmatrix} \middle| a \in Z_{10} \right\}$$

be a biset bivector space over the biset $S = 3Z^+ \cup Z_{10}$ ie $V_1$ is a set vector space over the set $3Z^+$ and $V_2$ is a set vector space over $Z_{10}$.

$$\eta = \eta_1 \cup \eta_2 : V = V_1 \cup V_2 \to [0, 1]$$

i.e., $\eta_1 : V_1 \to [0, 1]$ and $\eta_2 : V_2 \to [0, 1]$

$$\eta_1 (x\ y\ z) = \frac{1}{x+y+z}.$$

$$\eta_2 \begin{pmatrix} a & a & a \\ a & a & a \end{pmatrix} = \begin{cases} \dfrac{1}{6a} & \text{if } a \neq 0 \\ 1 & \text{if } a = 0 \end{cases}$$

$V\eta$ is a biset fuzzy bivector space.

We give yet another example.

***Example 6.23:*** Let

$$V = V_1 \cup V_2$$
$$= \{(1111), (0000), (0011), (1100), (1001), (0110)\}$$
$$\cup \left\{ \begin{pmatrix} a & b \\ c & d \end{pmatrix} \middle| a, b, c, d \in Z_{12} \right\}$$

be a biset vector bispace over the biset $S = S_1 \cup S_2$ where $S_1 = \{0, 1\}$ and $S_2 = Z_{12}$.

Define $\eta = \eta_1 \cup \eta_2 : V = V_1 \cup V_2 \to [0, 1]$ by $\eta_1 : V_1 \to [0, 1]$, $\eta_2 : V_2 \to [0, 1]$; with

$$\eta_1 (a\ b\ c\ d) = \begin{cases} \dfrac{1}{4(a+b+c+d)} & \text{if } a+b+c+d \neq 0 \\ 1 & \text{if } a+b+c+d = 0 \end{cases}.$$

$$\eta_2 \begin{pmatrix} a & b \\ c & d \end{pmatrix} = \begin{cases} \frac{1}{3(ad-bc)} & \text{if } ad - bc \neq 0 \\ 0 & \text{if } ad - bc = 0 \end{cases}.$$

Now we proceed onto define the notion of biset fuzzy vector bisubspaces.

**DEFINITION 6.10:** *Let $V = V_1 \cup V_2$ be a biset vector bispace over the biset $S = S_1 \cup S_2$. Let $W = W_1 \cup W_2 \subseteq V_1 \cup V_2 = V$ be a proper bisubset of V which is a biset vector bispace over the same biset $S = S_1 \cup S_2$.*

*Define $\eta = \eta_1 \cup \eta_2 : W \to [0, 1]$ by $\eta_1 : W_1 \to [0, 1]$ and $\eta_2 : W_2 \to [0, 1]$ then we call $W_\eta = (W_1 \cup W_2)_{\eta_1 \cup \eta_2}$ to be the biset fuzzy vector subbispace.*

We illustrate this by some examples.

***Example 6.24:*** Let $V = V_1 \cup V_2 = \{Z_{10} \times Z_{10} \times Z_{10}\} \cup \{Z_{12} \times Z_{12} \times Z_{12} \times Z_{12} \times Z_{12}\}$ be a biset vector bispace over the biset $S = Z_{10} \cup Z_{12}$. Let $W = \{S_1 \times S_1 \times S_1$ where $S_1 = \{0, 2, 6, 4, 8\}\} \cup \{S_2 \times S_2 \times S_2 \times \{0\} \times S_2$ where $S_2 = \{0, 3, 6, 9\}\} = W_1 \cup W_2 \subset V_1 \cup V_2$ be a biset vector subbispace of V over the biset $S = Z_{10} \cup Z_{12}$. Define the bimap $\eta = \eta_1 \cup \eta_2 : W_1 \cup W_2 \to [0, 1]$ : by $\eta_1 : W_1 \to [0 , 1]$ and

$\eta_2 : W_2 \to [0 , 1]$ given by

$$\eta_1 (x\ y\ z) = \begin{cases} \frac{1}{x+y+z} & \text{if } x + y + z \not\equiv 0 (\text{mod} 10) \\ 1 & \text{if } x + y + z \equiv 0 (\text{mod} 10) \end{cases}$$

$$\eta_2(x\ y\ z\ 0\ t) = \begin{cases} \frac{1}{x+y+z+0+t} & \text{if } x + y + z + 0 + t \not\equiv 0 (\text{mod} 12) \\ 1 & \text{if } x + y + z + 0 + t \equiv 0 (\text{mod} 12) \end{cases}$$

$W_\eta = (W_1 \cup W_2)_{\eta_1 \cup \eta_2}$ is a biset fuzzy vector bispace.

***Example 6.25:*** Let

$$V = V_1 \cup V_2$$
$$= \{Z \times Z \times Z \times Z\} \cup \left\{ \begin{pmatrix} a & b & c \\ d & e & f \end{pmatrix} \middle| \, a, b, c, d, e, f \in Z_{15} \right\}$$

be a biset vector bispace over the biset $S = Z \cup Z_{15}$. Consider

$$W = W_1 \cup W_2$$
$$= \{2Z \times 2Z \times \{0\} \times \{0\}\} \cup \left\{ \begin{pmatrix} a & a & a \\ a & a & a \end{pmatrix} \middle| \, a \in Z_{15} \right\}$$
$$\subseteq V_1 \cup V_2 = V$$

be a biset vector subbispace of V over the same biset $S = Z \cup Z_{15}$. Now consider $\eta$ the fuzzy set where $\eta = \eta_1 \cup \eta_2$ with $\eta_1 : W_1 \to [0, 1]$ and $\eta_2 : W_2 \to [0, 1]$ defined by

$$\eta_1 (a\ b\ 0\ 0) = \begin{cases} \dfrac{1}{a+b} & \text{if } a+b \neq 0 \\ 1 & \text{if } a+b = 0 \end{cases}$$

and

$$\eta_2 \begin{pmatrix} a & a & a \\ a & a & a \end{pmatrix} = \begin{cases} \dfrac{1}{a} & \text{if } a \neq 0 \\ 1 & \text{if } a = 0 \end{cases}.$$

Clearly $W\eta = (W_1 \cup W_2)_{\eta_1 \cup \eta_2}$ is a biset fuzzy bivector space.

**DEFINITION 6.11:** *Let $V = V_1 \cup V_2$ be a biset linear bialgebra over the biset $S = S_1 \cup S_2$. Define $\eta = \eta_1 \cup \eta_2 : V = V_1 \cup V_2 \to [0, 1]$ in the following way $\eta_i : V_i \to [0, 1]$, $i = 1, 2$ such that $\eta_i(a_i + b_i) \geq \min (\eta_i (a_i), \eta_i (b_i))$ for all $a_i, b_i \in V_i$, $i = 1, 2$.*

*$V_\eta = (V_1 \cup V_2)_{\eta_1 \cup \eta_2}$ is a biset fuzzy bilinear algebra.*

We illustrate this by some simple examples.

***Example 6.26:*** Let

$$V = V_1 \cup V_2 = \left\{\begin{pmatrix} a & b \\ c & d \end{pmatrix} \middle| a, b, c, d \in Z_{12}\right\} \cup \{Z_5 \times Z_5 \times Z_5 \times Z_5\}$$

be a biset linear bialgebra over the biset $S = Z_{12} \cup Z_5$. Define $\eta = \eta_1 \cup \eta_2 : V \to [0, 1]$ by $\eta_1 : V_1 \to [0, 1]$, $\eta_2 : V_2 \to [0, 1]$ such that

$$\eta_1 \begin{pmatrix} a & b \\ c & d \end{pmatrix} = \begin{cases} \dfrac{1}{a+b+c+d} & \text{if } a+b+c+d \not\equiv 0 (\text{mod} 12) \\ 1 & \text{if } a+b+c+d \equiv 0 (\text{mod} 12) \end{cases}.$$

$$\eta_2 \,(a\ b\ c\ d) = \begin{cases} \dfrac{1}{a+d} & \text{if } a+d \not\equiv 0 (\text{mod} 5) \\ 1 & \text{if } a+d \equiv 0 (\text{mod} 5) \end{cases}.$$

Clearly $V_\eta = (V_1 \cup V_2)_{\eta_1 \cup \eta_2}$ is a biset linear fuzzy bialgebra.

***Example 6.27:*** Let $V = V_1 \cup V_2 = \{Z^+[x]\}$ all polynomials in the variable x with coefficients from $Z^+\} \cup$

$$\left\{\begin{bmatrix} a_1 & a_2 & a_3 & a_4 \\ a_5 & a_6 & a_7 & a_8 \\ a_9 & a_{10} & a_{11} & a_{12} \end{bmatrix} \middle| a_i \in 2Z\right\}.$$

V is a biset linear bialgebra over the biset $Z^+ \cup 2Z$. Define the bimap $\eta = \eta_1 \cup \eta_2 : V = V_1 \cup V_2 \to [0, 1]$ by

$$\eta_1(p(x)) = \frac{1}{\deg\ p(x)};\ \eta_1(p(x)) = 1 \text{ if } p(x) \text{ is a constant}$$

and

$$\eta_2\begin{pmatrix} a_1 & a_2 & a_3 & a_4 \\ a_5 & a_6 & a_7 & a_8 \\ a_9 & a_{10} & a_{11} & a_{12} \end{pmatrix} =$$

$$\begin{cases} \dfrac{1}{a_1+a_2+a_3+a_4} & \text{if } a_1+a_2+a_3+a_4 \neq 0 \\ 1 & \text{if } a_1+a_2+a_3+a_4 = 0 \end{cases}$$

We see $V_\eta = (V_1 \cup V_2)_{\eta_1 \cup \eta_2}$ is a biset fuzzy linear bialgebra.

Now we proceed onto define the notion of biset fuzzy linear bialgebra.

**DEFINITION 6.12:** *Let $V = V_1 \cup V_2$ be a biset linear bialgebra over the biset $S = S_1 \cup S_2$. Let $W = W_1 \cup W_2 \subseteq V_1 \cup V_2 = V$ be a biset linear subbialgebra of W over the biset $S = S_1 \cup S_2$. Define $\eta = \eta_1 \cup \eta_2 : W_1 \cup W_2 \to [0, 1]$ by $\eta_1 : W_1 \to [0, 1]$ and $\eta_2 : W_2 \to [0, 1]$. We see $W\eta = (W_1 \cup W_2)_{\eta_1 \cup \eta_2}$ is a biset fuzzy linear subbialgebra.*

We illustrate this by some examples.

***Example 6.28:*** Let

$$\begin{aligned} V &= V_1 \cup V_2 \\ &= \left\{ \begin{pmatrix} a & b & c \\ d & e & f \end{pmatrix} \middle| \, a, b, c, d, e, f \in Z_5 \right\} \cup Z_4 \times Z_4 \times Z_4 \times Z_4 \end{aligned}$$

be the biset linear bialgebra over the biset $S = Z_5 \cup Z_4$. Let

$$\begin{aligned} W &= W_1 \cup W_2 \\ &= \left\{ \begin{pmatrix} a & a & a \\ a & a & a \end{pmatrix} \middle| \, a \in Z_5 \right\} \cup \{(a\ a\ a\ a) \mid a \in Z_4\} \\ &\subseteq V = V_1 \cup V_2 \end{aligned}$$

be the biset linear subbialgebra. Define $\eta = \eta_1 \cup \eta_2 : W \to [0, 1]$ as $\eta_i : W_1 \to [0, 1]$ and $\eta_2 : W_2 \to [0, 1]$ as

$$\eta_1\left(\begin{pmatrix} a & a & a \\ a & a & a \end{pmatrix}\right) = \begin{cases} \frac{1}{a} & \text{if } a \neq 0 \\ 1 & \text{if } a = 0 \end{cases}$$

$$\eta_2 = (a\ a\ a\ a) = \begin{cases} \frac{1}{4a} & \text{if } a \neq 0 \\ 1 & \text{if } a = 0 \end{cases}$$

$W_\eta = (W_1 \cup W_2)_{\eta_1 \cup \eta_2}$ is a biset fuzzy linear subbialgebra.

***Example 6.29:*** Let

$$V = V_1 \cup V_2 = \{Z[x]\} \cup \left\{\begin{pmatrix} a & b & c & g \\ d & e & f & h \end{pmatrix} \middle| a, b, c, d, e, f, g, h \in Z_{23}\right\}$$

be a biset linear bialgebra over the biset $S = Z \cup Z_{23}$. Take

$W = W_1 \cup W_2$
$= \{$all polynomials of even degree with coefficients from $Z\} \cup \left\{\begin{pmatrix} a & a & a & a \\ a & a & a & a \end{pmatrix} \middle| a \in Z_{23}\right\}$
$\subseteq V_1 \cup V_2$

which is clearly a biset linear bialgebra over the same biset $S = Z \cup Z_{23}$. Define the bimap $\eta = \eta_1 \cup \eta_2 : W_1 \cup W_2 \to [0, 1]$ as $\eta_1 : W_1 \to [0, 1]$, $\eta_2 : W_2 \to [0, 1]$ as,

$$\eta_1(p(x)) = \frac{1}{\deg p(x)}; \ \eta_1(p(x)) = 1 \text{ if } p(x) \text{ is a constant}$$

$$\eta_2\left(\begin{pmatrix} a & a & a & a \\ a & a & a & a \end{pmatrix}\right) = \begin{cases} \frac{1}{a} & \text{if } a \neq 0 \\ 1 & \text{if } a = 0 \end{cases}.$$

Now $W_\eta = (W_1 \cup W_2)_{\eta_1 \cup \eta_2}$ is a biset fuzzy linear subbialgebra.

***Example 6.30:*** Let

$$V = V_1 \cup V_2 = \left\{ \begin{pmatrix} a & b \\ c & d \\ e & f \\ g & h \end{pmatrix} \middle| \, a, b, c, d, e, f, g, h \in Z^+ \cup \{0\} = Z^o \right\} \cup \{2Z \times 2Z \times 2Z\}$$

be a biset linear bialgebra over the biset $S = Z^o \cup 2Z$. Take

$$W = W_1 \cup W_2 = \left\{ \begin{pmatrix} a & a \\ a & a \\ a & a \\ a & a \end{pmatrix} \middle| \, a \in Z^o \right\} \cup \{2Z \times 2Z \times \{0\}\} \subseteq V_1 \cup V_2 = V,$$

W is a biset linear subbialgebra of V. Define the bimap $\eta = \eta_1 \cup \eta_2 : W = W_1 \cup W_2 \to [0, 1]$ as $\eta_1 : W_1 \to [0, 1]$ and $\eta_2 : W_2 \to [0, 1]$ defined by

$$\eta_1 \begin{pmatrix} a & a \\ a & a \\ a & a \\ a & a \end{pmatrix} = \begin{cases} \dfrac{1}{a} & \text{if } a \neq 0 \\ 1 & \text{if } a = 0 \end{cases}$$

$$\eta_2 (a\ b\ 0) = \begin{cases} \dfrac{1}{a+b} & \text{if } a + b \neq 0 \\ 1 & \text{if } a + b = 0 \end{cases}.$$

We see $W_\eta = (W_1 \cup W_2)_{\eta_1 \cup \eta_2}$ is a biset fuzzy linear subbialgebra.

Now we proceed onto define the notion of quasi semigroup fuzzy bialgebra. However it is pertinent to mention though quasi semigroup bilinear algebra is different from semigroup linear bialgebra but semigroup fuzzy bilinear bialgebra and quasi semigroup fuzzy linear bialgebra are one and the same.

**DEFINITION 6.13:** *Let $V = V_1 \cup V_2$ be a quasi semigroup linear bialgebra over the set S. Define $\eta = \eta_1 \cup \eta_2 : V \to [0, 1]$ as $\eta_1 : V_1 \to [0, 1]$ and $\eta_2 : V_2 \to [0, 1]$ such that $V_\eta = (V_1 \cup V_2)_\eta$ is a quasi semigroup fuzzy linear bialgebra.*

We illustrate this by the following example.

***Example 6.31:*** Let

$$\begin{aligned} V &= V_1 \cup V_2 \\ &= \{(1\ 0\ 0\ 0), (0\ 0\ 0\ 0), (1\ 1\ 1\ 0\ 0), (0\ 0\ 0\ 0\ 0), (1\ 1\ 0), \\ &\quad (0\ 1\ 1), (0\ 0\ 0), (1\ 0\ 1)\} \cup \left\{ \begin{pmatrix} a & b \\ c & d \end{pmatrix} \middle|\ a, b\ c\ d \in Z_2 \right\} \end{aligned}$$

be a quasi semigroup, linear bialgebra over the set $Z_2$. Define $\eta = \eta_1 \cup \eta_2 : V = V_1 \cup V_2 \to [0, 1]$ by $\eta_1((1\ 0\ 0\ 0)) = 1$, $\eta_1((0\ 0\ 0\ 0)) = 0$, $\eta_1((0\ 0\ 0\ 0\ 0)) = 0$, $\eta_1((1\ 1\ 0)) = \frac{1}{2}$, $\eta_1((0\ 0\ 0)) = 0$, $\eta_1((0\ 1\ 1)) = \frac{1}{2}$, $\eta_1((1\ 0\ 1)) = \frac{1}{2}$ where $\eta_1: V_1 \to [0, 1]$ and $\eta_2: V_2 \to [0, 1]$;

$$\eta_2 \begin{pmatrix} a & b \\ c & d \end{pmatrix} = \begin{cases} \dfrac{1}{a+b+c+d} & \text{if } a+b+c+d \neq 0 (\text{mod } 2) \\ 1 & \text{if } a+b+c+d = 0 (\text{mod } 2) \end{cases}.$$

Clearly $V_\eta = V_{\eta_1 \cup \eta_2} = (V_1 \cup V_2)_{\eta_1 \cup \eta_2}$ is a quasi semigroup fuzzy linear bialgebra.

Now we proceed onto define the notion of bisemigroup fuzzy bivector space over a bisemigroup $S = S_1 \cup S_2$.

**DEFINITION 6.14:** *Let $V = V_1 \cup V_2$ be a bisemigroup bivector space over the bisemigroup $S = S_1 \cup S_2$. Let $\eta = \eta_1 \cup \eta_2$ be a bimap from $V = V_1 \cup V_2$ into [0, 1], if $\eta_i(a_i + b_i) \geq \min(\eta_i(a_i), \eta_i(b_i))$ and $\eta_i(r_i a_i) \geq r_i \eta_i(a_i)$ for all $a_i, b_i \in V_i$ and $r_i \in S_i$; i = 1, 2. Then we call $V_\eta = (V_1 \cup V_2)_{\eta_1 \cup \eta_2}$ to be a bisemigroup fuzzy vector bispace or bisemigroup fuzzy bivector space.*

We illustrate this by a simple example.

***Example 6.32:*** Let

$$\begin{aligned} V &= V_1 \cup V_2 \\ &= \{(a\ a\ a) \mid a \in Z_6\} \cup \left\{ \begin{pmatrix} a & b \\ c & d \end{pmatrix} \middle| \, a, b, c, d \in Z^+ \right\} \end{aligned}$$

be a bisemigroup vector bispace over the bisemigroup $S = Z_6 \cup Z^+$. Define a bimap $\eta = \eta_1 \cup \eta_2 : V = V_1 \cup V_2 \to [0, 1]$ by $\eta_1 : V_1 \to [0, 1]$ is such that

$$\eta_1 (a\ a\ a) = \begin{cases} \dfrac{1}{a} & \text{if } a \neq 0 \\ 1 & \text{if } a = 0 \end{cases}.$$

$\eta_2 : V_2 \to [0, 1]$ is such that

$$\eta_2 \left( \begin{pmatrix} a & b \\ c & d \end{pmatrix} \right) = \frac{1}{a + d}.$$

Clearly $V_\eta = (V_1 \cup V_2)_{\eta_1 \cup \eta_2}$ is a bisemigroup fuzzy vector bispace.

We define now the notion of bisemigroup fuzzy linear bialgebra.

**DEFINITION 6.15:** *Let $V = V_1 \cup V_2$ be a semigroup linear bialgebra over the bisemigroup $S = S_1 \cup S_2$. Suppose $\eta = \eta_1 \cup \eta_2$: $V = V_1 \cup V_2 \to [0, 1]$ is a bimap such that $V_1\eta_1$ is a semigroup fuzzy linear algebra and $V_2\eta_2$ is a semigroup fuzzy linear algebra then $V_\eta = (V_1 \cup V_2)_{\eta_1 \cup \eta_2}$ is a bisemigroup fuzzy linear bialgebra.*

We illustrate this by some examples.

***Example 6.33:*** Let

$$V = V_1 \cup V_2 = \left\{ \begin{pmatrix} a_1 & a_2 & a_3 \\ a_4 & a_5 & a_6 \\ a_7 & al & a_9 \end{pmatrix} \middle| a_i \in 3Z^+ \cup \{0\}; 1 \leq i \leq 9 \right\} \cup \{2Z^+ \times 2Z^+ \times 2Z^+ \times 2Z^+\}$$

be a bisemigroup linear bialgebra over the bisemigroup $S = (3Z^+ \cup \{0\}) \cup 2Z^+$. Define the bimap $\eta = \eta_1 \cup \eta_2 : V = V_1 \cup V_2 \to [0, 1]$ such that $\eta_1 : V_1 \to [0, 1]$

$$\eta_1 \left\{ \begin{pmatrix} a_1 & a_2 & a_3 \\ a_4 & a_5 & a_6 \\ a_7 & a_8 & a_9 \end{pmatrix} \right\} = \begin{cases} \dfrac{1}{a_1 + a_5 + a_9} & \text{if } a_1 + a_5 + a_9 \neq 0 \\ 1 & \text{if } a_1 + a_5 + a_9 = 0 \end{cases}$$

$$\eta_2 : V_2 \to [0, 1]; \qquad \eta_2 \,(a\ b\ c\ d) = \frac{1}{a+b+c+d}.$$

$V_\eta = (V_1 \cup V_2)_{\eta_1 \cup \eta_2}$ is a bisemigroup fuzzy linear bialgebra.

***Example 6.34:*** Let $V = V_1 \cup V_2 = \{Z^+[x]\} \cup \{3Z \times 3Z \times 3Z \times 3Z\}$ be a bisemigroup linear bialgebra over bisemigroup $S = Z^+$ $\cup$ 3Z. Define $\eta = \eta_1 \cup \eta_2 : V = V_1 \cup V_2 \to [0, 1]$; $\eta_1 : V_1 \to [0, 1]$ and $\eta_2 : V_2 \to [0, 1]$.

$$\eta_1(p(x)) = \frac{1}{\deg p(x)}; \ \eta_1(p(x)) = 1 \text{ if } p(x) \text{ is a constant}$$

$$\eta_2 (a\ b\ c\ d) = \begin{cases} \dfrac{1}{a+b+c+d} & \text{if } a+b+c+d \neq 0 \\ 1 & \text{if } a+b+c+d = 0 \end{cases}$$

$V_1\eta_1$ and $V_2\eta_2$ are semigroup, fuzzy linear algebras so $V_\eta = (V_1 \cup V_2)_{\eta_1 \cup \eta_2}$ is a bisemigroup fuzzy linear bialgebra.

Now we proceed onto define the notion of bisemigroup fuzzy linear subbialgebra.

**DEFINITION 6.16:** *Let $V = V_1 \cup V_2$ be a bisemigroup linear bialgebra over the bisemigroup $S = S_1 \cup S_2$. Let $W = W_1 \cup W_2 \subseteq V_1 \cup V_2 = V$ be a bisemigroup linear subbialgebra of V. Define a bimap $\eta = \eta_1 \cup \eta_2 : W_1 \cup W_2 \to [0, 1]$ such that $\eta_1 : W_1 \to [0, 1]$ and $\eta_2 : W_2 \to [0, 1]$ so that $W_1\ \eta_1$ and $W_2\ \eta_2$ are semigroup fuzzy linear subbialgebras. We call $W_\eta = (W_1 \cup W_2)_{\eta_1 \cup \eta_2}$ to be a bisemigroup fuzzy linear subbialgebra.*

We illustrate this by the following examples.

***Example 6.35:*** Let

$$V = V_1 \cup V_2 = \{3Z^+ \times 3Z^+ \times 3Z^+ \times 3Z^+ \times 3Z^+\} \cup \left\{ \begin{pmatrix} a & b \\ c & d \end{pmatrix} \middle| a, b, c, d \in Z_{20} \right\}$$

be a bisemigroup linear bialgebra over the bisemigroup $S = 3Z^+ \cup Z_{20}$.

Take

$$W = W_1 \cup W_2$$
$$= \{3Z^+ \times \{0\} \times \{0\} \times 3Z^+ \times 3Z^+\} \cup \left\{ \begin{pmatrix} a & a \\ a & a \end{pmatrix} \middle| \; a \in Z_{20} \right\}$$
$$\subseteq V_1 \cup V_2 = V,$$

W is a bisemigroup linear subbialgebra of V.

Define $\eta = \eta_1 \cup \eta_2 : W_1 \cup W_2 \to [0, 1]$ by

$$\eta_1 (x\ 0\ 0\ y\ z) = \frac{1}{x + y + z}$$

so that $W_1\eta_1$ is a semigroup fuzzy linear algebra.

$\eta_2 : W_2 \to [0, 1]$ defined by

$$\eta_2 \left( \begin{pmatrix} a & a \\ a & a \end{pmatrix} \right) = \begin{cases} \dfrac{1}{a} & \text{if } a \neq 0 \\ 1 & \text{if } a = 0 \end{cases},$$

$W_2\eta_2$ is a semigroup fuzzy linear algebra so that $W_\eta = W_1\eta_1 \cup W_2\eta_2$ $(W_1 \cup W_2)_{\eta_1 \cup \eta_2}$ is a bisemigroup fuzzy linear bialgebra.

Now we proceed onto define the notion of bigroup fuzzy bivector space.

**DEFINITION 6.17:** *Let $V = V_1 \cup V_2$ be a bigroup fuzzy bivector space over the bigroup $G = G_1 \cup G_2$. Let $\eta = \eta_1 \cup \eta_2 : V = V_1 \cup V_2 \to [0, 1]$ be such that $V_1\eta_1$ is a group fuzzy vector space and $V_2\eta_2$ is a group fuzzy vector space where $\eta_1 : V_1 \to [0, 1]$ and $\eta_2 : V_2 \to [0, 1]$, so that $V_\eta = (V_1 \cup V_2)_{\eta_1 \cup \eta_2}$ is a bigroup fuzzy vector bispace.*

We illustrate this by some simple examples.

***Example 6.36:*** Let $V = V_1 \cup V_2 = Z_3 [x] \cup Q \times Q$ be a bigroup vector bispace over the bigroup $G = Z_3 \cup Q$. Define $\eta = \eta_1 \cup \eta_2 : V = V_1 \cup V_2 \to [0, 1]$ where

$\eta_1 : V_1 \to [0, 1]$ is given by

$$\eta_1(p(x)) = \begin{cases} \dfrac{1}{\deg p(x)} & \\ 1 & \text{if } p(x) = 0 \end{cases}$$

$\eta_2 : V_2 \to [0, 1]$ is given by

$$\eta_2\,(a\;b) = \begin{cases} \dfrac{1}{a+b} & \text{if } a+b \neq 0, a+b \text{ is an integer} \\ 1 & \text{if } a+b = 0 \text{ or } a+b \text{ is a rational number} \end{cases}$$

$V_\eta = (V_1 \cup V_2)_{\eta_1 \cup \eta_2}$ is a bigroup fuzzy vector bispace.

Now we proceed onto define bigroup fuzzy vector subbispace.

**DEFINITION 6.18:** *Let $V = V_1 \cup V_2$ be a bigroup vector bispace over the bigroup $G = G_1 \cup G_2$. Take $W = W_1 \cup W_2 \subseteq V_1 \cup V_2$ be a bigroup vector subbispace of V over the same bigroup $G = G_1 \cup G_2$. Define $\eta = \eta_1 \cup \eta_2 : W_1 \cup W_2 \to [0, 1]$ such that $(W_1 \cup W_2)_{\eta_1 \cup \eta_2}$ is a bigroup fuzzy vector bispace then we call W to be a bigroup fuzzy vector subbispace.*

Now we define set n-fuzzy vector space.

**DEFINITION 6.19:** *Let $V = V_1 \cup V_2 \cup \ldots \cup V_n$ be a set n-vector space over the set S. If $\eta : \eta_1 \cup \eta_2 \cup \ldots \cup \eta_n : V = V_1 \cup V_2 \cup \ldots \cup V_n \to [0, 1]$ is such that $\eta_i : V_i \to [0, 1]$ so that $V_i\eta_i$ is a set fuzzy vector space for each i = 1, 2, …, n, then we call $V_\eta = (V_1 \cup V_2 \cup \ldots V_n)_{\eta_1 \cup \eta_2 \cup \ldots \cup \eta_n}$ to be a set fuzzy n-vector space.*

We illustrate this by an example.

***Example 6.37:*** Let

$V \;=\; V_1 \cup V_2 \cup V_3 \cup V_4$

$= \{(1\ 1\ 1), (0\ 0\ 0), (1\ 0\ 0), (0\ 1\ 0), (1\ 1\ 1), (0\ 0), (1\ 1\ 1\ 1),$
$(1\ 0\ 0\ 0), (0\ 0\ 0\ 0)\} \cup \left\{ \begin{pmatrix} a & a \\ a & a \end{pmatrix} \middle| a \in Z_2 \right\} \cup$
$\left\{ \begin{pmatrix} a & a & a & a \\ a & a & a & a \end{pmatrix} \middle| a \in Z_2 \right\} \cup \{Z_2[x]\}.$

V is set 4-vector space over the set S = {0, 1}. Define $\eta = \eta_1 \cup \eta_2 \cup \eta_3 \cup \eta_4 : V = V_1 \cup V_2 \cup V_3 \cup V_4 \to [0, 1]$; $\eta_1 : V_1 \to [0, 1]$ by

$$\eta_1 (a\ b\ c) = \begin{cases} \dfrac{1}{a+b+c} & \text{if } a+b+c \neq 0 \\ 1 & \text{if } a+b+c = 0 \end{cases}$$

$\eta_2 : V_2 \to [0, 1]$ by

$$\eta_2 \left( \begin{pmatrix} a & a \\ a & a \end{pmatrix} \right) = \begin{cases} \dfrac{1}{a} & \text{if } a \neq 0 \\ 1 & \text{if } a = 1 \end{cases}$$

$\eta_3 : V_3 \to [0, 1]$ defined by

$$\eta_3 \left( \begin{pmatrix} a & a & a & a \\ a & a & a & a \end{pmatrix} \right) = \begin{cases} \dfrac{1}{8a} & \text{if } a \neq 0 \\ 1 & \text{if } a = 1 \end{cases}$$

and $\eta_4 : V_4 \to [0, 1]$ given by

$\eta_4 (p(x) = \dfrac{1}{\deg p(x)}$; $\eta_4(p(x)) = 1$ if p(x) is a constant.

Clearly $V_i\eta_i$ 's are set fuzzy vector spaces for each i. Thus $V_\eta = (V_1 \cup V_2 \cup V_3 \cup V_4)_{\eta_1 \cup \eta_2 \cup \eta_3 \cup \eta_4}$ is a set fuzzy n-vector space.

Now we proceed onto define set fuzzy n-vector subspace.

**DEFINITION 6.20:** *Let $V = V_1 \cup V_2 \cup \ldots \cup V_n$ be a set n-vector space over the set S. If $W = W_1 \cup W_2 \cup \ldots \cup W_n \subseteq V_1 \cup V_2 \cup \ldots \cup V_n$ with $W_i \neq W_j$ if $i \neq j$; $W_i \not\subseteq W_j$, $W_j \not\subseteq W_i$; $1 \leq i, j \leq n$ and W is a set n-vector subspace of V over the set S. If $\eta = \eta_1 \cup \eta_2 \cup \ldots \cup \eta_n : W = W_1 \cup W_2 \cup \ldots \cup W_n \to [0, 1]$ where $\eta_i : W_i \to [0, 1]$; $1 \leq i \leq n$ are set fuzzy vector subspaces then we call $W_\eta = (W_1 \cup W_2 \cup \ldots \cup W_n)_{\eta_1 \cup \eta_2 \cup \ldots \cup \eta_n}$ to be the set fuzzy n-vector subspace.*

Now we illustrate this by a simple example.

***Example 6.38:*** Let

$$
\begin{aligned}
V &= V_1 \cup V_2 \cup V_3 \cup V_4 \\
&= \{(a\ a\ a), (a\ a) \mid a \in Z^+\} \cup \left\{ \begin{pmatrix} a & b \\ c & d \end{pmatrix} \middle| \ a, b, c, d \in Z^+ \right\} \cup \\
&\quad Z^+[x] \cup \left\{ \begin{pmatrix} a_1 \\ a_2 \\ a_3 \\ a_4 \end{pmatrix} \middle| \ a_i \in Z^+ \ 1 \leq i \leq 4 \right\}
\end{aligned}
$$

be a set 4-vector space over the set $S = Z^+$. Take

$$
\begin{aligned}
W &= W_1 \cup W_2 \cup W_3 \cup W_4 \\
&= \{(a\ a\ a) \mid a \in Z^+\} \cup \left\{ \begin{pmatrix} a & a \\ a & a \end{pmatrix} \middle| \ a \in Z^+ \right\} \cup \\
&\quad \{\text{all polynomials of even degree}\} \cup \left\{ \begin{pmatrix} a \\ a \\ a \\ a \end{pmatrix} \middle| \ a \in 2Z^+ \right\} \\
&\subseteq V_1 \cup V_2 \cup V_3 \cup V_4 = V
\end{aligned}
$$

is a set 4 vector subspace of V over the set $S = Z^+$. Define $\eta = \eta_1 \cup \eta_2 \cup \eta_3 \cup \eta_4 : W = W_1 \cup W_2 \cup W_3 \cup W_4 \to [0, 1]$ by $\eta_1 : W_1 \to [0, 1]$ with

$$\eta_1 (a\ a\ a) = \frac{1}{3a},$$

$\eta_2 : W_2 \to [0, 1]$ with

$$\eta_2 \left( \begin{pmatrix} a & a \\ a & a \end{pmatrix} \right) = \frac{1}{4a};$$

$\eta_3 : W_3 \to [0, 1]$;

$$\eta_3 (p(x)) = \frac{1}{\deg p(x)};$$

and $\eta_4 : W_4 \to [0, 1]$ by

$$\eta_4 \left( \begin{pmatrix} a \\ a \\ a \\ a \end{pmatrix} \right) = \frac{1}{3(a+a+a+a)}.$$

Clearly each $W_i\eta_i$ is a set fuzzy vector subspace. Thus $W_\eta = (W_1 \cup W_2 \cup W_3 \cup W_4)_{\eta_1 \cup \eta_2 \cup \eta_3 \cup \eta_4}$ is a set fuzzy 4-vector subspace.

Now we proceed onto define set fuzzy n-linear algebra. We see both set fuzzy n-vector space and set fuzzy n-linear algebra are fuzzy equivalent.

**DEFINITION 6.21:** *Let $V = V_1 \cup V_2 \cup \ldots \cup V_n$ be a set n-vector space over the set S. If each $V_i$ is closed under addition then we call V to be a n-set linear algebra over S; $1 \le i \le n$. If in a set n-linear algebra $V = V_1 \cup V_2 \cup \ldots \cup V_n$ we define a n-map $\eta = \eta_1 \cup \eta_2 \cup \ldots \cup \eta_n$ from $V = V_1 \cup V_2 \cup \ldots \cup V_n$ into [0, 1] such that each $V_i\eta_i$ is a set fuzzy linear algebra for i = 1, 2, …, n; then $V_\eta = (V_1 \cup V_2 \cup \ldots \cup V_n)_{\eta_1 \cup \eta_2 \cup \ldots \cup \eta_n}$, is called the set n fuzzy linear algebra.*

We illustrate this by the following example.

***Example 6.39:*** Let

$$V = V_1 \cup V_2 \cup V_3 \cup V_4$$

$$= \left\{ \begin{pmatrix} a & a & a & a \\ a & a & a & a \end{pmatrix} \middle| \, a \in Z^+ \cup \{0\} \right\} \cup \{S^o \times S^o \times S^o \times S^o$$

$\times S^o$ such that $S^o = Z^+ + \cup \{0\}\} \cup$ {all polynomials of degree less than or equal to 5 with coefficients from $S^o$ $= Z^+ \cup \{0\}\} \cup \left\{ \begin{pmatrix} a & b \\ c & d \end{pmatrix} \middle| \, a, b, c, d \text{ is in } S^o \right\}$ is a 4-set linear algebra over $S^o = Z^+ \cup \{0\}$.

Define $\eta = \eta_1 \cup \eta_2 \cup \eta_3 \cup \eta_4 : V = V_1 \cup V_2 \cup V_3 \cup V_4 \to [0, 1]$; such that $\eta_i : V_i \to [0, 1]$, i = 1, 2, 3, 4 with

$$\eta_1 \left( \begin{pmatrix} a & a & a & a \\ a & a & a & a \end{pmatrix} \right) = \begin{cases} \frac{1}{a} & \text{if } a \neq 0 \\ 1 & \text{if } a = 0 \end{cases};$$

$$\eta_2 \, (a\ b\ c\ d\ e) = \begin{cases} \frac{1}{a+b+c+d+e} & \text{if } a+b+c+d+e \neq 0 \\ 1 & \text{if } a+b+c+d+e = 0 \end{cases}$$

$$\eta_3 \, (p(x)) = \frac{1}{\deg p(x)};$$

$\eta_3(p(x)) = 1$ if p(x) is a constant polynomial or 0

and

$$\eta_4 \begin{pmatrix} a & b \\ c & d \end{pmatrix} = \begin{cases} \frac{1}{a+d} & \text{if } a+d \neq 0 \\ 1 & \text{if } a+d = 0 \end{cases}$$

Clearly each $V_i\eta_i$ is a set fuzzy linear algebra, thus $V_\eta = (V_1 \cup V_2 \cup V_3 \cup V_4)_{\eta_1 \cup \eta_2 \cup \eta_3 \cup \eta_4}$ is a set 4-linear algebra.

Now we proceed onto define the notion of semigroup n set fuzzy vector space.

**DEFINITION 6.22:** *Let $V = V_1 \cup V_2 \cup \ldots \cup V_n$ be a semigroup n-set vector space over the semigroup S. If the n-map $\eta = \eta_1 \cup \eta_2 \cup \ldots \cup \eta_n : V = V_1 \cup V_2 \cup \ldots \cup V_n \to [0, 1]$ where $\eta_i : V_i \to [0, 1]$ and $V_i\eta_i$ is a semigroup set fuzzy vector space for each i, i = 1, 2, ..., n then we call $V_\eta = (V_1 \cup V_2 \cup \ldots \cup V_n)_{\eta_1 \cup \eta_2 \cup \ldots \cup \eta_n}$ to be a semigroup n set fuzzy vector space.*

We illustrate this by a simple example.

***Example 6.40:*** Let

$$
\begin{aligned}
V &= V_1 \cup V_2 \cup V_3 \cup V_4 \\
&= \left\{ \begin{pmatrix} a & b \\ c & d \end{pmatrix} \middle| a, b, c, d \in Z_3 = \{0, 1\} \right\} \cup \{Z_3 \times Z_3 \times Z_3 \times Z_3\} \cup \{ (1\ 1\ 0\ 1), (2\ 0\ 2\ 0), (2\ 2\ 0\ 2), (1\ 0\ 1\ 0), (0\ 0\ 0\ 0), (1\ 0\ 1), (2\ 0\ 2), (0\ 0\ 0), (1\ 1\ 1), (2\ 2\ 2)\} \cup \\
&\quad \left\{ \begin{pmatrix} a & a & a \\ a & a & a \end{pmatrix}, \begin{pmatrix} a & a & a & a & a \\ a & a & a & a & a \end{pmatrix} \middle| a \in Z_3 \right\}
\end{aligned}
$$

be a semigroup 4-set vector space over the semigroup $Z_3$.

Define $\eta = \eta_1 \cup \eta_2 \cup \eta_3 \cup \eta_4 : V = V_1 \cup V_2 \cup V_3 \cup V_4 \to [0, 1]$ (with $\eta_i : V_i \to [0, 1]$; i = 1, 2, 3, 4) such that

$$
\eta_1 \left( \begin{pmatrix} a & b \\ c & d \end{pmatrix} \right) = \begin{cases} \dfrac{1}{a+b} & \text{if } a + b \neq 0 \\ 1 & \text{if } a + b = 0 \end{cases}
$$

$$
\eta_2 (a\ b\ c\ d) = \begin{cases} \dfrac{1}{a} & \text{if } a \neq 0 \\ 1 & \text{if } a = 0 \end{cases}
$$

$$\eta_3\,(1\ 1\ 0\ 1) = 0,\ \eta_3\,(2\ 0\ 2\ 0) = 1,$$

$$\eta_3\,(2\ 2\ 0\ 2) = 0,\ \eta_3\,(1\ 0\ 1\ 0) = \frac{1}{2},$$

$$\eta_3\,(0\ 0\ 0\ 0) = 0,\ \eta_3\,(1\ 0\ 1) = \frac{1}{2},$$

$$\eta_3\,(1\ 1\ 1) = 0,\ \eta_3\,(2\ 2\ 2) = 0$$

and

$$\eta_4\left(\begin{pmatrix} a & a & a \\ a & a & a \end{pmatrix}\right) = \begin{cases} \dfrac{1}{6a} & \text{if } a \neq 0 \\ 1 & \text{if } a = 0 \end{cases}$$

$$\eta_4\left(\begin{pmatrix} a & a & a & a & a \\ a & a & a & a & a \end{pmatrix}\right) = \begin{cases} \dfrac{1}{10a} & \text{if } a \neq 0 \\ 1 & \text{if } a = 0 \end{cases}$$

Clearly each $V_i\eta_i$ is a semigroup set fuzzy vector space hence $V_\eta = (V_1 \cup V_2 \cup V_3 \cup V_4)_{\eta_1 \cup \eta_2 \cup \eta_3 \cup \eta_4}$ is a semigroup 4-set fuzzy vector space.

We can as in case of other set fuzzy n-vector subspaces define the notion of semigroup n set fuzzy vector subspaces.

Now we proceed onto define the notion of semigroup n set fuzzy linear algebras.

**DEFINITION 6.23:** *Let $V = V_1 \cup V_2 \cup \ldots \cup V_n$ be a semigroup n-set linear algebra over the semigroup S. $\eta = \eta_1 \cup \eta_2 \cup \ldots \cup \eta_n : V = V_1 \cup V_2 \cup \ldots \cup V_n \to [0, 1]$ where each $\eta_i : V_i \to [0, 1]$ is such that $V_i\ \eta_i$ is a semigroup set fuzzy linear algebra for each i, i = 1, 2, …, n. $V_\eta = (V_1 \cup V_2 \cup \ldots \cup V_n)_{\eta_1 \cup \eta_2 \cup \ldots \cup \eta_n}$ is a semigroup n set fuzzy linear algebra.*

We illustrate this situation by an example.

***Example 6.41:*** Let

$$V = V_1 \cup V_2 \cup V_3 \cup V_4 \cup V_5$$
$$= \{S^o[x] \text{ where } S^o = Z^+ \cup \{0\}\} \cup \{S^o \times S^o \times S^o \times S^o \times S^o \times S^o\} \cup \left\{ \begin{pmatrix} a & b & c \\ d & e & f \\ g & h & i \end{pmatrix} \middle| a, b, c, d, e, f, g, h, i \in S^o \right\} \cup$$

$$\left\{ \begin{pmatrix} a & a & a & a & a \\ a & a & a & a & a \\ a & a & a & a & a \end{pmatrix} \middle| a \in S^o \right\} \cup$$

$$\left\{ \begin{pmatrix} a & b & c & d \\ 0 & e & f & g \\ 0 & 0 & h & i \\ 0 & 0 & 0 & j \end{pmatrix} \middle| a, b, c, d, e, f, g, h, i, j \in S^o \right\}$$

be a semigroup 5-set linear algebra over the semigroup $S^o$.

Define $\eta = \eta_1 \cup \eta_2 \cup \eta_3 \cup \eta_4 \cup \eta_5 : V = V_1 \cup V_2 \cup V_3 \cup V_4 \cup V_5 \to [0, 1]$ by $\eta_i : V_i \to [0, 1]$ ; $1 \leq i \leq 5$ by

$$\eta_1(p(x)) = \begin{cases} \dfrac{1}{\deg p(x)} & \\ 1 & \text{if } p(x) = \text{constant} \end{cases}$$

$$\eta_2((a_1, a_2, a_3, a_4, a_5, a_6)) = \begin{cases} \dfrac{1}{a_1} & \text{if } a_1 \neq 0 \\ 1 & \text{if } a_1 = 0 \end{cases}$$

$$\eta_3\left(\begin{pmatrix} a & b & c \\ d & e & f \\ g & h & i \end{pmatrix}\right) = \begin{cases} \dfrac{1}{a+e+i} & \text{if } a+e+i \neq 0 \\ 1 & \text{if } a+e+i = 0 \end{cases}$$

$$\eta_4\left(\begin{pmatrix} a & a & a & a & a \\ a & a & a & a & a \\ a & a & a & a & a \end{pmatrix}\right) = \begin{cases} \dfrac{1}{15a} & \text{if } a \neq 0 \\ 1 & \text{if } a = 0 \end{cases}$$

and

$$\eta_5\left(\begin{pmatrix} a & b & c & d \\ 0 & e & f & g \\ 0 & 0 & h & i \\ 0 & 0 & 0 & j \end{pmatrix}\right) = \begin{cases} \dfrac{1}{a+e+h+j} & \text{if } a+e+h+j \neq 0 \\ 1 & \text{if } a+e+h+j = 0 \end{cases}$$

Thus $V\eta$ is a semigroup 5 set fuzzy linear algebra.

Now we proceed onto define the notion of semigroup n-set fuzzy linear subalgebra.

**DEFINITION 6.24:** *Let $V = V_1 \cup V_2 \cup \ldots \cup V_n$ be a semigroup n-set linear algebra over the semigroup S. Let $W = W_1 \cup W_2 \cup \ldots \cup W_n \subseteq V_1 \cup \ldots \cup V_n = V$ be a n-subset of V such that W is a semigroup n set linear subalgebra over the semigroup S.*

*Let $\eta = \eta_1 \cup \eta_2 \cup \ldots \cup \eta_n : W = W_1 \cup W_2 \cup \ldots \cup W_n \to [0, 1]$ is such that each $W_i\eta_i$ is a semigroup set fuzzy linear subalgebra for i = 1, 2, …, n; then we call $W_\eta = (W_1 \cup W_2 \cup \ldots \cup W_n)_{\eta_1 \cup \eta_2 \cup \ldots \cup \eta_n}$ to be a semigroup n set fuzzy linear subalgebra.*

We illustrate this by an example.

***Example 6.42:*** Let

$$V = V_1 \cup V_2 \cup V_3 \cup V_4$$

$$= \{Z_{24} \times Z_{24} \times Z_{24} \times Z_{24}\} \cup \{Z_{24}[x]\} \cup$$

$$\left\{ \begin{pmatrix} a & b & e \\ c & d & f \end{pmatrix} \middle| \; a, b, c, d, e, f \in Z_{24} \right\} \cup \left\{ \begin{pmatrix} a & a \\ a & a \\ a & a \\ a & a \\ a & a \end{pmatrix} \middle| \; a \in Z_{24} \right\}$$

be a semigroup 4 set linear algebra over the semigroup $Z_{24}$.
Let

$W = W_1 \cup W_2 \cup W_3 \cup W_4$
$= \{S_1 \times S_1 \times S_1 \times S_1$ where $S_1 = \{0, 2, 4, 6, 8, \ldots, 22\}\}$
$\cup$ $\{Z_{24}'$ [x] polynomials of only even degree with coefficients from $Z_{24}\} \cup$

$$\left\{ \begin{pmatrix} a & a & a \\ a & a & a \end{pmatrix} \middle| \; a \in Z_{24} \right\} \cup \left\{ \begin{pmatrix} a & a \\ a & a \\ a & a \\ a & a \end{pmatrix} \middle| \; a \in S_1 \right\}$$

$\subseteq V_1 \cup V_2 \cup V_3 \cup V_4$
$= V$

is a semigroup 4-set linear subalgebra of V over the same semigroup $Z_{24}$. Let $\eta = \eta_1 \cup \eta_2 \cup \eta_3 \cup \eta_4 : W = W_1 \cup W_2 \cup W_3 \cup W_4 \to [0, 1]$ where each $\eta_i : W_i \to [0\ 1]$ is such that $W_i\eta_i$ is a semigroup set fuzzy linear subalgebra for i = 1, 2, 3, 4. Now $W_\eta = (W_1 \cup W_2 \cup W_3 \cup W_4)_{\eta_1 \cup \eta_2 \cup \eta_3 \cup \eta_4}$ is a semigroup 4-set fuzzy linear subalgebra.

Now we proceed onto define the notion of group n set fuzzy vector space.

**DEFINITION 6.25:** *Let $V = V_1 \cup V_2 \cup \ldots \cup V_n$ be a group n set vector space over the group G. Let $\eta = \eta_1 \cup \eta_2 \cup \ldots \cup \eta_n$. $V = V_1 \cup V_2 \cup \ldots \cup V_n \to [0, 1]$ such that $\eta_i : V_i \to$ with $V_i\eta_i$ is a* *group set fuzzy vector space for each i = 1, 2, … , n;* $V_\eta = (V_1 \cup V_2 \cup \ldots \cup V_n)_{\eta_1 \cup \eta_2 \cup \ldots \cup \eta_n}$ *is a group n-set fuzzy vector space.*

We illustrate this an example.

***Example 6.43:*** Let

$$V = V_1 \cup V_2 \cup V_3 \cup V_4 \cup V_5$$

$$= \left\{ \begin{pmatrix} a_1 & a_2 & a_3 \\ a_4 & a_5 & a_6 \end{pmatrix} \middle| \, a_i \in Z_{11}; 1 \leq i \leq 6 \right\} \cup \{Z_{11}[x]\} \cup$$

$$\{Z_{11} \times Z_{11} \times Z_{11} \times Z_{11} \times Z_{11}\} \cup$$

$$\left\{ \begin{pmatrix} a_1 & a_2 \\ a_3 & a_4 \\ a_5 & a_6 \\ a_7 & a_8 \\ a_9 & a_{10} \end{pmatrix} \middle| \, a_i \in Z_{11}; 1 \leq i \leq 10 \right\} \cup$$

$$\left\{ \begin{pmatrix} a & a & a \\ a & a & a \\ a & a & a \end{pmatrix}, \begin{pmatrix} a & a \\ a & a \end{pmatrix}, \begin{pmatrix} a & a & a & a & a \\ a & a & a & a & a \\ a & a & a & a & a \\ a & a & a & a & a \end{pmatrix} \middle| \, a \in Z_{11} \right\}$$

be a group 5 set vector space over the group $Z_{11}$.
Take

$$W = W_1 \cup W_2 \cup W_3 \cup W_4 \cup W_5$$

$$= \left\{ \begin{pmatrix} a & a & a \\ a & a & a \end{pmatrix} \middle| \, a \cup Z_{11} \right\} \cup \{Z_{11}[x]; \text{all polynomials of}$$

even degree with coefficients from $Z_{11}\} \cup \{(a\ a\ a\ a\ a) \mid a \in Z_{11}\}$

$$\cup \left\{ \left.\begin{pmatrix} a & a \\ a & a \\ a & a \\ a & a \\ a & a \end{pmatrix} \right| a \in Z_{11} \right\} \cup \left\{ \begin{pmatrix} a & a & a \\ a & a & a \\ a & a & a \end{pmatrix}, \left.\begin{pmatrix} a & a \\ a & a \end{pmatrix} \right| a \in Z_{11} \right\}$$

$$\subseteq \quad V_1 \cup V_2 \cup V_3 \cup V_4 \cup V_5$$
$$= \quad V,$$

to be a group set 5-vector subspace over the group $Z_{11}$. Define $\eta = \eta_1 \cup \eta_2 \cup \eta_3 \cup \eta_4 \cup \eta_5 : W = W_1 \cup W_2 \cup W_3 \cup W_4 \cup W_5 \to [0, 1]$ such that $\eta_i : W_i \to [0, 1]$ with $W_i\eta_i$ a group set fuzzy vector subspace for each i = 1, 2, …, 5. Clearly $W_\eta = (W_1 \cup W_2 \cup \ldots \cup W_5)_{\eta_1 \cup \eta_2 \cup \ldots \cup \eta_5}$ is a group set fuzzy 5-vector subspace.

Now we proceed onto define the notion of group n set fuzzy linear algebra. It is important to mention here that the notion of group n set fuzzy vector space and group n-set fuzzy linear algebra are fuzzy equivalent.

**DEFINITION 6.26:** *Let $V = V_1 \cup V_2 \cup \ldots \cup V_n$ be a group n set linear algebra over the group G. Let $\eta = \eta_1 \cup \eta_2 \cup \ldots \cup \eta_n : V \to [0, 1]$ where $\eta_i : V_i \to [0, 1]$ such that $V_i\eta_i$ is a group set fuzzy linear algebra for each i = 1, 2, …, n, then we call $V_\eta = (V_1 \cup V_2 \cup \ldots \cup V_n)_{\eta_1 \cup \eta_2 \cup \ldots \cup \eta_n}$ to be a group n set fuzzy linear algebra.*

Now we proceed onto indicate the applications of these new algebraic structures before which we give few differences.

1. If we use a set vector space over the set [0, 1] working can easily be carried out for we do not demand 1 + 1 = 0. We can use them in fuzzy models.

2. If we use n set vector space then we can use several sets to model the problem that too when the elements are not compatible with any of the operation, that is even the closure axiom under any operation cannot be expected.

3. When we give algebraic structures with many abstract operations on them, the use of these algebraic structures by students/ researchers other than mathematicians is absent. But as the world at present is in a computer age the mathematical models should be constructed with least abstractness. These algebraic structures given in this book will certainly cater to this need.

Now we give a sketch where these structures can be applied.

Applications to Markov Chains

Suppose we have a n set of states which is a time dependent model which has many transition matrices representing the transition possibilities. Suppose we have k such transition matrices given by $T_1, \ldots, T_k$, where each $T_i$ is a $n \times n$ matrix, i = 1, 2, …, k. Suppose $G(T_i)$ denotes the graph having vertex set {1, 2, …, n}. The arc set $(i, j) : t_{ij} > 0$. For each $G(T_i)$, we get the four classes

$C_1^i = \{n_{i_1}, \ldots, n_{i_k}\}$ essential class

$C_2^i = \{n_{j_1}, n_{j_2}, \ldots, n_{j_p}\}$ inessential class (self communicating)

$C_3^i = \{n_{r_1}, \ldots, n_{r_s}\}$ inessential non self communicating

$C_4^i = \{n_{m_1}, \ldots, n_{m_t}\}$ essential class (non self communicating)

such that

(1) $k \geq 1, p \geq 1, s \geq 1$ and $t \geq 1$.

(2) $C^i_p \cap C^i_j = \phi$ if $p \neq j$

(3) $C^i_1 \cup C^i_2 \cup C^i_3 \cup C^i_4 = \{1, 2, ..., n\}$.

Now we see $G(T_i)$, varies; for suppose $V = \{T_1, \ldots, T_k\}$ then we can take V as a set vector space over [0, 1] or as a semigroup set vector space over [0, 1].

Further the transition matrices vary from expert to expert for the $i^{th}$ expert may start his study from the $j^{th}$ state and say a $k^{th}$ expert ($i \neq k$) may start his study from the $p^{th}$ state and so on. Thus the very transition matrix may be different. One can use the semigroup set vector space to represent the collection of all transition matrices and work with them. If

$V_o = \{T_1, \ldots, T_n\}$ is a set vector space over the set [0, 1].

$V_1 = \{T^2_1, ..., T^2_n\}$ is again a set vector space over the set [0, 1], …,

$V_k = \{T^{k-1}_1, ..., T^{k-2}_n\}$ is a set vector space.

The best limiting case would the best chosen set vector space.

It is pertinent at this juncture to mention that the entries of $T_i$ can be any value lying between [0, 1] further even the assumption the row sums adds to one can be chosen.

Their is no loss of generality in assuming some of the columns or rows in $T_i$ have only zero entries.

These types of applications like studying the behaviour of children, war-peace situation and the related behaviour of military personals like stress, strain, relaxation, taking to drugs, taking to drink, taking to women, taking leave on basis of false ailments, taking medical leave and so on or like the study of students tensions during exams.

Another practical use of fuzzy Markov chains can be in statistical software testing.

It can be used in web monitoring, testing or counting. Thus when we use in web monitoring or testing and if the web pages or sometimes 'n' entries are taken as the set of states one has to keep in mind the states vary so when we use the set fuzzy n-vector space or semigroup fuzzy n-vector space or group fuzzy n-vector space it would be easy to handle the bulk problems at a time!

- These algebraic structures can be used by cryptologists and coding theorists.
- They can also be used in fuzzy models, by social scientists and doctors.
- The set fuzzy vector spaces in which the elements form $n_i \times n_i$ matrices can be used by computer scientists.

**Chapter Seven**

# SUGGESTED PROBLEMS

In this book we suggest over 300 problems so that the reader becomes familiar with the new notion of set vector spaces and their particularization. All vectors spaces and semivector spaces are set vector spaces, thus set vector spaces are the most generalized concept.

1. Let V = {–1, 0, 2, 3, …, ∞} and S = {1, 2, …, ∞}. Is V a set vector over S. Justify your claim.

2. Let V = R, the set of reals. Let S = {1, 2, …, ∞}. Is R a set vector space over the set S?

3. Let V = R, the set of reals. Let S = $Q^+$ (the set of positive rationals). What is the difference between the set vector space defined in problem (2) and in this problem.

4. Let V = {1, 2, ..., ∞}. Can V be a set vector space over the set S = $Q^+$?

5. Let V = {1, 2, …, ∞}. Can V be a set vector space over the set S = {–1, 1}? Justify your claim.

6. Is Z the set of integers a set vector space over the set S = {–1, 1}?

7. Is Z the set of integers a set vector space over the set S = {1, 1/2, 1/3, …, 1/n}?

8. Prove the set of positive rationals is not a set vector space over the set of positive reals.

9. Let V = {0, 2, 4, …, ∞}; prove V is a set vector space over the set S = {0, 1, 2, …, ∞} = $Z^+ \cup \{0\}$.

10. Can V = $\{2, 2^2, 2^3, \ldots\}$ be a set vector space over the set S = $Z^+$?

11. Prove V = {0, 3, 6, …, ∞} is a set vector space over the set S = $Z^+$ = {1, 2, …, ∞}.

12. Prove V = {1, 3, 6, …, ∞} is not a set vector space over the set S = $Z^+$ = {1, 2, …, ∞}.

13. Prove V = {1, p, 2p, …, ∞} is not a set vector space over the set S = {1, 2, ..., ∞}.

14. Can we say for V to be a set vector space over S, S must be a proper subset of V?

15. Is it true if $S \subset V$; V is a set vector space over S?

16. Can $M = \left\{\left.\begin{pmatrix} a & b & c \\ d & e & f \end{pmatrix}\right| a, b, c, d, e, f \in Z^+\right\}$ be a set vector space over $Q^+$ (the set of positive rationals)?

17. Prove $M = \left\{\left.\begin{pmatrix} a & b \\ c & d \end{pmatrix}\right| a, b, c, d \in Z^+\right\}$ is a set vector space over the set S = $Z^+$.

18. Can the set $V = \left\{\begin{pmatrix} 1 & 2 \\ 3 & 4 \end{pmatrix}, \begin{pmatrix} 0 & 1 \\ 2 & 3 \end{pmatrix}, \begin{pmatrix} 5 & 7 \\ 8 & 10 \end{pmatrix}, \begin{pmatrix} 1 & 2 \\ 3 & 0 \end{pmatrix}\right\}$ be a set vector space over the set S = {0, 1, 2, …, ∞}?

19. Find a set of linearly independent elements in the set vector space V = {1, 2, ..., ∞} defined over the set S = {2, 4, …, ∞}.

20. Find a set basis for the set vector space V = {1, 2, ..., ∞} over the set S = {2, 4, …, ∞}.

21. Is the set basis for the set vector space $V = Z^+ \times Z^+$ over the set $S = Z^+$ unique?

22. Does there exist a set vector space V over the set S which has more than one set basis?

23. Let $V = \left\{ \begin{pmatrix} a & b \\ c & d \end{pmatrix} \middle| \; a, b, c, d \in Z^+ \right\}$ be a set vector space over the set $S = Z^+$. Find a set basis for V? Is it unique?

24. Find a set basis for the set vector space V = {–1, 0, 1, 2, 3, …, ∞} over the set S = {0, 1}. Is the set basis unique? Is it finite dimensional?

25. Give an example of a set vector space of dimension 7?

26. Give an example of a set vector space V over the set S which is also a set vector space S over V.

27. Give an example of a set vector space V over S which is not a set vector space of S over V (i.e., show if V is a set vector space over S then S is not a set vector space over V).

28. Let $V = Q^+ \times Z^+ \times 2Z^+ \times R^+$ be a set vector space over the set S = {1, 2, …, ∞}. Find a basis for V? What is the dimension of V? Does V have a unique set basis? Give a finite set independent set in V.

29. Characterize those set vector spaces V on the set S of finite dimension such that $V \cong S \times \ldots \times S$.

30. Give an example of set vector spaces V of set finite dimension over the set S which cannot be written as $V \cong S \times \ldots \times S$.

31. Let $V = Z^+ \cup \{0\} = \{0, 1, \ldots, \infty\}$ be a set vector space over the set $S = \{0, 1, 3, 5, 7, \ldots\}$. Is $V \cong \underbrace{S \times \ldots \times S}_{n-\text{times}}$ ? , (n, finite).

32. Give some interesting properties of set vector spaces.

33. Does there exist a set vector space V over the set S which has no set vector subspaces?

34. Let $V = Z^+$ be a set vector space over $S = \{2, 4, \ldots, \infty\}$. How many set vector subspaces does V have?

35. Find at least two set vector subspaces of $V = \left\{\begin{pmatrix} a & b & c \\ d & e & f \end{pmatrix}\right|$ a, b, c, d, e and $f \in Z^+\}$; where V is a set vector space over the set $S = \{1, 2, \ldots, \infty\}$.

36. Give an example of a set vector space V over the set S which has no set vector subspace.

37. Characterize those set vector spaces V over the set S which has no set vector subspaces.

38. Does there exist for any set vector space V over the set S, a set vector subspace over S such that it is also a subset of S?

39. Find a generating set of the set vector space $V = Z^+ \times 2Z^+$ over $2Z^+ = S = \{2, 4, \ldots, \infty\}$.

40. Find a set basis for V given in problem 39. Does the V have a unique set basis?

41. How many generating sets does the set vector space V = $\left\{\begin{pmatrix} a & b \\ c & d \end{pmatrix} \middle| a, b, c, d \in Z^+ \right\}$ have, V defined over the set S = $Z^+$?

42. Find at least two set vector subspaces of V = $Z^+$ over the set S = $2Z^+$ = {2, 4, …, ∞}, which have non empty intersection.

43. Let V = $Z^+ \times Z^+ \times Z^+$ be a set vector space over the set S = $Z^+$. Find three set vector subspaces $W_1$, $W_2$ and $W_3$ of V which have a non empty intersection. If W = $W_1 \cap W_2 \cap W_3 \subset V$; is W a set vector subspace of V?

44. Let V be a set vector subspace over the set S. Suppose $W_1$, $W_2$, …, $W_n$ be set vector subspaces of V. Will W = $\bigcap_{i=1}^{n} W_i \subseteq V$ be always a proper set vector subspace of V?

45. Let V = $\left\{\begin{pmatrix} a & b & c \\ d & e & f \end{pmatrix} \middle|\right.$ a, b, c, d, e and f $\in Z^+$} be a set vector space over S. Find 3 subset vector subspaces over subsets of S. Does there exists $W_1$, $W_2$ and $W_3$, subset vector subspaces of V defined over the subsets $T_1$, $T_2$ and $T_3$ of S such that

(1) $\bigcap_{i=1}^{3} W_i = \phi$ but $\bigcap_{i=1}^{3} T_i \neq \phi$ or

(2) $\bigcap_{i=1}^{3} W_i \neq \phi$ but $\bigcap_{i=1}^{3} T_i = \phi$ or

(3) $W = \bigcap_{i=1}^{3} W_i \neq \phi$ and $T = \bigcap_{i=1}^{3} T_i \neq \phi$ and W a subset vector subspace over T.

46. Can V = {–1, 0, 2, 4, …, ∞} be set vector space over the set S = {1}?

i. If V is a set vector space over S; what is the dimension of V?

ii. Find a generating set G of V over S. Is G = V?

iii. Can V have subset vector subspaces over subsets of S?

iv. Does V contain set vector subspaces?

v. Is every singleton in V a set vector subspace of V?

47. Let $V = Z^+ \times Z^+ \times Z^+ \times Z^+$ be a set vector space over the set $S = Z^+$. Does V contain subsets which are subset vector subspaces over the subsets; $T_1 = \{2, 2^2, \ldots\}$, $T_2 = \{3, 3^2, \ldots\}$, $T_3 = \{p, p^2, \ldots, p^n \ldots$, p a prime$\}$ and $T_4 = \{n, n^2, \ldots;$ n a non prime$\}$ in S?

48. Let $V = 2Z^+ \times 3Z^+ \times 5Z^+ \times 2Z^+$ be a set vector space over the set $S = Z^+$. Find a set basis for V. Is the basis of V unique? Can V have more than one generating subset? Justify you claim. Find at least 4 subset vector subspaces of V such that their intersection is empty!

49. Let $V = 3Z^+$ and $W = 2Z^+$ be two set vector spaces over the set $Z^+$. Prove the collection of all set linear transformation $L_S$ (V, W) from V to W is a set vector space over the set $Z^+$.

50. Find a generating subset of $L_S$ (V, W). What is the set dimension of $L_S$ (V, W) given in problem 49? What is the dimension of $V = 3Z^+$ over $Z^+$? What is the dimension of $W = 2Z^+$ over $Z^+$? Find generating subsets of V and W over $Z^+$. Does their exist two subset vector subspaces of V and W which are defined over the same subset of $Z^+$? If so find a set linear transformation between them.

51. Give an example of a finite set linear algebra.

52. Give an example of an infinite set linear algebra.

53. Let $V = Z^+$ be the set under addition. V a set linear algebra over the set $S = 3\,Z^+$. What is the set dimension of V?

Can V have set linear subalgebra defined over S?

Can V have subset linear subalgebra defined over the set $T = 3^2Z^+$?

54. Prove $V = \left\{ \begin{pmatrix} a & b & c \\ d & e & f \end{pmatrix} \middle| a, b, c, d, e \text{ and } f \in Z^+ \cup \{0\} \right\}$ is a set linear algebra over the set $S = Z^+$. (V, +), a set with closure operation. Find a subset linear subalgebra over the subset $T = 3Z^+$. Give an example of a set linear subalgebra of V.

Find three subset linear subalgebras of V such that their intersection is non empty.

Suppose $W = \left\{ \begin{pmatrix} a & 0 & e \\ d & 0 & f \end{pmatrix} \middle| a, d, e \text{ and } f \text{ are in } Z^+ \right\} \subseteq V$; is W a subset linear subalgebra over the subset $T = 2\,Z^+$? If W is a subset linear subalgebra find its set dimension over T.

55. Does there exist a set linear algebra which has no set linear subalgebra?

56. Does their exist a set linear algebra which has no subset linear subalgebra?

57. Give an example of a set linear algebra which has only a finite number of subset linear subalgebras.

58. Let $A_1 = \{3\,Z^+$ under the operation $+\}$, be a set linear algebra over the set $S = Z^+$. Let $A_2 = \{5\,Z^+$ under the operation $+\}$ be a set linear algebra over the set $Z^+$. Let T be a set linear transformation from $A_1$ to $A_2$. Is the set of all set linear transformations from $A_1$ to $A_2$ a set linear algebra over $Z^+$? Give a set basis for $A_1$ and $A_2$ over $Z^+$.

59. Let V be the set vector space over the set S. Prove the set of all set linear operators from V to V is a set linear vector space over the set S.

60. Let V = $Z^+ \times 3Z^+ \times 5Z^+$ be a set vector space over $Z^+$. Is T(x, y, z) = (5x, 6y, 10z) for all x, y, z ∈ V a set linear operator on V. What is the set dimension of the set vector space of all set linear operators of V on the set $Z^+$.

61. Let V = $3Z^+ \times Z^+ \times Q^+$ be a set vector space over the set S = $3Z^+$. What is the dimension of V over S? Is the basis of V unique over S? Is T(x, y, z) = (x, y, 2z) a set linear operator? What is the dimension of the set of all set linear operators on V? Find a generating set for the set vector space of set linear operators on the set $Z^+$.

62. Let V = $\left\{ \begin{pmatrix} a & b & c \\ d & e & f \end{pmatrix} \middle| a, b, c, d, e, f \in Z^+ \right\}$ be a set vector space over the set S = $Z^+$. Is $T\left\{\begin{pmatrix} a & b & c \\ d & e & f \end{pmatrix}\right\}$ = a + b + c + d + e + f a set linear operator on V? Justify your claim.

Suppose $T\left\{\begin{pmatrix} a & b & c \\ d & e & f \end{pmatrix}\right\} = \begin{pmatrix} a+b & b+c & d \\ e+b & f & e+f \end{pmatrix}$; is T a set linear operator on V? What is the dimension of V over $Z^+$?

Can $W = \left\{ \begin{pmatrix} 2a & 3b & 5c \\ 2d & 5e & 3f \end{pmatrix} \middle| a, b, c, d, e, f \in Z^+ \right\}$ be a set vector subspace of V over $Z^+$?

63. Let V = {$Z_8^+$ [x] be the set of all polynomials of degree less than or equal to 8 with coefficients from the set $Z^+$}. Is V is a set vector space over $Z^+$? What is the set dimension of V over $Z^+$?

If V is endowed with polynomial addition, can V be a set linear algebra over $Z^+$. What is the dimension of V over $Z^+$ as a set linear algebra?

Suppose W = {All polynomials of degree less than or equal to four} $\subseteq$ V; will W be a set vector subspace of V over $Z^+$?

64. Give an example of a set linear algebra of dimension 6 which has a subset linear subalgebra of dimension 3.

65. Give an example of a set linear algebra of dimension 4 which has a set linear subalgebra of dimension 2.

66. Give an example of an infinite dimensional set linear algebra which has both a subset linear subalgebra and a set linear subalgebra.

67. Does there exist a set linear algebra of infinite dimension which has no subset linear subalgebra and no set linear subalgebra?

68. Does there exist a set linear algebra of infinite dimension which has no subset linear subalgebra?

69. Does there exists a set vector space which does not have a unique set basis?

70. Does there exists a set linear algebra with a unique set basis?

71. Give any interesting result on set linear algebra.

72. Enumerate at least 3 differences between a set vector space and the usual vector space.

73. Define a set linear operator T on a set linear algebra V = $2Z^+$ (set under addition) over the set S = $Z^+$ which is invertible. Give a set linear operator U on V which is non invertible.

74. Let $V = \{\{0\} \cup Z^+\} \times \{Z^+ \cup \{0\}\}$ and $U = \left\{\begin{pmatrix} a & b \\ 0 & c \end{pmatrix} \middle| a, b, c \in Z^+ \cup \{0\}\right\}$ be two set vector spaces defined over $Z^+$. Can we define a set linear transformation from V to U which is invertible?

Does their exists a set linear operator from U to V which is invertible? Find a set basis for the set vector space of all set linear operators from V to U. Is the set basis unique? What is the dimension of this set vector space?

75. Does there exists a set vector space V defined over the set S such that every subset of V is a set vector subspace of V over S?

76. Does there exists a set vector space V over the set S such that for every subset W of V there exists a subset T of S such that W is a subset vector subspace over T?

77. Does there exists a set vector space V over the set S such that V has a subset W which is a subset vector space over every subset of S?

78. Does there exists a set vector space V over the set S such that every subset of V is a subset vector subspace over a unique subset T of S?

79. Let $V = Z^+ \times Z^+ \times Z^+$ be a set vector space over $Z^+ = S$. Show that the set of all set linear functionals from V to S is a set vector space over S. What is the dimension of this set vector space? Is the set basis unique?

80. Let $V = \left\{\begin{pmatrix} x_1 & x_2 & x_3 \\ x_4 & x_5 & x_6 \end{pmatrix} \middle| x_i \in Z^+ \cup \{0\}; 1 \le i \le 6\right\}$ be the set vector space over the set $S = Z^+ \cup \{0\}$. Define $f_i : V \to S$ by $f_i\left[\begin{pmatrix} x_1 & x_2 & x_3 \\ x_4 & x_5 & x_6 \end{pmatrix}\right] = x_i$ for i = 1, 2, …, 6.

Prove $f_i$'s are set linear functionals of V. Is the set $\{f_1, \ldots, f_6\}$ a set linearly independent set?

81. Find a basis of V for the V given in the above problem. Is it unique?

82. Let $V = Z^+ \cup \{0\} \times Z^+ \cup \{0\}$ a set vector space with zero defined over the set $S = 2Z^+ \cup \{0\}$. Let T be a set linear operator on V given by $T(x, y) = (x - y, x + y)$. Find the set nullity of T! Is it a set vector subspace of V over S?

83. If V is a set vector space with zero can we say every set linear operator on V which has a non trivial nullity T is a set vector subspace of V?

84. If nullity $T \neq \phi$ can we say T is non invertible?

85. Let $V = 3Z^+ \times 2Z^+ \times Z^+$ be a set vector space over the set $S = Z^+$. Is the set $\{(6, 4, 10), (3, 8, 11), (9, 2, 11)\}$ a set linearly independent subset of V?

86. Let $V = (3Z^+ \cup \{0\}) \times (5Z^+ \cup \{0\}) \times (2Z^+ \cup \{0\})$ be a set vector space over the set $S = Z^+ \cup \{0\}$. Is the set $\{(3, 0, 0)$ $(0, 5, 0), (0, 0, 2)\}$ a set basis of V over S?

    If V has a set basis over S, is it unique? What is the set dimension of V?

87. Let $V = \left\{ \begin{pmatrix} a & b \\ c & d \end{pmatrix} \middle| a, b, c, d \in Z^+ \cup \{0\} \right\}$, V a set vector space over $Z^+ \cup \{0\}$. Is the set $G = \left\{ \begin{pmatrix} 1 & 0 \\ 0 & 0 \end{pmatrix}, \begin{pmatrix} 0 & 1 \\ 0 & 0 \end{pmatrix}, \begin{pmatrix} 0 & 0 \\ 1 & 0 \end{pmatrix}, \begin{pmatrix} 0 & 0 \\ 0 & 1 \end{pmatrix} \right\}$ a set basis of V?

    Is the set G a set linearly dependent subset of V? What is the dimension of V?

88. If V is a set vector space over the set S. Can we find set vector subspaces $W_1, \ldots, W_n$ of V such that $V = W_1 \cup W_2 \cup \ldots \cup W_n$ such that $W_i \cap W_j = \phi$, $i \neq j$, $1 \leq i, j \leq n$?

89. If A is a set linear algebra over the set $S = Z^+ \cup \{0\}$ where
$$A = \left\{ \begin{pmatrix} a & b \\ c & d \end{pmatrix} \middle| a, b, c, d \in Z^+ \cup \{0\} \right\}.$$

Can $\left\{ \begin{pmatrix} 1 & 0 \\ 0 & 0 \end{pmatrix}, \begin{pmatrix} 0 & 1 \\ 0 & 0 \end{pmatrix}, \begin{pmatrix} 0 & 0 \\ 1 & 0 \end{pmatrix}, \begin{pmatrix} 0 & 0 \\ 0 & 1 \end{pmatrix} \right\}$ generate V? What is the set dimension of this set linear algebra over the set S?

90. Let $V = Z^+ \cup \{0\} \times 2Z^+ \cup \{0\} \times 3Z^+ \cup \{0\}$ be a set vector space over the set $S = Z^+ \cup \{0\}$.

Can V be represented uniquely as $W_1 \cup W_2 \cup W_3$ where $W_i$'s are set vector subspace of V; i = 1, 2, 3 defined over S; with $W_i \cap W_j = \{0\}$, if $i \neq j$; $1 \leq i, j \leq 3$?

91. Give examples of set vector spaces V defined over the set S such that $V = \bigcup_{i=1}^{n} W_i$ ; $W_i$ set vector subspaces of V. ($W_i \cap W_j = \{0\}$ or $\phi$, $i \neq j$; $1 \leq i, j \leq n$).

92. Give an example of a set linear algebra A defined over the set A such that

$A = A_1 + A_2 + \ldots + A_n$ where $A_i$'s are set linear subalgebras of A defined over the set A; i = 1, 2, …, n with $A_i \cap A_j = 0$ or $\phi$ if $i \neq j$; $1 \leq i, j \leq n$.

93. Is every set linear algebra A representable as $A = A_1 + \ldots + A_n$; $A_i$ set linear subalgebra of A ; i = 1, 2, …, n. $A_i \cap A_j = 0$ or $\phi$; $1 \leq i, j \leq n$?

94. Let V be a set vector space over the set S. Let $W_1, W_2, ..., W_n$ be subset vector subspaces over the subsets $T_1, T_2, \ldots, T_n$ respectively of S.

Can $V = \bigcup_{i=1}^{n} W_i$ and $S = \bigcup_{i=1}^{n} T_i$ ?

95. In problem 94 can $V = \bigcup_{i=1}^{n} W_i$ and $S \neq \bigcup_{i=1}^{n} T_i$ ?

96. Give an example of a set vector space in which conditions of problem 94 is true. (Does such an example exist?)

97. Can $V = Z_{15} = \{0, 1, 2, \ldots, 14\}$ integers modulo 15 be a set vector space over the set $S = \{0, 1\}$?

98. What is dimension of V? Can V have more than one set basis (V given in problem 97)?

99. Let $V = Z_n = \{0, 1, 2, \ldots, n\}$ where $n = p_1 p_2 \ldots p_t$ ; $p_i$'s are distinct primes be a set vector space over the set $S = \{1, 0\}$.

What is the set dimension of V over S. How many set basis has V?

| | | | | | |
|---|---|---|---|---|---|
| Study when | n | = | 2.3.5 | = | 30 |
| | n | = | 3.5.11.13 | = | 2145 |
| and | n | = | 2.5.17 | = | 170. |

100. Let $V = Z_p = \{0, 1, 2, \ldots, p - 1\}$, p a prime, $Z_p$ a set vector space over the set $S = \{1\}$. What is the set dimension of V?

101. Let $V = Z_7 = \{0, 1, 2, \ldots, 6\}$ be a set vector space over the set $S = \{1\}$. What is the set dimension of V? Is the set basis for V unique?

102. Let $V = Z_{11} = \{0, 1, 2, \ldots, 10\}$ be a set vector space over the set $S = \{1, 0\}$. What is the set dimension of V? Is the set basis unique? Can V have set vector subspaces?

103. Let $V = Z_5 = \{0, 1, 2, 3, 4\}$ be a set vector space over the set $S = \{0, 1\}$. What is set dimension of V over S? Is the set basis unique? Is $W = \{2, 3\}$ a set vector subspace of V over $\{1, 0\}$? Is $W = \{2, 3\}$ a subset vector subspace of V over $\{1\}$?

104. Let $V = Z_5 \times Z_3$ be a set vector space over the set $S = \{(0, 1), (0, 0) (1, 0), (1, 1)\}$. What is the set dimension of V over S? Find a set basis of V. Is the set basis of V unique?

105. Let $V = Z_7 \times Z_5 \times Z_3$ be a set vector space over the set $S = \{0, 1\}$. Does their exists set vector subspaces $W_1, \ldots, W_n$ such that $V = \bigcup_{i=1}^{n} W_i$ ?

106. Prove if V is a finite set which is a set vector space over another finite set S; the set dimension of V is finite.

107. If V is an infinite set which is a set vector space over a set S. Can V have a finite set dimension? Justify your claim !

108. Let $V = Z_7 \times Z_{12}$ be a set vector space over the set $S = \{0, 1\}$. What is the dimension of V? Can V have subset vector subspaces? Can V have set vector subspaces? Is the set basis of V unique over S?

109. Let $V = Z_4 \times Z_{10} \times Z_{16} \times Z_2$ be the set vector space over the $S = \{0, 1\}$. What is the set dimension of V? Can we write V as a union of set vector subspaces? Can V have more than one subset vector subspace?

110. Let $V = Z_7 [x]$ be a set vector space defined over the set $\{0, 1\}$. What is the set dimension of V?

111. If $V = Z_6 [x]$ is a set linear algebra over the set $S = \{0, 1\}$. Find a set basis of V. Find a set linear subalgebra of V?

112. If $V = Z_5^2[x]$ = {all polynomials of degree less than or equal to two} be a set linear algebra over the set $S = \{0, 1\}$. What is the set dimension of V? Is the set basis of V unique? Can V have a subset linear subalgebra over the subset {1}? Can V have set linear subalgebra?

113. Give some interesting properties of semigroup vector spaces.

114. Obtain some properties enjoyed by semigroup linear algebra and not by semigroup vector spaces.

115. Find pseudo semigroup vector subspaces of the semigroup linear algebra V = $\left\{ \begin{pmatrix} a_1 & a_2 & a_3 \\ a_4 & a_5 & a_6 \\ a_7 & a_8 & a_9 \end{pmatrix} \middle| a_i \in 2Z^+ \cup \{0\}, 1 \leq i \leq 9 \right\}$ over the semigroup S = $Z^+ \cup \{0\}$. Does V have semigroup linear subalgebras?

116. Does there exist any semigroup linear algebra without a pseudo semigroup vector subspace?

117. Obtain any interesting result about semigroup linear algebra and its substructures.

118. Does there exist a semigroup linear algebra which has no pseudo subsemigroup vector subspace?

119. Let F[x} be the semigroup under polynomial addition over the prime field F = $Z_3$ = {0, 1, 2} of characteristic 3. Let S = $Z_3$ = {0, 1, 2} be the semigroup under addition. F[x] is a semigroup linear algebra over $Z_3$ = S. Can F[x] have a pseudo subsemigroup vector subspace?

120. Give an example of a simple semigroup linear algebra defined over a semigroup which has proper subsemigroups.

121. Can a semigroup linear algebra V defined over a non simple semigroup S be a pseudo simple semigroup linear algebra?

122. Define a semigroup linear operator on the semigroup linear algebra V = P × P × P × P × P where P = $2Z^+ \cup \{0\}$ over the

semigroup $S = Z^+ \cup \{0\}$ which is invertible. Prove the set of all semigroup linear operators on V is a semigroup linear algebra over S. What is the dimension of $Hom_S(V, V)$ = {all semigroup linear operators on V}?

123. Let $V = \left\{\begin{pmatrix} a_1 & a_2 & \ldots & a_7 \\ a_8 & a_9 & \ldots & a_{14} \end{pmatrix} \middle| a_i \in Z^+ \cup \{0\}; 1 \leq i \leq 14\right\}$ be a semigroup linear algebra over the semigroup $S = Z^+ \cup \{0\}$. Define a semigroup linear operator on V which is non invertible. What is the dimension of $Hom_S$ (V, V)? Is $Hom_S$ (V, V) a semigroup linear algebra over S? Suppose

$$W = \left\{\begin{pmatrix} a_1 & a_2 & \ldots & a_7 \\ 0 & 0 & \ldots & 0 \end{pmatrix} \middle| a_i \in Z^+ \cup \{0\}; 1 \leq i \leq 7\right\} \subseteq V$$

be a semigroup linear subalgebra of V. Define a semigroup linear operator T on V which is a semigroup projection of V into W. Is T o T = T?

124. Let $V = \left\{\begin{pmatrix} a_1 & a_2 & a_3 & a_4 & a_5 \\ a_6 & a_7 & a_8 & a_9 & a_{10} \end{pmatrix} \middle| a_i \in Z^+ \cup \{0\}, 1 \leq i \leq 10\right\}$ be a semigroup linear algebra over the semigroup $S = Z^+ \cup \{0\}$.

$$W = \left\{\begin{pmatrix} a_1 & a_2 & 0 & 0 & 0 \\ a_6 & a_7 & 0 & 0 & 0 \end{pmatrix} \middle| a_1, a_2, a_6, a_7 \in Z^+ \cup \{0\}\right\} \subseteq V \text{ be}$$

a semigroup linear subalgebra of V over the subsemigroup $T = 2Z^+ \cup \{0\}$. Define pseudo projection $U : V \to W$. Is U o U = U?

125. Obtain some interesting properties about

1. semigroup projection operator of a semigroup linear algebra

2. pseudo semigroup projection operator of a subsemigroup linear subalgebra.

Find the difference between them.

126. Is every semigroup linear projection operator on a semigroup linear subalgebra an idempotent map?

127. Is every pseudo semigroup linear projection operator of V an idempotent map?

128. Let V be a semigroup vector space over the semigroup S. Is it always possible to represent V as a direct union of the semigroup vector subspaces of V? Justify your claim.

129. Let V = $\left\{\begin{pmatrix} a_1 & a_2 \\ a_3 & a_4 \end{pmatrix}, (a_1\ a_2\ a_3\ a_4), \begin{pmatrix} 0 & 0 & 0 & 0 \\ a_1 & a_2 & a_3 & a_4 \end{pmatrix},\right.$ $\left.\begin{pmatrix} a_1 & 0 & a_2 & 0 \\ 0 & a_3 & 0 & a_4 \end{pmatrix} \middle| a_i \in Q^+ \cup \{0\}; 1 \leq i \leq 4\right\}$ be a semigroup vector space over the semigroup $S = Z^+ \cup \{0\}$.

Is V representable as the direct union of semigroup vector subspaces of V? Does V contain subsemigroup vector subspaces?

130. Does there exist a semigroup vector space which cannot be represented as a direct union of semigroup vector subspaces?

131. Is it possible to represent a semigroup vector space over the semigroup S both as a direct union of semigroup vector subspaces and pseudo direct union of semigroup vector subspaces of V? Justify your claim.

132. Let V = $\left\{\begin{pmatrix} a_1 & a_2 & a_3 \\ a_4 & a_5 & a_6 \end{pmatrix} \middle| a_i \in Q^+ \cup \{0\};\ 1 \leq i \leq 6\right\}$ be a semigroup linear algebra over the semigroup $S = Q^+ \cup \{0\}$. Is the set of all semigroup linear operators on V a semigroup linear algebra over S? Justify your answer.

133. Let V be as in problem 132 suppose W = $\left\{\begin{pmatrix} a_1 & a_2 & a_3 \\ a_4 & a_5 & a_6 \end{pmatrix} \middle| a_i \in Q^+ \cup \{0\}; 1 \leq i \leq 6\right\}$ be a subsemigroup linear subalgebra of V over the subsemigroup $T = Z^+ \cup \{0\}$. Define a projection from V to W. Represent V as a direct union of semigroup linear subalgebras of V over S.

Can V be represented as a pseudo direct union of subsemigroup linear subalgebras?

134. Does there exist a semigroup vector space which cannot be written as a direct union of semigroup subvector spaces?

135. Does there exist a semigroup vector space V which cannot be written as a pseudo direct union of semigroup subvector spaces of V?

136. Can the semigroup vector space $V = \{(a_1, a_2, \ldots, a_5), (a_1, a_2, a_3), (a_1\ a_2 \ldots a_7), \begin{pmatrix} a_1 & a_2 \\ a_3 & a_4 \end{pmatrix} \middle| a_i \in Z^+ \cup \{0\}; 1 \leq i \leq 7$ over the semigroup $S = Z^+ \cup \{0\}$ be represented as a direct union of semigroup vector subspaces of V?

137. Can the semigroup vector space V = {0, 1, 3, 5, 7, …} over the semigroup S = {0, 1} be such that 1 + 1 = 1 be represented as the direct union of semigroup subvector spaces of V?

138. Let P[x] = {all polynomials in the variable x with coefficients from $Z^+ \cup \{0\}$}, be a semigroup linear algebra over the semigroup $S = Z^+ \cup \{0\}$. W[x] = {all polynomials of even degree with coefficients from $S = Z^+ \cup \{0\}$} be a semigroup linear subalgebra of P[x].

Define a projection from P[x] into W[x].

Is that projection semigroup operator idempotent semigroup linear operator on P[x]?

139. Let $V = \{0, 1, 3, 3^2, 3^3, \ldots\}$ be a semigroup vector space over the semigroup $\{0, 1\}$ with $1 + 1 = 1$. Can V be represented as a pseudo direct union of semigroup subvector spaces of V? Is V representable as the direct union of semigroup subvector spaces of V?

140. Does there exists semigroup linear algebra V which has a unique direct sum?

141. Let $V = P \times P$ where $P = Z^+ \cup \{0\}$ be a semigroup linear algebra over P. Can V have many direct sum representations?

142. Suppose $V = \{(a, b, c) \mid a, b, c \in Z^+ \cup \{0\}\}$, be a semigroup linear algebra over $S = Z^+ \cup \{0\}$. How many ways can V be represented as a direct sum?

143. Does their exists a semigroup linear algebra V which can be represented as a direct sum in infinite number of ways?

144. We know $V = Z^+ \cup \{0\}$ is a semigroup linear algebra over $S = Z^+ \cup \{0\} = V$. Clearly dimension of V is one but V has infinitely many semigroup linear subalgebras over V given by $W_n = \{nZ^+ \cup \{0\} \mid n \in Z^+\} \subseteq V$, as n can take infinite number of values. V has infinite number of semigroup linear subalgebras. Can V be written as a finite direct sum of semigroup linear subalgebras?

145. Let $V = Q^+ \cup \{0\}$ be a semigroup linear algebra over the semigroup $S = Z^+ \cup \{0\}$. What is the dimension of V? Can V be represented as a finite direct sum of semigroup linear subalgebras of V?

146. Characterize all those semigroup linear algebras V which can be written as a finite direct sum of semigroup subalgebras of V.

147. Characterize all those semigroup linear algebras V of dimension one which is an infinite direct sum of semigroup linear subalgebras of V.

148. Let $V = Z_7 = \{0, 1, 2, \ldots, 6\}$ be a semigroup under addition modulo 7. $S = Z_7 = V$ the semigroup. V is a semigroup linear algebra of dimension one over S. Can V have any semigroup linear subalgebras? Can V be written as a non trivial direct sum of semigroup linear subalgebras?

149. Find all semigroup linear subalgebras of the semigroup linear algebra $Z_p[x]$ = {all polynomials in the variable x with coefficients from $Z_p$} over the semigroup $Z_p$, p any prime. Can $Z_p [x]$ be a direct sum of semigroup linear subalgebras?

150. Let $V = Z_p \times Z_p \times Z_p \times Z_p$ be the semigroup linear algebra over the semigroup $S = Z_p$ under addition. Can V be a direct sum of semigroup linear subalgebras, p be any prime?

151. Let $V = Z_6 \times Z_6 \times Z_6$ be the semigroup linear algebra over the semigroup $S = Z_6$. Can V be a direct sum of semigroup linear subalgebras?

152. Let $V = Z_6 \times Z_6 \times Z_6$ be the semigroup linear algebra over the semigroup $S = \{0, 1\}$. Can V be a direct sum of semigroup linear subalgebras?

153. Let $V = Z_6 \times Z_6 \times Z_6$ be the semigroup linear algebra over the semigroup $S = \{0, 1\}$? What can we say about its direct sum?

154. Let $V = Z_9 \times Z_9 \times Z_9 \times Z_9$ be a group vector space over the group $H = \{0, 3, 6\}$. Find group vector subspaces of V. Can V have subgroup vector subspaces? Justify your claim.

155. Give some interesting results about group vector spaces.

156. Let

$$V = \left\{ \begin{pmatrix} 1 & 1 & 1 & 1 \\ 0 & 0 & 0 & 1 \\ 0 & 0 & 0 & 1 \end{pmatrix}, \begin{pmatrix} 0 & 0 & 0 & 0 \\ 0 & 0 & 0 & 0 \\ 0 & 0 & 0 & 0 \end{pmatrix}, \begin{pmatrix} 1 & 1 & 0 \\ 1 & 1 & 0 \\ 1 & 1 & 0 \end{pmatrix}, \begin{pmatrix} 0 & 0 & 0 \\ 0 & 0 & 0 \\ 0 & 0 & 0 \end{pmatrix}, \right.$$

(1 1 1 1 1), (0 0 0 0 0)},

Is V a group vector space over the group $Z_2 = G$ (group under addition modulo 2)? Can V have proper group vector subspaces? Is W = {(1 1 1 1 1), (0 0 0 0 0)} ⊂ V a proper group vector subspace of V?

157. Does their exists a group vector space V over the group G; G has proper subgroups, yet V is a simple group vector space?

158. Let $V = Z_5 \times Z_5 \times Z_5 \times Z_5$ be a group vector space over the group $G = Z_5$. What is the dimension of V over $Z_5$? Find a generating subset of V. Is T = {(1 1 1 1 1)} a generating subset of V?

159. Let $V = Z \times Z$ be a group vector space over the group Z. What is the dimension of V? Can a finite subset of V be a generating subset of V?

160. Let $V = \left\{\begin{pmatrix} a_1 & a_2 \\ a_3 & a_4 \end{pmatrix} \middle| a_i \in Z, 1 \leq i \leq 4\right\}$ be a group vector space over Z.

Can $\left\{\begin{pmatrix} 1 & 0 \\ 0 & 0 \end{pmatrix}, \begin{pmatrix} 0 & 1 \\ 0 & 0 \end{pmatrix}, \begin{pmatrix} 0 & 0 \\ 1 & 0 \end{pmatrix}, \begin{pmatrix} 0 & 0 \\ 0 & 1 \end{pmatrix}, \begin{pmatrix} 1 & 1 \\ 0 & 0 \end{pmatrix}, \begin{pmatrix} 1 & 0 \\ 1 & 0 \end{pmatrix}, \begin{pmatrix} 0 & 1 \\ 0 & 1 \end{pmatrix}, \begin{pmatrix} 1 & 0 \\ 0 & 1 \end{pmatrix}, \begin{pmatrix} 0 & 1 \\ 1 & 0 \end{pmatrix}, \begin{pmatrix} 0 & 0 \\ 1 & 1 \end{pmatrix}, \begin{pmatrix} 1 & 1 \\ 1 & 1 \end{pmatrix}\right\} \subset V$ be the generating subset of V?

What is the dimension of V over Z? Find proper subgroup vector subspaces of V? Give a proper group vector subspace of V? Is $W = \left\{\begin{pmatrix} a & a \\ 0 & 0 \end{pmatrix} \middle| a \in Z\right\} \subset V$ a group vector subspace of V? What is the dimension of W over Z? Is $\left\{\begin{pmatrix} 1 & 1 \\ 0 & 0 \end{pmatrix}\right\} \subset W$ a generating subset of W? Can $\left\{\begin{pmatrix} 5 & 5 \\ 0 & 0 \end{pmatrix}\right\}$ be a generating subset of V? Can W have more than one generating set?

161. Let $V = \{(a_1\ a_2, \ldots, a_n) \mid a_i \in R\}$ be a group vector space over the group 2Z. Find proper subgroup vector subspaces of V? Does V have subgroup vector subspaces? Is $W = \{(a, a, \ldots, a) \mid a \in Z\}$, a group vector subspace of V? What is the dimension of W? What is the dimension of V? Find a generating subset of W? Find a generating subset of V? Is V finite dimensional over 2Z? Suppose V is defined over the group R will V be finite dimensional?

162. Let $V = Z_5 \times Z_5 \times Z_5 \times Z_5$ be a group vector space over $Z_5$. Can V have subgroup vector subspaces? Find all group vector subspaces of V. What is dimension of V over $Z_5$? Find the finite set which generates V over $Z_5$.

163. Let $V = Z_2 \times Z_2 \times Z_2$ be a group vector space over $Z_2$. Prove V is a simple group vector space. Find all group vector subspaces of V. Is B = {(0 0 0), (1 0 0), (0 1 0), (0 0 1), (1 1 0), (0 1 1), (1 1 1), (1 0 1)} a basis of V over $Z_2$?

164. Let $V = Z_{16} \times Z_{16} \times Z_{16} \times Z_{16}$ be a group vector space over the group G = {0, 4, 8, 12} under addition modulo 16. Find

   a. Subgroup vector subspaces of V using the subgroup H = {0, 8} contained in G.

   b. $W = Z_{16} \times Z_{16} \times \{0\} \times \{0\}$, a group vector subspace of V. What is dimension of W?

   c. How many group vector subspaces can V have?

   d. What is dimension of V over G?

165. Give examples of a one dimensional group vector space which has a group vector subspace of infinite dimension?

166. Let $V = Z \times Z \times Z$ be a group vector space over Z. Let $W = \{(p, p, p) \mid p \text{ is a prime}\} \subset V$. Is W a group vector subspace of V? What is dimension of V over Z? What is the dimension of W over Z? Suppose S = $\{(x_1, x_2, x_3) \mid x_i \in 2Z; 1 \leq i \leq 3\} \subseteq V$; is S a group vector subspace of V? What is dimension of S over Z?

167. Let V = $\{(a, a, a, a) \mid a \in Z\}$ a group vector space over Z. What is the dimension of V over Z? Can V have group vector subspaces of higher dimension than V? Suppose S = $\{(7, 7, 7, 7), (2, 2, 2, 2), (5, 5, 5, 5), (1, 1, 1, 1), (0, 0, 0, 0)\} \subseteq V$. Can S be a group vector subspace of V? Justify your claim?

168. Let V = $\left\{\begin{pmatrix} 1 & 1 & 1 \\ 0 & 1 & 1 \\ 0 & 0 & 1 \end{pmatrix}, \begin{pmatrix} 0 & 0 & 0 \\ 0 & 0 & 0 \\ 0 & 0 & 0 \end{pmatrix}, \begin{pmatrix} 1 & 1 \\ 1 & 1 \end{pmatrix}, \begin{pmatrix} 0 & 0 \\ 0 & 0 \end{pmatrix}\right\}$ be a group vector space over the group $Z_2 = \{0, 1\}$. Prove V is a simple group vector space. Prove V has only two group vector subspaces. Is dimension of V over $Z_2$ four?

169. Let V = {1 1 1 1 1 1 1), (1 1 1 0 0 0 0), (0 0 0 0 0 0 0), (0 1 0 1 0 1 0), (1 0 1 0 1 0 1)} be a group vector space over the group $Z_2 = \{0, 1\}$. Find all group vector subspaces of V. What is dimension of V over $Z_2$?

170. Let V = Q × Q be a group vector space over Q. What is dimension of V over Q? Find group vector subspaces of V. Does V contain subgroup vector subspaces of V? Does V contain subgroup vector subspaces? Let $W = Z \times Z \subset Q \times Q$, is W a subgroup vector subspace of V over the subgroup Z? What is dimension of W? Suppose $V_1 = Q \times Q$ is a group vector space over the group 3Z. What is the dimension of $V_1$ over 3Z. Find generating subsets of V and $V_1$ over Q and 3Z respectively.

171. Give examples of group vector spaces V such that dimension of V is the same as cardinality of V?

172. Let $V = \mathbb{C} \times \mathbb{C} \times \mathbb{C}$ be a group vector space over $\mathbb{C}$. What is dimension V over $\mathbb{C}$?

If $V_1 = \mathbb{C} \times \mathbb{C} \times \mathbb{C}$ is a group vector space over R. What is dimension of $V_1$ over R?
If $V_2 = \mathbb{C} \times \mathbb{C} \times \mathbb{C}$ is a group vector space over Q. What is dimension of $V_2$ over Q?
If $V_3 = \mathbb{C} \times \mathbb{C} \times \mathbb{C}$ is a group vector space over Z what is dimension of $V_3$ over Z?
If $V_4 = \mathbb{C} \times \mathbb{C} \times \mathbb{C}$ is a group vector space over 5Z what is dimension of $V_4$ over 5Z?
Find group vector subspaces of V, $V_1$, $V_2$, $V_3$ and $V_4$.
Can V, $V_1$, $V_2$, $V_3$ and $V_4$ have duo subgroup vector subspaces?

173. Give an example of a group vector space which has no duo subgroup vector subspaces?

174. Let $V = Z_5 \times Z_5 \times Z_5$ be a group vector space over the group $Z_5$. Can V have duo subgroup vector subspaces?

175. Does their exists a group vector space V over a group G in which every group vector subspace is a subgroup vector subspace and vice versa?

176. Let $V = Z_{13} \times Z_{13} \times Z_{13} \times Z_{13}$ be group vector space over $Z_{13}$. Can V have duo subgroup vector subspace?

177. Let $V = Z_{17} \times Z_{17} \times Z_{17}$ be a group vector space over the group $Z_{17}$. Prove V is a simple group vector space. Can V have group vector subspace? What is dimension of V?

178. If V = $\underbrace{Z_p \times ... \times Z_p}_{n-\text{times}}$, $Z_p$ the group under addition modulo p. V the group vector space over $Z_p$. What is dimension of V over $Z_p$? How many generating subsets can V have?

179. Let $V = Z_6 \times Z_6 \times Z_6$ be a group vector space over $Z_6$, what is dimension of V. Find one subgroup vector subspace and group vector subspace of V?

180. Let V = $\underbrace{Z_n \times ... \times Z_n}_{t-times}$ where n is not a prime. Prove V is not simple as a group vector space over $Z_n$. What is dimension of V over $Z_n$? Find a generating subset of V. Does V contain group vector subspaces?

181. Let V = $\left\{\begin{pmatrix} a & a \\ a & a \end{pmatrix} \middle| a \in Z_5\right\}$ be a group vector space over the group $Z_5$. Prove V is simple? What is the dimension of V over $Z_5$? Is S = $\left\{\begin{pmatrix} 1 & 1 \\ 1 & 1 \end{pmatrix}, \begin{pmatrix} 0 & 0 \\ 0 & 0 \end{pmatrix}\right\}$ the generating subset of V? Can V contain group vector subspaces? Justify your claim.

182. Let V = $\left\{\begin{pmatrix} a & a & ... & a \\ \vdots & \vdots & & \vdots \\ a & a & ... & a \end{pmatrix}_{m \times n} \middle| a \in Z_p\right\}$ be a group vector space over $Z_p$, p a prime; prove V is simple. What is the dimension of V over $Z_p$? Can V have group vector subspace?

183. Let V = $\left\{\begin{pmatrix} a & a & a \\ a & a & a \end{pmatrix} \middle| a \in Z_{12}\right\}$ be a group vector space over $Z_{12}$. Is V simple? Find group vector subspaces of V. What is dimension of V over $Z_{12}$?

184. What is the dimension of V = $Z_p$ as a group vector space over $Z_p$? Is V simple? Can V have group vector subspaces?

185. Let V = $\left\{\begin{pmatrix} a & a & ... & a \\ a & a & ... & a \\ \vdots & \vdots & & \vdots \\ a & a & ... & a \end{pmatrix}_{m \times n} \middle| a \in Z_n\right.$ ; n not a prime} be a group vector space over the group $Z_n$. Prove V is not

simple! What is the dimension of V over $Z_n$. Find proper group vector subspaces of V.

186. Let V = $\{(a_{ij})_{m\times n} \mid a_{ij} \in Z_n;\ 1 \le i \le m \text{ and } 1 \le j \le n\}$ be a group vector space over $Z_n$. What is dimension of V over $Z_n$? When is V simple? When is V not simple? Find group vector subspaces of V.

187. Let V = $\left\{\begin{pmatrix} a_1 & a_2 & a_3 \\ a_4 & a_5 & a_6 \end{pmatrix} \middle| a_i \in Z_9;\ 1 \le i \le 6\right\}$ be the group vector space over $Z_9$. What is the dimension of V over $Z_9$? Is V simple? If V is defined over G = {0, 3, 6} what can one say about subgroup vector subspaces of V? What is the dimension of V over G = {0, 3, 6}?

188. Let V = $\left\{\begin{pmatrix} a_1 & a_2 & a_3 \\ a_4 & a_5 & a_6 \end{pmatrix} \middle| a_i \in Z_{3^n};\ n \ge 2;\ 1 \le i \le 6\right\}$ be a group vector space over $Z_{3^n}$. What is the dimension of V over $Z_{3^n}$? Prove V is not simple. Find subgroup vector subspace W of V which is simple.

189. Let V = $\left\{\begin{pmatrix} a & a \\ a & a \end{pmatrix}, \begin{pmatrix} b & b & b \\ b & b & b \end{pmatrix}, \begin{pmatrix} c & c \\ c & c \\ c & c \end{pmatrix} \middle| a, b, c \in Z\right\}$ be a group vector space over Z. Find the subgroup vector subspace over 5Z. Prove V is not simple. What is the dimension of V over Z? Find some group vector subspaces of V.

190. Let V = $\left\{\begin{pmatrix} a & b \\ a & b \\ a & b \end{pmatrix}, \begin{pmatrix} a & c & c \\ a & c & c \\ a & c & c \\ a & c & c \end{pmatrix} \middle| a, b, c \in Z\right\}$ be the group vector space over Z. Is V simple? What is dimension of V

over Z? Prove W = $\left\{\begin{pmatrix} a & b \\ a & b \\ a & b \end{pmatrix} \middle| a, b \in Z\right\}$ is a group subspace of V.

191. Let V = $\left\{\begin{pmatrix} a \\ a \\ a \\ b \end{pmatrix}, \begin{pmatrix} a \\ b \\ a \end{pmatrix}, \begin{pmatrix} a \\ b \end{pmatrix} \middle| a, b \in Z\right\}$ be a group vector space over the group Z. What is dimension of V over Z? Find proper group vector subspaces of V. Prove V is not simple. What is the dimension of the group vector subspace $W_1$ = $\left\{\begin{pmatrix} a \\ a \\ a \\ b \end{pmatrix} \middle| a, b \in Z\right\}$? Let $W_2$ = $\left\{\begin{pmatrix} a \\ b \\ a \end{pmatrix} \middle| a, b \in Z\right\}$ be a group vector subspace of V. What is the dimension of $W_2$? Let $W_3$ = $\left\{\begin{pmatrix} a \\ b \end{pmatrix} \middle| a, b \in Z\right\}$ be the group vector subspace of V. What dimension of $W_3$ over Z? Is dimension of $W_1$ = dimension of $W_2$ = dimension of $W_3$? W' = $\left\{\begin{pmatrix} a \\ a \\ a \\ a \end{pmatrix} \middle| a \in 2Z\right\}$ is a subgroup vector subspace over the subgroup H = 2Z. What is dimension of W' over 2Z? Is the generating subset of W' unique?

192. Does there exist a group vector space which has more than one generating set?

193. Let V = {(x, x, x), (y, y, y, y) | x $\in$ 2Z and y $\in$ 3Z} be the group vector space over Z. What is dimension of V over Z? Is V simple? What is the dimension of W = {(x x x) | x $\in$ 2Z} over Z. What is dimension of W over 2Z?

194. Let X = {(a, a) (a, b) | a $\in$ Z, b $\in$ 4Z} be a group vector space over Z. What is dimension of X over Z? Find a generating subset of X. Is it unique?

195. Let Y = {(a, a, a), (b, b, b, b), ($x_1$, $x_2$, $x_3$, $x_4$, $x_5$) | a, b $x_i \in$ Z; $1 \leq i \leq 5$} be a group vector space over the group Z. Is Y simple? What is the dimension of Y over Z? Find a generating subset S of Y. Is S unique? T = {(a, a, a), (b, b, b, b) (c, c, c, c, c)} $\subseteq$ Y, is a group vector subspace of Y over Z. Find dimension of T over Z?

196. Let V = $Z_{12} \times Z_{12}$ be a group linear algebra over $Z_{12}$. What is the dimension of V? Find a generating subset X of V.

197. Let V = $Z_{12} \times Z_{12}$ be a group linear algebra over G = {0, 2, 4, 6, 8, 10}. What is the dimension of V over G? Give a generating subset of V.

198. Let V = $Z_{12} \times Z_{12}$ be a group linear algebra over the group G = {0, 3, 6, 9}. What is dimension of V over G'? Give a generating subset of V?

199. Let V = $Z_{12} \times Z_{12}$ be the group linear algebra over the group S = {0, 6}. What is the dimension of V over S? Give a generating subset of V over S.

200. Let V = $Z_{12} \times Z_{12}$ be the group linear algebra over the group P = {0, 4, 8}. What is the dimension of V over P? Give a generating subset of V over P.

201. Does the dimension of a group linear algebra dependent on the group over which it is defined? Find the relation between the dimension of V in problems (196-200) and the cardinality of the group over which they are defined. Can one say the dimension increases with the decrease of the cardinality of the group over which V is defined?

202. Let $V = Z_{12} \times Z_{12} \times Z_{12} \times Z_{12}$ be the group linear algebra over the group $Z_{12}$. Does V have pseudo semigroup linear subalgebra? Find in V a pseudo group vector subspace over the group $Z_{12}$.

203. Let $V = Z_5 \times Z_5 \times Z_5$ be the group linear algebra over $Z_5$. Can V have pseudo semigroup linear subalgebra over $Z_5$?

204. Let $V = Z \times Z \times Z$ be the group linear algebra over Z. $W = 2Z \times Z \times Z \times 2Z$ be the group linear algebra over Z. Define a group linear transformation from V to W so that T is invertible. Define a group linear transformation from V to W, which is not invertible.

205. Let $V = Z_7 \times Z_7 \times Z_7$ be a group linear algebra over $G = Z_7$. Define a group linear operator T on V which has $T^{-1}$ such that $T \circ T^{-1} = T^{-1} \circ T$; is the identity group linear operator on V. Is T (x y z) = T (0 y x) from V to V an invertible group linear operator on V? If T (x y z) = (x + y x-y 0) from V to V, is T an invertible group linear operator on V? What is the structure of $M_G (V, V) = M_{Z_7} (V, V)$? Is $M_{Z_7} (V, V)$ a group linear algebra over $Z_7$ with composition of maps as an operation on it? Prove your claim.

206. Let V be a group linear algebra over a group G. Let $M_G$ (V, V) denote the set of all group linear operators from V to V. What is the algebraic structure of $M_G$ (V, V) ? Can $M_G$ (V, V) be a group linear algebra over G?

207. Let V and W be two group linear algebras over the group G. Let $M_G$ (V, W) denote the collection of all group linear transformations from V to W. What is the algebraic structure of $M_G$ (V, W)? Is $M_G$ (V, W) a group linear algebra over G? Is $M_G$ (V, W) a group vector space over G?

208. Let $V = Z_2 \times Z_2$ and $W = Z_2 \times Z_2 \times Z_2 \times Z_2$ be group linear algebras over the group $Z_2$. Find $M_G$ (V, W). What is the dimension of $M_G$ (V, W)? Find $M_G$ (V, V) and $M_G$ (W, W) and determine their dimensions. If T (x, y) = (x, x + y, x – y, y) is the map from V to W; is T a group linear transformation? Is T a invertible group linear transformation?

209. Let $V = \left\{ \begin{pmatrix} a & a & a \\ a & a & a \end{pmatrix} \middle| a \in Z \right\}$ and $W = \{Z \times Z \times Z\}$ be group linear algebras over Z. Find dimensions of V and W. What is the dimension of $M_z$ (V, W)?

210. Let $V = \left\{ \begin{pmatrix} a & b & c \\ d & e & f \end{pmatrix} \middle| a, b, c, d, e, f \in Z_{15} \right\}$ be a group linear algebra over G = {0, 5, 10}, group under addition modulo 15. What is dimension of V over G? Can V have pseudo semigroup linear subalgebras? Justify your claim. Is $W = \left\{ \begin{pmatrix} a & a & a \\ a & a & a \end{pmatrix} \middle| a \in Z_{15} \right\}$ a proper group linear subalgebra of V? Find a proper linear subalgebra of V with dimension 3.

211. Let $V = Z_2 [x]$ be a group linear algebra over $Z_2$. What is the dimension of V over $Z_2$?
Find $M_{Z_2}$ [V, V]. What is the dimension of $M_{Z_2}$ [V, V] over $Z_2$?

212. Let $V = Z_8 \times Z_8 \times Z_8$ be a group linear algebra over the group G = {0, 4} under addition modulo 8. Find the dimension of V over G. Find $M_G$ [V, V]. What is the dimension of $M_G$ [V, V]? Does their exist any relation between dimension of V and $M_G$ [V, V]? Which has higher dimension V or $M_G$ [V, V]?

213. Let V be a group vector space over the group G. Suppose H is a subgroup of G, V be also a group vector space over H.

Will the dimension of V increase with the subgroups on which V is defined or decrease when V is defined on the subgroups of G.

214. Let $V = Z_{12} \times Z_{12} \times Z_{12} \times Z_{12}$ be a group vector space over the group $Z_{12}$. What is dimension of V over $Z_{12}$? Let $V_1 = Z_{12} \times Z_{12} \times Z_{12} \times Z_{12}$ be the group vector space over the subgroup $H = \{0, 2, 4, 6, 8, 10\} \subseteq Z_{12}$. What is the dimension of $V_1$ over H? Suppose $V_2 = Z_{12} \times Z_{12} \times Z_{12} \times Z_{12}$ group vector space over the subgroup $P = \{0, 3, 6, 9\}$, find the dimension of $V_2$ over P. If $V_3 = Z_{12} \times Z_{12} \times Z_{12} \times Z_{12}$ group vector space over the subgroup $K = \{0, 6\}$ find the dimension of $V_3$ over K. Compare the dimensions of group vector spaces V, $V_1$, $V_2$ and $V_3$.

215. Let $V = \left\{\begin{pmatrix} a_1 & a_2 & a_3 \\ a_4 & a_5 & a_6 \end{pmatrix} \middle| a_i \in Z_8\right\}$ be the group linear algebra over $Z_8$. Find dimension of V over $Z_8$. If $V_1 = \left\{\begin{pmatrix} a_1 & a_2 & a_3 \\ a_4 & a_5 & a_6 \end{pmatrix} \middle| a_i \in Z_8; 1 \leq i \leq 6\right\}$, what is the dimension of $V_1$ over the subgroup $H = \{0, 4\} \subset Z_8$ as a group linear algebra over $H \subset Z_8$. Let $V_2 = \left\{\begin{pmatrix} a_1 & a_2 & a_3 \\ a_4 & a_5 & a_6 \end{pmatrix} \middle| a_i \in Z_8; 1 \leq i \leq 6\right\}$ be a group linear algebra over the subgroup $K = \{0, 2, 4, 6\} \subseteq Z_8$. What is the dimension of $V_2$ over K. Compare the dimensions of V, $V_1$ and $V_2$.

216. Let V be a group linear algebra over the group G, and also proper subgroups $H_1, H_2, \ldots, H_n$ of G. Find dimensions of the group linear algebra V over the subgroup $H_1, \ldots, H_n$ and compare them. Compare dimension of the group linear algebra V over the group G with the group linear algebra over $H_i$'s; $i = 1, 2, \ldots, n$.

217. Let V = $\left\{\begin{pmatrix} a & b \\ c & d \end{pmatrix} \middle| a,b,c,d \in Z_4 \right\}$ be the group linear algebra over $Z_4$. What is dimension of V over $Z_4$. Suppose V = $\left\{\begin{pmatrix} a & b \\ c & d \end{pmatrix} \middle| a,b,c,d \in Z_4 \right\}$ be the group linear algebra over H = {0, 2} $\subseteq Z_4$. What is dimension of V over H? Let V = $\left\{\begin{pmatrix} a & b \\ c & d \end{pmatrix} \middle| a,b,c,d \in Z_4 \right\}$ be a group vector space over $Z_4$ what is dimension of V as a group vector space over $Z_4$? Does dimension of V as a group vector space over $Z_4$ different from V as a group linear algebra over $Z_4$? What is the dimension of V as a vector space over H = {0, 2}?

218. Let V = $Z_5 \times Z_5 \times Z_5$ be a group vector space over $Z_5$. What is dimension of V over $Z_5$? Let V be the group linear algebra over $Z_5$. What is the dimension of V as a group linear algebra over $Z_5$?

219. Let V = Z × Z × Z × Z × Z be a group linear algebra over Z. What is dimension of V over Z? If W = 11Z × 2Z × 3Z × 5Z × 7Z $\subset$ V be the group linear subalgebra of V over Z. What is dimension of W over Z? Let P = Z × {0} × {0} × {0} × {0} $\subset$ V be a group linear subalgebra of V over Z. What is dimension of P over Z?

220. Give an example of a fuzzy vector space over the field of reals.

221. Given V = $\{Z_{12} \times Z_{12} \times Z_{12}\}$ is the set vector space over the set S = $Z_{12}$. Define a map $\eta : V \to [0, 1]$ so that $V\eta$ is a fuzzy set vector space.

222. Let $V = \left\{\begin{pmatrix} a & b \\ c & d \end{pmatrix} \middle| a, b, c, d \in Z_4\right\}$ be a set vector space over the set {0, 2}. Define a map $\eta : V \to [0, 1]$ so that $V_\eta$ is a fuzzy set vector space.

223. Let V = {(0001), (1100), (000), (0000), (111), (110), (010), (00000), (10000), (10101), (01010)} be the set vector space over the set S = {0, 1}. Define $\eta : V \to [0, 1]$ so that $V_\eta$ is a fuzzy set vector space.

224. Let $V = \left\{\begin{pmatrix} a & b \\ c & d \end{pmatrix} \middle| a, b, c, d \in Z\right\}$ be a set linear algebra over the set Z. Define the map $\eta : V \to [0, 1]$ so that $V_\eta$ is a fuzzy set linear algebra.

225. Let $V = \left\{\begin{pmatrix} a & a \\ a & a \end{pmatrix} \middle| a \in Z_{12}\right\}$ be a set linear algebra over the set S = {0, 2, 5, 7, 3}. Define $\eta : V \to [0, 1]$ so that $V\eta$ is a fuzzy set linear algebra.

226. Let $V = \left\{\begin{pmatrix} a & b \\ c & d \end{pmatrix} \middle| a, b, c, d \in Z_8\right\}$ be a semigroup vector space over the semigroup S = {0, 2, 4, 6} under addition modulo 8. Define $\eta : V \to [0, 1]$ so that $V_\eta$ is a semigroup fuzzy vector space.

227. Let $V = \left\{\begin{pmatrix} a & b & c & d \\ e & f & g & h \end{pmatrix} \middle| a, b, c, d, e, f, g, h \in Z^+\right\}$ be the semigroup vector space over the semigroup $Z^+$. Does there exist a $\eta : V \to [0, 1]$ so that $V_\eta$ is a semigroup fuzzy vector space?

228. Let $V = \left\{ \begin{pmatrix} a & b \\ c & d \end{pmatrix} \middle| a, b, c, d \in 3Z \right\}$ be a group vector space over G = Z the group under addition. Define $\eta : V \to [0, 1]$ so that $V_\eta$ is a group fuzzy vector space?

229. Let V = {(a b c d e) | a, b, c, d, e $\in Z_{20}$} be a group vector space over $Z_{20}$. Define $\eta : V \to [0, 1]$ so that $V_\eta$ is a group fuzzy vector space.

230. Let $V = Z_{20} \times Z_{20} \times Z_{20} \times Z_{20}$ be a group vector space over the group G = {0, 5, 10, 15}. Define $\eta : V \to [0, 1]$ so that $V_\eta$ is a group fuzzy vector space.

231. Let V = {(00111), (01000), (11000), (10001), (00001), (00000)} be a group vector space over the group $Z_2$ = {0, 1}. Define $\eta$: $V \to [0, 1]$ so that $V\eta$ is a group fuzzy vector space.

232. Given any group vector space V over the group G. Does there always exists a group fuzzy vector space $V_\eta$?

233. Let $V = \left\{ \begin{pmatrix} a & b \\ c & d \end{pmatrix} \middle| a, b, c, d \in Z \right\}$ be a group linear algebra over the group Z. Define map $\eta : V \to [0, 1]$ so that $V_\eta$ is a group fuzzy linear algebra. Find a group fuzzy linear subalgebra of V.

234. Let V = {(a a a a a) | a $\in$ Z} be a group linear algebra over the group Z. Define $\eta : V \to [0, 1]$ so that $V_\eta$ is a group fuzzy linear algebra. Find a group linear fuzzy subalgebra of V.

235. Obtain some interesting properties about set bivector spaces.

236. Let $V = V_1 \cup V_2 = \{(a \ a \ a \ a) \mid a \in Z^+\} \cup \left\{\begin{pmatrix} a & a & a \\ a & a & a \end{pmatrix} \middle| a \in 2Z^+\right\}$ be the set bivector space over the set $Z^+$. Find a generating bisubset of V. What is the bidimension of V? Find for V a set bivector subbispace.

237. Let $V = V_1 \cup V_2 = \left\{\begin{pmatrix} a & 0 & b \\ 0 & c & 0 \end{pmatrix} \middle| a, b, c \in Z^+\right\} \cup \left\{\begin{pmatrix} a & b & 0 \\ 0 & 0 & d \end{pmatrix} \middle| a, b, d \in Z^+\right\}$ be a set bivector space over the set $Z^+$. Find the generating bisubset of V. What the bidimension of V?

238. Given $V_1 = \left\{\begin{pmatrix} a & b \\ c & d \end{pmatrix} \middle| a, b, c, d \in Z\right\}$ the group vector space over the group G = Z. $V_2 = \{Z \times Z \times Z\}$ the group vector space over the same group G = Z. $V = V_1 \cup V_2$ is the group bivector space over G. Find the bidimension of V. Determine the generating bisubset of V.

239. Let $V = V_1 \cup V_2$ be a semigroup bivector space over the semigroup $Z_{12}$ where $V_1 = \{(a \ a \ a) \mid a \in Z_{12}\}$ and $V_2 = \left\{\begin{pmatrix} a & a \\ a & a \end{pmatrix} \middle| a \in Z_{12}\right\}$. Find a bigenerating subset $X = X_1 \cup X_2$ of $V_1 \cup V_2$. What is the bidimension of V?

240. Let $V = \{Z_{12} \times Z_{12} \times Z_{12}\} \cup Z_{12}[x]$ be the semigroup bivector space over the semigroup S = {0, 6} under addition modulo 12. Find the bidimension of V. Find a semigroup bivector subspace of V.

241. Let $V = V_1 \cup V_2 = \{(1111), (333), (456), (000)\} \cup \{Z_{12} \times Z_{12}\}$. Is V a semigroup bivector subspace? Justify your claim.

242. Let $V = \left\{\begin{pmatrix} a & a \\ a & a \end{pmatrix} \middle| a \in Z_3\right\} \cup \{(a\ a\ a\ a) \mid a \in Z_3\}$ be a semigroup bivector space over $Z_3$. What is bidimension of V over $Z_3$ as a semigroup bivector space. Suppose V is viewed as a semigroup bilinear algebra over $Z_3$, what is the bidimension of V over $Z_3$ as a semigroup bilinear algebra over $Z_3$. Hence or other wise can one prove if V is a semigroup bivector space over the semigroup S, if V be viewed as a semigroup bilinear algebra them bidimension V as a semigroup bivector space is always greater than the bidimension of V as a semigroup bilinear algebra over the semigroup S.

Does their exists a semigroup bivector space V whose bidimension is same as that of V viewed as a semigroup bilinear algebra?

243. Let $V = \{Z_4 \times Z_4 \times Z_4\} \cup \left\{\begin{pmatrix} a & b \\ c & d \end{pmatrix} \middle| a, b, c, d \in Z_4\right\} = V_1 \cup V_2$ be a semigroup bivector space over $Z_4$. If $V = V_1 \cup V_2$ is a semigroup bilinear algebra over $Z_4$. What are the bidimension of V as a semigroup bivector space over $Z_4$ and as a semigroup bilinear algebra over $Z_4$?

244. Let $V = V_1 \cup V_2 = \{(1110), (1111), (0000), (011), (000), (100), (101)\} \cup \{(1111111), (1101100), (0000000)\}$ be a semigroup bivector space over the semigroup $S = \{0, 1\} = Z_2$. Find the bidimension of V. Find a bigenerating subset of V. Can this V be made into a semigroup bilinear algebra?

245. Let $V = V_1 \cup V_2 = \left\{\begin{pmatrix} a & a \\ a & a \end{pmatrix}, \begin{pmatrix} b \\ b \end{pmatrix}, \begin{pmatrix} a & a & a \\ a & a & a \end{pmatrix} \middle| a \in Z^+\right\}$ $\left\{\begin{pmatrix} a & b \\ c & d \end{pmatrix}, (a\ \ b) \middle| a, b, c, d \in Z^+\right\}$ be the semigroup bivector

space over the semigroup $Z^+$. What is the bidimension V? Find a nontrivial semigroup bivector subspace of V?

246. Obtain some interesting results about semigroup bivector subspaces of a semigroup bivector space V.

247. Given $V = V_1 \cup V_2$ is a quasi semigroup bilinear algebra over the semigroup $S = Z^+$, here $V_1 = \{(a\ a\ a) \mid a \in Z^+\}$ is a semigroup linear algebra over S and $V_2 = \left\{\begin{pmatrix} a & a \\ a & 0 \end{pmatrix}, \begin{pmatrix} a & a \\ 0 & a \end{pmatrix} \middle| a \in Z^+\right\}$ is only a semigroup vector space.

(1) Find the bigenerating subset of V.
(2) What is the bidimension of V?
(3) Find a proper quasi semi linear subalgebra of V.

248. Let $V = V_1 \cup V_2$ be a quasi semigroup linear algebra over $Z_8$ where $V_1 = Z_8[x]$ is the semigroup linear algebra over $Z_8$ and $V_2$ = {(1110), (2220), (3330), (4440), (5550), (6660), (7770), (0000), (000a) | $a \in Z_8$} is the semigroup vector space over $Z_8$. Find a quasi semigroup linear subalgebra of V. What is bidimension of V? Is $W = \{ax + b / a, b \in Z_8\} \cup$ {(2220), (4440), (0000), (0002), (0004), (6660), (0006)} a quasi semigroup subalgebra of V. What is the bidimension of W? Find the generating bisubset of V and W.

249. Give some interesting properties about group bivector spaces defined over a group G.

250. Define bisemigroup bivector subspace of a bisemigroup bivector space V over the bisemigroup $S = S_1 \cup S_2$. Illustrate this by an example.

251. Let $V = V_1 \cup V_2 = \{Z_6[x]\} \cup \{Z^+[x]\}$ be the bisemigroup bivector space defined over the bisemigroup $S = S_1 \cup S_2 = Z_6 \cup Z^+$ . Find a bisemigroup bivector subspace of V. Find the bigenerating biset of V over S. What the bidimension of V?

252. Give an example of a (7,4) bidimensional bisemigroup bivector space $V = V_1 \cup V_2$ over the bisemigroup $S = S_1 \cup S_2$.

253. Let V = {(1110), (0001), (110), (000), (0000), (101), (011), (0110)} ∪ {(000), (210), (320), (020), (200)} = $V_1 \cup V_2$ be a bisemigroup bivector space over the bisemigroup S = S = $Z_2 \cup Z_3$. What is the bidimension of V over the bisemigroup $S = S_1 \cup S_2$ ?

254. Find the bidimension of the group bivector space $V = V_1 \cup V_2 = \{Z[x]\} \cup \left\{\begin{pmatrix} a & b \\ c & d \end{pmatrix} \middle| a, b, c, d \in Z\right\}$ over the group G = Z.

255. Obtain some interesting properties about bigroup bivector spaces over the bigroup $G = G_1 \cup G_2$.

256. Let $V = \{2Z\ [x]\} \cup \left\{\begin{pmatrix} a & a \\ a & a \end{pmatrix} \middle| a \in 3Z\right\} = V_1 \cup V_2$ be the bigroup bivector space over the bigroup G = 2Z ∪ 3Z. Find the bidimension of V over G. What the bigenerating subset of V?

257. Let $V = \{Z_7\ [x]\} \cup \{(a, a, a) \mid a \in Z_9\}$ be the bigroup bivector space over the bigroup $Z_7 \cup Z_9$. What is the bidimension of V? Find a bigenerating bisubset $X = X_1 \cup X_2$ of V. Does V have a bigroup bivector subspace?

258. Let $V = V_1 \cup V_2 = \{2Z \times 2Z \times 2Z\} \cup \{9Z \times 9Z \times 9Z \times 9Z\}$ be a bigroup bivector space over the bigroup G = 2Z ∪ 9Z. Find a bigenerating subset of $X_1 \cup X_2 \subseteq V = V_1 \cup V_2$ . Suppose $V = V_1 \cup V_2$ is just considered as a group bivector space over the group G = Z. What is the bidimension of V? Find a bigenerating subset of V. Compare the properties of these spaces.

259. Let V = $V_1 \cup V_2 = \left\{\begin{pmatrix} a & a \\ a & a \end{pmatrix} \middle| a \in Z\right\} \cup \left\{\begin{pmatrix} a & b \\ 0 & c \end{pmatrix} \middle| a, b, c \in Z\right\}$ be the group bivector space over the group G = Z. What is the bidimension of V? Find a bigenerating bisubset X of V.

260. Let V = $V_1 \cup V_2 = \left\{\begin{pmatrix} a & b \\ c & d \end{pmatrix} \middle| a, b, c, d \in 3Z\right\} \cup \left\{\begin{pmatrix} a & b \\ c & d \end{pmatrix} \middle| a, b, c, d \in 5Z\right\}$ be the bigroup bivector space over the bigroup G = 3Z ∪ 5Z. Find at least two proper bigroup bivector subspaces of V. If = $X_1 \cup X_2 \subseteq V_1 \cup V_2$ is a proper subset of V which is the bigenerator V. What is bidimension of V over the bigroup? Is W = $\left\{\begin{pmatrix} a & a \\ 0 & 0 \end{pmatrix} \middle| a \in 3Z\right\} \cup \left\{\begin{pmatrix} a & a \\ 0 & 0 \end{pmatrix} \middle| a \in 5Z\right\}$ a bigroup bivector subspace of V? What is the bidimension of W over G = 3Z ∪ 5Z?

261. Let V = $V_1 \cup V_2$ = {(a, a, a, a, a) / a ∈ $Z_{11}$} ∪ {(a a a) / a ∈ $Z_{31}$} be a bigroup bivector space over the bigroup G = $Z_{11} \cup Z_{31}$.What is the bidimension of V over the bigroup G = $Z_{11} \cup Z_{31}$? Does V have bigroup bivector subbispaces? Find a bigenerating bisubset X = $X_1 \cup X_2 \subseteq V_1 \cup V_2$ over G? Is X unique or can V have several generating bisubsets? Justify your claim.

262. Let V = $V_1 \cup V_2 = \{Z_5^7\} \cup \{Z_6^8\}$ be a bigroup bivector space over the bigroup G = = $Z_5 \cup Z_6$ . What is the bidimension of V over G? Does these exist a pseudo bisemigroup bivector subspace of W of V?

263. Can every bigroup bivector space over a bigroup have pseudo bisemigroup bivector subspace? Justify your claim.

264. Can every bigroup bivector space over a bigroup have bigroup bivector subspace?
Substantiate your claim!

265. Let $V = V_1 \cup V_2 = \{Z_{11} [x]\} \cup \{Z_{14} \times Z_{14} \times Z_{14} \times Z_{14}\}$ be a bigroup bivector space over the bigroup $G = Z_{11} \cup Z_{14}$ . Can V have a pseudo bisemigroup bivector subspace over a bisemigroup $H = H_1 \cup H_2 \subseteq Z_{11} \cup Z_{14}$ ? Justify your claim.

266. Let $V = V_1 \cup V_2 = \left\{\begin{pmatrix} a & b \\ c & d \end{pmatrix} \middle| a,b,c,d \in Z_7\right\} \cup \left\{\begin{pmatrix} a & a & a \\ a & a & a \end{pmatrix} \middle| a \in Z_5\right\}$ be a bigroup bivector space over the bigroup $G = Z_7 \cup Z_5$. Find bigroup bivector subspace of V over the bigroup G.

267. Does the bigroup $G = Z_7 \cup Z_3$ (addition the operator on $Z_7$ and $Z_3$) have nontrivial bisubsemigroups?

268. Let $V = V_1 \cup V_2 = \{Z_3[x]\} \cup \left\{\begin{pmatrix} a & b \\ c & d \end{pmatrix} \middle| a,b,c,d \in Z_7\right\}$ be a bigroup bivector space over the bigroup $G = Z_3 \cup Z_7$ . Can V have pseudo bisemigroup bivector space over any bisubsemigroup $H = H_1 \cup H_2$ of the bigroup $G = Z_3 \cup Z_7$?

269. Define a bigroup bilinear transformation T of the bigroup bivector spaces V and W defined over the same bigroup $G = G_1 \cup G_2$. Hint: $V = V_1 \cup V_2$ and $W = W_1 \cup W_2$ ; $V_i$ and $W_i$ are group vector spaces over the group $G_i$ ; i = 1, 2; $T : V \to W$ is such that $T = T_1 \cup T_2 : V_1 \cup V_2 \to W_1 \cup W_2$. $T_i : V_i \to W_i$, i = 1, 2.

270. Let $\text{Hom}_{G_1 \cup G_2}$ (V, W) be the collection of all bigroup bilinear transformations from V to W. Is $\text{Hom}_{G_1 \cup G_2}$ (V, W) a bigroup bilinear bivector space over $G = G_1 \cup G_2$?

271. Define a bigroup bilinear operator on a bigroup bivector space $V = V_1 \cup V_2$ over the bigroup $G = G_1 \cup G_2$. Will $\text{Hom}_{G_1 \cup G_2}$ (V, V) the collection of all bigroup bilinear operators form a bigroup bivector space over the bigroup G $= G_1 \cup G_2$?

272. Let $V = \{Z_5\ [x]\} \cup \{Z_7 \times Z_7\}$ and $W = \left\{\begin{pmatrix} a & b \\ c & d \end{pmatrix} \middle| a, b, c, d \in Z_5\right\} \cup \left\{\begin{pmatrix} a & a \\ a & a \end{pmatrix} \middle| a \in Z_7\right\}$ be bigroup bivector bispaces over the bigroup $G = Z_5 \cup Z_7$ . Define a bigroup bilinear transformation from V to W. What is the bidimension of V and W over G? Find $\text{Hom}_{G_1 \cup G_2}$ (V, W). Is $\text{Hom}_{G_1 \cup G_2}$ (V, W) a bigroup bivector bispace over $G = G_1 \cup G_2 = Z_5 \cup Z_7$? Justify your claim!

273. Let $V = V_1 \cup V_2$ be a bigroup bivector space over the bigroup $G = Z_7 \cup Z_{11}$ where $V_1 = Z_7 \times Z_7$ and $Z_{11} \times Z_{11} \times Z_{11} = V_2$. Find $\text{Hom}_{G_1 \cup G_2}$ (V, V). Is $\text{Hom}_{G_1 \cup G_2}$ (V, V) a bigroup linear bialgebra over $G = G_1 \cup G_2 = Z_7 \cup Z_{11}$ under composition of maps?

274. Let $V = V_1 \cup V_2 = \{Z_5 \times Z_5 \times Z_5\} \cup \{Z_7 \times Z_7 \times Z_7\}$ be a bigroup bivector bispace over the bigroup $G = Z_5 \cup Z_7$. Can we define for $W = W_1 \cup W_2 = \{Z_5 \times \{0\} \times Z_5\} \cup \{\{0\} \times Z_7 \times Z_7\}$ a bigroup bivector subspace of V a biprojection E of V to V? Is E o E = E? Define E explicitly E from V to V.

275. Let $V = \{Z_7 \times Z_7 \times Z_7\} \cup \{Z_5 \times Z_5 \times Z_5\}$ a bigroup bivector space over the bigroup $G = Z_7 \cup Z_5$. Write V as a direct sum of bigroup bivector subspaces of V. Suppose $V = \bigoplus_i W_i$ .

Can we for each i define $E_i$ a biprojector? Will $E_i \circ E_j = 0$ if $i \neq j$ and $E_i \circ E_j = E_i$ if $i = j$ ? Justify your claim.

276. Let $V = V_1 \cup V_2 = \left\{\begin{pmatrix} a & b \\ c & d \end{pmatrix} \middle| a, b, c, d \in Z_{12}\right\} \cup \left\{\begin{pmatrix} a & b \\ c & d \end{pmatrix} \middle| a, b, c, d \in Z_{18}\right\}$ be a bigroup bivector space over the bigroup $G = Z_{12} \cup Z_{18}$. Define biprojection operators on V. Suppose $W = W_1 \cup W_2 = \left\{\begin{pmatrix} a & a \\ a & a \end{pmatrix} \middle| a \in Z_{12}\right\} \cup \left\{\begin{pmatrix} a & a \\ a & a \end{pmatrix} \middle| a \in Z_{18}\right\} \subseteq V_1 \cup V_2$ be a bigroup bivector subspace of V over $G = Z_{12} \cup Z_{18}$ can we have biprojection of V to W?

277. Let $V = V_1 \cup V_2 = \{Z^+ [x]\} \cup \left\{\begin{pmatrix} a & b \\ c & d \end{pmatrix} \middle| a, b, c, d \in 5Z^+\right\}$ be a bisemigroup bivector space over the bisemigroup $G = G_1 \cup G_2 = 3Z^+ \cup 5Z^+$. Find the bidimension of V. Is V finite bidimensional? Find proper bisemigroup bivector subspaces of V over the bisemigroup $G = 3Z^+ \cup 5Z^+$. Define a bisemigroup bilinear operator on V. Find $\text{Hom}_{G_1 \cup G_2}$ (V, V). Is $\text{Hom}_{3Z^+ \cup 5Z^+}$ (V, V) again a bisemigroup bivector space over the bisemigroup $G = 3Z^+ \cup 5Z^+$.

278. Let $V = V_1 \cup V_2$ and $W = W_1 \cup W_2$ be bisemigroup bivector space over the bisemigroup $S = S_1 \cup S_2$. Find $P = \text{Hom}_{S_1 \cup S_2}$ (V, W). Is P again a bisemigroup bivector space over $S_1 \cup S_2$? Find bidimension of P if V is $(n_1, n_2)$ bidimension and W is of $(m_1, m_2)$ bidimension over $S = S_1 \cup S_2$.

279. Let $V = \left\{\begin{pmatrix} a & b \\ c & d \end{pmatrix} \middle| a, b, c, d \in 3Z^+\right\} \cup \{(a, a) \mid a \in 5Z^+\}$ be a bisemigroup bivector space over the bisemigroup $S = S_1 \cup S_2 = 3Z^+ \cup 5Z^+$. What the bidimension of V?
If $P = \text{Hom}_{S_1 \cup S_2}(V, V)$ where $S = S_1 \cup S_2 = 3Z^+ \cup 5Z^+$, what is the bidimension of P over $S_1 \cup S_2$?

280. Let $V = V_1 \cup V_2 = \{Z_2 \times Z_2 \times Z_2\} \cup \{3Z^+ \times 3Z^+ \times 3Z^+\}$ a bisemigroup bivector space over the bisemigroup $S = S_1 \cup S_2 = Z_2 \cup 3Z^+$.
1. Find bidimension of V.
2. Find $\text{Hom}_{Z_2 \cup 3Z^+}(V, V)$ over $S = Z_2 \cup 3Z^+$.
3. Does V have proper bisemigroup bivector subspace? Find at least two distinct proper bisemigroup bivector subspaces of V of same bidimension.

281. Let $V = V_1 \cup ... \cup V_5 = \{Z_7 \times Z_7\} \cup \{Z_7[x]\} \cup \left\{\begin{pmatrix} a & a & a \\ a & a & a \end{pmatrix} \middle| a \in Z_7\right\} \cup \left\{\begin{pmatrix} a & b \\ c & d \end{pmatrix} \middle| a, b, c, d \in Z_7\right\} \cup \left\{\begin{pmatrix} a & a \\ a & a \\ a & a \end{pmatrix} \middle| a \in Z_7\right\}$ be a 5-set vector space over $S = Z_7$. Find a generating 5-set subset of $Z_7$. What is the 5-dimension of 5 set vector space over $Z_7$, $Z_7$ is taken just as a set and not as a semigroup?
If the same V is taken as a semigroup 5-set vector space over $Z_7$ what the 5-dimension of V over $Z_7$ as a semigroup? Is V is a semigroup 5 set linear algebra over $Z_7$?

282. Let $V = V_1 \cup V_2 \cup V_3 \cup V_4 = \{Z^+ \times Z^+ \times Z^+\} \cup \left\{\begin{pmatrix} a & a & a \\ a & a & a \end{pmatrix} \middle| a \in Z^+\right\} \cup \left\{\begin{pmatrix} a & a \\ a & a \end{pmatrix} \middle| a \in Z^+\right\} \cup$

$\left\{\begin{pmatrix} a & a & a \\ a & a & a \\ a & a & a \end{pmatrix} \middle| a \in Z^+ \right\}$ be a semigroup 4 set linear algebra over the semigroup $Z^+ = S$. Find semigroup 4 set linear operator on V. Is $\text{Hom}_{Z^+}$ (V, V) a semigroup set linear algebra over $Z^+$ where $\text{Hom}_{Z^+}$ (V, V) = {T : V → V where $T = T_1 \cup T_2 \cup T_3 \cup T_4$}?

283. Let $V = V_1 \cup V_2 \cup V_3 \cup V_4 \cup V_5 = \{Z^+ \times Z^+\} \cup \left\{\begin{pmatrix} a & b \\ c & d \end{pmatrix} \middle| a,b,c,d \in Z^+ \right\} \cup \left\{\begin{pmatrix} a & a & a \\ a & a & a \end{pmatrix} \middle| a \in Z^+ \right\} \cup \left\{\begin{pmatrix} a & 0 & 0 \\ b & c & 0 \\ d & e & f \end{pmatrix} \middle| a,b,c,d,e,f \in Z^+ \right\} \cup \left\{\begin{pmatrix} a & a \\ a & a \\ a & a \\ a & a \end{pmatrix} \middle| a \in Z^+ \right\}$ be a semigroup 5 set linear algebra over $Z^+$. Let $W = W_1 \cup W_2 \cup W_3 \cup W_4 \cup W_5 = \{Z^+ \times Z^+ \times Z^+ \times Z^+\} \cup \left\{\begin{pmatrix} a & b & c \\ 0 & d & e \\ 0 & 0 & f \end{pmatrix} \middle| a,b,c,d,e,f \in Z^+ \right\} \cup \left\{\begin{pmatrix} a & a & a & a & a \\ a & a & a & a & a \end{pmatrix} \middle| a \in Z^+ \right\} \cup \left\{\begin{pmatrix} a & a \\ a & a \\ a & a \\ a & a \end{pmatrix} \middle| a \in Z^+ \right\} \cup \left\{\begin{pmatrix} a & a \\ a & a \end{pmatrix} \middle| a \in Z^+ \right\}$. $V_1 \cup V_2 \cup \ldots \cup V_5 = V$ be a semigroup 5-set linear algebra over $Z^+$. Let $\text{Hom}_{Z^+}$ (V, W) = {$T = T_1 \cup T_2 \cup T_3 \cup T_4 \cup T_5$: V → W /$T_i$ : $V_i \to W_i$ be the collection of all semigroup linear algebra transformation; $1 \leq i \leq 5$}. Is $\text{Hom}_{Z^+}$ (V, W) a semigroup 5 set linear algebra over $Z^+$?

284. Let $V = V_1 \cup V_2 \cup V_3 \cup V_4 \cup V_5 = \left\{ \begin{pmatrix} a & a \\ 0 & a \end{pmatrix} \middle| a \in Z_2 = \{0,1\} \right\} \cup (Z_2 \times Z_2 \times Z_2\} \cup \left\{ \begin{pmatrix} a & a \\ a & a \\ a & a \\ a & a \end{pmatrix} \middle| a \in Z_2 = \{0,1\} \right\} \cup \left\{ \begin{pmatrix} a & a & a \\ a & a & a \\ a & a & a \end{pmatrix} \middle| a \in Z_2 = \{0,1\} \right\} \cup \left\{ \begin{pmatrix} a & a & a \\ b & b & b \end{pmatrix} \middle| a \in Z_2 = \{0,1\} \right\}$ be a semigroup 5 set linear algebra over $Z_2$. Find

(a) semigroup 5 set linear subalgebra of V over $Z_2$.

(b) find a pseudo semigroup 5-set linear subalgebra of V over $Z_2$.

(c) Find a semigroup 5 linear operator on V.

(d) Is $\text{Hom}_{Z_2}$ (V, V) a semigroup set 5-linear algebra over $Z_2$?

285. Let $V = V_1 \cup V_2 \cup V_3 \cup V_4 \cup V_5 = \{Z_{12} \times Z_{12}\} \cup \{Z_{12}[x]$ | all polynomials of degree less than or equal to 4$\} \cup \{(a,a,a,a) | a \in Z_{12}\} \cup \left\{ \begin{pmatrix} a & a & a \\ a & a & a \\ a & a & a \end{pmatrix} \middle| a \in Z_{12} \right\} \cup \left\{ \begin{pmatrix} a & b \\ c & d \end{pmatrix} \middle| a,b,c,d \in Z_{12} \right\}$ be a semigroup set 5 linear algebra over $Z_{12}$. What is the 5-dimension of V over $Z_{12}$? Find a semigroup set 5 linear subalgebra of V over $Z_{12}$. Is $\text{Hom}_{Z_{12}}$ (V, V) a semigroup set 4 -linear algebra over $Z_{12}$?

286. Let $V = V_1 \cup V_2 \cup V_3 \cup V_4 = \{3Z^+ \times 3Z^+ \times 3Z^+\} \cup Z^+[x] \cup \left\{ \begin{pmatrix} a & a & a \\ a & a & a \end{pmatrix} \middle| a \in Z^+ \right\} \cup \left\{ \begin{pmatrix} a & b & c \\ 0 & d & e \\ 0 & 0 & f \end{pmatrix} \middle| a,b,c,d,e,f \in Z^+ \right\}$ be a semigroup 4-set linear algebra over $Z^+$.

(a) Find a semigroup 4-linear subalgebra of V over $Z^+$.
(b) Find a pseudo semigroup 4-linear subalgebra of V over $Z^+$.
(c) Find $\text{Hom}_{Z^+}$ (V, V).
(d) What is the 4-dimension of V over $Z^+$?

287. Let $V = V_1 \cup V_2 \cup V_3 \cup V_4 = \{Z_9 \times Z_9 \times Z_9\} \cup \{Z_9[x] \mid$ all polynomials of degree less than or equal to three$\}$ $\cup \left\{ \begin{pmatrix} a & a & a \\ a & a & a \\ a & a & a \end{pmatrix} \middle| a \in Z_9 \right\} \cup \{(a, a, a, a, a, a), (a\ a\ a\ a) / a \in Z_9\}$ be a semigroup 4-set vector space over $Z_9 = \{0, 1, 2, \ldots, 8\}$. Find a semigroup 4-set vector subspace of V. What is the 4-dimension of V over $Z_9$?
Find a 4-generating subset of V. Is $\text{Hom}_{Z_9}$ (V, V) a semigroup set vector space over $Z_9$?

288. Let $V = V_1 \cup V_2 \cup V_3 \cup V_4 = \{Z_{15} \times Z_{15} \times Z_{15}\}$ $\cup \left\{ \begin{pmatrix} a & a & a \\ a & a & a \\ a & a & a \end{pmatrix} \middle| a \in Z_{15} \right\} \cup \{Z_{15}[x]$ / all polynomials of degree less than or equal to $10\} \cup \left\{ \begin{pmatrix} a & a \\ a & a \\ a & a \\ a & a \end{pmatrix} \middle| a \in Z_{15} \right\}$ be a semigroup 4-linear algebra over the semigroup $Z_{15}$. Find a semigroup 4-linear subalgebra of V over $Z_{15}$. What is the 4-dimension of V over $Z_{15}$? Find a 4-generating subset X of V over $Z_{15}$. Is $\text{Hom}_{Z_{15}}$ (V, V) a semigroup 4-linear algebra over the semigroup $Z_{15}$?

289. Let $V = V_1 \cup V_2 \cup V_3 \cup V_4 \cup V_5 = \{Z_{10} \times Z_{10} \times Z_{10}\} \cup \left\{\begin{pmatrix} a & a & a \\ a & a & a \end{pmatrix}, \begin{pmatrix} a & a \\ a & a \end{pmatrix} \middle| a \in Z_{10}\right\} \cup \left\{\begin{pmatrix} a & b \\ c & d \end{pmatrix} \middle| a, b, c, d \in Z_{10}\right\} \cup \left\{\begin{pmatrix} a & a & a \\ a & a & a \\ a & a & a \\ a & a & a \end{pmatrix} \middle| a \in \{0, 5\}\right\} \cup \{(a\ a\ a\ a\ a\ a), (a\ a), (a\ a\ a\ a\ a) \mid a \in Z_{10}\}$ be a semigroup 5-set vector space over the semigroup $S = Z_{10}$. Find a semigroup 5-set vector subspace of V. Find a semigroup 5-set linear operator on V? Can V be made into a semigroup 5-set linear algebra over $Z_{10}$?

290. Let $V = V_1 \cup V_2 \cup V_3 \cup V_4 \cup V_5 = \{Z^+ [x]\} \cup \left\{\begin{pmatrix} a & b \\ c & d \end{pmatrix} \middle| a, b, c, d \in Z^+\right\} \cup \{Z^+ \times Z^+ \times Z^+\} \cup \left\{\begin{pmatrix} a & a & a & a \\ a & a & a & a \end{pmatrix}, \begin{pmatrix} a & a & a & a & a & a \\ a & a & a & a & a & a \end{pmatrix} \middle| a \in Z^+\right\} \cup \left\{\begin{pmatrix} a & a & a \\ a & a & a \\ a & a & a \end{pmatrix}, (a \quad a) \middle| a \in Z^+\right\}$ be a semigroup 5 set vector space over $Z^+$. Find a semigroup 5-set vector subspace of V. Find a 5-dimension generating subset of V.

291. Give some interesting properties about semigroup n set – vector space over a semigroup S.

292. Obtain some interesting properties about semigroup n-set linear algebra over a semigroup S.

293. Characterize those semigroup n-set linear algebras over a semigroup S which has pseudo semigroup n-linear algebras.

294. Let $V = V_1 \cup ... \cup V_n$ be a group n-set vector space over the group G and $W = W_1 \cup W_2 \cup ... \cup W_m$ $(m > n)$ be a group m set vector space over the same group G. Let $\text{Hom}_G^q$ (V, W) denote the collection of all quasi n-set linear transformation of V into W. What is the algebraic structure of $\text{Hom}_G^q$ (V, W)? Is $\text{Hom}_G^q$ (V, W) a group n- set vector space over the group G?

295. Obtain some interesting properties about $\text{Hom}_G^q$ (V, W).

296. Let $V = V_1 \cup V_2 \cup \ldots \cup V_n$ be a group n-set vector space over the group G. $W = W_1 \cup W_2 \cup \ldots \cup W_n$ $(n > m)$ be a group m-set vector space over the group G. Let $\text{Hom}_G^{pq}$ (V, W) denote the collection of all pseudo quasi group n-set linear transformation of V into W. What is the algebraic structure of $\text{Hom}_G^{pq}$ (V,W) ? Will $\text{Hom}_G^{pq}$ (V, W) be a group n-set vector space over G?

297. Find some interesting results about $\text{Hom}_G^{pq}$ (V, W).

298. If $V = V_1 \cup ... \cup V_n$ be a group n-set vector space over the group G. If the n-dimension of V is $(t_1, \ldots, t_n)$. Suppose $W = W_1 \cup ... \cup W_m$ $(m > n)$ be the group m set vector space over the group G. If the m-dimension of W is $(r_1, r_2, \ldots, r_m)$. What are the dimension of $\text{Hom}_G^q$ (V, W)? Likewise if $(m < n)$ what is the dimension of $\text{Hom}_G^{pq}$ (V, W)?

299. Find some good applications of set fuzzy vector spaces to other fields other than the ones mentioned in this book.

300. Let $V = V_1 \cup ... \cup V_n$ be a set n-vector space over the set S. How many set fuzzy n-vector spaces can be constructed?

301. Let $V = V_1 \cup V_2 \cup V_3 = \{Z^+[x]\} \cup \{Z^+ \times Z^+ \times Z^+ \times Z^+ \times Z^+\} \cup \left\{\begin{pmatrix} a & b & c \\ d & e & f \end{pmatrix} \middle| a,b,c,d,e,f \in Z^+\right\}$ be a set 3-vector space over the set $Z^+$. Find at least 5-distinct set fuzzy 3 vector spaces.

302. Let $V = \{A_i = (a_{ij}) \mid (a_{ij}) \in [0, 1]$ be a $n \times n$ matrix with sum of the rows equal to one $i = 1, 2, \ldots, k\}$ be a set vector space over the set $[0, 1]$. i.e., $V = (A_1, \ldots, A_k)$, $V_1 = \left(A_1^1,...,A_k^1\right)$, $V_2 = \left(A_1^2,...,A_k^2\right)$, … , $V_t = \left(A_1^t,...,A_k^t\right)$. Construct a fuzzy markov chain using V and apply it to web related problems.

303. Let $V = V_1 \cup V_2$ be a group vector bispace where $V_1 = \{A = [a_{ij}];\ m \times n$ matrices with entries from $Z\}$ $\cup \left\{\begin{pmatrix} a & a & a & a & a & a \\ a & a & a & a & a & a \end{pmatrix} \middle| a \in Z\right\}$ be a group vector bispace over the group $G = Z$.
(1) Find at least four group fuzzy vector bispace
(2) Find at least 4 distinct group fuzzy vector subbispace.

304. Let $V = V_1 \cup V_2 \cup V_3 \cup V_4 \cup V_5 = \{Z_5[x]\} \cup \left\{\begin{pmatrix} a & b \\ c & d \end{pmatrix} \middle| a,b,c,d \in Z_5\right\} \cup \{Z_5 \times Z_5\} \cup \left\{\begin{pmatrix} a & b \\ c & d \\ e & f \end{pmatrix} \middle| a,b,c,d,e,f \in Z_5\right\} \cup \left\{\begin{pmatrix} a & a & a & a \\ a & a & a & a \\ a & a & a & a \end{pmatrix} \middle| a \in Z_5\right\}$ be a group 5 vector space over the group $Z_5$. Find a group fuzzy vector space.

# INDEX

**D**

**F**

**G**

**N**

**P**

## Q

## S

## V

## W

# About the Authors

**Dr.W.B.Vasantha Kandasamy** is an Associate Professor in the Department of Mathematics, Indian Institute of Technology Madras, Chennai. In the past decade she has guided 12 Ph.D. scholars in the different fields of non-associative algebras, algebraic coding theory, transportation theory, fuzzy groups, and applications of fuzzy theory of the problems faced in chemical industries and cement industries.

She has to her credit 640 research papers. She has guided over 64 M.Sc. and M.Tech. projects. She has worked in collaboration projects with the Indian Space Research Organization and with the Tamil Nadu State AIDS Control Society. This is her 34rd book.

On India's 60th Independence Day, Dr.Vasantha was conferred the Kalpana Chawla Award for Courage and Daring Enterprise by the State Government of Tamil Nadu in recognition of her sustained fight for social justice in the Indian Institute of Technology (IIT) Madras and for her contribution to mathematics. (The award, instituted in the memory of Indian-American astronaut Kalpana Chawla who died aboard Space Shuttle Columbia). The award carried a cash prize of five lakh rupees (the highest prize-money for any Indian award) and a gold medal. She can be contacted at vasanthakandasamy@gmail.com You can visit her on the web at: http://mat.iitm.ac.in/~wbv or: http://www.vasantha.net

**Dr. Florentin Smarandache** is a Professor of Mathematics and Chair of Math & Sciences Department at the University of New Mexico in USA. He published over 75 books and 150 articles and notes in mathematics, physics, philosophy, psychology, rebus, literature.
In mathematics his research is in number theory, non-Euclidean geometry, synthetic geometry, algebraic structures, statistics, neutrosophic logic and set (generalizations of fuzzy logic and set respectively), neutrosophic probability (generalization of classical and imprecise probability). Also, small contributions to nuclear and particle physics, information fusion, neutrosophy (a generalization of dialectics), law of sensations and stimuli, etc. He can be contacted at smarand@unm.edu

**K. Ilanthenral** is the editor of The Maths Tiger, Quarterly Journal of Maths. She can be contacted at ilanthenral@gmail.com